\documentclass[11pt]{amsart}

\usepackage[T1]{fontenc}
\usepackage[utf8]{inputenc}
\usepackage{lmodern}
\usepackage{microtype}
\usepackage{amsmath,amssymb,amsthm,mathtools,mathrsfs}
\usepackage{enumitem}
\usepackage{graphicx}
\usepackage{float}
\usepackage[margin=1.15in]{geometry}
\usepackage[hidelinks]{hyperref}
\usepackage{bookmark}
\hypersetup{
  pdftitle={Factorial Calculi and the Canonical Stirling Defect of the Prime Bhargava Factorial},
  pdfauthor={Brian Diaz}
}

\allowdisplaybreaks[2]
\newtheorem{theorem}{Theorem}[section]
\newtheorem{conjecture}[theorem]{Conjecture}
\newtheorem{proposition}[theorem]{Proposition}
\newtheorem{lemma}[theorem]{Lemma}
\newtheorem{corollary}[theorem]{Corollary}
\theoremstyle{definition}
\newtheorem{definition}[theorem]{Definition}
\theoremstyle{remark}
\newtheorem{remark}[theorem]{Remark}
\newtheorem{example}[theorem]{Example}

\newcommand{\PP}{\mathbb P}
\newcommand{\NN}{\mathbb N}
\newcommand{\ZZ}{\mathbb Z}
\newcommand{\QQ}{\mathbb Q}
\newcommand{\RR}{\mathbb R}
\newcommand{\CC}{\mathbb C}
\newcommand{\Zp}{\mathbb Z_p}
\newcommand{\cA}{\mathcal A}
\newcommand{\cB}{\mathcal B}
\newcommand{\cC}{\mathcal C}
\newcommand{\cD}{\mathcal D}
\newcommand{\cE}{\mathcal E}

\newcommand{\cG}{\mathcal G}
\newcommand{\cH}{\mathcal H}
\newcommand{\cI}{\mathcal I}
\newcommand{\cK}{\mathcal K}
\newcommand{\cL}{\mathcal L}
\newcommand{\cM}{\mathcal M}
\newcommand{\cO}{\mathcal O}
\newcommand{\cP}{\mathcal P}
\newcommand{\cQ}{\mathcal Q}
\newcommand{\cR}{\mathcal R}
\newcommand{\cS}{\mathcal S}
\newcommand{\cT}{\mathcal T}
\newcommand{\cU}{\mathcal U}
\newcommand{\cV}{\mathcal V}
\newcommand{\cW}{\mathcal W}
\newcommand{\cX}{\mathcal X}
\newcommand{\cZ}{\mathcal Z}
\newcommand{\FactCalc}{\mathbf{FactCalc}}
\newcommand{\LayP}{\mathbf{Lay}_{\PP}}
\newcommand{\Tors}{\mathbf{Tors}}
\newcommand{\Per}{\operatorname{Per}}
\newcommand{\Lift}{\operatorname{Lift}}
\newcommand{\Rec}{\operatorname{Rec}}
\newcommand{\Ext}{\operatorname{Ext}}
\newcommand{\Sel}{\operatorname{Sel}}
\newcommand{\St}{\operatorname{St}}
\newcommand{\Ob}{\operatorname{Ob}}
\newcommand{\Int}{\operatorname{Int}}

\newcommand{\Norm}{\operatorname{N}}
\newcommand{\Id}{\operatorname{Id}}
\newcommand{\Dom}{\operatorname{Dom}}
\newcommand{\Fin}{\operatorname{Fin}}
\newcommand{\Logp}{\operatorname{Log}}
\newcommand{\Log}{\operatorname{Log}}
\newcommand{\res}{\operatorname{res}}
\newcommand{\supp}{\operatorname{supp}}

\newcommand{\cond}{\operatorname{cond}}
\newcommand{\dist}{\operatorname{dist}}
\newcommand{\lcm}{\operatorname{lcm}}
\newcommand{\sinc}{\operatorname{sinc}}
\newcommand{\vp}{v_p}
\newcommand{\e}{\mathrm e}
\newcommand{\ee}[1]{\e\!\left(#1\right)}
\newcommand{\dd}{\,\mathrm d}
\newcommand{\Harm}{\mathrm H}
\newcommand{\ind}{\mathbf 1}
\newcommand{\eps}{\varepsilon}
\newcommand{\dyad}{\mathrel{\sim}}
\newcommand{\pos}[1]{\left[#1\right]_{+}}
\newcommand{\norm}[1]{\lVert #1\rVert}
\newcommand{\widebar}[1]{\overline{#1}}
\newcommand{\boxt}{\boxtimes}
\newcommand{\boxplusT}{\mathbin{\boxplus}}
\newcommand{\projlimb}{\varprojlim^{\!b}}
\newcommand{\indlim}{\varinjlim}
\newcommand{\Deltaop}{\Delta}
\newcommand{\Mean}{\mathsf A}
\newcommand{\Ces}{\mathsf C}

\title[Factorial Calculi and the Prime Bhargava Factorial]{Factorial Calculi and the Canonical Stirling Defect\\of the Prime Bhargava Factorial}
\author{Brian Diaz}
\date{}
\begin{document}
\begin{abstract}
The classical Euler gamma function is distinguished among positive solutions of its recurrence by the Bohr--Mollerup characterization.  We develop an analogous normalization theory for the prime Bhargava factorial attached to the rational primes, constructing a canonical prime gamma function and placing it in an analytic category of factorial calculi.  The prime factorial is realized as a filtered symmetric-monoidal product of atomic valuation layers, while its reciprocal coefficients produce both Bhargava's unshifted prime exponential and the shifted eigenfunction intrinsic to the completed differential calculus.  The logarithmic recurrence admits a torsor of entire lifts under normalized entire $1$-periodic functions.  Centered cyclotomic quadrature supplies a normally convergent lift, and an orbitwise Stirling normalization theorem selects its unique admissible representative.  Thus orbitwise normalization plays for the prime recurrence the same conceptual role that Bohr--Mollerup plays for Euler's recurrence, without asserting a direct extension of the classical convexity theorem.

The resulting zero-free entire function $\Gamma_{\mathbb P}^{\mathrm{cyc}}$ interpolates the prime Bhargava factorial, satisfies its exact recurrence, and obeys the Euler-type reflection law
\[
 \Gamma_{\mathbb P}^{\mathrm{cyc}}(z)
 \Gamma_{\mathbb P}^{\mathrm{cyc}}(1-z)=1,
 \qquad
 \Gamma_{\mathbb P}^{\mathrm{cyc}}\!\left(\frac12\right)=1.
\]
We include a conductor-truncated real-axis graph of this distinguished function.  At the positive integers the construction yields the explicit unsmoothed Stirling formula
\[
 \log (n+1)!_{\mathbb P}
 =n\log n+(C_{\mathbb P}-1)n
 +\frac12\log(2\pi n)
 +\mathcal R_{\mathbb P}(n)
 +O\!\left(\frac1n\right),
 \qquad
 C_{\mathbb P}=\sum_p\frac{\log p}{(p-1)^2},
\]
where $\mathcal R_{\mathbb P}(n)=o(n)$ unconditionally and admits an exact absolutely convergent prime-local fractional-part expansion.  This gives a prime-case response to Questions~31 and~33 of Bhargava's 2000 paper \emph{The Factorial Function and Generalizations}.

No positive Euler--Mellin representation can exist for either $\Gamma_{\mathbb P}^{\mathrm{cyc}}(z+1)$ or its reciprocal, as the corresponding factorial sequences fail the necessary Hankel positivity conditions.  Nevertheless, the cyclotomic construction gives an exact logarithmic distributional integral and motivates a prime Hankel conjecture asserting a canonical signed, complex, Jackson-type, distributional, or contour representation of the reciprocal.

Finally, we study the sharp unsmoothed second moment of $\mathcal R_{\mathbb P}(n)$.  All higher-conductor, projective, affine long-line, Type~I, and long-free Type~II sectors are controlled unconditionally, leaving one explicit centered subcritical Vaughan character moment.  Assuming the stated square-root bound for that moment, we prove
\[
 \sum_{n\le X}\mathcal R_{\mathbb P}(n)^2
 \sim -\frac43\zeta\!\left(-\frac12\right)X^{3/2}\log X.
\]
Thus the categorical construction, canonical continuation, reflection law, and Stirling expansion are unconditional, while the variance asymptotic has all of its conditionality isolated in a single spectral estimate.
\end{abstract}
\maketitle
\clearpage
\tableofcontents
\clearpage
\part{Factorial calculi}
% ===== Source: section1_introduction_main_results.tex =====
\section{Introduction and main results}
\label{sec:introduction}

Euler's recurrence
\[
 \Gamma(x+1)=x\Gamma(x)
\]
does not by itself determine the classical gamma function.  The Bohr--Mollerup theorem~\cite{BohrMollerup1922,Artin1964} supplies the missing normalization: positivity, log-convexity, and the value $\Gamma(1)=1$ single out the classical solution on the positive real axis.  This suggests a general question for arithmetic factorials: when a factorial datum determines a recurrence with many analytic lifts, what principle identifies the correct lift?  The purpose of this paper is to answer that question for the Bhargava factorial attached to the rational primes.  Our replacement for log-convexity is an orbitwise Stirling normalization inside the periodic torsor of entire recurrence lifts.

The classical factorial $n!$ supports several parallel structures.  It gives the weights of ordinary differentiation on monomials, supplies the reciprocal coefficients of the exponential series, and extends through Euler's gamma function to a holomorphic object governed by a recurrence and a normalization principle.  Jackson's $q$-gamma functions~\cite{Jackson1905,GasperRahman2004,DeSoleKac2005} exhibit the same general pattern in a deformed setting: discrete factorial data determine a recurrence, but an additional analytic condition is required to select a distinguished continuation.  We develop an analogous framework for the prime Bhargava factorial and explain how its arithmetic, operator-theoretic, gamma-theoretic, and asymptotic structures fit into one categorical picture.

Our starting point is the prime Bhargava factorial~\cite{Bhargava2000,DiazNormanyo2026}
\begin{equation}
 (n+1)!_{\PP}
 =
 \prod_p
 p^{\displaystyle\sum_{j\ge 0}
 \left\lfloor \frac{n}{(p-1)p^j}\right\rfloor},
 \qquad n\ge 0,
 \label{eq:prime-bhargava-factorial}
\end{equation}
where $\PP$ denotes the set of rational primes.  The expression in \eqref{eq:prime-bhargava-factorial} is assembled from atomic valuation layers indexed by pairs $(p,j)$, each contributing a step factor of period $(p-1)p^j$.  Every finite collection of layers has exponential growth, whereas the filtered completion has factorial growth.  A central theme of the paper is therefore that Stirling behavior is not visible layer by layer: it emerges only after completion.

To formalize this phenomenon, we introduce a category of \emph{factorial calculi}.  A positive factorial datum
\[
 F=(F_n)_{n\ge 0},
 \qquad F_0=1,
\]
determines the weighted Hilbert space
\[
 \cH_F
 =
 \left\{
 f(z)=\sum_{n\ge 0}a_nz^n:
 \sum_{n\ge 0}|a_n|^2F_n<\infty
 \right\}
\]
and the associated generalized derivative
\[
 D_Fz^n=\frac{F_n}{F_{n-1}}z^{n-1}.
\]
Its canonical exponential is
\begin{equation}
 \mathcal E_F(z)=\sum_{n\ge 0}\frac{z^n}{F_n},
 \qquad
 D_F\mathcal E_F=\mathcal E_F,
 \label{eq:intro-generalized-exponential}
\end{equation}
whenever the series converges.  If two factorial data $F$ and $G$ satisfy
\[
 F\preccurlyeq G
 \qquad\Longleftrightarrow\qquad
 \sup_{n\ge 0}\frac{F_n}{G_n}<\infty,
\]
then there is a canonical bounded diagonal intertwiner
\[
 J_{F,G}:\cH_F\longrightarrow\cH_G,
 \qquad
 J_{F,G}z^n=\frac{F_n}{G_n}z^n,
\]
with
\[
 D_GJ_{F,G}=J_{F,G}D_F.
\]
The adjoint $J_{F,G}^{*}$ is the canonical inclusion $\cH_G\hookrightarrow\cH_F$.  Thus the same diagonal mechanism yields an inductive structure on generalized differential calculi and a projective structure on coefficient spaces.

This construction gives a concrete prime-case response to two questions posed by Bhargava in his 2000 paper \emph{The Factorial Function and Generalizations}~\cite{Bhargava2000}.  Question~31 asks for analogues of Stirling's formula for generalized factorials, while Question~33 asks for the $S$-analogue of the exponential function and its properties.  With the convention $0!_{\PP}=1$, the unshifted Bhargava prime exponential is
\begin{equation}
 \operatorname{Exp}_{\PP}(z)
 =
 \sum_{n\ge 0}\frac{z^n}{n!_{\PP}}.
 \label{eq:intro-unshifted-prime-exponential}
\end{equation}
The completed factorial calculus is indexed instead by $A_n=(n+1)!_{\PP}$, and therefore carries the shifted eigenfunction
\begin{equation}
 \mathcal E_{\PP}(z)
 =
 \sum_{n\ge 0}\frac{z^n}{(n+1)!_{\PP}},
 \qquad
 D_A\mathcal E_{\PP}=\mathcal E_{\PP}.
 \label{eq:intro-shifted-prime-exponential}
\end{equation}
The two normalizations are related by the exact identity
\begin{equation}
 \operatorname{Exp}_{\PP}(z)
 =1+z\mathcal E_{\PP}(z).
 \label{eq:intro-exponential-shift-relation}
\end{equation}
Thus \eqref{eq:intro-unshifted-prime-exponential} is the natural answer to Bhargava's Question~33, while \eqref{eq:intro-shifted-prime-exponential} is the canonical exponential intrinsic to the operator calculus used in this paper.  The gamma and Stirling constructions developed below provide the corresponding analytic continuation and asymptotic theory.  The paper does not propose a universal answer for every subset $S$; rather, it gives a complete structural answer for the prime factorial and a categorical language designed to isolate which ingredients generalize.

The prime Bhargava factorial fits naturally into this framework.  Each valuation layer defines an atomic factorial datum; finite products of layers form a filtered symmetric-monoidal system; and the completed norm realization recovers \eqref{eq:prime-bhargava-factorial}.  The first objective of the paper is to make this structure precise.

The second objective is to construct and characterize a distinguished \emph{prime gamma object}.  An entire recurrence symbol $L$ may admit many entire primitives $H$ satisfying
\[
 H(z+1)-H(z)=L(z),
 \qquad H(1)=0.
\]
The normalized lifts form a torsor under the additive group of normalized entire $1$-periodic functions.  Thus recurrence and interpolation alone do not determine the continuation.  This is the prime analogue of the familiar need for a normalization principle in the classical and $q$-gamma settings, but the ambiguity here is explicitly periodic.  We show that the weak Stirling asymptotic fixes the correct affine comparison model
\[
 M_{\PP}(z)
 =
 \log\Gamma(z)+C_{\PP}(z-1),
 \qquad
 C_{\PP}=\sum_p\frac{\log p}{(p-1)^2},
\]
and that an orbitwise Stirling normalization theorem selects a unique normalized lift $H_{\PP}^{\mathrm{cyc}}$.  Exponentiation gives the distinguished zero-free entire prime gamma function
\[
 \Gamma_{\PP}^{\mathrm{cyc}}(z)
 =
 \exp\!\bigl(H_{\PP}^{\mathrm{cyc}}(z)\bigr).
\]
It satisfies an Euler-type reflection law
\begin{equation}
 \Gamma_{\PP}^{\mathrm{cyc}}(z)
 \Gamma_{\PP}^{\mathrm{cyc}}(1-z)=1,
 \label{eq:intro-prime-reflection}
\end{equation}
and hence the distinguished midpoint evaluation
\begin{equation}
 \Gamma_{\PP}^{\mathrm{cyc}}\!\left(\frac12\right)=1.
 \label{eq:intro-prime-midpoint}
\end{equation}

The associated analytic defect is
\[
 \mathscr R_{\PP}(z)
 =
 H_{\PP}^{\mathrm{cyc}}(z+1)-M_{\PP}(z+1).
\]
At the positive integers, set
\begin{equation}
 \cR_{\PP}(n)
 =
 \log (n+1)!_{\PP}-\log n!-C_{\PP}n.
 \label{eq:prime-stirling-defect}
\end{equation}
Combining this identity with the classical Stirling expansion gives the explicit unsmoothed formula
\begin{equation}
 \boxed{
 \log (n+1)!_{\PP}
 =
 n\log n+(C_{\PP}-1)n
 +\frac12\log(2\pi n)
 +\cR_{\PP}(n)
 +O\!\left(\frac1n\right).
 }
 \label{eq:intro-full-stirling}
\end{equation}
The term $\cR_{\PP}(n)$ is not introduced by smoothing or interpolation: it is the exact arithmetic discrepancy in \eqref{eq:prime-stirling-defect}.  Smoothing enters only in the analytic theorem that selects the unique continuation $H_{\PP}^{\mathrm{cyc}}$ from its periodic torsor.

More importantly, \eqref{eq:intro-full-stirling} is not merely a formal rewriting with an uncontrolled remainder.  The discrepancy admits the exact prime-local expansion
\begin{equation}
 \cR_{\PP}(n)
 =
 \sum_p (\log p)
 \left(
  \sum_{k\ge 1}\left\{\frac{n}{p^k}\right\}
  -
  \sum_{j\ge 0}
  \left\{\frac{n}{(p-1)p^j}\right\}
 \right),
 \label{eq:intro-prime-local-discrepancy}
\end{equation}
whose displayed prime series is absolutely convergent.  In particular, one has the unconditional pointwise control
\begin{equation}
 \cR_{\PP}(n)=o(n).
 \label{eq:intro-sublinear-discrepancy}
\end{equation}
Thus the new term is a genuinely arithmetic, sublinear discrepancy rather than an unspecified error absorbed into the main expansion.  Under Conjecture~\ref{conj:intro-vaughan}, its averaged size is much sharper:
\begin{equation}
 \sum_{n\le X}\cR_{\PP}(n)^2
 \sim
 -\frac43\zeta\!\left(-\frac12\right)
 X^{3/2}\log X.
 \label{eq:intro-second-moment-discrepancy}
\end{equation}
This is an unsmoothed second-moment asymptotic, not a pointwise estimate; equivalently, the typical quadratic fluctuation scale is of order $X^{1/4}(\log X)^{1/2}$.

The third objective is analytic.  We study the second moment of $\cR_{\PP}(n)$.  The resulting variance problem has a rigid internal geometry.  Its long affine-line sector is controlled by an affine-lattice dispersion argument in the tradition of~\cite{BFI1986,DeshouillersIwaniec1982}, while the short balanced sector is reduced to a centered Vaughan-type spectral moment.  This reduction is unconditional.  The final asymptotic is conditional on a square-root cancellation estimate for the remaining centered character moment.

We now state the principal results.

\subsection{Factorial calculi and prime completion}
\label{subsec:factorial-calculi}

For factorial data $F$ and $G$, define
\[
 (F\boxt G)_n:=F_nG_n.
\]
The unit is the constant datum $\mathbf 1_n=1$.  For a prime layer $\lambda=(p,j)$, write
\[
 d_\lambda=(p-1)p^j,
 \qquad
 w_\lambda=\log p,
\]
and set
\[
 F_n^{(\lambda)}
 =
 \exp\!\left(
 w_\lambda\left\lfloor\frac{n}{d_\lambda}\right\rfloor
 \right).
\]
For a finite layer set $\Lambda$, define
\[
 A^\Lambda
 =
 \boxt_{\lambda\in\Lambda}F^{(\lambda)}.
\]
These objects form a filtered system as $\Lambda$ ranges over finite subsets of the layer set.

\begin{theorem}[Categorical prime realization]
\label{thm:intro-prime-realization}
The finite prime-layer data $A^\Lambda$ form a filtered symmetric-monoidal system in $\FactCalc$.  Their generalized derivatives assemble through bounded diagonal intertwiners, while their coefficient Hilbert spaces assemble through the adjoint inclusions.  The completed factorial datum
\[
 A_n=\sup_\Lambda A_n^\Lambda
\]
is the norm realization of the prime Bhargava factorial:
\[
 A_n=(n+1)!_{\PP}.
\]
Moreover,
\[
 \cH_A\cong\varprojlim_\Lambda\cH_{A^\Lambda}
\]
in the natural bounded projective-limit sense, and the entire function
\[
 \mathcal E_{\PP}(z)=\sum_{n\ge0}\frac{z^n}{A_n}
\]
satisfies
\[
 D_A\mathcal E_{\PP}=\mathcal E_{\PP},
 \qquad
 \mathcal E_{\PP}(0)=1.
\]
The unshifted Bhargava exponential
\[
 \operatorname{Exp}_{\PP}(z)
 =\sum_{n\ge0}\frac{z^n}{n!_{\PP}}
\]
is entire and obeys
\[
 \operatorname{Exp}_{\PP}(z)=1+z\mathcal E_{\PP}(z).
\]
\end{theorem}

The theorem makes precise both halves of the completed calculus: the prime Bhargava factorial is a filtered monoidal completion of atomic valuation calculi, while its reciprocal coefficients produce both the unshifted Bhargava exponential and the shifted eigenfunction intrinsic to the completed derivative.

\subsection{Gamma-lift torsors, canonical normalization, and reflection}
\label{subsec:gamma-lifts}

For an entire recurrence symbol $L$, define
\[
 \Lift(L)
 =
 \left\{
 H\in\cO(\mathbb C):
 H(z+1)-H(z)=L(z),\ H(1)=0
 \right\}.
\]
When nonempty, this set is a torsor under
\[
 \Per_0
 =
 \left\{
 P\in\cO(\mathbb C):
 P(z+1)=P(z),\ P(1)=0
 \right\}.
\]
For the completed prime recurrence, a symmetric cyclotomic construction produces a normalized lift, and the orbitwise Stirling theorem proves that this lift is the unique admissible point of the torsor.

\begin{theorem}[Canonical prime gamma object and Stirling defect]
\label{thm:intro-canonical-defect}
Let $L_{\PP}^{\mathrm{cyc}}$ denote the completed prime recurrence symbol.  Then the following statements hold.
\begin{enumerate}[label=\textup{(\roman*)}]
 \item The weak first-order asymptotic
 \[
  \log (n+1)!_{\PP}
  =
  \log n!+C_{\PP}n+o(n)
 \]
 uniquely determines the affine comparison model
 \[
  M_{\PP}(z)=\log\Gamma(z)+C_{\PP}(z-1).
 \]
 \item Among the normalized entire lifts of $L_{\PP}^{\mathrm{cyc}}$, the orbitwise Stirling condition selects a unique lift $H_{\PP}^{\mathrm{cyc}}$.
 \item The function
 \[
  \Gamma_{\PP}^{\mathrm{cyc}}(z)
  =e^{H_{\PP}^{\mathrm{cyc}}(z)}
 \]
 is zero-free and entire, interpolates the prime factorial,
 \[
  \Gamma_{\PP}^{\mathrm{cyc}}(m)=m!_{\PP}
  \qquad(m\ge1),
 \]
 and satisfies the recurrence
 \[
  \Gamma_{\PP}^{\mathrm{cyc}}(z+1)
  =e^{L_{\PP}^{\mathrm{cyc}}(z)}
   \Gamma_{\PP}^{\mathrm{cyc}}(z).
 \]
 \item The Euler-type reflection law and midpoint value are
 \[
  \Gamma_{\PP}^{\mathrm{cyc}}(z)
  \Gamma_{\PP}^{\mathrm{cyc}}(1-z)=1,
  \qquad
  \Gamma_{\PP}^{\mathrm{cyc}}\!\left(\frac12\right)=1.
 \]
 \item The difference
 \[
  \mathscr R_{\PP}(z)
  =
  H_{\PP}^{\mathrm{cyc}}(z+1)-M_{\PP}(z+1)
 \]
 is canonically determined, satisfies
 \[
  \mathscr R_{\PP}(n)=\cR_{\PP}(n)
 \]
 for every positive integer $n$, and yields the unsmoothed Stirling expansion \eqref{eq:intro-full-stirling}.
\end{enumerate}
\end{theorem}

Every normalized lift of the same recurrence agrees at the positive integers, because a normalized $1$-periodic perturbation vanishes there.  Orbitwise rigidity therefore does not alter the discrete interpolation data; it makes their analytic continuation and transverse behavior across translation orbits unique.  In this sense it plays for the prime gamma torsor the selecting role that convexity-type normalizations play for classical and $q$-gamma functions.

\begin{remark}[A norm-one entire unit]
Let $\iota(z)=1-z$.  Since $\Gamma_{\PP}^{\mathrm{cyc}}$ is zero-free and entire, it lies in the unit group $\cO(\mathbb C)^\times$.  The reflection law can be written
\[
 \Gamma_{\PP}^{\mathrm{cyc}}\,\iota^*\Gamma_{\PP}^{\mathrm{cyc}}=1.
\]
Thus the prime gamma function is a norm-one unit for the reflection involution.  The identity is structurally analogous to Euler's classical reflection formula, but its right-hand side is the trivial unit rather than $\pi/\sin(\pi z)$.
\end{remark}

\begin{figure}[H]
 \centering
 \includegraphics[width=0.96\textwidth]{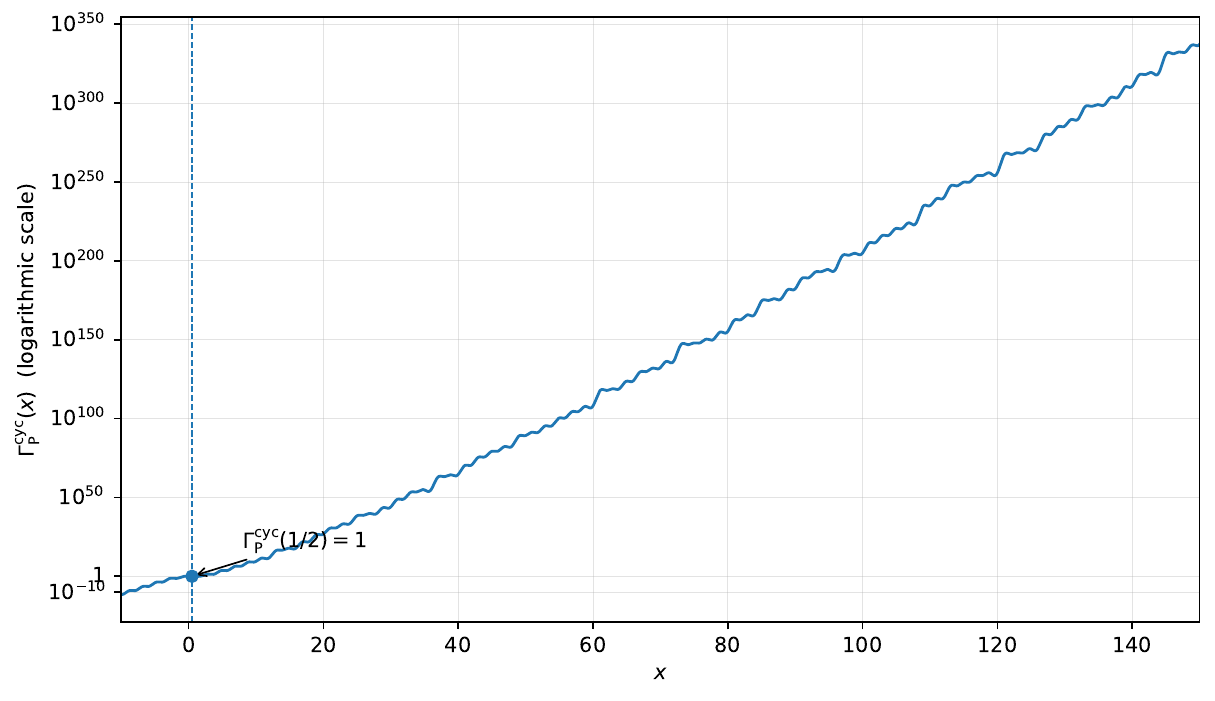}
 \caption{A conductor-truncated numerical plot of the distinguished prime gamma function on the real axis.  The ordinate is logarithmic, the horizontal range is $-10\le x\le150$, and the curve is computed from the truncation $H_{[10^4]}^{\mathrm{cyc}}$.  This truncation is exact at every displayed positive integer, and the conductor approximants converge normally on compact subsets as the truncation grows.  Each finite truncation already satisfies the reflection law, and the marked point records $\Gamma_{\PP}^{\mathrm{cyc}}(1/2)=1$.}
 \label{fig:intro-prime-gamma-real-axis}
\end{figure}

The classical Euler integral does not have a positive-measure analogue for the prime gamma function.  Indeed, the values
\[
 1!_{\PP}=1,\qquad 2!_{\PP}=2,\qquad
 3!_{\PP}=24,\qquad 4!_{\PP}=48
\]
violate the shifted Hankel positivity required of a Stieltjes moment sequence; the reciprocal values violate ordinary Hankel positivity already at order two.  On the other hand, the cyclotomic construction gives an exact logarithmic integral representation
\[
 H_{\PP}^{\mathrm{cyc}}(z)
 =\left\langle\mathfrak N_{\PP},K_z\right\rangle,
 \qquad
 \frac{1}{\Gamma_{\PP}^{\mathrm{cyc}}(z)}
 =\exp\!\left(-\left\langle\mathfrak N_{\PP},K_z\right\rangle\right),
\]
for a canonical order-two prime-cyclotomic distribution $\mathfrak N_{\PP}$ on $[-1/2,1/2]$.  Section~11 proves these assertions.  They suggest a linear integral representation for the reciprocal which cannot be positive, but may still exist in a signed, complex, distributional, Jackson-type, or contour form; compare the $q$-integral models of~\cite{DeSoleKac2005}.

\begin{conjecture}[Prime Hankel representation]
\label{conj:intro-prime-hankel}
There exist an oriented contour $\mathscr C\subset\mathbb C^\times$, a specified branch of $w^{-z}=\exp(-z\log w)$ along $\mathscr C$, and a canonical prime-cyclotomic kernel $\mathcal K_{\PP}(w)$, independent of $z$, such that on a nonempty vertical strip
\begin{equation}
 \boxed{
 \frac1{\Gamma_{\PP}^{\mathrm{cyc}}(z)}
 =\frac1{2\pi i}\int_{\mathscr C}
   \mathcal K_{\PP}(w)w^{-z}\,\mathrm dw.
 }
 \label{eq:intro-prime-hankel-target}
\end{equation}
The right-hand side admits entire continuation to $\CC$, and the recurrence
\[
 \Gamma_{\PP}^{\mathrm{cyc}}(z+1)
 =e^{L_{\PP}^{\mathrm{cyc}}(z)}\Gamma_{\PP}^{\mathrm{cyc}}(z)
\]
and reflection law
\[
 \Gamma_{\PP}^{\mathrm{cyc}}(z)
 \Gamma_{\PP}^{\mathrm{cyc}}(1-z)=1
\]
are induced by explicit transformations of the contour kernel.
\end{conjecture}

The conjecture asks for more than an abstract inverse Mellin transform.  The adjective \emph{canonical} requires the kernel to arise functorially from the prime-cyclotomic layer data and to expose, rather than merely reproduce after analytic continuation, the defining recurrence and Euler-type reflection symmetry.

\subsection{Unsmoothed fluctuations}
\label{subsec:unsmoothed-fluctuations}

We analyze
\[
 \sum_{n\le X}\cR_{\PP}(n)^2.
\]
The proof begins with an exact arithmetic expansion and covariance decomposition.  It then passes through quotient geometry, a long-line/short-line decomposition for the relevant affine congruences, and a complete Type~I/Type~II analysis of the short sector.  The geometric and multiplicative-energy obstructions can all be removed unconditionally.  The final residual term is a centered small-free-variable Vaughan moment.

\begin{theorem}[Unconditional reduction]
\label{thm:intro-unconditional-reduction}
The second moment of $\cR_{\PP}(n)$ admits an unconditional decomposition with the following properties.
\begin{enumerate}[label=\textup{(\roman*)}]
 \item The diagonal and projective-diagonal contributions are evaluated.
 \item The affine long-line sector is controlled at the required power scale.
 \item The relevant numerator multiplicative energies satisfy divisor-bounded $L^2$ estimates.
 \item Every packet containing a long coefficient-$1$ Vaughan variable is negligible by Poisson summation.
 \item The remaining short balanced sector is reduced to a centered subcritical Vaughan character moment.
\end{enumerate}
Thus the only unresolved unconditional input is a centered spectral estimate for the residual short-free-variable packet.
\end{theorem}

The remaining conjecture retains the exact provenance of the Vaughan decomposition.  Let $K_s(\chi)$ denote the short coefficient-$1$ character polynomial, let $D_s(\chi)$ denote the localized M\"obius polynomial, and let $P_{s,L}^{\circ}(\chi)$ denote the centered prime character sum.

\begin{conjecture}[Subcritical centered Vaughan moment]
\label{conj:intro-vaughan}
For dyadic parameters satisfying
\[
 K\le SX^\eta,
 \qquad
 DKL\asymp P,
 \qquad
 S^2\ll L,
\]
one has
\[
 \sum_{s\asymp S}
 \frac{1}{\varphi(s)}
 \sum_{\chi\ne\chi_0}
 |K_s(\chi)|^2
 |D_s(\chi)|^2
 |P_{s,L}^{\circ}(\chi)|^2
 \ll_\varepsilon
 SDKL\,X^\varepsilon.
\]
The character family is understood with the principal and explicitly separated exceptional or major-arc modes removed in the manner specified in the body of the paper.
\end{conjecture}

This is substantially narrower than a generic large-sieve statement: the conjecture uses the shortness of the free variable, the truncated M\"obius provenance, the centering of the prime major arcs, and the factorization-sensitive short-line mask.

\subsection{Conditional variance law}
\label{subsec:conditional-variance}

Under Conjecture~\ref{conj:intro-vaughan}, or under a square-root character hypothesis strong enough to imply it, the complete variance asymptotic follows.

\begin{theorem}[Conditional unsmoothed variance law]
\label{thm:intro-conditional-variance}
Assume Conjecture~\ref{conj:intro-vaughan}.  Then
\[
 \sum_{n\le X}\cR_{\PP}(n)^2
 \sim
 -\frac{4}{3}\zeta\!\left(-\frac12\right)
 X^{3/2}\log X.
\]
Equivalently,
\[
 \sum_{n\le X}
 \left(
 \log (n+1)!_{\PP}
 -\log n!
 -C_{\PP}n
 \right)^2
 \sim
 -\frac{4}{3}\zeta\!\left(-\frac12\right)
 X^{3/2}\log X.
\]
\end{theorem}

The sequence in Theorem~\ref{thm:intro-conditional-variance} is completely unsmoothed.  Smoothing enters only in the analytic selection and construction of the distinguished lift $H_{\PP}^{\mathrm{cyc}}$, not in the fluctuation sequence itself.

\subsection{Conceptual summary}
\label{subsec:conceptual-summary}

The architecture of the paper can be summarized by the chain
\[
 \boxed{
 \begin{gathered}
 \text{arithmetic factorial datum}
 \longrightarrow
 \text{categorical calculus and exponential}
 \longrightarrow
 \text{gamma-lift torsor}
 \\
 \longrightarrow
 \text{orbitwise Stirling selection}
 \longrightarrow
 \text{unsmoothed fluctuation theory}.
 \end{gathered}
 }
\]
More concretely:
\begin{enumerate}[label=\textup{(\roman*)}]
 \item the prime Bhargava factorial is realized as a filtered monoidal completion in a category of factorial calculi;
 \item its reciprocal coefficients define the unshifted prime exponential $\operatorname{Exp}_{\PP}$ and the shifted calculus eigenfunction $\mathcal E_{\PP}$;
 \item its associated gamma object is naturally a torsor of normalized lifts rather than a single function;
 \item the weak Stirling coefficient selects the affine comparison model, and orbitwise rigidity selects the distinguished lift;
 \item the distinguished prime gamma function obeys the Euler-type reflection law \eqref{eq:intro-prime-reflection}, with fixed-point value \eqref{eq:intro-prime-midpoint};
 \item its positive-integer values yield the explicit unsmoothed expansion \eqref{eq:intro-full-stirling};
 \item the second moment of the resulting fluctuation is reduced unconditionally to Conjecture~\ref{conj:intro-vaughan} and evaluated conditionally once that final input is assumed.
\end{enumerate}

The paper therefore has two intertwined themes.  The first is structural: factorial calculi, generalized exponentials, filtered completion, gamma-lift torsors, reflection, and Stirling selection.  The second is analytic: the dispersion theory of the resulting canonical unsmoothed defect.

\subsection{Organization of the paper}
\label{subsec:organization}

The manuscript is organized into four parts.  Part~I develops the category of factorial calculi, the generalized exponential, the symmetric-monoidal layer formalism, and the filtered completion realizing the prime Bhargava factorial.  Part~II studies recurrence symbols and gamma-lift torsors, culminating in orbitwise Stirling selection, the Euler-type reflection law, and the construction of the canonical analytic defect $\mathscr R_{\PP}$.  Part~III analyzes the sharp second moment of $\cR_{\PP}(n)$, beginning with the exact arithmetic expansion and ending with the reduction to the centered subcritical Vaughan moment and the conditional variance law.  Part~IV contains the functional-analytic, cyclotomic, divisor-energy, gcd-descent, affine-dispersion, and Vaughan--Poisson appendices needed to keep the main argument auditable.

\section{Bhargava factorial ideals and norm realization}
\label{sec:bhargava-factorials}

This section records the arithmetic object from which the later factorial calculus is built.  We begin with the local ordering construction for an arbitrary subset of the integers, specialize it to the set of rational primes, and then pass from the resulting factorial ideals to their positive norms.  The shift
\[
 A_n=(n+1)!_{\PP},\qquad n\ge 0,
\]
is chosen so that the prime factorial admits the layer expansion
\[
 \log A_n
 =
 \sum_p\sum_{j\ge 0}
 \left\lfloor\frac{n}{(p-1)p^j}\right\rfloor\log p.
\]
This expansion will later become the filtered symmetric-monoidal decomposition of the prime factorial calculus.

\subsection{Characteristic ideals and local orderings}
\label{subsec:characteristic-ideals}

Let $S\subseteq\ZZ$ be infinite.  Its ring of integer-valued polynomials is
\[
 \Int(S,\ZZ)
 =
 \{f\in\QQ[x]:f(S)\subseteq\ZZ\}.
\]
For $k\ge 0$, define the $k$th characteristic module by
\[
 \cI_k(S)
 =
 \left\{
 a_k\in\QQ:
 \text{there exists }
 f(x)=a_kx^k+\cdots+a_0\in\Int(S,\ZZ)
 \right\}.
\]
It is a fractional ideal of $\ZZ$.  Since every fractional ideal of $\ZZ$ is principal, its inverse has a unique positive generator.

\begin{definition}[Factorial ideal]
\label{def:factorial-ideal}
The $k$th Bhargava factorial ideal of $S$ is
\[
 \mathfrak F_k(S)=\cI_k(S)^{-1}.
\]
Its positive generator is denoted by $k!_S$.  Thus
\[
 \mathfrak F_k(S)=k!_S\ZZ.
\]
\end{definition}

The local description of $\mathfrak F_k(S)$ is obtained from $p$-orderings.  Fix a rational prime $p$ and write $\vp$ for the normalized $p$-adic valuation.

\begin{definition}[$p$-ordering]
\label{def:p-ordering}
A $p$-ordering of $S$ is a sequence $(a_k)_{k\ge 0}$ in $S$ obtained recursively as follows.  Choose $a_0\in S$.  Having chosen $a_0,\dots,a_{k-1}$, choose $a_k\in S$ so that
\[
 \vp\!\left(\prod_{i<k}(a_k-a_i)\right)
\]
is minimal among all choices in $S$.  The minimum is denoted by
\[
 \nu_k(S,p)
 =
 \vp\!\left(\prod_{i<k}(a_k-a_i)\right).
\]
We put $\nu_0(S,p)=0$.
\end{definition}

Bhargava's local ordering theorem~\cite{Bhargava1997,Bhargava2000,CahenChabert1997} asserts that $\nu_k(S,p)$ is independent of the chosen $p$-ordering.  It also identifies these local exponents with the characteristic ideals.

\begin{theorem}[Local factorial formula]
\label{thm:local-factorial-formula}
For every infinite $S\subseteq\ZZ$ and every $k\ge 0$,
\[
 \mathfrak F_k(S)
 =
 \prod_p (p)^{\nu_k(S,p)}.
\]
Equivalently,
\[
 k!_S
 =
 \prod_p p^{\nu_k(S,p)}.
\]
For fixed $k$, only finitely many local exponents are nonzero.
\end{theorem}

The theorem may be viewed as a local-to-global denominator statement.  A $p$-ordering supplies a Newton basis whose $k$th denominator has valuation $\nu_k(S,p)$; intersecting the resulting local leading-coefficient modules over all primes recovers $\cI_k(S)$ and hence $\mathfrak F_k(S)$.

We shall also use the following elementary invariance principle.

\begin{lemma}[$p$-adic closure invariance]
\label{lem:closure-invariance}
Let $S,T\subseteq\Zp$ be dense in the same $p$-adically closed subset.  Then
\[
 \nu_k(S,p)=\nu_k(T,p)
\]
for every $k\ge 0$.
\end{lemma}

\begin{proof}
After finitely many elements have been chosen, the function
\[
 x\longmapsto
 \vp\!\left(\prod_{i<k}(x-a_i)\right)
\]
takes values in $\ZZ_{\ge 0}\cup\{\infty\}$ and is locally constant on the complement of its finite zero set.  Its finite minimum on a closed $p$-adic set is therefore attained on a nonempty open subset.  Every dense subset meets that open set.  Induction on $k$ shows that the same sequence of minima is obtained from either dense subset.
\end{proof}

\subsection{The local ordering of the prime set}
\label{subsec:prime-local-ordering}

Fix a rational prime $p$.  Dirichlet's theorem on primes in arithmetic progressions implies that the rational primes distinct from $p$ are dense in $\Zp^{\times}$.  The prime $p$ itself is isolated from those unit classes.  Consequently,
\begin{equation}
 \overline{\PP}^{\,p}
 =
 \{p\}\sqcup\Zp^{\times},
 \label{eq:p-adic-closure-primes}
\end{equation}
where the closure is taken inside $\Zp$.

We first compute the $p$-sequence of the unit group.  Let
\[
 1=u_0<u_1<u_2<\cdots
\]
be the increasing sequence of positive integers not divisible by $p$.

\begin{lemma}[Balanced ordering of the $p$-adic units]
\label{lem:unit-p-ordering}
The sequence $(u_k)_{k\ge 0}$ is a $p$-ordering of $\Zp^{\times}$, and
\begin{equation}
 \nu_k(\Zp^{\times},p)
 =
 \sum_{r\ge 1}
 \left\lfloor\frac{k}{\varphi(p^r)}\right\rfloor.
 \label{eq:unit-p-sequence}
\end{equation}
\end{lemma}

\begin{proof}
For $x\in\Zp^{\times}$ distinct from $u_0,\dots,u_{k-1}$, one has
\[
 \vp\!\left(\prod_{i<k}(x-u_i)\right)
 =
 \sum_{r\ge 1}
 \#\{i<k:u_i\equiv x\pmod{p^r}\}.
\]
Among the first $k$ positive $p$-units, each reduced residue class modulo $p^r$ occurs either
\[
 \left\lfloor\frac{k}{\varphi(p^r)}\right\rfloor
 \quad\text{or}\quad
 \left\lceil\frac{k}{\varphi(p^r)}\right\rceil
\]
times.  Hence every candidate $x$ gives the lower bound
\[
 \vp\!\left(\prod_{i<k}(x-u_i)\right)
 \ge
 \sum_{r\ge 1}
 \left\lfloor\frac{k}{\varphi(p^r)}\right\rfloor.
\]
For $x=u_k$, the residue class of $u_k$ is precisely the next class in the periodic ordered list of reduced residues modulo every $p^r$.  It has occurred exactly
\[
 \left\lfloor\frac{k}{\varphi(p^r)}\right\rfloor
\]
times among $u_0,\dots,u_{k-1}$.  Equality therefore holds simultaneously for all $r$, proving both assertions.
\end{proof}

The isolated point $p$ accounts for the shift appearing in the prime factorial.

\begin{proposition}[Local prime factorial exponent]
\label{prop:local-prime-exponent}
For every rational prime $p$,
\[
 \nu_0(\PP,p)=0,
\]
and for $k\ge 1$,
\begin{equation}
 \nu_k(\PP,p)
 =
 \sum_{r\ge 1}
 \left\lfloor\frac{k-1}{\varphi(p^r)}\right\rfloor
 =
 \sum_{j\ge 0}
 \left\lfloor
 \frac{k-1}{(p-1)p^j}
 \right\rfloor.
 \label{eq:local-prime-factorial-exponent}
\end{equation}
\end{proposition}

\begin{proof}
By Lemma~\ref{lem:closure-invariance} and \eqref{eq:p-adic-closure-primes}, it is enough to order
\[
 \{p\}\sqcup\Zp^{\times}.
\]
Choose $p$ first and then use the unit ordering of Lemma~\ref{lem:unit-p-ordering}.  For every unit $u$, one has $\vp(u-p)=0$.  The isolated initial point therefore contributes no valuation to any later stage, and the $k$th minimum for the prime set is the $(k-1)$st minimum for the unit group.  Formula \eqref{eq:local-prime-factorial-exponent} follows from
\[
 \varphi(p^r)=(p-1)p^{r-1}.
\]
\end{proof}

\subsection{The prime factorial ideal and its norm}
\label{subsec:prime-factorial-norm}

Combining Theorem~\ref{thm:local-factorial-formula} with Proposition~\ref{prop:local-prime-exponent} gives the explicit prime factorial.

\begin{theorem}[Prime Bhargava factorial]
\label{thm:prime-bhargava-factorial}
For every $n\ge 0$, the $(n+1)$st factorial ideal of the prime set is
\begin{equation}
 \mathfrak A_n
 :=
 \mathfrak F_{n+1}(\PP)
 =
 \prod_p
 (p)^{\displaystyle
 \sum_{j\ge 0}
 \left\lfloor\frac{n}{(p-1)p^j}\right\rfloor}.
 \label{eq:prime-factorial-ideal}
\end{equation}
It is principal, and its positive generator is
\begin{equation}
 A_n
 :=
 (n+1)!_{\PP}
 =
 \prod_p
 p^{\displaystyle
 \sum_{j\ge 0}
 \left\lfloor\frac{n}{(p-1)p^j}\right\rfloor}.
 \label{eq:prime-factorial-positive-generator}
\end{equation}
For fixed $n$, both products are finite.
\end{theorem}

\begin{proof}
Only primes with $p-1\le n$ can contribute a nonzero exponent.  For each such $p$, only those $j$ satisfying $(p-1)p^j\le n$ contribute.  The claimed formulas now follow directly from \eqref{eq:local-prime-factorial-exponent}.
\end{proof}

The terminology ``norm realization'' is literal.  For a nonzero ideal $m\ZZ\subseteq\ZZ$, its absolute norm is
\[
 \Norm(m\ZZ)=|\ZZ/m\ZZ|=m.
\]
Thus Theorem~\ref{thm:prime-bhargava-factorial} gives
\begin{equation}
 \boxed{
 A_n
 =
 \Norm(\mathfrak A_n)
 =
 (n+1)!_{\PP}.
 }
 \label{eq:norm-realization}
\end{equation}
The analytic factorial datum used in the remainder of the paper is therefore not introduced independently of the arithmetic factorial ideal: it is its absolute norm.

\begin{remark}[Indexing]
The conventional factorial ideal satisfies $0!_{\PP}=1!_{\PP}=1$.  The first nontrivial value is
\[
 2!_{\PP}=2.
\]
The shifted datum $A_n=(n+1)!_{\PP}$ has $A_0=1$ and places every local exponent in the uniform form $\lfloor n/((p-1)p^j)\rfloor$.
\end{remark}

\subsection{Atomic valuation layers}
\label{subsec:atomic-layers}

Let
\[
 \cA
 =
 \{(p,j):p\in\PP,\ j\in\ZZ_{\ge 0}\}
\]
be the set of valuation layers.  For $\lambda=(p,j)$, define
\[
 d_\lambda=(p-1)p^j,
 \qquad
 w_\lambda=\log p.
\]
The corresponding atomic ideal and its norm are
\begin{equation}
 \mathfrak A_n^{(\lambda)}
 =
 (p)^{\lfloor n/d_\lambda\rfloor},
 \qquad
 F_n^{(\lambda)}
 =
 \Norm\!\left(\mathfrak A_n^{(\lambda)}\right)
 =
 p^{\lfloor n/d_\lambda\rfloor}.
 \label{eq:atomic-layer}
\end{equation}

For a finite subset $\Lambda\subset\cA$, put
\begin{equation}
 \mathfrak A_n^{\Lambda}
 =
 \prod_{\lambda\in\Lambda}
 \mathfrak A_n^{(\lambda)},
 \qquad
 A_n^{\Lambda}
 =
 \Norm(\mathfrak A_n^{\Lambda})
 =
 \prod_{\lambda\in\Lambda}
 F_n^{(\lambda)}.
 \label{eq:finite-layer-data}
\end{equation}
The norm turns ideal multiplication into pointwise multiplication of factorial data.  This is the arithmetic source of the monoidal operation used in the next section.

\begin{proposition}[Filtered layer factorization]
\label{prop:filtered-layer-factorization}
The finite-layer ideals form a filtered system under inclusion of layer sets.  For every fixed $n$,
\[
 \mathfrak A_n
 =
 \prod_{\lambda\in\cA}
 \mathfrak A_n^{(\lambda)},
 \qquad
 A_n
 =
 \prod_{\lambda\in\cA}
 F_n^{(\lambda)},
\]
where only finitely many factors differ from the unit ideal or from $1$.  Equivalently,
\[
 A_n=\sup_{\Lambda\Subset\cA}A_n^{\Lambda}.
\]
Moreover,
\begin{equation}
 \log A_n
 =
 \sum_{\lambda\in\cA}
 w_\lambda
 \left\lfloor\frac{n}{d_\lambda}\right\rfloor.
 \label{eq:log-layer-factorization}
\end{equation}
\end{proposition}

\begin{proof}
The formulas are restatements of \eqref{eq:prime-factorial-ideal} and \eqref{eq:prime-factorial-positive-generator}.  Local finiteness follows from $d_\lambda\le n$ for every nontrivial factor.  If $\Lambda\subseteq\Lambda'$, then $A_n^\Lambda$ divides $A_n^{\Lambda'}$ and hence $A_n^\Lambda\le A_n^{\Lambda'}$, which proves the supremum statement.
\end{proof}

Every finite layer system has at most exponential growth:
\[
 \log A_n^\Lambda
 =
 n\sum_{\lambda\in\Lambda}\frac{w_\lambda}{d_\lambda}
 +O_\Lambda(1).
\]
The completed object has factorial-scale growth.  Thus completion changes the asymptotic category of the sequence, a phenomenon that will be central to the later Stirling theory.

\subsection{The exact factorial increment}
\label{subsec:factorial-increment}

The quotient of consecutive norms supplies the weight of the generalized derivative associated with the prime factorial calculus.  Define
\[
 \rho_{\PP}(n)=\frac{A_n}{A_{n-1}},
 \qquad n\ge 1.
\]

\begin{corollary}[Prime factorial recurrence weight]
\label{cor:prime-factorial-increment}
For every $n\ge 1$,
\begin{align}
 \rho_{\PP}(n)
 &=
 \prod_{\substack{p\in\PP,\ j\ge 0\\(p-1)p^j\mid n}}p
 \label{eq:increment-layer-form}\\
 &=
 \prod_{p-1\mid n}p^{\vp(n)+1}.
 \label{eq:increment-prime-form}
\end{align}
Consequently,
\begin{equation}
 \log\rho_{\PP}(n)
 =
 \sum_{p-1\mid n}(\vp(n)+1)\log p.
 \label{eq:log-increment}
\end{equation}
\end{corollary}

\begin{proof}
For every positive integer $d$,
\[
 \left\lfloor\frac{n}{d}\right\rfloor
 -
 \left\lfloor\frac{n-1}{d}\right\rfloor
 =
 \mathbf 1_{d\mid n}.
\]
Applying this identity to each layer in \eqref{eq:prime-factorial-positive-generator} gives \eqref{eq:increment-layer-form}.  If $p-1\mid n$, then $p\nmid p-1$, and the number of integers $j\ge 0$ satisfying $(p-1)p^j\mid n$ is $\vp(n)+1$.  If $p-1\nmid n$, there is no such layer.  This proves \eqref{eq:increment-prime-form} and \eqref{eq:log-increment}.
\end{proof}

The formula shows that the prime factorial is governed by a locally finite divisor incidence structure.  Each layer $(p,j)$ fires exactly at multiples of $(p-1)p^j$.  In the language developed in Section~3, the atomic data \eqref{eq:atomic-layer} will be monoidal generators, while \eqref{eq:increment-layer-form} will be the diagonal weight of the completed generalized derivative.

The arithmetic construction has therefore produced three compatible levels of structure:
\[
 \boxed{
 \begin{gathered}
 \text{factorial ideals }
 \mathfrak A_n^{\Lambda}
 \quad\xrightarrow{\ \Norm\ }\quad
 \text{positive factorial data }
 A_n^{\Lambda},
 \\
 \Lambda\Subset\cA
 \quad\longrightarrow\quad
 \mathfrak A_n
 \quad\xrightarrow{\ \Norm\ }\quad
 A_n=(n+1)!_{\PP},
 \\
 \frac{A_n}{A_{n-1}}
 =
 \prod_{d_\lambda\mid n}e^{w_\lambda}.
 \end{gathered}
 }
\]

% ===== Source: section3_analytic_category_factorial_calculi.tex =====
\section{The analytic category of factorial calculi}
\label{sec:analytic-factorial-category}

The arithmetic construction of Section~2 produces a positive sequence
\[
 A_n=(n+1)!_{\PP},\qquad A_0=1,
\]
together with its atomic layer factorizations. We now explain how any such positive sequence determines a Hilbert-space calculus and how comparison of factorial data produces canonical intertwiners. The resulting category, formulated with the standard conventions of~\cite{MacLane1998}, is deliberately more rigid than the algebraic category obtained by working only with polynomials: on the polynomial level every two positive factorial data are diagonally conjugate, whereas after Hilbert completion the conjugacy is bounded only under a genuine growth condition.

This distinction is essential. It is the Hilbert geometry, rather than the formal recurrence alone, that remembers the asymptotic size of a factorial datum.

\subsection{Admissible factorial data and their Hilbert realizations}
\label{subsec:admissible-factorial-data}

\begin{definition}[Factorial datum]
\label{def:factorial-datum}
A \emph{factorial datum} is a sequence
\[
 F=(F_n)_{n\ge 0},\qquad F_0=1,\qquad F_n>0.
\]
Its recurrence weights are
\[
 \rho_F(n)=\frac{F_n}{F_{n-1}},\qquad n\ge 1.
\]
The datum is called \emph{analytic} if
\begin{equation}
 R_F
 :=
 \left(\liminf_{n\to\infty}F_n^{1/n}\right)^{1/2}>0.
 \label{eq:analytic-radius}
\end{equation}
We allow $R_F=+\infty$.
\end{definition}

To an analytic factorial datum we associate the coefficient Hilbert space
\begin{equation}
 \cH_F
 =
 \left\{
 f(z)=\sum_{n\ge 0}a_nz^n:
 \|f\|_F^2:=\sum_{n\ge 0}|a_n|^2F_n<\infty
 \right\}.
 \label{eq:factorial-hilbert-space}
\end{equation}
The power series in \eqref{eq:factorial-hilbert-space} converge on the disk
\[
 \mathbb D_F=\{z\in\CC:|z|<R_F\}.
\]
Indeed, Cauchy--Schwarz gives
\[
 |f(z)|
 \le
 \|f\|_F
 \left(\sum_{n\ge 0}\frac{|z|^{2n}}{F_n}\right)^{1/2}.
\]

\begin{proposition}[Reproducing-kernel realization]
\label{prop:rkhs-realization}
The space $\cH_F$ is a reproducing-kernel Hilbert space on $\mathbb D_F$ with kernel
\begin{equation}
 K_F(z,w)
 =
 \sum_{n\ge 0}\frac{(z\overline w)^n}{F_n}.
 \label{eq:factorial-kernel}
\end{equation}
The normalized monomials
\[
 e_n^F(z)=\frac{z^n}{\sqrt{F_n}},\qquad n\ge 0,
\]
form an orthonormal basis.
\end{proposition}

\begin{proof}
The monomial statement follows directly from the definition of the norm. For $w\in\mathbb D_F$, the coefficient sequence of $K_F(\,\cdot\,,w)$ is $(\overline w^{\,n}/F_n)_{n\ge 0}$, and
\[
 \|K_F(\,\cdot\,,w)\|_F^2
 =
 \sum_{n\ge 0}\frac{|w|^{2n}}{F_n}<\infty.
\]
Hence
\[
 \langle f,K_F(\,\cdot\,,w)\rangle_F
 =
 \sum_{n\ge 0}a_nw^n=f(w).
\]
\end{proof}

The factorial datum determines two mutually adjoint weighted shifts. On polynomials define
\begin{equation}
 D_Fz^n=\rho_F(n)z^{n-1}\quad(n\ge 1),
 \qquad
 D_F1=0,
 \label{eq:generalized-derivative}
\end{equation}
and let $M_F$ denote multiplication by $z$.

\begin{proposition}[Creation and annihilation operators]
\label{prop:creation-annihilation}
The operators $D_F$ and $M_F$ extend to closed densely defined operators on $\cH_F$ with domains
\begin{align}
 \Dom(D_F)
 &=
 \left\{
 \sum_{n\ge 0}a_nz^n:
 \sum_{n\ge 1}|a_n|^2\frac{F_n^2}{F_{n-1}}<\infty
 \right\},
 \label{eq:domain-D}\\
 \Dom(M_F)
 &=
 \left\{
 \sum_{n\ge 0}a_nz^n:
 \sum_{n\ge 0}|a_n|^2F_{n+1}<\infty
 \right\}.
 \label{eq:domain-M}
\end{align}
They satisfy
\begin{equation}
 D_F=M_F^*.
 \label{eq:adjoint-shifts}
\end{equation}
In the orthonormal monomial basis,
\begin{equation}
 D_Fe_n^F=\sqrt{\rho_F(n)}\,e_{n-1}^F,
 \qquad
 M_Fe_n^F=\sqrt{\rho_F(n+1)}\,e_{n+1}^F.
 \label{eq:weighted-shift-form}
\end{equation}
On polynomials their commutator is diagonal:
\begin{equation}
 [D_F,M_F]e_n^F
 =
 \bigl(\rho_F(n+1)-\rho_F(n)\bigr)e_n^F,
 \label{eq:factorial-commutator}
\end{equation}
where $\rho_F(0):=0$.
\end{proposition}

\begin{proof}
The domain formulas and closedness are the standard coordinate description of unilateral weighted shifts. Formula \eqref{eq:weighted-shift-form} follows from $F_n=\rho_F(n)F_{n-1}$. The two displayed shifts have transposed weights, proving \eqref{eq:adjoint-shifts}. The commutator formula follows by applying $D_FM_F-M_FD_F$ to $e_n^F$.
\end{proof}

The formal generalized exponential attached to $F$ is
\begin{equation}
 E_F(z)=\sum_{n\ge 0}\frac{z^n}{F_n}.
 \label{eq:generalized-exponential}
\end{equation}
On its disk of convergence it satisfies
\begin{equation}
 D_FE_F=E_F.
 \label{eq:generalized-exponential-eigenfunction}
\end{equation}
More generally, the kernel vectors are eigenvectors of the annihilation operator:
\begin{equation}
 D_FK_F(\,\cdot\,,w)=\overline w\,K_F(\,\cdot\,,w).
 \label{eq:kernel-eigenvector}
\end{equation}

\begin{example}[Classical factorial calculus]
For $F_n=n!$, the space $\cH_F$ is the standard Bargmann--Fock space in coefficient normalization,
\[
 K_F(z,w)=e^{z\overline w},
 \qquad
 D_F=\frac{d}{dz},
\]
and \eqref{eq:factorial-commutator} becomes $[D_F,M_F]=\Id$.
\end{example}

\begin{example}[The monoidal unit]
For the constant datum $\mathbf 1_n=1$, the space $\cH_{\mathbf 1}$ is the Hardy coefficient space $H^2(\mathbb D)$, $M_{\mathbf 1}$ is the unilateral shift, and $D_{\mathbf 1}$ is the backward shift.
\end{example}

\subsection{Canonical diagonal intertwiners}
\label{subsec:canonical-intertwiners}

On the polynomial algebra, every two factorial derivatives are conjugate by a diagonal map. The Hilbert-space question is whether this map is bounded.

\begin{definition}[Analytic domination]
\label{def:analytic-domination}
For factorial data $F$ and $G$, write
\begin{equation}
 F\preccurlyeq G
 \quad\Longleftrightarrow\quad
 C(F,G):=\sup_{n\ge 0}\frac{F_n}{G_n}<\infty.
 \label{eq:factorial-domination}
\end{equation}
When this holds, define the canonical diagonal map
\begin{equation}
 J_{F,G}z^n=\frac{F_n}{G_n}z^n.
 \label{eq:canonical-intertwiner}
\end{equation}
\end{definition}

The relation $\preccurlyeq$ is a preorder. It records continuous comparison of the coefficient norms and is strictly finer than formal equivalence of the underlying polynomial calculi.

\begin{theorem}[Bounded intertwiner criterion]
\label{thm:bounded-intertwiner-criterion}
Let $F$ and $G$ be analytic factorial data. The following are equivalent.
\begin{enumerate}[label=\textup{(\roman*)}]
 \item $F\preccurlyeq G$.
 \item The map $J_{F,G}$ extends boundedly from polynomials to $\cH_F\to\cH_G$.
 \item There exists a bounded diagonal operator $T:\cH_F\to\cH_G$ such that
 \[
  T1=1,
  \qquad
  D_GT=TD_F
 \]
 on polynomials.
\end{enumerate}
When these conditions hold, $T=J_{F,G}$ and
\begin{equation}
 \|J_{F,G}\|^2=C(F,G).
 \label{eq:intertwiner-norm}
\end{equation}
\end{theorem}

\begin{proof}
For a polynomial $f(z)=\sum_na_nz^n$,
\[
 \|J_{F,G}f\|_G^2
 =
 \sum_n|a_n|^2\frac{F_n^2}{G_n}
 \le
 C(F,G)\sum_n|a_n|^2F_n.
\]
Testing on a monomial shows that the best constant is exactly $C(F,G)$, proving the equivalence of (i) and (ii) and formula \eqref{eq:intertwiner-norm}.

Now let $Tz^n=t_nz^n$. The normalization gives $t_0=1$. The relation $D_GT=TD_F$ yields
\[
 t_n\frac{G_n}{G_{n-1}}
 =
 \frac{F_n}{F_{n-1}}t_{n-1}.
\]
Induction gives $t_n=F_n/G_n$, so $T=J_{F,G}$. This proves the remaining assertion.
\end{proof}

The adjoint of the canonical intertwiner has an especially simple form.

\begin{proposition}[Projective--inductive adjunction]
\label{prop:projective-inductive-adjunction}
If $F\preccurlyeq G$, then the coefficientwise identity map
\begin{equation}
 \iota_{G,F}:\cH_G\hookrightarrow\cH_F,
 \qquad
 \iota_{G,F}f=f,
 \label{eq:canonical-inclusion}
\end{equation}
is bounded and
\begin{equation}
 J_{F,G}^*=\iota_{G,F}.
 \label{eq:intertwiner-adjoint}
\end{equation}
Moreover,
\begin{equation}
 J_{G,H}J_{F,G}=J_{F,H},
 \qquad
 \iota_{G,F}\iota_{H,G}=\iota_{H,F}
 \label{eq:composition-laws}
\end{equation}
whenever $F\preccurlyeq G\preccurlyeq H$.
\end{proposition}

\begin{proof}
The boundedness of $\iota_{G,F}$ follows from
\[
 \|f\|_F^2\le C(F,G)\|f\|_G^2.
\]
For monomials,
\[
 \langle J_{F,G}z^n,z^m\rangle_G
 =
 \delta_{nm}F_n
 =
 \langle z^n,\iota_{G,F}z^m\rangle_F,
\]
which proves \eqref{eq:intertwiner-adjoint}. The composition identities are immediate from the diagonal formulas.
\end{proof}

\begin{definition}[The category of factorial calculi]
\label{def:factorial-calculus-category}
The category $\FactCalc$ is the thin category whose objects are analytic factorial data, regarded together with their realizations $(\cH_F,D_F)$, and for which
\[
 \operatorname{Hom}_{\FactCalc}(F,G)
 =
 \begin{cases}
  \{J_{F,G}\},&F\preccurlyeq G,\\
  \varnothing,&\text{otherwise}.
 \end{cases}
\]
\end{definition}

Thus the assignment
\[
 F\longmapsto(\cH_F,D_F),
 \qquad
 (F\preccurlyeq G)\longmapsto J_{F,G},
\]
is a faithful realization of the factorial preorder by Hilbert spaces and bounded intertwiners. Taking adjoints reverses every arrow and produces the coefficient-space diagram
\[
 \cH_G\xhookrightarrow{\ J_{F,G}^*\ }\cH_F.
\]

\begin{corollary}[Isomorphism classes]
\label{cor:isomorphism-classes}
Two objects $F$ and $G$ are isomorphic in $\FactCalc$ if and only if their weights are uniformly equivalent:
\begin{equation}
 0<\inf_n\frac{F_n}{G_n}
 \le
 \sup_n\frac{F_n}{G_n}<\infty.
 \label{eq:uniform-equivalence}
\end{equation}
The skeletal quotient of $\FactCalc$ is therefore the partially ordered set of uniform-equivalence classes of analytic factorial data.
\end{corollary}

\begin{remark}[Why completion is necessary]
On $\CC[z]$, the map $J_{F,G}$ is always an invertible diagonal intertwiner, with inverse $J_{G,F}$. If one defined morphisms only algebraically, all positive factorial data would become canonically isomorphic and their growth would disappear from the category. Boundedness after Hilbert completion is what makes the category nontrivial.
\end{remark}

\subsection{The symmetric-monoidal structure}
\label{subsec:monoidal-structure}

The multiplication of factorial ideals becomes pointwise multiplication after taking norms. This motivates the following operation.

\begin{definition}[Monoidal product]
\label{def:monoidal-product}
For factorial data $F$ and $G$, define
\begin{equation}
 (F\boxt G)_n=F_nG_n.
 \label{eq:monoidal-product}
\end{equation}
The monoidal unit is the constant datum $\mathbf 1_n=1$.
\end{definition}

If $F$ and $G$ are analytic, then so is $F\boxt G$. At the level of recurrence weights,
\begin{equation}
 \rho_{F\boxt G}(n)=\rho_F(n)\rho_G(n).
 \label{eq:monoidal-recurrence}
\end{equation}
The generalized exponentials and reproducing kernels multiply coefficientwise:
\begin{align}
 E_{F\boxt G}
 &=E_F\odot E_G,
 \label{eq:hadamard-exponential}\\
 K_{F\boxt G}(z,w)
 &=K_F(z,w)\odot K_G(z,w),
 \label{eq:hadamard-kernel}
\end{align}
where $\odot$ denotes the Hadamard product in the monomial variable.

\begin{theorem}[Symmetric-monoidal factorial category]
\label{thm:symmetric-monoidal-category}
The operation $\boxt$ makes $\FactCalc$ a strict symmetric-monoidal thin category. More precisely, if
\[
 F\preccurlyeq F',
 \qquad
 G\preccurlyeq G',
\]
then
\[
 F\boxt G\preccurlyeq F'\boxt G'
\]
and
\begin{equation}
 J_{F\boxt G,F'\boxt G'}z^n
 =
 \left(\frac{F_n}{F'_n}\right)
 \left(\frac{G_n}{G'_n}\right)z^n.
 \label{eq:monoidal-intertwiner}
\end{equation}
Furthermore,
\begin{equation}
 \|J_{F\boxt G,F'\boxt G'}\|
 \le
 \|J_{F,F'}\|\,\|J_{G,G'}\|.
 \label{eq:monoidal-norm-bound}
\end{equation}
\end{theorem}

\begin{proof}
The ratio of the product data factors as
\[
 \frac{F_nG_n}{F'_nG'_n}
 =
 \frac{F_n}{F'_n}\frac{G_n}{G'_n}.
\]
Taking suprema proves domination and \eqref{eq:monoidal-norm-bound}; the diagonal formula gives \eqref{eq:monoidal-intertwiner}. Associativity, commutativity, and the unit law hold identically at the sequence level.
\end{proof}

The operation $\boxt$ is not the usual Hilbert tensor product. It is a diagonal or Hadamard tensor operation: it combines the $n$th factorial weights at the same degree. This is exactly the operation required by valuation-layer factorization.

\subsection{The filtered prime-layer diagram}
\label{subsec:prime-layer-diagram}

Let
\[
 \cA=\{(p,j):p\in\PP,\ j\ge 0\}
\]
be the prime-layer set, and for $\lambda=(p,j)$ write
\[
 d_\lambda=(p-1)p^j,
 \qquad
 w_\lambda=\log p.
\]
The atomic factorial datum is
\begin{equation}
 F_n^{(\lambda)}
 =
 \exp\!\left(
 w_\lambda\left\lfloor\frac{n}{d_\lambda}\right\rfloor
 \right)
 =
 p^{\lfloor n/d_\lambda\rfloor}.
 \label{eq:atomic-factorial-datum-section3}
\end{equation}
Its recurrence weight is the periodic step function
\begin{equation}
 \rho_\lambda(n)
 =
 \begin{cases}
  p,&d_\lambda\mid n,\\
  1,&d_\lambda\nmid n.
 \end{cases}
 \label{eq:atomic-recurrence-weight}
\end{equation}

For a finite layer set $\Lambda\Subset\cA$, define
\begin{equation}
 A^\Lambda
 =
 \boxt_{\lambda\in\Lambda}F^{(\lambda)}.
 \label{eq:finite-layer-monoidal-product}
\end{equation}
Then
\begin{align}
 A_n^\Lambda
 &=
 \prod_{\lambda\in\Lambda}
 e^{w_\lambda\lfloor n/d_\lambda\rfloor},
 \label{eq:finite-layer-weight}\\
 \rho_\Lambda(n)
 &=
 \prod_{\substack{\lambda\in\Lambda\\d_\lambda\mid n}}
 e^{w_\lambda}.
 \label{eq:finite-layer-recurrence}
\end{align}

Let $\Fin(\cA)$ denote the directed poset of finite subsets of $\cA$. If $\Lambda\subseteq\Lambda'$, then
\[
 A_n^\Lambda\le A_n^{\Lambda'}
\]
for every $n$. Hence $A^\Lambda\preccurlyeq A^{\Lambda'}$, and the transition map
\begin{equation}
 J_{\Lambda,\Lambda'}
 :=J_{A^\Lambda,A^{\Lambda'}}
 :\cH_{A^\Lambda}\longrightarrow\cH_{A^{\Lambda'}}
 \label{eq:layer-transition-map}
\end{equation}
is a contraction. Since both data equal $1$ at $n=0$, its norm is exactly $1$.

\begin{proposition}[Prime-layer functors]
\label{prop:prime-layer-functors}
The assignments
\begin{align}
 \Lambda
 &\longmapsto
 (\cH_{A^\Lambda},D_{A^\Lambda}),
 &
 (\Lambda\subseteq\Lambda')
 &\longmapsto J_{\Lambda,\Lambda'},
 \label{eq:covariant-layer-functor}\\
 \Lambda
 &\longmapsto
 \cH_{A^\Lambda},
 &
 (\Lambda\subseteq\Lambda')
 &\longmapsto J_{\Lambda,\Lambda'}^*
 \label{eq:contravariant-layer-functor}
\end{align}
define respectively a covariant diagram of generalized differential calculi and a contravariant diagram of coefficient Hilbert spaces. The two diagrams are adjoint degree by degree.
\end{proposition}

The first diagram remembers how the generalized derivatives change when valuation layers are added. The second records the reverse nesting
\begin{equation}
 \Lambda\subseteq\Lambda'
 \quad\Longrightarrow\quad
 \cH_{A^{\Lambda'}}\hookrightarrow\cH_{A^\Lambda}.
 \label{eq:reverse-hilbert-nesting}
\end{equation}
Thus adding arithmetic layers enlarges the factorial weights, strengthens coefficient decay, and produces a smaller function space.

\subsection{Bounded projective completion}
\label{subsec:bounded-projective-completion}

Define the completed prime datum by
\begin{equation}
 A_n
 =
 \sup_{\Lambda\Subset\cA}A_n^\Lambda
 =
 (n+1)!_{\PP}.
 \label{eq:completed-prime-datum}
\end{equation}
For each $n$, only finitely many layers contribute, so the supremum is attained once $\Lambda$ contains all layers with $d_\lambda\le n$.

The ordinary locally convex projective limit of the spaces $\cH_{A^\Lambda}$ is their intersection equipped with the family of seminorms $\|\cdot\|_{A^\Lambda}$. The completed factorial datum corresponds to its bounded part.

\begin{definition}[Bounded projective limit]
\label{def:bounded-projective-limit}
For the contravariant layer diagram, define
\begin{equation}
 \varprojlim{}^{\!b}_{\Lambda}\cH_{A^\Lambda}
 =
 \left\{
 f\in\bigcap_{\Lambda\Subset\cA}\cH_{A^\Lambda}:
 \sup_{\Lambda}\|f\|_{A^\Lambda}<\infty
 \right\},
 \label{eq:bounded-projective-limit}
\end{equation}
with norm
\[
 \|f\|_{\mathrm{proj}}
 =
 \sup_{\Lambda}\|f\|_{A^\Lambda}.
\]
\end{definition}

\begin{theorem}[Prime completion theorem]
\label{thm:prime-completion}
Coefficientwise identification induces an isometric isomorphism
\begin{equation}
 \boxed{
 \cH_A
 \cong
 \varprojlim{}^{\!b}_{\Lambda\Subset\cA}
 \cH_{A^\Lambda}.
 }
 \label{eq:prime-projective-limit}
\end{equation}
More precisely, for every formal series $f(z)=\sum_{n\ge 0}a_nz^n$,
\begin{equation}
 \|f\|_A^2
 =
 \sup_{\Lambda\Subset\cA}
 \|f\|_{A^\Lambda}^2.
 \label{eq:projective-norm-identity}
\end{equation}
The transition maps to the completed calculus,
\[
 J_{\Lambda,A}:\cH_{A^\Lambda}\longrightarrow\cH_A,
 \qquad
 J_{\Lambda,A}z^n=\frac{A_n^\Lambda}{A_n}z^n,
\]
are contractive intertwiners, and their adjoints are the canonical inclusions
\[
 \cH_A\hookrightarrow\cH_{A^\Lambda}.
\]
\end{theorem}

\begin{proof}
For a fixed coefficient sequence, the quantities
\[
 \|f\|_{A^\Lambda}^2
 =
 \sum_{n\ge 0}|a_n|^2A_n^\Lambda
\]
form an increasing net as $\Lambda$ grows. Since $A_n^\Lambda\uparrow A_n$ for every $n$, monotone convergence gives
\[
 \sup_\Lambda\sum_{n\ge 0}|a_n|^2A_n^\Lambda
 =
 \sum_{n\ge 0}|a_n|^2A_n.
\]
This proves \eqref{eq:projective-norm-identity} and the isometric identification. The assertions about $J_{\Lambda,A}$ and its adjoint follow from Theorem~\ref{thm:bounded-intertwiner-criterion} and Proposition~\ref{prop:projective-inductive-adjunction}.
\end{proof}

The completed recurrence weight is the locally finite product
\begin{equation}
 \rho_A(n)
 =
 \frac{A_n}{A_{n-1}}
 =
 \prod_{\substack{\lambda\in\cA\\d_\lambda\mid n}}
 e^{w_\lambda}
 =
 \prod_{p-1\mid n}p^{v_p(n)+1}.
 \label{eq:completed-prime-recurrence}
\end{equation}
Accordingly,
\begin{equation}
 D_Az^n
 =
 \left(\prod_{p-1\mid n}p^{v_p(n)+1}\right)z^{n-1}.
 \label{eq:completed-prime-derivative}
\end{equation}
Every finite-layer derivative maps to $D_A$ through the canonical intertwiner:
\begin{equation}
 D_AJ_{\Lambda,A}=J_{\Lambda,A}D_{A^\Lambda}
 \qquad\text{on }\CC[z].
 \label{eq:finite-to-complete-intertwining}
\end{equation}

\begin{remark}[Completion changes the growth regime]
For fixed finite $\Lambda$,
\[
 \log A_n^\Lambda
 =
 n\sum_{\lambda\in\Lambda}\frac{w_\lambda}{d_\lambda}
 +O_\Lambda(1),
\]
so $A^\Lambda$ has exponential growth and $\cH_{A^\Lambda}$ is naturally realized on a disk of finite radius. The completed datum has factorial-scale growth and, after the weak Stirling estimate is established, $R_A=+\infty$. Thus the passage from the filtered diagram to its bounded projective completion changes not only the weights but the analytic type of the realization.
\end{remark}

\subsection{Categorical interpretation}
\label{subsec:categorical-interpretation}

The constructions of this section may be summarized in the diagram
\begin{equation}
 \boxed{
 \begin{gathered}
 \text{factorial data}
 \quad
 A^\Lambda\preccurlyeq A^{\Lambda'}\preccurlyeq A,
 \\
 \text{differential calculi}
 \quad
 \cH_{A^\Lambda}
 \xrightarrow{\ J_{\Lambda,\Lambda'}\ }
 \cH_{A^{\Lambda'}}
 \xrightarrow{\ J_{\Lambda',A}\ }
 \cH_A,
 \\
 \text{coefficient spaces}
 \quad
 \cH_A
 \hookrightarrow
 \cH_{A^{\Lambda'}}
 \hookrightarrow
 \cH_{A^\Lambda}.
 \end{gathered}
 }
 \label{eq:projective-inductive-summary}
\end{equation}
The upper arrows are canonical normalized intertwiners of generalized derivatives. The lower arrows are their adjoints. The monoidal product records multiplication of valuation layers, while the bounded projective limit records their completion.

This framework separates three notions that would otherwise be conflated:
\begin{enumerate}[label=\textup{(\roman*)}]
 \item the formal polynomial calculus, in which all positive factorial data are diagonally conjugate;
 \item the analytic factorial category, in which boundedness detects relative growth;
 \item the filtered completion, in which infinitely many exponential layers combine to produce the prime factorial calculus.
\end{enumerate}

The category $\FactCalc$ therefore supplies the natural domain for the gamma-lift construction developed in the next part of the paper. The recurrence weights determine the discrete factorial calculus, but they do not by themselves choose a unique analytic primitive of the associated logarithmic recurrence. That remaining ambiguity is encoded by a torsor, and it is Stirling rigidity that will select its distinguished point.

% ===== Source: section4_symmetric_monoidal_layer_decompositions.tex =====
\section{Symmetric-monoidal layer decompositions}
\label{sec:symmetric-monoidal-layers}

The preceding section introduced pointwise multiplication of factorial data as the monoidal operation
\[
 (F\boxt G)_n=F_nG_n.
\]
For the prime Bhargava factorial this operation is not auxiliary: it is the exact mechanism by which the local valuation layers assemble. We now isolate the free commutative monoid generated by those layers, linearize its monoidal law through logarithmic recurrence increments, and prove a unique-factorization statement for the prime-layer realization.

Two consequences will be used repeatedly. First, every finite layer system has an exact affine-periodic normal form and therefore purely exponential growth. Second, the natural filtration by layer conductor is degreewise stationary: once all layers of conductor at most $N$ have been inserted, the resulting factorial norm and generalized derivative agree exactly with the completed prime calculus through degree $N$.

\subsection{The free commutative monoid of weighted layers}
\label{subsec:free-layer-monoid}

Let
\begin{equation}
 \cA=\{(p,j):p\in\PP,\ j\ge 0\}
 \label{eq:prime-layer-index-set-section4}
\end{equation}
be the prime-layer index set. For $\lambda=(p,j)$ put
\begin{equation}
 d_\lambda=(p-1)p^j,
 \qquad
 w_\lambda=\log p.
 \label{eq:prime-layer-conductor-weight}
\end{equation}
The associated atomic factorial datum is
\begin{equation}
 E_n^\lambda
 =
 \exp\!\left(w_\lambda\left\lfloor\frac{n}{d_\lambda}\right\rfloor\right)
 =
 p^{\lfloor n/((p-1)p^j)\rfloor}.
 \label{eq:atomic-layer-datum-section4}
\end{equation}

It is useful to allow multiplicities rather than only finite subsets. Let
\begin{equation}
 \NN^{(\cA)}
 =
 \left\{
 \mathbf m=(m_\lambda)_{\lambda\in\cA}:
 m_\lambda\in\NN\cup\{0\},\ \supp(\mathbf m)\text{ finite}
 \right\}.
 \label{eq:finite-multiplicity-monoid}
\end{equation}
This is the free commutative monoid on $\cA$, with addition defined coordinatewise.

\begin{definition}[Layer realization]
\label{def:layer-realization}
For $\mathbf m\in\NN^{(\cA)}$, define
\begin{equation}
 A^{\mathbf m}
 =
 \mathop{\boxt}_{\lambda\in\cA}
 (E^\lambda)^{\boxt m_\lambda}.
 \label{eq:layer-realization}
\end{equation}
Equivalently,
\begin{equation}
 A_n^{\mathbf m}
 =
 \prod_{(p,j)\in\cA}
 p^{m_{p,j}\lfloor n/((p-1)p^j)\rfloor}.
 \label{eq:layer-realization-explicit}
\end{equation}
The product is finite because $\mathbf m$ has finite support.
\end{definition}

The passage from finite subsets to multiplicity vectors removes a minor defect in set notation. If two finite layer sets overlap, their monoidal product counts the common layers twice, and is therefore represented by the sum of their indicator vectors rather than by their union.

Let $\LayP$ denote the thin category whose objects are $\NN^{(\cA)}$ and in which there is a unique arrow
\[
 \mathbf m\longrightarrow\mathbf m'
\]
when $m_\lambda\le m'_\lambda$ for every $\lambda$. The monoidal product is addition and the unit is the zero vector.

\begin{proposition}[Monoidal realization functor]
\label{prop:monoidal-realization-functor}
The assignment
\begin{equation}
 \Phi_{\PP}:\LayP\longrightarrow\FactCalc,
 \qquad
 \Phi_{\PP}(\mathbf m)=A^{\mathbf m},
 \label{eq:prime-layer-realization-functor}
\end{equation}
defines a strict symmetric-monoidal functor. Explicitly,
\begin{equation}
 A^{\mathbf m+\mathbf n}
 =
 A^{\mathbf m}\boxt A^{\mathbf n},
 \qquad
 A^{\mathbf 0}=\mathbf 1.
 \label{eq:strict-monoidal-realization}
\end{equation}
If $\mathbf m\le\mathbf m'$, then $A^{\mathbf m}\preccurlyeq A^{\mathbf m'}$, and the corresponding morphism is
\begin{equation}
 J_{\mathbf m,\mathbf m'}z^n
 =
 \frac{A_n^{\mathbf m}}{A_n^{\mathbf m'}}z^n
 =
 \frac{1}{A_n^{\mathbf m'-\mathbf m}}z^n.
 \label{eq:layer-complement-intertwiner}
\end{equation}
\end{proposition}

\begin{proof}
The monoidal identities follow directly from addition of the exponents in \eqref{eq:layer-realization-explicit}. If $\mathbf m\le\mathbf m'$, then
\[
 A_n^{\mathbf m'}
 =
 A_n^{\mathbf m}A_n^{\mathbf m'-\mathbf m}
 \ge A_n^{\mathbf m},
\]
so $A^{\mathbf m}\preccurlyeq A^{\mathbf m'}$. The diagonal formula is the canonical intertwiner of the factorial category.
\end{proof}

Thus every transition map is determined by the complementary layer system. This exact complementarity will later reappear for logarithmic recurrence symbols and their gamma-lift torsors.

\subsection{Logarithmic linearization}
\label{subsec:logarithmic-linearization}

The monoidal product is multiplicative at the level of factorial data but additive after taking the logarithmic first difference.

\begin{definition}[Logarithmic recurrence transform]
\label{def:logarithmic-recurrence-transform}
For a factorial datum $F$, define
\begin{equation}
 \cL_F(n)
 =
 \log\frac{F_n}{F_{n-1}},
 \qquad n\ge 1.
 \label{eq:logarithmic-recurrence-transform}
\end{equation}
\end{definition}

Then
\begin{equation}
 \cL_{F\boxt G}=\cL_F+\cL_G.
 \label{eq:log-transform-monoidal-additivity}
\end{equation}
For an atomic layer,
\begin{equation}
 \cL_{E^\lambda}(n)
 =
 w_\lambda\mathbf 1_{d_\lambda\mid n}.
 \label{eq:atomic-log-recurrence}
\end{equation}
Consequently,
\begin{equation}
 \boxed{
 \cL_{A^{\mathbf m}}(n)
 =
 \sum_{\lambda\in\cA}
 m_\lambda w_\lambda\mathbf 1_{d_\lambda\mid n}.
 }
 \label{eq:finite-layer-log-recurrence}
\end{equation}

Formula \eqref{eq:finite-layer-log-recurrence} is the divisor-poset zeta transform of the layer weights. To state its scalar inversion, aggregate all layers having the same conductor:
\begin{equation}
 \Omega_{\mathbf m}(d)
 =
 \sum_{\substack{\lambda\in\cA\\d_\lambda=d}}
 m_\lambda w_\lambda.
 \label{eq:aggregated-conductor-weight}
\end{equation}
Only finitely many $d$ occur, and
\begin{equation}
 \cL_{A^{\mathbf m}}(n)
 =
 \sum_{d\mid n}\Omega_{\mathbf m}(d).
 \label{eq:divisor-zeta-transform}
\end{equation}

\begin{proposition}[Incidence-algebra inversion]
\label{prop:incidence-algebra-inversion}
For every $d\ge 1$,
\begin{equation}
 \boxed{
 \Omega_{\mathbf m}(d)
 =
 \sum_{e\mid d}
 \mu\!\left(\frac{d}{e}\right)
 \cL_{A^{\mathbf m}}(e).
 }
 \label{eq:mobius-layer-inversion}
\end{equation}
Hence the logarithmic recurrence determines the total weight carried by each conductor.
\end{proposition}

\begin{proof}
This is ordinary M\"obius inversion applied to \eqref{eq:divisor-zeta-transform}.
\end{proof}

Different prime layers may share the same conductor. For example,
\[
 d_{2,1}=d_{3,0}=2.
\]
Thus scalar incidence inversion recovers the total logarithmic weight at conductor $2$, but not yet the labels of its two prime components. The prime factorization of the recurrence ratio supplies the missing refinement.

For $p\in\PP$, put
\begin{equation}
 \cL_{\mathbf m,p}(n)
 =
 v_p\!\left(\frac{A_n^{\mathbf m}}{A_{n-1}^{\mathbf m}}\right).
 \label{eq:p-component-log-recurrence}
\end{equation}
By \eqref{eq:finite-layer-log-recurrence},
\begin{equation}
 \cL_{\mathbf m,p}(n)
 =
 \sum_{\substack{j\ge 0\\(p-1)p^j\mid n}}m_{p,j}.
 \label{eq:p-component-layer-sum}
\end{equation}

\begin{theorem}[Unique prime-layer factorization]
\label{thm:unique-prime-layer-factorization}
The realization map
\[
 \mathbf m\longmapsto A^{\mathbf m}
\]
is injective. More precisely, with the convention
\[
 \cL_{\mathbf m,p}(d_{p,-1})=0,
\]
one has
\begin{equation}
 \boxed{
 m_{p,j}
 =
 \cL_{\mathbf m,p}(d_{p,j})
 -
 \cL_{\mathbf m,p}(d_{p,j-1}),
 \qquad j\ge 0.
 }
 \label{eq:prime-chain-inversion}
\end{equation}
Consequently, the image of $\Phi_{\PP}$ is a free commutative monoid under $\boxt$, with irreducible generators $E^{(p,j)}$.
\end{theorem}

\begin{proof}
For fixed $p$ and $j$, the divisibility relation
\[
 (p-1)p^k\mid(p-1)p^j
\]
holds exactly when $k\le j$. Therefore \eqref{eq:p-component-layer-sum} gives
\[
 \cL_{\mathbf m,p}(d_{p,j})
 =
 \sum_{k=0}^{j}m_{p,k}.
\]
Taking consecutive differences proves \eqref{eq:prime-chain-inversion}. Thus the factorial datum determines every multiplicity $m_{p,j}$, proving injectivity. Since $\NN^{(\cA)}$ is the free commutative monoid on $\cA$, its injective image has the asserted unique factorization.
\end{proof}

\begin{remark}
The theorem uses both pieces of arithmetic structure: divisibility of the conductors along the $p$-power chain and unique factorization of the recurrence ratios into rational primes. If the bases $e^{w_\lambda}$ were arbitrary positive real numbers, the scalar transform \eqref{eq:mobius-layer-inversion} would remain valid, but the labeled factorization need not be unique.
\end{remark}

\subsection{Affine-periodic normal form of finite systems}
\label{subsec:affine-periodic-normal-form}

Every finite layer system has an exact normal form consisting of a linear term and a bounded periodic fluctuation. This is the finite-level growth law that completion will eventually destroy.

For $\mathbf m\in\NN^{(\cA)}$, define its slope and total layer weight by
\begin{align}
 \sigma(\mathbf m)
 &=
 \sum_{\lambda\in\cA}
 \frac{m_\lambda w_\lambda}{d_\lambda},
 \label{eq:finite-layer-slope}\\
 W(\mathbf m)
 &=
 \sum_{\lambda\in\cA}m_\lambda w_\lambda.
 \label{eq:finite-layer-total-weight}
\end{align}
Both sums are finite.

\begin{theorem}[Finite-layer affine-periodic normal form]
\label{thm:finite-layer-normal-form}
Let $\mathbf m\in\NN^{(\cA)}$. Then
\begin{equation}
 \boxed{
 \log A_n^{\mathbf m}
 =
 n\sigma(\mathbf m)+\psi_{\mathbf m}(n),
 }
 \label{eq:affine-periodic-normal-form}
\end{equation}
where
\begin{equation}
 \psi_{\mathbf m}(n)
 =
 -\sum_{\lambda\in\cA}
 m_\lambda w_\lambda
 \left\{\frac{n}{d_\lambda}\right\}
 \label{eq:finite-layer-periodic-error}
\end{equation}
is periodic and satisfies
\begin{equation}
 -W(\mathbf m)<\psi_{\mathbf m}(n)\le 0.
 \label{eq:finite-layer-error-bound}
\end{equation}
One may take
\begin{equation}
 T(\mathbf m)
 =
 \lcm\{d_\lambda:\lambda\in\supp(\mathbf m)\}
 \label{eq:finite-layer-period}
\end{equation}
as a period. In particular,
\begin{equation}
 \lim_{n\to\infty}(A_n^{\mathbf m})^{1/n}
 =e^{\sigma(\mathbf m)}
 \label{eq:finite-layer-root-growth}
\end{equation}
and the analytic radius of its Hilbert realization is
\begin{equation}
 \boxed{
 R_{A^{\mathbf m}}=e^{\sigma(\mathbf m)/2}.
 }
 \label{eq:finite-layer-analytic-radius}
\end{equation}
\end{theorem}

\begin{proof}
Using $\lfloor x\rfloor=x-\{x\}$ in \eqref{eq:layer-realization-explicit} gives \eqref{eq:affine-periodic-normal-form} and \eqref{eq:finite-layer-periodic-error}. Each fractional-part function is periodic modulo its denominator, so their finite sum has period \eqref{eq:finite-layer-period}. The bounds follow from $0\le\{x\}<1$. Dividing by $n$ and passing to the limit proves \eqref{eq:finite-layer-root-growth}; the radius formula follows from the definition of the analytic radius of a factorial datum.
\end{proof}

The slope is a monoidal character:
\begin{equation}
 \sigma(\mathbf m+\mathbf n)
 =
 \sigma(\mathbf m)+\sigma(\mathbf n).
 \label{eq:slope-monoidal-character}
\end{equation}
The periodic defect is additive as well. Thus finite layer systems are completely linearized by the pair
\[
 \bigl(\sigma(\mathbf m),\psi_{\mathbf m}\bigr).
\]

\begin{corollary}[Finite systems are never factorial in scale]
\label{cor:finite-systems-exponential}
For every finite layer system $A^{\mathbf m}$,
\begin{equation}
 e^{-W(\mathbf m)}e^{n\sigma(\mathbf m)}
 <
 A_n^{\mathbf m}
 \le
 e^{n\sigma(\mathbf m)}.
 \label{eq:finite-layer-two-sided-growth}
\end{equation}
In particular, no fixed finite layer system has $n\log n$ growth in its logarithm.
\end{corollary}

The eventual Stirling law of the completed prime datum is therefore a completion phenomenon. It cannot be obtained by assigning a Stirling expansion to each finite object and then passing formally to a limit, because the asymptotic type changes from affine-periodic to factorial.

\subsection{The conductor filtration}
\label{subsec:conductor-filtration}

The correct sequential exhaustion of the prime layers is by conductor rather than by the prime coordinate alone.

\begin{definition}[Conductor truncation]
\label{def:conductor-truncation}
For $T\ge 1$, let
\begin{equation}
 \cA_{\le T}
 =
 \{(p,j)\in\cA:(p-1)p^j\le T\},
 \label{eq:conductor-truncated-layer-set}
\end{equation}
and let $\mathbf 1_{\le T}$ be its indicator vector. Define
\begin{equation}
 A^{[T]}:=A^{\mathbf 1_{\le T}}.
 \label{eq:conductor-truncated-factorial-datum}
\end{equation}
\end{definition}

The set $\cA_{\le T}$ is finite. Indeed, $p\le T+1$, and for each such prime only finitely many $j$ satisfy $(p-1)p^j\le T$. Explicitly,
\begin{equation}
 \#\cA_{\le T}
 =
 \sum_{p\le T+1}
 \left(
 1+
 \left\lfloor
 \frac{\log(T/(p-1))}{\log p}
 \right\rfloor
 \right),
 \label{eq:number-of-active-layers}
\end{equation}
where the summand is interpreted as zero when $p-1>T$.

\begin{proposition}[Cofinality and exact finite detection]
\label{prop:cofinality-exact-detection}
The sequence $(\cA_{\le T})_{T\ge 1}$ is cofinal in $\Fin(\cA)$. Moreover, for every $N\ge 0$ and every $T\ge N$,
\begin{equation}
 A_n^{[T]}=A_n,
 \qquad
 \frac{A_n^{[T]}}{A_{n-1}^{[T]}}
 =
 \frac{A_n}{A_{n-1}},
 \qquad
 0\le n\le N,
 \label{eq:degreewise-stationarity-weights}
\end{equation}
where $A_n=(n+1)!_{\PP}$ is the completed prime datum. Consequently, on the polynomial space $\CC[z]_{\le N}$,
\begin{equation}
 \langle f,g\rangle_{A^{[T]}}
 =
 \langle f,g\rangle_A,
 \qquad
 D_{A^{[T]}}f=D_Af.
 \label{eq:degreewise-stationarity-calculus}
\end{equation}
\end{proposition}

\begin{proof}
Every finite layer set has a largest conductor, proving cofinality. If $n\le N\le T$, then a layer contributes to $A_n$ only when $d_\lambda\le n$, hence it already belongs to $\cA_{\le T}$. This proves equality of the factorial weights. A layer contributes to the recurrence ratio at $n$ only when $d_\lambda\mid n$, which again implies $d_\lambda\le n\le T$. The Hilbert inner product and generalized derivative through degree $N$ depend only on these weights and recurrence ratios.
\end{proof}

The proposition is stronger than pointwise convergence. The system is eventually constant in every fixed degree. Thus the completed polynomial calculus is a genuine degreewise colimit of finite layer calculi, even though its Hilbert completion requires the bounded projective construction of the next section.

\begin{corollary}[Stabilization of generalized exponentials]
\label{cor:stabilization-generalized-exponentials}
Let
\[
 E_{[T]}(z)=\sum_{n\ge 0}\frac{z^n}{A_n^{[T]}},
 \qquad
 E_{\PP}(z)=\sum_{n\ge 0}\frac{z^n}{A_n}.
\]
For each $N$ and $T\ge N$, the Taylor polynomials of degree $N$ agree:
\begin{equation}
 [z^{\le N}]E_{[T]}=[z^{\le N}]E_{\PP}.
 \label{eq:generalized-exponential-stabilization}
\end{equation}
The same holds for the reproducing kernels in each fixed bidegree.
\end{corollary}

\subsection{Complementary layers and transition maps}
\label{subsec:complementary-layers}

For $T\le T'$, define the complementary factorial datum
\begin{equation}
 C^{(T,T')}
 =
 \mathop{\boxt}_{\substack{\lambda\in\cA\\T<d_\lambda\le T'}}E^\lambda.
 \label{eq:complementary-layer-datum}
\end{equation}
Then
\begin{equation}
 A^{[T']}=A^{[T]}\boxt C^{(T,T')}
 \label{eq:truncation-complement-factorization}
\end{equation}
and
\begin{equation}
 C^{(T,T'')}=C^{(T,T')}\boxt C^{(T',T'')}
 \qquad(T\le T'\le T'').
 \label{eq:complement-cocycle}
\end{equation}
The canonical transition map is therefore
\begin{equation}
 J_{T,T'}z^n
 =
 \frac{1}{C_n^{(T,T')}}z^n.
 \label{eq:truncation-transition-map}
\end{equation}

\begin{proposition}[Local identity of transition maps]
\label{prop:local-identity-transition-maps}
If $T\le T'$, then
\begin{equation}
 J_{T,T'}f=f
 \qquad
 \text{for every }f\in\CC[z]_{\le T}.
 \label{eq:transition-map-local-identity}
\end{equation}
Moreover,
\begin{equation}
 J_{T',T''}J_{T,T'}=J_{T,T''},
 \label{eq:transition-map-cocycle}
\end{equation}
and the cocycle identity is the operator realization of the complementary factorization \eqref{eq:complement-cocycle}.
\end{proposition}

\begin{proof}
If $n\le T$, every layer in $C^{(T,T')}$ has conductor greater than $n$, so its $n$th factorial weight is $1$. Hence $C_n^{(T,T')}=1$ and $J_{T,T'}z^n=z^n$. The cocycle identity follows either from the general composition law for canonical intertwiners or directly from \eqref{eq:complement-cocycle}.
\end{proof}

This local identity is an important compatibility between the arithmetic filtration and polynomial degree. New layers alter only the degrees at which their conductors become visible. In particular, the transition system has no hidden renormalization on previously stabilized finite-degree data.

\subsection{The completed prime layer profile}
\label{subsec:completed-prime-layer-profile}

Let $\mathbf 1_{\cA}$ denote the formal multiplicity vector with every coordinate equal to $1$. It is not finitely supported, but its realization is degreewise locally finite:
\begin{equation}
 A_n
 =
 \prod_{(p,j)\in\cA}
 p^{\lfloor n/((p-1)p^j)\rfloor}
 =
 (n+1)!_{\PP}.
 \label{eq:completed-prime-layer-product-section4}
\end{equation}
Its logarithmic recurrence is
\begin{align}
 \cL_A(n)
 &=
 \sum_{\substack{p\in\PP,\ j\ge 0\\(p-1)p^j\mid n}}
 \log p
 \label{eq:completed-log-recurrence-layer-sum}\\
 &=
 \sum_{p-1\mid n}(v_p(n)+1)\log p,
 \label{eq:completed-log-recurrence-prime-sum}
\end{align}
and hence
\begin{equation}
 \frac{A_n}{A_{n-1}}
 =
 \prod_{p-1\mid n}p^{v_p(n)+1}.
 \label{eq:completed-recurrence-ratio-section4}
\end{equation}

For the conductor truncations, the finite-layer slopes are
\begin{equation}
 \sigma_T
 =
 \sum_{\substack{p,j\\(p-1)p^j\le T}}
 \frac{\log p}{(p-1)p^j}.
 \label{eq:truncated-layer-slope}
\end{equation}
They increase without bound. Indeed, already the $j=0$ layers contribute
\[
 \sum_{p\le T+1}\frac{\log p}{p-1},
\]
which diverges. Thus the radii
\begin{equation}
 R_{A^{[T]}}=e^{\sigma_T/2}
 \label{eq:truncated-radii-diverge}
\end{equation}
escape to infinity as more layers are inserted.

This divergence is the first categorical signal of the change in growth regime. Each finite object has a finite exponential slope and a bounded periodic defect, while no finite slope survives completion. The completed object must instead be compared with the classical factorial scale. The categorical completion of the coefficient spaces and differential calculi will be formalized in the next section; the later gamma-lift construction will then show how Stirling normalization selects a distinguished analytic primitive of the completed recurrence.

\subsection{Structural summary}
\label{subsec:layer-structural-summary}

The layer theory developed above may be condensed into the sequence
\begin{equation}
 \boxed{
 \begin{gathered}
 \NN^{(\cA)}
 \xrightarrow{\ \Phi_{\PP}\ }
 \FactCalc,
 \qquad
 \mathbf m+\mathbf n
 \longmapsto
 A^{\mathbf m}\boxt A^{\mathbf n},
 \\
 \cL_{A^{\mathbf m}}(n)
 =
 \sum_{\lambda}m_\lambda w_\lambda\mathbf 1_{d_\lambda\mid n},
 \\
 m_{p,j}
 =
 v_p(\rho_{\mathbf m}(d_{p,j}))
 -v_p(\rho_{\mathbf m}(d_{p,j-1})),
 \\
 \log A_n^{\mathbf m}
 =
 n\sigma(\mathbf m)+\psi_{\mathbf m}(n),
 \qquad
 \psi_{\mathbf m}\text{ bounded and periodic},
 \\
 A^{[T]}|_{\{0,\dots,N\}}
 =
 A|_{\{0,\dots,N\}}
 \qquad(T\ge N).
 \end{gathered}
 }
 \label{eq:layer-structural-summary}
\end{equation}

The first line records the free symmetric-monoidal generation of the finite calculi. The second and third express the arithmetic recoverability of the layers. The fourth identifies the exact finite-level asymptotic type. The final line shows that the infinite prime object is assembled without ambiguity in each finite degree.

These statements distinguish two limits that must not be conflated. The polynomial calculus is obtained degreewise, where the conductor filtration eventually stabilizes exactly. The Hilbert-space realization is obtained by completing infinitely many coefficient constraints simultaneously. The projective--inductive duality governing that second limit is the subject of the next section.

% ===== Source: section5_filtered_completion_projective_inductive_duality.tex =====
\section{Filtered completion and projective--inductive duality}
\label{sec:filtered-completion-duality}

The prime-layer realization has two complementary limiting structures.  On each fixed polynomial degree, the conductor filtration becomes stationary and therefore produces an ordinary inductive limit of finite-dimensional calculi.  On the full coefficient spaces, adding layers strengthens the norm and reverses the inclusions; the completed prime space is consequently a bounded projective limit.  The purpose of this section is to place these two constructions in one exact diagram and to prove that they recover the same Hilbert realization.

This distinction is not terminological.  The ordinary intersection of the finite-layer Hilbert spaces is strictly larger than the completed prime space.  Uniform boundedness of the finite-stage norms is the condition that detects the infinite product of valuation layers.  At the same time, the stable polynomial cores already know the completed generalized derivative, and their graph closure recovers the full unbounded operator.

For $T\ge 1$, let
\begin{equation}
 A^{[T]}_n
 =
 \prod_{\substack{p,j\\(p-1)p^j\le T}}
 p^{\lfloor n/((p-1)p^j)\rfloor},
 \qquad
 A_n=(n+1)!_{\PP},
 \label{eq:section5-truncated-completed-data}
\end{equation}
and abbreviate
\begin{equation}
 \cH_T:=\cH_{A^{[T]}},
 \qquad
 D_T:=D_{A^{[T]}},
 \qquad
 \cH_\PP:=\cH_A,
 \qquad
 D_\PP:=D_A.
 \label{eq:section5-space-operator-notation}
\end{equation}
The conductor filtration satisfies the exact finite-detection property
\begin{equation}
 A^{[T]}_n=A_n,
 \qquad
 \frac{A^{[T]}_n}{A^{[T]}_{n-1}}
 =
 \frac{A_n}{A_{n-1}}
 \qquad(0\le n\le T).
 \label{eq:section5-exact-finite-detection}
\end{equation}

\subsection{Stable polynomial cores}
\label{subsec:section5-stable-polynomial-cores}

Let
\begin{equation}
 \cE_T:=\CC[z]_{\le T}\subseteq\cH_T.
 \label{eq:section5-stable-core}
\end{equation}
We equip $\cE_T$ with the norm inherited from $\cH_T$.  Formula \eqref{eq:section5-exact-finite-detection} shows that this norm is already the completed norm:
\begin{equation}
 \left\|\sum_{n=0}^{T}a_nz^n\right\|_{\cH_T}^2
 =
 \sum_{n=0}^{T}|a_n|^2A_n
 =
 \left\|\sum_{n=0}^{T}a_nz^n\right\|_{\cH_\PP}^2.
 \label{eq:section5-stable-core-norm}
\end{equation}
Consequently, for $T\le T'$ the identity map
\begin{equation}
 \iota^{\mathrm{alg}}_{T,T'}:\cE_T\longrightarrow\cE_{T'}
 \label{eq:section5-algebraic-core-map}
\end{equation}
is an isometry.

\begin{theorem}[Hilbert inductive completion]
\label{thm:section5-hilbert-inductive-completion}
The stable cores form an isometric inductive system, and coefficientwise identification gives
\begin{equation}
 \boxed{
 \cH_\PP
 \cong
 \overline{\indlim_{T}\cE_T}.
 }
 \label{eq:section5-hilbert-inductive-limit}
\end{equation}
Equivalently, $\CC[z]=\bigcup_T\cE_T$ carries a well-defined pre-Hilbert norm, independent of the stage at which a polynomial is represented, and its Hilbert completion is $\cH_\PP$.
\end{theorem}

\begin{proof}
The isometric compatibility follows from \eqref{eq:section5-stable-core-norm}.  The algebraic inductive limit is therefore the polynomial algebra equipped with
\[
 \left\|\sum_{n=0}^{N}a_nz^n\right\|^2
 =
 \sum_{n=0}^{N}|a_n|^2A_n.
\]
Polynomials are dense in the weighted coefficient space $\cH_\PP$, because Taylor truncation converges in its defining $\ell^2$ norm.  Completing the inductive limit gives $\cH_\PP$.
\end{proof}

The exactness also holds at the operator level.  If $T\ge N$, then
\begin{equation}
 D_Tp=D_\PP p
 \qquad(p\in\CC[z]_{\le N}).
 \label{eq:section5-operator-stationarity}
\end{equation}
Thus the finite derivatives define an operator on the algebraic inductive limit:
\begin{equation}
 D_{\mathrm{alg}}z^n
 =
 \left(\prod_{p-1\mid n}p^{v_p(n)+1}\right)z^{n-1}.
 \label{eq:section5-algebraic-limit-derivative}
\end{equation}

\begin{proposition}[Graph-inductive realization]
\label{prop:section5-graph-inductive-realization}
The operator $D_{\mathrm{alg}}$ is closable in $\cH_\PP$, and
\begin{equation}
 \boxed{
 \overline{D_{\mathrm{alg}}}=D_\PP.
 }
 \label{eq:section5-closure-algebraic-derivative}
\end{equation}
More precisely, if $P_N$ denotes Taylor truncation through degree $N$, then for every $f\in\Dom(D_\PP)$,
\begin{equation}
 P_Nf\longrightarrow f,
 \qquad
 D_\PP P_Nf\longrightarrow D_\PP f
 \quad\text{in }\cH_\PP.
 \label{eq:section5-graph-core-convergence}
\end{equation}
Hence the stable polynomial system is a common core for the completed annihilation operator.
\end{proposition}

\begin{proof}
Write $f(z)=\sum_{n\ge0}a_nz^n$.  The domain condition for the weighted backward shift is
\[
 \sum_{n\ge1}|a_n|^2\frac{A_n^2}{A_{n-1}}<\infty.
\]
The tails of both this series and $\sum_n|a_n|^2A_n$ tend to zero.  This proves \eqref{eq:section5-graph-core-convergence}.  Since $D_\PP$ is closed and extends $D_{\mathrm{alg}}$, it is the closure of the latter.
\end{proof}

The same argument applies to the creation operator $M_z$.  In particular, the completed prime calculus is not obtained by guessing an infinite weighted shift after the fact: it is the graph closure of the degreewise stationary finite-layer calculi.

\subsection{The contravariant coefficient diagram}
\label{subsec:section5-contravariant-diagram}

When $T\le T'$, the factorial weights satisfy
\[
 A^{[T]}_n\le A^{[T']}_n.
\]
Accordingly, the identity on coefficients defines a contractive inclusion
\begin{equation}
 I_{T',T}:\cH_{T'}\hookrightarrow\cH_T.
 \label{eq:section5-projective-inclusion}
\end{equation}
These maps satisfy
\begin{equation}
 I_{T'',T}=I_{T',T}I_{T'',T'}
 \qquad(T\le T'\le T''),
 \label{eq:section5-projective-cocycle}
\end{equation}
and therefore form a projective system of Hilbert spaces.

An element of the set-theoretic projective limit is a compatible family $(f_T)_T$ with
\[
 I_{T',T}f_{T'}=f_T.
\]
Because each $I_{T',T}$ is the identity on Taylor coefficients, compatibility forces all $f_T$ to be represented by one formal power series.  The projective limit is therefore the intersection $\bigcap_T\cH_T$.  The completed prime space is its bounded part.

\begin{definition}[Bounded projective limit]
\label{def:section5-bounded-projective-limit}
Define
\begin{equation}
 \projlimb_T\cH_T
 =
 \left\{
 (f_T)_T:
 I_{T',T}f_{T'}=f_T,
 \ \sup_T\|f_T\|_{\cH_T}<\infty
 \right\},
 \label{eq:section5-bounded-projective-limit}
\end{equation}
with norm
\begin{equation}
 \|(f_T)_T\|_{\mathrm{proj}}
 =
 \sup_T\|f_T\|_{\cH_T}.
 \label{eq:section5-projective-limit-norm}
\end{equation}
\end{definition}

\begin{theorem}[Sequential prime completion]
\label{thm:section5-sequential-prime-completion}
The diagonal map
\begin{equation}
 \Delta:\cH_\PP\longrightarrow\projlimb_T\cH_T,
 \qquad
 f\longmapsto(f)_T,
 \label{eq:section5-diagonal-projective-map}
\end{equation}
is an isometric isomorphism.  For every formal power series $f(z)=\sum_{n\ge0}a_nz^n$,
\begin{equation}
 \boxed{
 \|f\|_{\cH_\PP}^2
 =
 \sup_T\|f\|_{\cH_T}^2
 =
 \lim_{T\to\infty}\|f\|_{\cH_T}^2.
 }
 \label{eq:section5-projective-norm-identity}
\end{equation}
\end{theorem}

\begin{proof}
For every $n$, the sequence $A^{[T]}_n$ increases and is eventually equal to $A_n$.  Therefore
\[
 \|f\|_{\cH_T}^2
 =
 \sum_{n\ge0}|a_n|^2A^{[T]}_n
\]
is increasing in $T$.  Monotone convergence gives
\[
 \lim_{T\to\infty}\sum_{n\ge0}|a_n|^2A^{[T]}_n
 =
 \sum_{n\ge0}|a_n|^2A_n.
\]
This proves the norm identity and both injectivity and surjectivity of $\Delta$.
\end{proof}

The preceding theorem has a useful universal formulation.

\begin{corollary}[Bounded-cone universal property]
\label{cor:section5-bounded-cone-universal-property}
Let $X$ be a normed space, and suppose that bounded maps $u_T:X\to\cH_T$ satisfy
\begin{equation}
 I_{T',T}u_{T'}=u_T
 \qquad(T\le T')
 \label{eq:section5-compatible-bounded-cone}
\end{equation}
and
\begin{equation}
 \sup_T\|u_Tx\|_{\cH_T}\le C\|x\|_X
 \qquad(x\in X).
 \label{eq:section5-uniform-cone-bound}
\end{equation}
Then there is a unique bounded map $u:X\to\cH_\PP$ such that $u_T=I_{\PP,T}u$ for every $T$, and $\|u\|\le C$.
\end{corollary}

\begin{proof}
For each $x$, the compatible family $(u_Tx)_T$ belongs to the bounded projective limit.  Apply Theorem~\ref{thm:section5-sequential-prime-completion} pointwise.  Linearity, uniqueness, and the norm estimate follow from the isometric identification.
\end{proof}

Here $I_{\PP,T}:\cH_\PP\hookrightarrow\cH_T$ denotes the canonical inclusion.  The corollary states that $\cH_\PP$ is the projective limit in the category whose cones are required to be uniformly bounded.  Without this boundedness condition one obtains a larger locally convex object.

\subsection{Why the ordinary intersection is too large}
\label{subsec:section5-intersection-too-large}

The bounded modifier in $\projlimb$ is essential.  It cannot be omitted even though every coefficient stabilizes degreewise.

\begin{proposition}[Strictness of the unbounded projective limit]
\label{prop:section5-strict-projective-intersection}
One has a strict inclusion
\begin{equation}
 \boxed{
 \cH_\PP
 \subsetneq
 \bigcap_{T\ge1}\cH_T.
 }
 \label{eq:section5-strict-intersection}
\end{equation}
In particular, membership in every finite-layer space does not imply membership in the completed prime space.
\end{proposition}

\begin{proof}
Fix $T$.  Choose a prime layer $\lambda=(p,j)$ with conductor $d_\lambda>T$.  This layer is absent from $A^{[T]}$ but present in $A$, so
\[
 \frac{A^{[T]}_n}{A_n}
 \le
 p^{-\lfloor n/d_\lambda\rfloor}
 \longrightarrow0
 \qquad(n\to\infty).
\]
Choose inductively integers $n_j$ such that
\begin{equation}
 \frac{A^{[j]}_{n_j}}{A_{n_j}}\le2^{-j}.
 \label{eq:section5-sparse-diagonal-choice}
\end{equation}
Define the formal series
\begin{equation}
 f(z)=\sum_{j\ge1}\frac{z^{n_j}}{\sqrt{A_{n_j}}}.
 \label{eq:section5-intersection-counterexample}
\end{equation}
Its completed norm diverges:
\[
 \|f\|_{\cH_\PP}^2=\sum_{j\ge1}1=\infty.
\]
On the other hand, for fixed $T$ and $j\ge T$, monotonicity gives
\[
 \frac{A^{[T]}_{n_j}}{A_{n_j}}
 \le
 \frac{A^{[j]}_{n_j}}{A_{n_j}}
 \le2^{-j}.
\]
The terms with $j<T$ are finite in number, while the remaining terms are summable.  Hence $f\in\cH_T$ for every $T$.
\end{proof}

\begin{remark}
The ordinary projective limit $\bigcap_T\cH_T$, equipped with the seminorms $\|\cdot\|_{\cH_T}$, remains a natural Fr\'echet-type coefficient space.  It records simultaneous finite-layer regularity.  The prime Hilbert space records the stronger requirement that this regularity remain uniformly bounded as all arithmetic layers are inserted.
\end{remark}

There is a parallel categorical asymmetry.  For every finite $T$ there is a contraction from the truncated calculus to the completed calculus, but no bounded inverse.

\begin{proposition}[The completion boundary is one-sided]
\label{prop:section5-one-sided-completion-boundary}
For each $T$,
\begin{equation}
 A^{[T]}\preccurlyeq A,
 \qquad
 A\not\preccurlyeq A^{[T]}.
 \label{eq:section5-strict-domination}
\end{equation}
Thus the canonical intertwiner $J_{T,\PP}:\cH_T\to\cH_\PP$ is contractive but not boundedly invertible.
\end{proposition}

\begin{proof}
The first assertion follows from $A^{[T]}_n\le A_n$.  For an omitted layer $\lambda=(p,j)$ with $d_\lambda>T$,
\[
 \frac{A_n}{A^{[T]}_n}
 \ge p^{\lfloor n/d_\lambda\rfloor},
\]
which is unbounded.  Hence $A\not\preccurlyeq A^{[T]}$.
\end{proof}

This is the analytic signature of completion: the infinite prime object lies strictly beyond every finite exponential layer system in the factorial preorder.

\subsection{Adjoint transition maps and the dual diagram}
\label{subsec:section5-adjoint-transition-maps}

The covariant transition map between differential calculi is
\begin{equation}
 J_{T,T'}:\cH_T\longrightarrow\cH_{T'},
 \qquad
 J_{T,T'}z^n
 =
 \frac{A^{[T]}_n}{A^{[T']}_n}z^n
 \qquad(T\le T').
 \label{eq:section5-covariant-transition-map}
\end{equation}
It is a contraction and satisfies
\begin{equation}
 D_{T'}J_{T,T'}=J_{T,T'}D_T
 \quad\text{on }\CC[z].
 \label{eq:section5-transition-intertwining}
\end{equation}
Its adjoint is exactly the contravariant coefficient inclusion:
\begin{equation}
 \boxed{
 J_{T,T'}^*=I_{T',T}.
 }
 \label{eq:section5-adjoint-inclusion-identity}
\end{equation}
Indeed, for monomials,
\[
 \left\langle
 J_{T,T'}z^n,z^m
 \right\rangle_{\cH_{T'}}
 =
 \delta_{m,n}A^{[T]}_n
 =
 \left\langle
 z^n,I_{T',T}z^m
 \right\rangle_{\cH_T}.
\]

The transition maps become the identity on the stable core:
\begin{equation}
 J_{T,T'}|_{\cE_T}=\iota^{\mathrm{alg}}_{T,T'}.
 \label{eq:section5-transition-restriction-core}
\end{equation}
Thus the isometric inductive system from Theorem~\ref{thm:section5-hilbert-inductive-completion} is not imposed separately; it is the stationary finite-degree part of the canonical factorial intertwiners.

\begin{theorem}[Projective--inductive bicompletion]
\label{thm:section5-projective-inductive-bicompletion}
The completed prime Hilbert space admits the two canonical realizations
\begin{equation}
 \boxed{
 \overline{\indlim_T\cE_T}
 \ \cong\ 
 \cH_\PP
 \ \cong\ 
 \projlimb_T\cH_T.
 }
 \label{eq:section5-bicompletion}
\end{equation}
Under these identifications:
\begin{enumerate}[label=\textup{(\roman*)}]
 \item the left-hand system is covariant and isometric;
 \item the right-hand system is contravariant and contractive;
 \item the two transition systems are adjoint through \eqref{eq:section5-adjoint-inclusion-identity};
 \item the completed derivative $D_\PP$ is the graph closure of the operator induced on the algebraic inductive limit.
\end{enumerate}
\end{theorem}

\begin{proof}
The two Hilbert-space identifications are Theorems~\ref{thm:section5-hilbert-inductive-completion} and \ref{thm:section5-sequential-prime-completion}.  The adjoint relation was proved above, and the operator assertion is Proposition~\ref{prop:section5-graph-inductive-realization}.
\end{proof}

The theorem can be summarized by the diagram
\begin{equation}
 \boxed{
 \begin{gathered}
 \cE_T
 \xrightarrow{\ \Id\ }
 \cE_{T'}
 \xrightarrow{\ \Id\ }
 \CC[z]
 \xrightarrow{\ \mathrm{completion}\ }
 \cH_\PP,
 \\
 \cH_T
 \xrightarrow{\ J_{T,T'}\ }
 \cH_{T'}
 \xrightarrow{\ J_{T',\PP}\ }
 \cH_\PP,
 \\
 \cH_\PP
 \xhookrightarrow{\ I_{\PP,T'}\ }
 \cH_{T'}
 \xhookrightarrow{\ I_{T',T}\ }
 \cH_T,
 \\
 J_{T,T'}^*=I_{T',T}.
 \end{gathered}
 }
 \label{eq:section5-projective-inductive-diagram}
\end{equation}
The upper line records degreewise stabilization, the middle line the covariant derivative intertwiners, and the lower line the projective coefficient constraints.

\subsection{Cofinal invariance}
\label{subsec:section5-cofinal-invariance}

The conductor truncation is convenient but not intrinsic to the limit.  Let $(\Lambda_i)_{i\in I}$ be any directed cofinal family of finite prime-layer sets, and write $A^{\Lambda_i}$ for the corresponding finite data.

\begin{proposition}[Independence of the exhaustion]
\label{prop:section5-cofinal-invariance}
For every cofinal layer exhaustion,
\begin{equation}
 \cH_\PP
 \cong
 \projlimb_{i\in I}\cH_{A^{\Lambda_i}}
 \label{eq:section5-cofinal-projective-limit}
\end{equation}
isometrically.  If $\cE_i$ denotes the polynomial degrees already detected by $\Lambda_i$, then the completed inductive union of the stable cores is again $\cH_\PP$.
\end{proposition}

\begin{proof}
Cofinality implies that for each fixed $n$, the net $A_n^{\Lambda_i}$ is eventually equal to $A_n$.  Monotone convergence proves the projective norm identity exactly as before.  Every polynomial involves finitely many degrees, so cofinality also places it in an eventually stationary core.  Their union is $\CC[z]$, whose completion is $\cH_\PP$.
\end{proof}

Thus the completed object belongs to the filtered prime-layer system itself, not to a particular enumeration of the layers.

\subsection{Kernel convergence and change of analytic type}
\label{subsec:section5-kernel-convergence}

The projective limit also has a direct analytic manifestation.  Let
\begin{equation}
 K_T(z,w)
 =
 \sum_{n\ge0}\frac{(z\overline w)^n}{A^{[T]}_n},
 \qquad
 K_\PP(z,w)
 =
 \sum_{n\ge0}\frac{(z\overline w)^n}{A_n}
 \label{eq:section5-kernels}
\end{equation}
be the reproducing kernels.  The finite-layer datum has an affine-periodic logarithm
\[
 \log A^{[T]}_n=n\sigma_T+O_T(1),
\]
so its kernel is defined on a disk of radius $e^{\sigma_T/2}$.  The slopes $\sigma_T$ tend to infinity as $T\to\infty$.

\begin{proposition}[Locally uniform kernel completion]
\label{prop:section5-local-uniform-kernel-completion}
The completed kernel $K_\PP$ is entire in $z$ and anti-entire in $w$.  Moreover,
\begin{equation}
 K_T(z,w)\longrightarrow K_\PP(z,w)
 \label{eq:section5-kernel-local-uniform-convergence}
\end{equation}
locally uniformly on $\CC^2$.  The same conclusion holds for the generalized exponentials
\begin{equation}
 E_T(z)=\sum_{n\ge0}\frac{z^n}{A^{[T]}_n},
 \qquad
 E_\PP(z)=\sum_{n\ge0}\frac{z^n}{A_n}.
 \label{eq:section5-generalized-exponentials}
\end{equation}
\end{proposition}

\begin{proof}
Let $C\subset\CC^2$ be compact and choose $R$ with $|z\overline w|\le R$ on $C$.  Since $\sigma_T\to\infty$, there is $T_0$ such that
\[
 \sum_{n\ge0}\frac{R^n}{A^{[T_0]}_n}<\infty.
\]
For $T\ge T_0$,
\[
 0<\frac{1}{A_n}\le\frac{1}{A^{[T]}_n}\le\frac{1}{A^{[T_0]}_n}.
\]
For each fixed $n$, the middle term is eventually equal to $1/A_n$.  Dominated convergence for the power series proves local uniform convergence.  Since $R$ was arbitrary, the limiting kernel is entire in the stated sense.  The proof for $E_T$ is identical.
\end{proof}

This result exhibits the change of analytic type without invoking a Stirling expansion.  Every finite layer system naturally lives on a disk of finite radius, while the completed prime Hilbert space is an entire-function space.  The divergence of the finite exponential slopes supplies larger and larger domains, and the projective coefficient limit selects the entire kernel.

\begin{corollary}[Finite-jet stabilization]
\label{cor:section5-finite-jet-stabilization}
For $T\ge N$, the Taylor coefficients of $K_T$ and $K_\PP$ agree in every bidegree at most $N$, and those of $E_T$ and $E_\PP$ agree through degree $N$.
\end{corollary}

The local uniform convergence and finite-jet stabilization are complementary.  The latter is algebraic and exact; the former is analytic and global.

\subsection{Completion does not commute with finite-layer asymptotics}
\label{subsec:section5-noncommuting-asymptotics}

For fixed $T$, the finite-layer logarithm has the exact form
\begin{equation}
 \log A^{[T]}_n=n\sigma_T+\psi_T(n),
 \qquad
 \psi_T\text{ bounded and periodic}.
 \label{eq:section5-finite-affine-periodic-form}
\end{equation}
No estimate in \eqref{eq:section5-finite-affine-periodic-form} is uniform in $T$: the slopes diverge and the periods acquire new prime-power conductors.  The completed datum therefore cannot be assigned a Stirling expansion by taking a termwise limit of the finite affine-periodic expansions.

The categorical completion records exactly where the obstruction lies.  Degreewise evaluation commutes with the filtered limit, because each degree sees only finitely many layers.  Large-$n$ asymptotics do not commute with the filtered limit, because the set of visible layers grows with $n$.  Symbolically,
\begin{equation}
 \boxed{
 \text{finite degree first: exact stabilization},
 \qquad
 \text{large degree first: changing asymptotic regime}.
 }
 \label{eq:section5-two-orders-of-limit}
\end{equation}

This is the point at which a new analytic object is required.  The completed logarithmic recurrence is canonically determined by the filtered factorial calculus, but a holomorphic solution of its difference equation is not.  The set of normalized solutions forms a torsor under entire periodic functions.  The next part of the paper develops that torsor and introduces the orbitwise Stirling condition that selects its distinguished point.

\subsection{Structural summary}
\label{subsec:section5-structural-summary}

The completion theory of Part~I is summarized by
\begin{equation}
 \boxed{
 \begin{gathered}
 A^{[T]}_n=A_n\quad(n\le T),
 \\
 \overline{\indlim_T\CC[z]_{\le T}}
 \cong\cH_\PP
 \cong\projlimb_T\cH_T,
 \\
 J_{T,T'}^*=I_{T',T},
 \qquad
 \overline{D_{\mathrm{alg}}}=D_\PP,
 \\
 \cH_\PP\subsetneq\bigcap_T\cH_T,
 \\
 K_T\longrightarrow K_\PP
 \quad\text{locally uniformly on }\CC^2.
 \end{gathered}
 }
 \label{eq:section5-structural-summary}
\end{equation}
The completed prime factorial calculus is therefore simultaneously an inductive completion of its stationary polynomial cores and a bounded projective completion of its full coefficient spaces.  This projective--inductive duality closes the categorical construction of the discrete calculus.  The remaining problem is to pass from the completed recurrence to a distinguished analytic gamma object.

\clearpage
\part{Gamma objects and Stirling selection}
% ===== Source: section6_recurrence_symbols_gamma_lift_torsors.tex =====
\section{Recurrence symbols and gamma-lift torsors}
\label{sec:recurrence-symbols-gamma-lift-torsors}

The completed factorial calculus constructed in Part~I determines a discrete logarithmic recurrence exactly.  It does not, by itself, determine a preferred holomorphic solution of the corresponding difference equation.  This is the first genuinely analytic ambiguity in the theory.

The purpose of this section is to separate two levels of nonuniqueness.  A discrete factorial skeleton admits many entire recurrence symbols, and a fixed recurrence symbol admits many normalized entire primitives.  The second ambiguity is controlled by the additive group of normalized entire $1$-periodic functions.  Consequently, the normalized primitives of a fixed symbol form a torsor rather than a vector space with a preferred origin.

This torsor structure is compatible with the monoidal product of factorial data.  Addition of logarithmic recurrence symbols corresponds to a contracted sum of torsors, while compatible complementary primitives turn a filtered layer system into an inverse system of torsors.  The resulting formalism is the analytic counterpart of the projective--inductive duality developed in the preceding section.

\subsection{The logarithmic skeleton of a factorial datum}
\label{subsec:section6-logarithmic-skeleton}

Let
\begin{equation}
 F=(F_n)_{n\ge 0},
 \qquad
 F_0=1,
 \qquad
 F_n>0,
 \label{eq:section6-factorial-datum}
\end{equation}
be a positive factorial datum.  Its logarithmic factorial skeleton is the sequence
\begin{equation}
 h_F(m)=\log F_{m-1},
 \qquad m\ge 1,
 \label{eq:section6-logarithmic-skeleton}
\end{equation}
and its discrete logarithmic recurrence is
\begin{equation}
 \ell_F(m)
 =h_F(m+1)-h_F(m)
 =\log\frac{F_m}{F_{m-1}},
 \qquad m\ge 1.
 \label{eq:section6-discrete-logarithmic-recurrence}
\end{equation}
The indexing is chosen so that $h_F(1)=0$ and, in the prime case,
\[
 h_{\PP}(m)=\log m!_{\PP}.
\]

\begin{definition}[Entire factorial extension]
\label{def:section6-entire-factorial-extension}
An \emph{entire extension} of $F$ is a function $H\in\cO(\CC)$ satisfying
\begin{equation}
 H(m)=h_F(m)
 \qquad(m\ge 1).
 \label{eq:section6-entire-extension-condition}
\end{equation}
We write $\Ext(F)$ for the set of all such functions.
\end{definition}

The positive integers form a discrete subset of $\CC$.  The classical entire interpolation theorem therefore implies that $\Ext(F)$ is nonempty for every positive factorial datum.  The interpolation is highly nonunique.  Put
\begin{equation}
 \cZ_+
 =
 \left\{
 Q\in\cO(\CC):Q(m)=0\text{ for every }m\ge 1
 \right\}.
 \label{eq:section6-discrete-null-space}
\end{equation}

\begin{proposition}[The extension torsor]
\label{prop:section6-extension-torsor}
The space $\Ext(F)$ is a torsor under the additive group $\cZ_+$.  Explicitly, if $H_0\in\Ext(F)$, then
\begin{equation}
 \Ext(F)=H_0+\cZ_+.
 \label{eq:section6-extension-affine-space}
\end{equation}
\end{proposition}

\begin{proof}
The difference of two entire extensions vanishes at every positive integer and therefore belongs to $\cZ_+$.  Conversely, adding an element of $\cZ_+$ preserves all prescribed values.
\end{proof}

Thus the integer factorial values determine an affine space of entire functions, not a distinguished analytic continuation.  We now refine this space by fixing the entire function that encodes the one-step recurrence.

\subsection{Entire recurrence symbols}
\label{subsec:section6-entire-recurrence-symbols}

Let
\begin{equation}
 (\Deltaop H)(z)=H(z+1)-H(z)
 \label{eq:section6-difference-operator}
\end{equation}
be the forward difference operator on $\cO(\CC)$.

\begin{definition}[Recurrence symbol]
\label{def:section6-recurrence-symbol}
An \emph{entire recurrence symbol} for $F$ is an entire function $L$ satisfying
\begin{equation}
 L(m)=\ell_F(m)
 \qquad(m\ge 1).
 \label{eq:section6-recurrence-symbol-interpolation}
\end{equation}
It is called \emph{liftable} if there exists $H\in\cO(\CC)$ such that
\begin{equation}
 \Deltaop H=L,
 \qquad
 H(1)=0.
 \label{eq:section6-normalized-difference-equation}
\end{equation}
We denote the set of liftable recurrence symbols by $\Rec(F)$.
\end{definition}

Every entire extension produces a liftable recurrence symbol:
\begin{equation}
 \Deltaop:\Ext(F)\longrightarrow\Rec(F).
 \label{eq:section6-extension-to-symbol-map}
\end{equation}
Conversely, a normalized lift of an interpolating recurrence symbol necessarily recovers the factorial skeleton.

\begin{lemma}[Discrete recovery]
\label{lem:section6-discrete-recovery}
Suppose $L$ satisfies \eqref{eq:section6-recurrence-symbol-interpolation} and $H$ satisfies \eqref{eq:section6-normalized-difference-equation}.  Then
\begin{equation}
 H(m)=h_F(m)=\log F_{m-1}
 \qquad(m\ge 1).
 \label{eq:section6-discrete-recovery}
\end{equation}
In particular, $H\in\Ext(F)$.
\end{lemma}

\begin{proof}
The assertion holds at $m=1$.  If it holds at $m$, then
\[
 H(m+1)
 =H(m)+L(m)
 =h_F(m)+\ell_F(m)
 =h_F(m+1).
\]
Induction proves the result.
\end{proof}

Let
\begin{equation}
 \cI_+
 =
 \left\{
 E\in\cO(\CC):E(m)=0\text{ for every }m\ge 1
 \right\}.
 \label{eq:section6-symbol-null-space}
\end{equation}
As sets, $\cI_+=\cZ_+$, but the notation records two different roles: elements of $\cZ_+$ perturb factorial extensions, while elements of $\cI_+$ perturb recurrence symbols.  Only those symbol perturbations lying in $\Deltaop\cZ_+$ preserve liftability automatically.

\begin{proposition}[Two levels of analytic gauge]
\label{prop:section6-two-level-gauge}
For every positive factorial datum $F$,
\begin{equation}
 \Rec(F)=L_0+\Deltaop\cZ_+
 \label{eq:section6-recurrence-symbol-affine-space}
\end{equation}
for any $L_0\in\Rec(F)$.  Moreover, the map
\begin{equation}
 \Deltaop:\Ext(F)\longrightarrow\Rec(F)
 \label{eq:section6-extension-symbol-quotient-map}
\end{equation}
is an affine surjection whose fibers are translates of the normalized periodic group introduced below.
\end{proposition}

\begin{proof}
Choose $H_0\in\Ext(F)$ with $\Deltaop H_0=L_0$.  By Proposition~\ref{prop:section6-extension-torsor}, every extension has the form $H_0+Q$ with $Q\in\cZ_+$, and its recurrence symbol is $L_0+\Deltaop Q$.  Conversely, every such symbol is lifted by $H_0+Q$.  The description of the fibers follows from the kernel of $\Deltaop$ together with the normalization at $1$.
\end{proof}

The distinction in Proposition~\ref{prop:section6-two-level-gauge} is important.  Choosing an entire recurrence symbol already imposes analytic information not contained in the discrete factorial sequence.  Once a symbol is fixed, a second and more structured ambiguity remains.

\subsection{Normalized gamma-lift torsors}
\label{subsec:section6-gamma-lift-torsors}

Define
\begin{equation}
 \Per
 =
 \ker(\Deltaop)
 =
 \left\{
 P\in\cO(\CC):P(z+1)=P(z)
 \right\}
 \label{eq:section6-periodic-group}
\end{equation}
and its normalized subgroup
\begin{equation}
 \Per_0
 =
 \left\{
 P\in\Per:P(1)=0
 \right\}.
 \label{eq:section6-normalized-periodic-group}
\end{equation}
This is an infinite-dimensional additive group.  For example,
\[
 \sin(2\pi z)\in\Per_0.
\]

\begin{definition}[Gamma-lift torsor]
\label{def:section6-gamma-lift-torsor}
For a liftable entire symbol $L$, define
\begin{equation}
 \Lift(L)
 =
 \left\{
 H\in\cO(\CC):
 \Deltaop H=L,
 \ H(1)=0
 \right\}.
 \label{eq:section6-lift-torsor-definition}
\end{equation}
\end{definition}

\begin{theorem}[Torsor theorem]
\label{thm:section6-torsor-theorem}
For every liftable entire symbol $L$, the set $\Lift(L)$ is a torsor under $\Per_0$.  If $H_0\in\Lift(L)$, then
\begin{equation}
 \boxed{
 \Lift(L)=H_0+\Per_0.
 }
 \label{eq:section6-lift-torsor-affine-form}
\end{equation}
In particular, the action of $\Per_0$ is free and transitive.
\end{theorem}

\begin{proof}
If $H_1,H_2\in\Lift(L)$, then
\[
 \Deltaop(H_1-H_2)=0
\]
and $(H_1-H_2)(1)=0$, so $H_1-H_2\in\Per_0$.  Conversely, if $P\in\Per_0$, then $H_0+P$ has the same difference $L$ and the same normalization.
\end{proof}

\begin{corollary}[Integer trace is gauge invariant]
\label{cor:section6-integer-trace-gauge-invariant}
Let $L\in\Rec(F)$.  Every $H\in\Lift(L)$ satisfies
\begin{equation}
 H(m)=\log F_{m-1}
 \qquad(m\ge 1).
 \label{eq:section6-all-lifts-same-integer-trace}
\end{equation}
Thus the action of $\Per_0$ is invisible on the positive integers.
\end{corollary}

\begin{proof}
The first assertion is Lemma~\ref{lem:section6-discrete-recovery}.  Alternatively, if $P\in\Per_0$, then periodicity gives $P(m)=P(1)=0$ for every integer $m$.
\end{proof}

The corollary isolates the role of the later Stirling condition.  It will not alter the discrete prime factorial values; those are already fixed.  It will select one analytic extension among functions that agree at every positive integer.

\begin{remark}[No preferred origin]
\label{rem:section6-no-preferred-origin}
Although $\Per_0$ is a vector space, $\Lift(L)$ is not canonically one.  Declaring a lift to be the zero element is precisely the additional choice that the orbitwise Stirling theory must justify.  Until such a choice is made, all statements should be formulated equivariantly on the torsor.
\end{remark}

\subsection{Monoidal addition and contracted sums}
\label{subsec:section6-contracted-sums}

The logarithm converts the monoidal product of factorial data into addition.  If $F$ and $G$ are positive factorial data, then
\begin{equation}
 h_{F\boxtimes G}=h_F+h_G,
 \qquad
 \ell_{F\boxtimes G}=\ell_F+\ell_G.
 \label{eq:section6-logarithmic-monoidal-additivity}
\end{equation}
Accordingly, $L_F+L_G$ is a recurrence symbol for $F\boxtimes G$ whenever $L_F$ and $L_G$ are symbols for the two factors.

The correct addition operation on lift torsors is the contracted sum.  If $T_1$ and $T_2$ are torsors under the additive group $\Per_0$, define
\begin{equation}
 T_1\boxplusT T_2
 =
 (T_1\times T_2)/\Per_0,
 \label{eq:section6-contracted-sum-definition}
\end{equation}
where $P\in\Per_0$ acts anti-diagonally by
\begin{equation}
 P\cdot(t_1,t_2)=(t_1+P,t_2-P).
 \label{eq:section6-antidiagonal-action}
\end{equation}
The quotient inherits a free transitive $\Per_0$-action from translation in either coordinate.

\begin{proposition}[Monoidal law for lift torsors]
\label{prop:section6-monoidal-law-lifts}
For liftable entire symbols $L_1,L_2$, addition of functions induces a canonical isomorphism of $\Per_0$-torsors
\begin{equation}
 \boxed{
 \Lift(L_1)\boxplusT\Lift(L_2)
 \xrightarrow{\ \sim\ }
 \Lift(L_1+L_2),
 \qquad
 [H_1,H_2]\longmapsto H_1+H_2.
 }
 \label{eq:section6-contracted-sum-isomorphism}
\end{equation}
These isomorphisms are associative, symmetric, and unital, with unit torsor
\begin{equation}
 \Lift(0)=\Per_0.
 \label{eq:section6-unit-torsor}
\end{equation}
\end{proposition}

\begin{proof}
The sum $H_1+H_2$ is normalized and has difference $L_1+L_2$.  It is unchanged by the anti-diagonal action, so the map descends to the contracted sum.  If $H\in\Lift(L_1+L_2)$ and $H_1\in\Lift(L_1)$ is chosen, then $H-H_1\in\Lift(L_2)$, proving surjectivity.  Two pairs have the same sum exactly when they differ by the anti-diagonal action, proving injectivity.  Associativity, symmetry, and the unit law follow from addition in $\cO(\CC)$.
\end{proof}

Let $\Tors(\Per_0)$ denote the groupoid of $\Per_0$-torsors and equivariant maps.  Proposition~\ref{prop:section6-monoidal-law-lifts} says that liftable logarithmic symbols take values naturally in the symmetric monoidal groupoid
\[
 (\Tors(\Per_0),\boxplusT).
\]
This is the torsor-valued analogue of the symmetric-monoidal factorial-calculus realization from Part~I.

\subsection{Complementary primitives and filtered systems}
\label{subsec:section6-filtered-torsor-systems}

Suppose $I$ is a directed poset with a least element $0$.  Let $(L_i)_{i\in I}$ be liftable entire symbols.  Assume that for every $i\le j$ there is a complementary symbol $L_{j/i}$ satisfying
\begin{equation}
 L_j=L_i+L_{j/i}.
 \label{eq:section6-complementary-symbol}
\end{equation}
A choice of normalized complementary primitives is a family
\begin{equation}
 C_{i,j}\in\Lift(L_{j/i})
 \qquad(i\le j)
 \label{eq:section6-complementary-primitives}
\end{equation}
satisfying
\begin{equation}
 C_{i,i}=0,
 \qquad
 C_{i,k}=C_{i,j}+C_{j,k}
 \qquad(i\le j\le k).
 \label{eq:section6-complementary-cocycle}
\end{equation}
The cocycle law is the analytic counterpart of multiplicativity of complementary layer products.

Define
\begin{equation}
 r_{i,j}:\Lift(L_j)\longrightarrow\Lift(L_i),
 \qquad
 r_{i,j}(H)=H-C_{i,j}.
 \label{eq:section6-torsor-transition-map}
\end{equation}

\begin{proposition}[Filtered torsor diagram]
\label{prop:section6-filtered-torsor-diagram}
The maps $r_{i,j}$ are $\Per_0$-equivariant bijections and satisfy
\begin{equation}
 r_{i,k}=r_{i,j}\circ r_{j,k}.
 \label{eq:section6-torsor-transition-cocycle}
\end{equation}
Hence
\begin{equation}
 i\longmapsto\Lift(L_i)
 \label{eq:section6-torsor-valued-functor}
\end{equation}
defines a contravariant functor from $I$ to $\Tors(\Per_0)$.
\end{proposition}

\begin{proof}
Subtracting $C_{i,j}$ changes the difference from $L_j$ to $L_i$ and preserves the normalization.  The inverse is $H\mapsto H+C_{i,j}$.  Equivariance is immediate, and the cocycle identity follows from \eqref{eq:section6-complementary-cocycle}.
\end{proof}

\begin{theorem}[Compatible lifts form a torsor]
\label{thm:section6-compatible-lifts-torsor}
For a filtered system with compatible complementary primitives, the inverse limit
\begin{equation}
 \varprojlim_{i\in I}\Lift(L_i)
 =
 \left\{
 (H_i)_i:r_{i,j}(H_j)=H_i\text{ whenever }i\le j
 \right\}
 \label{eq:section6-inverse-limit-lifts}
\end{equation}
is a nonempty $\Per_0$-torsor under diagonal translation
\begin{equation}
 (H_i)_i+P=(H_i+P)_i.
 \label{eq:section6-diagonal-periodic-translation}
\end{equation}
More precisely, evaluation at the least stage gives an equivariant bijection
\begin{equation}
 \boxed{
 \varprojlim_{i\in I}\Lift(L_i)
 \xrightarrow{\ \sim\ }
 \Lift(L_0).
 }
 \label{eq:section6-inverse-limit-base-torsor}
\end{equation}
\end{theorem}

\begin{proof}
Given $H_0\in\Lift(L_0)$, define
\[
 H_i=H_0+C_{0,i}.
\]
The cocycle law gives
\[
 r_{i,j}(H_j)
 =H_0+C_{0,j}-C_{i,j}
 =H_0+C_{0,i}
 =H_i.
\]
Thus every element at the base stage extends uniquely to a compatible family.  The construction commutes with translation by $\Per_0$.
\end{proof}

The theorem expresses a central principle:
\begin{equation}
 \boxed{
 \text{filtered compatibility transports the ambiguity; it does not remove it.}
 }
 \label{eq:section6-completion-preserves-ambiguity}
\end{equation}
In the prime-layer system, the least object is the empty layer datum and $L_0=0$.  Therefore the compatible finite-layer lift system remains a torsor modeled on $\Per_0$.  Later sections construct a canonical family of complementary cyclotomic primitives and study its analytic passage to the completed prime symbol.

\subsection{Exponentiation and multiplicative gamma objects}
\label{subsec:section6-exponentiation-gamma-objects}

The logarithmic formulation is primary because it is additive and has no branch ambiguity.  Nevertheless, every lift has a multiplicative realization.

For $H\in\Lift(L)$ define
\begin{equation}
 \Gamma_H(z)=\exp(H(z)),
 \qquad
 R_L(z)=\exp(L(z)).
 \label{eq:section6-gamma-and-multiplier}
\end{equation}
Then
\begin{equation}
 \Gamma_H(1)=1,
 \qquad
 \Gamma_H(z+1)=R_L(z)\Gamma_H(z).
 \label{eq:section6-multiplicative-functional-equation}
\end{equation}
If $L\in\Rec(F)$, discrete recovery gives
\begin{equation}
 \Gamma_H(m)=F_{m-1}
 \qquad(m\ge 1).
 \label{eq:section6-gamma-integer-values}
\end{equation}

Define the zero-winding periodic unit group
\begin{equation}
 \Per_0^{\exp}
 =
 \exp(\Per_0)
 =
 \left\{
 e^P:P\in\Per_0
 \right\}.
 \label{eq:section6-zero-winding-periodic-units}
\end{equation}

\begin{proposition}[Multiplicative gamma torsor]
\label{prop:section6-multiplicative-gamma-torsor}
Exponentiation identifies $\Lift(L)$ with a torsor under $\Per_0^{\exp}$.  Two gamma objects associated with the same logarithmic symbol satisfy
\begin{equation}
 \frac{\Gamma_{H_1}(z)}{\Gamma_{H_2}(z)}=e^{P(z)}
 \label{eq:section6-gamma-ratio-periodic}
\end{equation}
for a unique $P\in\Per_0$.
\end{proposition}

\begin{proof}
The ratio is $e^{H_1-H_2}$, and Theorem~\ref{thm:section6-torsor-theorem} gives $H_1-H_2\in\Per_0$.  Uniqueness follows because if $e^{P_1}=e^{P_2}$, then $P_1-P_2$ is an entire function with values in $2\pi i\ZZ$, hence is constant; its value at $1$ is zero.
\end{proof}

\begin{remark}[Winding and the need for logarithmic symbols]
\label{rem:section6-winding-logarithmic-symbols}
The full group of normalized nonvanishing entire $1$-periodic functions is larger than $\Per_0^{\exp}$.  For example,
\[
 e^{2\pi i(z-1)}
\]
is normalized and periodic but has nonzero winding on the cylinder $\CC/\ZZ$ and is not the exponential of a periodic entire function.  Fixing the additive symbol $L$, rather than only its multiplier $e^L$, removes this integral winding ambiguity.  This is one reason the logarithmic gamma-lift torsor is the natural primary object.
\end{remark}

\subsection{The classical factorial as a model}
\label{subsec:section6-classical-factorial-model}

For the ordinary factorial datum $F_n=n!$, the discrete logarithmic recurrence is
\begin{equation}
 \ell_F(m)=\log m.
 \label{eq:section6-classical-discrete-recurrence}
\end{equation}
No entire function agrees globally with a fixed branch of $\log z$, so the classical gamma function does not fit the entire-symbol category on all of $\CC$ without modifying the ambient domain.  On a simply connected slit domain containing the positive axis, however, one may take
\begin{equation}
 L(z)=\Log z,
 \qquad
 H(z)=\Log\Gamma(z),
 \label{eq:section6-classical-log-gamma-lift}
\end{equation}
and the same torsor statement holds with holomorphic $1$-periodic functions on that domain.

The prime construction developed later differs in an important respect.  Its cyclotomic layer symbols and their normalized primitives are entire.  The completed prime gamma-lift problem can therefore be formulated globally on $\CC$, even though the comparison model used in Stirling theory will involve the classical logarithmic gamma function on a right half-plane.

\begin{remark}[Domain-relative version]
\label{rem:section6-domain-relative-version}
Every construction in this section has a domain-relative analogue.  If $\Omega\subseteq\CC$ satisfies $\Omega+1\subseteq\Omega$, replace $\cO(\CC)$ by $\cO(\Omega)$ and define periodicity wherever both $z$ and $z+1$ lie in $\Omega$.  The torsor proofs are purely algebraic and remain unchanged.  We use the entire version for the prime recurrence and the half-plane version when comparing it with $\log\Gamma$.
\end{remark}

\subsection{What remains to be selected}
\label{subsec:section6-what-remains-selected}

For the completed prime factorial datum $A_n=(n+1)!_{\PP}$, the discrete skeleton
\begin{equation}
 h_{\PP}(m)=\log m!_{\PP}
 \label{eq:section6-prime-discrete-skeleton}
\end{equation}
is canonical.  The prime-layer construction will also provide a canonical entire recurrence symbol
\begin{equation}
 L_{\PP}^{\mathrm{cyc}}.
 \label{eq:section6-prime-cyclotomic-symbol-placeholder}
\end{equation}
Once that symbol is fixed, the normalized analytic extensions form the torsor
\begin{equation}
 \Lift(L_{\PP}^{\mathrm{cyc}})
 =H_0+\Per_0.
 \label{eq:section6-prime-lift-torsor-placeholder}
\end{equation}
Every element has the same integer values and the same exact difference.  The remaining distinction is transverse to the integer translation orbits.

A selection principle must therefore satisfy two requirements.  First, it must be equivariant under the $\Per_0$-action so that its dependence on the lift ambiguity is explicit.  Second, it must detect a periodic perturbation even though that perturbation vanishes on every positive integer.  Ordinary pointwise asymptotics along the integers cannot do this.

The next section introduces orbitwise Stirling profiles.  These profiles compare a lift with a prescribed asymptotic model simultaneously along all translates $N+u$, $u\in[0,1]$.  Under periodic translation they transform by addition of the translating function itself.  Their vanishing therefore provides exactly the torsor-trivializing condition needed to select a unique prime gamma lift.

\subsection{Structural summary}
\label{subsec:section6-structural-summary}

The analytic ambiguity of a factorial datum is organized by the diagram
\begin{equation}
 \boxed{
 \begin{gathered}
 \Ext(F)=H_0+\cZ_+,
 \qquad
 \Rec(F)=L_0+\Deltaop\cZ_+,
 \\
 \Deltaop:\Ext(F)\twoheadrightarrow\Rec(F),
 \qquad
 \ker(\Deltaop|_{H(1)=0})=\Per_0,
 \\
 \Lift(L)=H+\Per_0,
 \qquad
 \Lift(L_1)\boxplusT\Lift(L_2)\cong\Lift(L_1+L_2),
 \\
 \varprojlim_i\Lift(L_i)\text{ is again a }\Per_0\text{-torsor},
 \\
 \Gamma_H(z+1)=e^{L(z)}\Gamma_H(z),
 \qquad
 \Gamma_H(m)=F_{m-1}.
 \end{gathered}
 }
 \label{eq:section6-structural-summary}
\end{equation}
The discrete factorial values, the entire recurrence symbol, and the normalized gamma lift are thus three distinct levels of structure.  The first determines neither of the latter two.  The prime-layer theory will canonically determine the second, while orbitwise Stirling rigidity will determine the third.

% ===== Source: section7_stirling_data_equivariant_torsor_profiles.tex =====
\section{Stirling data as equivariant torsor profiles}
\label{sec:stirling-data-equivariant-torsor-profiles}

A normalized gamma lift is invisible on the positive integers: every element of a fixed torsor $\Lift(L)$ has the same discrete trace.  A selection principle must therefore examine the behavior of a lift away from the integer lattice.  The natural coordinates are the translation orbits
\[
 N+[0,1],
 \qquad N\in\NN,
\]
and the natural observable is the variation of the asymptotic defect across one such orbit.

This section develops that observable abstractly.  Given a holomorphic comparison model $M$, we average the defect $H-M$ along the integer direction and then remove its value at the base point of the orbit.  The resulting transverse profile is a function of $u\in[0,1]$.  Translation of $H$ by a normalized entire periodic function translates the profile by the restriction of that same function.  Thus the profile construction is equivariant for the torsor action, rather than invariant under it.

The distinction is essential.  A gauge-invariant quantity could not distinguish two lifts in the same torsor.  An equivariant profile records the ambiguity faithfully and makes it possible to select a preferred lift by requiring the profile to vanish asymptotically.  The present section formulates this mechanism, isolates the affine comparison parameter, and introduces the twice-averaged profile used for the completed prime calculus.  The rigidity and existence consequences are proved in the next section.

\subsection{Translation orbits and transverse defect}
\label{subsec:section7-translation-orbits}

Let
\begin{equation}
 \Omega_\sigma
 =\{z\in\CC:\Re z>\sigma\}
 \label{eq:section7-right-half-plane}
\end{equation}
be a right half-plane.  Suppose that $L\in\cO(\CC)$ is liftable, that
\[
 H\in\Lift(L),
\]
and that $M\in\cO(\Omega_\sigma)$ is a comparison model.  For $N$ sufficiently large, define the defect on the $N$th translation orbit by
\begin{equation}
 E_{H,M}(N+u)=H(N+u)-M(N+u),
 \qquad 0\le u\le 1.
 \label{eq:section7-orbit-defect}
\end{equation}

The scalar value $E_{H,M}(N)$ contains the longitudinal behavior of the defect along the positive axis.  The difference
\begin{equation}
 E_{H,M}(N+u)-E_{H,M}(N)
 \label{eq:section7-raw-transverse-defect}
\end{equation}
contains its transverse behavior across a fundamental translation interval.  The latter is the component that detects the periodic gauge.

Let
\begin{equation}
 \cX=C([0,1];\CC),
 \qquad
 \cX_0=\{\varphi\in\cX:\varphi(0)=0\},
 \label{eq:section7-profile-spaces}
\end{equation}
with the uniform norm.  Restriction to the fundamental interval defines a linear map
\begin{equation}
 \res:\Per_0\longrightarrow\cX_0,
 \qquad
 (\res P)(u)=P(u).
 \label{eq:section7-periodic-restriction}
\end{equation}
Since $P$ is $1$-periodic and $P(1)=0$, one also has $P(0)=0$.  The map \eqref{eq:section7-periodic-restriction} is injective: an entire function vanishing on the interval $[0,1]$ vanishes identically.

\begin{definition}[Raw orbitwise profile]
\label{def:section7-raw-profile}
The raw orbitwise profile of $H$ relative to $M$ at height $N$ is
\begin{equation}
 \cT_N^M(H)(u)
 =E_{H,M}(N+u)-E_{H,M}(N),
 \qquad 0\le u\le 1.
 \label{eq:section7-raw-profile}
\end{equation}
Thus $\cT_N^M(H)\in\cX_0$.
\end{definition}

The profile has an exact endpoint interpretation.  Since $\Deltaop H=L$,
\begin{equation}
 \cT_N^M(H)(1)
 =L(N)-\Deltaop M(N).
 \label{eq:section7-endpoint-recurrence-mismatch}
\end{equation}
Thus the endpoint measures the discrepancy between the exact recurrence symbol and the one-step increment of the comparison model.

\begin{proposition}[Exact gauge covariance]
\label{prop:section7-raw-gauge-covariance}
For every $P\in\Per_0$,
\begin{equation}
 \boxed{
 \cT_N^M(H+P)
 =\cT_N^M(H)+\res P.
 }
 \label{eq:section7-raw-gauge-covariance}
\end{equation}
\end{proposition}

\begin{proof}
Periodicity gives $P(N+u)=P(u)$ and $P(N)=P(0)=0$.  Substitution into \eqref{eq:section7-raw-profile} proves the identity.
\end{proof}

The equality is exact for every $N$; it is not an asymptotic statement.  The torsor ambiguity has been converted into translation by the fixed vector $\res P$ in the profile space.

\subsection{Regular orbitwise means}
\label{subsec:section7-regular-orbitwise-means}

The raw profile is appropriate when $H-M$ has a pointwise asymptotic expansion along every translate.  In the prime problem, however, the integer defect is deliberately left unsmoothed and exhibits arithmetic fluctuations.  We therefore permit averaging in the longitudinal variable while retaining the full transverse coordinate $u$.

\begin{definition}[Regular orbitwise mean]
\label{def:section7-regular-orbitwise-mean}
A \emph{regular orbitwise mean} is a family of nonnegative weights
\begin{equation}
 \omega_{N,n}\ge 0,
 \qquad N,n\ge 1,
 \label{eq:section7-mean-weights}
\end{equation}
such that:
\begin{enumerate}[label=\textnormal{(\roman*)}]
 \item for every $N$, only finitely many $\omega_{N,n}$ are nonzero;
 \item $\sum_{n\ge 1}\omega_{N,n}=1$;
 \item for every fixed $R\ge 1$,
 \begin{equation}
  \sum_{n\le R}\omega_{N,n}\longrightarrow 0
  \qquad(N\to\infty).
  \label{eq:section7-weights-escape}
 \end{equation}
\end{enumerate}
For a function $f$ defined on a sufficiently far right half-plane, put
\begin{equation}
 (\Mean_N f)(u)
 =\sum_{n\ge 1}\omega_{N,n}f(n+u).
 \label{eq:section7-orbitwise-mean}
\end{equation}
\end{definition}

The escape condition is the usual Toeplitz regularity requirement: bounded initial segments receive asymptotically no mass.  The normalization has a second consequence particularly important here.  If $P$ is $1$-periodic, then
\begin{equation}
 \Mean_NP(u)=P(u)
 \label{eq:section7-mean-fixes-periodic}
\end{equation}
for every $N$.  Longitudinal averaging regularizes the nonperiodic defect but leaves the periodic gauge untouched.

\begin{definition}[Averaged transverse profile]
\label{def:section7-averaged-transverse-profile}
For a regular orbitwise mean $\Mean=(\Mean_N)$, define
\begin{equation}
 \cR_{N,\Mean}^M(u;H)
 =\Mean_NE_{H,M}(u)
 =\sum_{n\ge1}\omega_{N,n}
   \bigl(H(n+u)-M(n+u)\bigr),
 \label{eq:section7-averaged-defect}
\end{equation}
and its transverse part
\begin{equation}
 \boxed{
 \cT_{N,\Mean}^M(H)(u)
 =\cR_{N,\Mean}^M(u;H)
  -\cR_{N,\Mean}^M(0;H).
 }
 \label{eq:section7-averaged-transverse-profile}
\end{equation}
\end{definition}

The notation separates the scalar averaged defect $\cR_{N,\Mean}^M(0;H)$ from the transverse profile $\cT_{N,\Mean}^M(H)$.  Only the second component is used to trivialize the torsor.

\begin{proposition}[Equivariance of averaged profiles]
\label{prop:section7-averaged-profile-equivariance}
For every regular orbitwise mean, every $H\in\Lift(L)$, and every $P\in\Per_0$,
\begin{equation}
 \boxed{
 \cT_{N,\Mean}^M(H+P)
 =\cT_{N,\Mean}^M(H)+\res P.
 }
 \label{eq:section7-averaged-profile-equivariance}
\end{equation}
Moreover,
\begin{equation}
 \cR_{N,\Mean}^M(0;H+P)
 =\cR_{N,\Mean}^M(0;H).
 \label{eq:section7-scalar-defect-gauge-invariant}
\end{equation}
\end{proposition}

\begin{proof}
By \eqref{eq:section7-mean-fixes-periodic},
\[
 \cR_{N,\Mean}^M(u;H+P)
 =\cR_{N,\Mean}^M(u;H)+P(u).
\]
At $u=0$, the added term is $P(0)=0$.  Subtracting the base value proves the first assertion, and the second follows directly.
\end{proof}

Thus the longitudinal scalar component is gauge invariant, while the transverse component is gauge equivariant.  This is the canonical splitting of the asymptotic data for the present problem.

The endpoint identity also survives averaging:
\begin{equation}
 \boxed{
 \cT_{N,\Mean}^M(H)(1)
 =\sum_{n\ge1}\omega_{N,n}
   \bigl(L(n)-\Deltaop M(n)\bigr).
 }
 \label{eq:section7-averaged-endpoint-identity}
\end{equation}
It follows that no periodic gauge transformation can repair a persistent averaged mismatch in the recurrence model, because every $P\in\Per_0$ vanishes at $u=1$.

\subsection{Stirling data and admissibility}
\label{subsec:section7-stirling-data-admissibility}

The profile sequence itself is the appropriate torsor-valued version of asymptotic data.

\begin{definition}[Orbitwise Stirling datum]
\label{def:section7-orbitwise-stirling-datum}
Fix a comparison model $M$ and a regular orbitwise mean $\Mean$.  The associated \emph{orbitwise Stirling datum} is the map
\begin{equation}
 \St_{M,\Mean}:\Lift(L)
 \longrightarrow (\cX_0)^{\NN},
 \qquad
 H\longmapsto
 \bigl(\cT_{N,\Mean}^M(H)\bigr)_{N\ge1}.
 \label{eq:section7-stirling-datum-map}
\end{equation}
The target is acted on by $\Per_0$ through diagonal translation by $\res P$.
\end{definition}

Proposition~\ref{prop:section7-averaged-profile-equivariance} says precisely that
\begin{equation}
 \St_{M,\Mean}(H+P)
 =\St_{M,\Mean}(H)+\res P.
 \label{eq:section7-stirling-datum-equivariance}
\end{equation}
In other words, $\St_{M,\Mean}$ is an equivariant map of $\Per_0$-sets, where the profile sequence space carries diagonal translation by $\res P$.

\begin{definition}[Stirling-admissible lift]
\label{def:section7-stirling-admissible-lift}
A lift $H\in\Lift(L)$ is \emph{$(M,\Mean)$-Stirling admissible} if
\begin{equation}
 \boxed{
 \left\|\cT_{N,\Mean}^M(H)\right\|_{C([0,1])}
 \longrightarrow 0.
 }
 \label{eq:section7-stirling-admissibility}
\end{equation}
Equivalently,
\begin{equation}
 \sup_{0\le u\le1}
 \left|
 \cR_{N,\Mean}^M(u;H)
 -\cR_{N,\Mean}^M(0;H)
 \right|
 \longrightarrow0.
 \label{eq:section7-stirling-admissibility-expanded}
\end{equation}
\end{definition}

The condition asks the averaged defect to become asymptotically constant across each translation orbit.  It does not require the scalar constant
\[
 \cR_{N,\Mean}^M(0;H)
\]
to tend to zero, or even to converge.  This deliberate asymmetry allows the same analytic lift to coexist with a nontrivial unsmoothed integer fluctuation sequence.

\begin{remark}[Restriction of the gauge action]
\label{rem:section7-restriction-gauge-action}
If two lifts differ by $P\in\Per_0$, their Stirling data differ by the fixed profile $\res P$ for every $N$.  Since restriction to $[0,1]$ is injective on $\Per_0$, a limiting profile determines at most one periodic gauge correction.  The next section turns this observation into the abstract orbitwise rigidity theorem.
\end{remark}

\subsection{Selection of the affine comparison coefficient}
\label{subsec:section7-affine-comparison-coefficient}

The periodic ambiguity of the lift is not the only normalization present in a Stirling problem.  The comparison model itself may contain an undetermined affine coefficient.  Orbitwise profiles separate these two ambiguities cleanly.

Fix a base model $M_0\in\cO(\Omega_\sigma)$ and consider
\begin{equation}
 M_c(z)=M_0(z)+c(z-1),
 \qquad c\in\CC.
 \label{eq:section7-affine-model-family}
\end{equation}

\begin{proposition}[Affine covariance]
\label{prop:section7-affine-covariance}
For every regular orbitwise mean,
\begin{equation}
 \boxed{
 \cT_{N,\Mean}^{M_c}(H)(u)
 =\cT_{N,\Mean}^{M_0}(H)(u)-cu.
 }
 \label{eq:section7-affine-covariance}
\end{equation}
In particular,
\begin{equation}
 \cT_{N,\Mean}^{M_c}(H)(1)
 =\cT_{N,\Mean}^{M_0}(H)(1)-c.
 \label{eq:section7-affine-endpoint-shift}
\end{equation}
\end{proposition}

\begin{proof}
For every $n$ and $u$,
\[
 M_c(n+u)-M_0(n+u)=c(n+u-1).
\]
The difference between the averaged values at $u$ and $0$ is therefore $cu$, independently of the weights.
\end{proof}

\begin{corollary}[Uniqueness of the affine slope]
\label{cor:section7-unique-affine-slope}
For a fixed lift torsor and regular mean, at most one value of $c$ can admit an $(M_c,\Mean)$-Stirling-admissible lift.
\end{corollary}

\begin{proof}
Suppose that $H_i$ is $(M_{c_i},\Mean)$-admissible for $i=1,2$.  Since both lifts belong to $\Lift(L)$, there is a $P\in\Per_0$ with $H_2=H_1+P$.  Propositions~\ref{prop:section7-averaged-profile-equivariance} and~\ref{prop:section7-affine-covariance} give
\[
 \cT_{N,\Mean}^{M_{c_2}}(H_2)
 -\cT_{N,\Mean}^{M_{c_1}}(H_1)
 =\res P+(c_1-c_2)u.
\]
The left-hand side tends uniformly to zero.  Hence $P(u)+(c_1-c_2)u=0$ on $[0,1]$.  Evaluating at $u=1$ and using $P(1)=0$ gives $c_1=c_2$.
\end{proof}

This is the first stage of the prime normalization.  The arithmetic weak Stirling law identifies the coefficient $C_{\PP}$ in the family
\begin{equation}
 M_c(z)=\log\Gamma(z)+c(z-1).
 \label{eq:section7-prime-affine-model-family}
\end{equation}
Once $c=C_{\PP}$ has been fixed, the remaining ambiguity is purely periodic and belongs to the gamma-lift torsor.

\subsection{Ces\`aro orbitwise profiles}
\label{subsec:section7-cesaro-orbitwise-profiles}

We now specify the averaging schemes used later.  For a function $f$ on a right half-plane, define the first orbitwise Ces\`aro mean by
\begin{equation}
 (\Ces_N^{(1)}f)(u)
 =\frac1N\sum_{n=1}^N f(n+u),
 \label{eq:section7-first-cesaro-mean}
\end{equation}
and the triangular second mean by
\begin{equation}
 \begin{aligned}
 (\Ces_N^{(2)}f)(u)
 &=\frac{2}{N(N+1)}
   \sum_{m=1}^N\sum_{n=1}^m f(n+u)\\
 &=\frac{2}{N(N+1)}
   \sum_{n=1}^N(N+1-n)f(n+u).
 \end{aligned}
 \label{eq:section7-second-cesaro-mean}
\end{equation}
Both are regular orbitwise means.  Their weights satisfy
\begin{equation}
 \omega_{N,n}^{(1)}=\frac1N\mathbf1_{n\le N},
 \qquad
 \omega_{N,n}^{(2)}
 =\frac{2(N+1-n)}{N(N+1)}\mathbf1_{n\le N}.
 \label{eq:section7-cesaro-weights}
\end{equation}

For $r\in\{1,2\}$, set
\begin{equation}
 \cR_N^{(r),M}(u;H)
 =\Ces_N^{(r)}(H-M)(u)
 \label{eq:section7-cesaro-averaged-defect}
\end{equation}
and
\begin{equation}
 \cT_N^{(r),M}(H)(u)
 =\cR_N^{(r),M}(u;H)
  -\cR_N^{(r),M}(0;H).
 \label{eq:section7-cesaro-transverse-profile}
\end{equation}
The completed prime lift will be selected using the second profile
\begin{equation}
 \boxed{
 \sup_{0\le u\le1}
 \left|
 \cR_N^{(2),M_{\PP}}(u;H)
 -\cR_N^{(2),M_{\PP}}(0;H)
 \right|
 \longrightarrow0.
 }
 \label{eq:section7-prime-second-cesaro-condition}
\end{equation}

The second averaging is an analytic regularization, not an alteration of the factorial data.  In particular, it does not replace the integer defect by a smoothed sequence in the later variance problem.

\begin{lemma}[Regularity of orbitwise means]
\label{lem:section7-toeplitz-regularity}
Let $\Mean$ be a regular orbitwise mean and suppose that $f$ satisfies
\begin{equation}
 f(n+u)-f(n)\longrightarrow\Phi(u)
 \label{eq:section7-uniform-orbit-limit}
\end{equation}
uniformly for $u\in[0,1]$.  Then
\begin{equation}
 \Mean_Nf(u)-\Mean_Nf(0)
 \longrightarrow\Phi(u)
 \label{eq:section7-mean-preserves-orbit-limit}
\end{equation}
uniformly on $[0,1]$.
\end{lemma}

\begin{proof}
Put
\[
 g_n(u)=f(n+u)-f(n)-\Phi(u).
\]
The sequence $g_n$ tends uniformly to zero.  Given $\varepsilon>0$, choose $R$ so that $\|g_n\|_\infty<\varepsilon$ for $n>R$.  Then
\[
 \left\|\sum_n\omega_{N,n}g_n\right\|_\infty
 \le
 \sum_{n\le R}\omega_{N,n}\|g_n\|_\infty+\varepsilon.
\]
The first term tends to zero by \eqref{eq:section7-weights-escape}.  This proves the assertion.
\end{proof}

Thus every regular profile recovers an existing uniform orbit limit.  The use of the second Ces\`aro mean becomes substantive only when the raw transverse defect fails to converge but its averaged transverse geometry does.

\subsection{Monoidal additivity of Stirling profiles}
\label{subsec:section7-monoidal-additivity}

The profile construction respects the logarithmic monoidal structure.  Let $L_1,L_2$ be liftable symbols, let $H_i\in\Lift(L_i)$, and let $M_i$ be comparison models on a common right half-plane.

\begin{proposition}[Additivity]
\label{prop:section7-profile-additivity}
For every regular orbitwise mean,
\begin{equation}
 \boxed{
 \cT_{N,\Mean}^{M_1+M_2}(H_1+H_2)
 =\cT_{N,\Mean}^{M_1}(H_1)
  +\cT_{N,\Mean}^{M_2}(H_2).
 }
 \label{eq:section7-profile-additivity}
\end{equation}
Consequently, the orbitwise Stirling datum is compatible with the contracted-sum isomorphism
\begin{equation}
 \Lift(L_1)\boxplusT\Lift(L_2)
 \cong\Lift(L_1+L_2).
 \label{eq:section7-profile-contracted-sum}
\end{equation}
\end{proposition}

\begin{proof}
Both the defect and the orbitwise mean are additive.  The transverse subtraction at $u=0$ is additive as well.  Compatibility with the contracted sum follows because anti-diagonal periodic translations cancel in $H_1+H_2$.
\end{proof}

This identity permits finite prime-layer profiles to be assembled layer by layer.  It also makes clear why the comparison model for a completed product must be chosen at the logarithmic level.

\subsection{Compatibility with filtered transition maps}
\label{subsec:section7-filtered-transition-compatibility}

Let $(L_i)_{i\in I}$ be a filtered system as in Section~6, with normalized complementary primitives
\[
 C_{i,j}\in\Lift(L_{j/i}).
\]
Suppose that the comparison models satisfy a parallel decomposition
\begin{equation}
 M_j=M_i+M_{j/i}.
 \label{eq:section7-comparison-model-decomposition}
\end{equation}
For $H_j\in\Lift(L_j)$, the transition map is
\[
 r_{i,j}(H_j)=H_j-C_{i,j}.
\]

\begin{proposition}[Naturality under complementary subtraction]
\label{prop:section7-filtered-profile-naturality}
For every regular mean,
\begin{equation}
 \boxed{
 \cT_{N,\Mean}^{M_i}(r_{i,j}H_j)
 =\cT_{N,\Mean}^{M_j}(H_j)
  -\cT_{N,\Mean}^{M_{j/i}}(C_{i,j}).
 }
 \label{eq:section7-filtered-profile-naturality}
\end{equation}
\end{proposition}

\begin{proof}
Since
\[
 r_{i,j}H_j-M_i
 =(H_j-M_j)-(C_{i,j}-M_{j/i}),
\]
the result follows from additivity.
\end{proof}

A complementary primitive whose transverse profile tends to zero is therefore asymptotically invisible to the Stirling datum.  Later, the cyclotomic primitives of finite prime layers will provide an explicit family for which this compatibility can be analyzed uniformly through the filtered completion.

\subsection{The unsmoothed integer trace}
\label{subsec:section7-unsmoothed-integer-trace}

The orbitwise regularization must not be confused with smoothing the arithmetic sequence ultimately studied in Part~III.  For a fixed model $M$, define the integer defect
\begin{equation}
 r_{H,M}(n)=H(n+1)-M(n+1).
 \label{eq:section7-integer-defect}
\end{equation}
If $H$ is replaced by $H+P$ with $P\in\Per_0$, then
\begin{equation}
 P(n+1)=P(1)=0,
\end{equation}
so
\begin{equation}
 \boxed{
 r_{H+P,M}(n)=r_{H,M}(n).
 }
 \label{eq:section7-integer-defect-gauge-invariant}
\end{equation}
The integer defect is already canonical at the level of the torsor.  Orbitwise Stirling theory selects its analytic extension, not its values on the integer lattice.

For the prime model
\begin{equation}
 M_{\PP}(z)=\log\Gamma(z)+C_{\PP}(z-1),
 \label{eq:section7-prime-comparison-model}
\end{equation}
the common integer trace is
\begin{equation}
 \boxed{
 r_{H,M_{\PP}}(n)
 =\log(n+1)!_{\PP}-\log n!-C_{\PP}n
 =\cR_{\PP}(n).
 }
 \label{eq:section7-prime-unsmoothed-trace}
\end{equation}
No Ces\`aro mean appears in \eqref{eq:section7-prime-unsmoothed-trace}.  The second-order mean in \eqref{eq:section7-prime-second-cesaro-condition} is used only to identify the distinguished analytic representative $H_{\PP}^{\mathrm{cyc}}$.

\begin{remark}[Longitudinal and transverse information]
\label{rem:section7-longitudinal-transverse-information}
The decomposition
\begin{equation}
 \cR_{N,\Mean}^M(u;H)
 =\cR_{N,\Mean}^M(0;H)
  +\cT_{N,\Mean}^M(H)(u)
 \label{eq:section7-longitudinal-transverse-decomposition}
\end{equation}
separates two logically independent questions.  The transverse term determines the analytic gauge.  The scalar term records averaged arithmetic growth and may retain substantial fluctuations.  Part~II controls the former; Part~III studies the unsmoothed integer sequence underlying the latter.
\end{remark}

\subsection{Structural summary}
\label{subsec:section7-structural-summary}

The orbitwise formalism may be summarized by
\begin{equation}
 \boxed{
 \begin{gathered}
 E_{H,M}=H-M,
 \qquad
 \cR_{N,\Mean}^M(u;H)=\Mean_NE_{H,M}(u),
 \\
 \cT_{N,\Mean}^M(H)(u)
 =\cR_{N,\Mean}^M(u;H)
  -\cR_{N,\Mean}^M(0;H),
 \\
 \cT_{N,\Mean}^M(H+P)
 =\cT_{N,\Mean}^M(H)+\res P,
 \qquad P\in\Per_0,
 \\
 \cT_{N,\Mean}^{M_c}(H)(u)
 =\cT_{N,\Mean}^{M_0}(H)(u)-cu,
 \\
 \cT_{N,\Mean}^{M_1+M_2}(H_1+H_2)
 =\cT_{N,\Mean}^{M_1}(H_1)
  +\cT_{N,\Mean}^{M_2}(H_2),
 \\
 H\text{ is Stirling admissible}
 \Longleftrightarrow
 \|\cT_{N,\Mean}^M(H)\|_\infty\to0.
 \end{gathered}
 }
 \label{eq:section7-structural-summary}
\end{equation}
The profile construction is therefore simultaneously transverse, equivariant, monoidal, and compatible with filtered complementary subtraction.  It fixes the affine comparison coefficient before addressing the periodic gauge and leaves the integer fluctuation sequence untouched.

The next section proves the abstract rigidity theorem.  Uniform convergence of the profile of one lift to the restriction of a periodic entire function is shown to trivialize the entire gamma-lift torsor.  The theorem also gives stability under changes of regular mean and under asymptotically orbit-flat perturbations of the comparison model.

% ===== Source: section8_abstract_orbitwise_rigidity.tex =====
\section{Abstract orbitwise rigidity}
\label{sec:abstract-orbitwise-rigidity}

The orbitwise profile introduced in the preceding section is equivariant under the periodic gauge group.  This elementary covariance has a strong consequence: once the comparison model and longitudinal averaging procedure have been fixed, a gamma-lift torsor contains at most one lift whose transverse profile tends to zero.  More generally, the limiting profile of any one lift records the unique periodic correction needed to reach the distinguished representative, whenever such a correction exists.

The purpose of this section is to formulate this principle independently of the prime construction.  We prove an exact separation identity, characterize existence and uniqueness of an admissible lift, and isolate the obstruction to trivializing a gamma-lift torsor.  We then treat the affine coefficient and the periodic gauge simultaneously, establish stability under changes of regular mean and comparison model, and record the monoidal and filtered functoriality needed later.

Throughout, let $L\in\cO(\CC)$ be liftable, let $\Lift(L)$ be its normalized gamma-lift torsor, and let
\begin{equation}
 \Per_0
 =\{P\in\cO(\CC):P(z+1)=P(z),\ P(1)=0\}
 \label{eq:section8-periodic-gauge-group}
\end{equation}
act by translation.  We write
\begin{equation}
 \cX_0=\{\varphi\in C([0,1];\CC):\varphi(0)=0\}
 \label{eq:section8-profile-space}
\end{equation}
with the uniform norm, and
\begin{equation}
 \res:\Per_0\longrightarrow\cX_0,
 \qquad
 (\res P)(u)=P(u).
 \label{eq:section8-restriction-map}
\end{equation}
The restriction map is injective.  Fix a holomorphic comparison model $M$ on a right half-plane and a regular orbitwise mean $\Mean$.  We abbreviate the associated transverse profile by
\begin{equation}
 \cT_N(H)=\cT_{N,\Mean}^{M}(H)\in\cX_0.
 \label{eq:section8-abbreviated-profile}
\end{equation}
The only structural fact used initially is the exact covariance law
\begin{equation}
 \cT_N(H+P)=\cT_N(H)+\res P.
 \label{eq:section8-profile-covariance}
\end{equation}

\subsection{Exact separation of lifts}
\label{subsec:section8-exact-separation}

The difference of the profiles of two lifts is independent of the orbit height, the comparison model, and the averaging weights.

\begin{proposition}[Exact separation identity]
\label{prop:section8-exact-separation}
Let $H_1,H_2\in\Lift(L)$ and put $P=H_2-H_1\in\Per_0$.  Then, for every $N$,
\begin{equation}
 \boxed{
 \cT_N(H_2)-\cT_N(H_1)=\res P.
 }
 \label{eq:section8-exact-separation}
\end{equation}
Consequently,
\begin{equation}
 \|\res(H_2-H_1)\|_\infty
 \le
 \|\cT_N(H_1)\|_\infty
 +\|\cT_N(H_2)\|_\infty.
 \label{eq:section8-quantitative-separation}
\end{equation}
\end{proposition}

\begin{proof}
The first assertion is the covariance law applied to $H_2=H_1+P$.  The second follows from the triangle inequality.
\end{proof}

The estimate is useful even before passage to the limit.  Two approximate selectors whose profiles are small at a single common height must already be close on a full fundamental interval.  Since their difference is entire and periodic, asymptotic closeness forces exact equality.

\begin{theorem}[Orbitwise rigidity]
\label{thm:section8-orbitwise-rigidity}
For fixed $(L,M,\Mean)$, there is at most one lift $H\in\Lift(L)$ satisfying
\begin{equation}
 \|\cT_N(H)\|_\infty\longrightarrow0.
 \label{eq:section8-admissibility-condition}
\end{equation}
Equivalently, the set of $(M,\Mean)$-Stirling-admissible lifts is either empty or a singleton.
\end{theorem}

\begin{proof}
If $H_1$ and $H_2$ are both admissible, Proposition~\ref{prop:section8-exact-separation} gives
\[
 \|\res(H_2-H_1)\|_\infty=0.
\]
Injectivity of $\res$ implies $H_1=H_2$.
\end{proof}

\begin{remark}[No compactness is required]
\label{rem:section8-no-compactness}
The proof uses neither normal families nor growth bounds on periodic entire functions.  Rigidity is a formal consequence of exact equivariance and the injectivity of restriction to one translation interval.  All analytic work in an application is therefore concentrated in proving existence.
\end{remark}

\subsection{Existence, reconstruction, and the obstruction class}
\label{subsec:section8-existence-reconstruction}

Fix a base lift $H_0\in\Lift(L)$.  If an admissible lift exists, the profile of $H_0$ must converge to the restriction of a periodic entire function.  This necessary condition is also sufficient.

\begin{theorem}[Reconstruction from one limiting profile]
\label{thm:section8-reconstruction}
The following are equivalent:
\begin{enumerate}[label=\textnormal{(\roman*)}]
 \item there exists an $(M,\Mean)$-Stirling-admissible lift $H_*$;
 \item the sequence $\cT_N(H_0)$ converges in $\cX_0$ to a limit
 \begin{equation}
  \Phi\in\res(\Per_0).
  \label{eq:section8-limit-in-periodic-image}
 \end{equation}
\end{enumerate}
If these conditions hold, there is a unique $P_\Phi\in\Per_0$ such that
\begin{equation}
 \res P_\Phi=\Phi,
 \label{eq:section8-periodic-extension-limit}
\end{equation}
and the unique admissible lift is
\begin{equation}
 \boxed{
 H_*=H_0-P_\Phi.
 }
 \label{eq:section8-reconstructed-selector}
\end{equation}
Moreover, for every $H=H_*+Q$ with $Q\in\Per_0$,
\begin{equation}
 \boxed{
 \cT_N(H)\longrightarrow\res Q
 }
 \label{eq:section8-all-profile-limits}
\end{equation}
uniformly on $[0,1]$.
\end{theorem}

\begin{proof}
Suppose first that $H_*$ is admissible.  Write $H_0=H_*+P$ with $P\in\Per_0$.  Exact separation gives
\[
 \cT_N(H_0)=\cT_N(H_*)+\res P,
\]
so $\cT_N(H_0)\to\res P$.  This proves (ii).

Conversely, suppose that $\cT_N(H_0)\to\Phi=\res P_\Phi$.  By covariance,
\[
 \cT_N(H_0-P_\Phi)
 =\cT_N(H_0)-\Phi\longrightarrow0.
\]
Thus $H_0-P_\Phi$ is admissible.  Its uniqueness follows from Theorem~\ref{thm:section8-orbitwise-rigidity}.  Finally, if $H=H_*+Q$, then
\[
 \cT_N(H)=\cT_N(H_*)+\res Q\longrightarrow\res Q.
\]
\end{proof}

The theorem gives a practical reconstruction formula: compute the limiting transverse profile of any convenient lift, recognize it as the restriction of an entire periodic function, and subtract that function.

\begin{corollary}[Base-lift independence of convergence]
\label{cor:section8-base-lift-independence}
If the profile of one lift converges in $\cX_0$, then the profile of every lift converges.  Their limits differ by the exact periodic restriction prescribed by the torsor action.
\end{corollary}

\begin{proof}
If $H=H_0+P$, then $\cT_N(H)=\cT_N(H_0)+\res P$ for every $N$.
\end{proof}

The image $\res(\Per_0)$ need not be treated as a closed subspace of $\cX_0$ for the preceding theorem.  It is useful, however, to record the algebraic obstruction carried by a limiting profile.

\begin{definition}[Asymptotic gauge obstruction]
\label{def:section8-gauge-obstruction}
Assume that $\cT_N(H_0)$ converges to $\Phi\in\cX_0$.  Define
\begin{equation}
 \Ob_{M,\Mean}(L)
 =\Phi+\res(\Per_0)
 \in
 \cX_0/\res(\Per_0),
 \label{eq:section8-obstruction-class}
\end{equation}
where the quotient is understood algebraically.
\end{definition}

\begin{proposition}[Intrinsic nature of the obstruction]
\label{prop:section8-obstruction-intrinsic}
The class \eqref{eq:section8-obstruction-class} is independent of the base lift.  It vanishes if and only if an admissible lift exists.
\end{proposition}

\begin{proof}
Replacing $H_0$ by $H_0+P$ replaces $\Phi$ by $\Phi+\res P$, leaving its quotient class unchanged.  Vanishing means precisely that $\Phi\in\res(\Per_0)$, and Theorem~\ref{thm:section8-reconstruction} applies.
\end{proof}

\begin{remark}[Failure of existence]
\label{rem:section8-failure-existence}
A convergent continuous profile need not extend to an entire periodic function.  Such a limit defines a nonzero obstruction class and cannot be removed by changing the lift.  Orbitwise rigidity therefore separates the existence problem into two distinct tasks: convergence in the profile space and analytic periodic extendability of the limit.
\end{remark}

\subsection{Simultaneous affine and periodic rigidity}
\label{subsec:section8-affine-periodic-rigidity}

The comparison model often belongs to an affine family
\begin{equation}
 M_c(z)=M_0(z)+c(z-1),
 \qquad c\in\CC.
 \label{eq:section8-affine-family}
\end{equation}
Let
\begin{equation}
 \iota(u)=u,
 \qquad 0\le u\le1.
 \label{eq:section8-identity-profile}
\end{equation}
The affine covariance from the preceding section reads
\begin{equation}
 \cT_{N,\Mean}^{M_c}(H)
 =\cT_{N,\Mean}^{M_0}(H)-c\iota.
 \label{eq:section8-affine-covariance}
\end{equation}
The affine direction and the periodic gauge direction are transverse.

\begin{lemma}[Directness of the affine-periodic decomposition]
\label{lem:section8-affine-periodic-directness}
One has
\begin{equation}
 \res(\Per_0)\cap\CC\iota=\{0\}.
 \label{eq:section8-direct-sum}
\end{equation}
Consequently,
\begin{equation}
 \res(\Per_0)\oplus\CC\iota
 \label{eq:section8-affine-periodic-space}
\end{equation}
is an algebraic direct sum in $\cX_0$.
\end{lemma}

\begin{proof}
Suppose that $\res P=c\iota$.  Evaluating at $u=1$ gives $0=P(1)=c$, hence $\res P=0$.  Injectivity of restriction gives $P=0$.
\end{proof}

\begin{theorem}[Joint affine-periodic selection]
\label{thm:section8-joint-selection}
Fix $H_0\in\Lift(L)$ and suppose that
\begin{equation}
 \cT_{N,\Mean}^{M_0}(H_0)\longrightarrow\Phi
 \qquad\text{in }\cX_0.
 \label{eq:section8-joint-limit}
\end{equation}
There exists a pair $(c_*,H_*)\in\CC\times\Lift(L)$ such that $H_*$ is $(M_{c_*},\Mean)$-admissible if and only if
\begin{equation}
 \Phi\in\res(\Per_0)\oplus\CC\iota.
 \label{eq:section8-joint-existence-criterion}
\end{equation}
When it exists, the pair is unique.  More precisely, write uniquely
\begin{equation}
 \Phi=\res Q+c_*\iota.
 \label{eq:section8-joint-decomposition}
\end{equation}
Then
\begin{equation}
 \boxed{
 H_*=H_0-Q.
 }
 \label{eq:section8-joint-selected-lift}
\end{equation}
The affine coefficient is recovered directly from the endpoint:
\begin{equation}
 \boxed{
 c_*=\Phi(1).
 }
 \label{eq:section8-endpoint-recovers-slope}
\end{equation}
\end{theorem}

\begin{proof}
If \eqref{eq:section8-joint-decomposition} holds, then
\[
 \cT_{N,\Mean}^{M_{c_*}}(H_0-Q)
 =\cT_{N,\Mean}^{M_0}(H_0)-\res Q-c_*\iota
 \longrightarrow0.
\]
Conversely, suppose that $H_*=H_0+P$ is $(M_c,\Mean)$-admissible.  Then
\[
 \cT_{N,\Mean}^{M_0}(H_0)+\res P-c\iota\longrightarrow0,
\]
so $\Phi=-\res P+c\iota$ belongs to the direct sum.  Uniqueness follows from Lemma~\ref{lem:section8-affine-periodic-directness} and orbitwise rigidity.  Finally, $\res Q(1)=0$, so evaluating \eqref{eq:section8-joint-decomposition} at $u=1$ gives $c_*=\Phi(1)$.
\end{proof}

The theorem formalizes the two-stage Stirling normalization.  The endpoint determines the affine slope; after subtracting the linear profile $\Phi(1)u$, the remaining transverse limit must extend to a normalized entire periodic function, which determines the gamma lift.

\begin{corollary}[Two-stage reconstruction]
\label{cor:section8-two-stage-reconstruction}
Under the hypotheses of Theorem~\ref{thm:section8-joint-selection}, put
\begin{equation}
 \Phi^{\mathrm{per}}(u)=\Phi(u)-\Phi(1)u.
 \label{eq:section8-periodic-residual-profile}
\end{equation}
Then a jointly admissible pair exists if and only if
\begin{equation}
 \Phi^{\mathrm{per}}\in\res(\Per_0).
 \label{eq:section8-periodic-residual-criterion}
\end{equation}
In that case the slope is $c_*=\Phi(1)$ and the periodic correction is the unique $Q\in\Per_0$ satisfying $\res Q=\Phi^{\mathrm{per}}$.
\end{corollary}

\subsection{Stability under the averaging procedure}
\label{subsec:section8-stability-means}

Rigidity holds for every regular mean separately.  To compare two means, it is enough to compare the corresponding profile sequences.

\begin{definition}[Asymptotically equivalent means for a defect]
\label{def:section8-equivalent-means}
Let $\Mean$ and $\Mean'$ be regular orbitwise means.  They are \emph{asymptotically equivalent for $(H,M)$} if
\begin{equation}
 \left\|
 \cT_{N,\Mean}^{M}(H)
 -\cT_{N,\Mean'}^{M}(H)
 \right\|_\infty
 \longrightarrow0.
 \label{eq:section8-equivalent-means}
\end{equation}
\end{definition}

\begin{proposition}[Mean stability]
\label{prop:section8-mean-stability}
Suppose that $\Mean$ and $\Mean'$ are asymptotically equivalent for one lift $H_0\in\Lift(L)$.  Then they are asymptotically equivalent for every lift.  Moreover, they have the same admissible lift whenever either selector exists.
\end{proposition}

\begin{proof}
For $H=H_0+P$, both profiles acquire the same term $\res P$, so their difference is unchanged.  If $H_*$ is admissible for one mean, asymptotic equivalence implies admissibility for the other.  Uniqueness then identifies the selectors.
\end{proof}

A useful sufficient condition is the existence of a raw orbitwise limit.  If
\begin{equation}
 H(n+u)-M(n+u)-H(n)+M(n)
 \longrightarrow\Phi(u)
 \label{eq:section8-raw-uniform-limit}
\end{equation}
uniformly for $u\in[0,1]$, every regular orbitwise mean has profile limit $\Phi$.  Thus the selector, when it exists, is independent of the regular summability method.

The first and second Ces\`aro profiles satisfy a one-way consistency relation even without a raw limit.

\begin{proposition}[Ces\`aro consistency]
\label{prop:section8-cesaro-consistency}
If the first Ces\`aro profiles satisfy
\begin{equation}
 \cT_N^{(1),M}(H)\longrightarrow\Phi
 \qquad\text{in }\cX_0,
 \label{eq:section8-first-cesaro-limit}
\end{equation}
then the triangular second Ces\`aro profiles satisfy
\begin{equation}
 \cT_N^{(2),M}(H)\longrightarrow\Phi.
 \label{eq:section8-second-cesaro-same-limit}
\end{equation}
Consequently, first-order admissibility implies second-order admissibility, and the selected lifts agree.
\end{proposition}

\begin{proof}
The triangular second mean is a regular weighted mean of the first means:
\begin{equation}
 \cT_N^{(2),M}(H)
 =\frac{2}{N(N+1)}
  \sum_{m=1}^N m\,\cT_m^{(1),M}(H).
 \label{eq:section8-second-as-mean-first}
\end{equation}
The weights are nonnegative, sum to one, and escape every fixed initial segment.  Toeplitz regularity therefore preserves the limit $\Phi$.
\end{proof}

\begin{remark}[Second-order selection is weaker]
\label{rem:section8-second-weaker}
The converse need not hold: triangular averaging may converge when the first Ces\`aro profiles do not.  This additional regularization is precisely why the second profile is suited to the completed prime system.  Rigidity is not weakened, because the gauge covariance remains exact under every normalized longitudinal average.
\end{remark}

\subsection{Stability under changes of comparison model}
\label{subsec:section8-stability-models}

Let $M$ and $\widetilde M=M+Q$ be two comparison models on a common right half-plane.  Define the transverse profile of the perturbation by
\begin{equation}
 \cU_{N,\Mean}(Q)(u)
 =\Mean_NQ(u)-\Mean_NQ(0).
 \label{eq:section8-model-perturbation-profile}
\end{equation}
Then
\begin{equation}
 \cT_{N,\Mean}^{\widetilde M}(H)
 =\cT_{N,\Mean}^{M}(H)-\cU_{N,\Mean}(Q).
 \label{eq:section8-model-change-identity}
\end{equation}

\begin{definition}[Orbit-flat perturbation]
\label{def:section8-orbit-flat-perturbation}
The perturbation $Q$ is \emph{$\Mean$-orbit-flat} if
\begin{equation}
 \|\cU_{N,\Mean}(Q)\|_\infty\longrightarrow0.
 \label{eq:section8-orbit-flat-condition}
\end{equation}
\end{definition}

\begin{proposition}[Model stability]
\label{prop:section8-model-stability}
If $Q$ is $\Mean$-orbit-flat, then $(M,\Mean)$ and $(M+Q,\Mean)$ have the same admissible lift.  More generally, if
\begin{equation}
 \cU_{N,\Mean}(Q)\longrightarrow\res P
 \label{eq:section8-model-periodic-limit}
\end{equation}
for some $P\in\Per_0$, then their selectors, when they exist, differ by $P$:
\begin{equation}
 \Sel_{M+Q,\Mean}(L)
 =\Sel_{M,\Mean}(L)+P.
 \label{eq:section8-model-selector-shift}
\end{equation}
\end{proposition}

\begin{proof}
The first statement follows immediately from \eqref{eq:section8-model-change-identity}.  For the second, let $H_*$ be admissible relative to $M$.  Then
\[
 \cT_{N,\Mean}^{M+Q}(H_*+P)
 =\cT_{N,\Mean}^{M}(H_*)+\res P-\cU_{N,\Mean}(Q)
 \longrightarrow0.
\]
Uniqueness identifies the selector.
\end{proof}

This motivates an equivalence relation on comparison models:
\begin{equation}
 M\sim_\Mean\widetilde M
 \quad\Longleftrightarrow\quad
 \widetilde M-M\text{ is }\Mean\text{-orbit-flat}.
 \label{eq:section8-model-equivalence}
\end{equation}
The selected lift depends only on the $\sim_\Mean$-class of the model.  Affine changes are deliberately excluded from this equivalence: for $Q(z)=c(z-1)$ one has $\cU_{N,\Mean}(Q)=c\iota$, so the endpoint records the nonzero change in slope.

\subsection{Monoidal rigidity and filtered coherence}
\label{subsec:section8-monoidal-filtered}

The additive profile law turns selection into a monoidal construction wherever it is defined.

\begin{proposition}[Monoidal closure of admissible lifts]
\label{prop:section8-monoidal-closure}
Let $L_1,L_2$ be liftable symbols with comparison models $M_1,M_2$, and use the same regular mean $\Mean$.  If $H_i\in\Lift(L_i)$ is $(M_i,\Mean)$-admissible for $i=1,2$, then
\begin{equation}
 \boxed{
 H_1+H_2
 }
 \label{eq:section8-monoidal-selected-sum}
\end{equation}
is $(M_1+M_2,\Mean)$-admissible in $\Lift(L_1+L_2)$.  Hence, under the contracted-sum identification,
\begin{equation}
 \Sel_{M_1,\Mean}(L_1)\boxplusT
 \Sel_{M_2,\Mean}(L_2)
 =
 \Sel_{M_1+M_2,\Mean}(L_1+L_2).
 \label{eq:section8-selector-monoidal}
\end{equation}
\end{proposition}

\begin{proof}
Additivity of profiles gives
\[
 \cT_{N,\Mean}^{M_1+M_2}(H_1+H_2)
 =\cT_{N,\Mean}^{M_1}(H_1)
  +\cT_{N,\Mean}^{M_2}(H_2),
\]
which tends uniformly to zero.  Uniqueness gives the final identity.
\end{proof}

There is also a cancellation principle.  If $H_{12}$ is selected for $L_1+L_2$ relative to $M_1+M_2$ and $H_1$ is selected for $L_1$ relative to $M_1$, then $H_{12}-H_1$ is selected for $L_2$ relative to $M_2$, provided it belongs to the normalized lift torsor of $L_2$.  Thus any two sides of a monoidal selection identity determine the third.

Now let $(L_i)_{i\in I}$ be a filtered system with complementary symbols
\begin{equation}
 L_j=L_i+L_{j/i}
 \qquad(i\le j),
 \label{eq:section8-filtered-symbol-decomposition}
\end{equation}
and models
\begin{equation}
 M_j=M_i+M_{j/i}.
 \label{eq:section8-filtered-model-decomposition}
\end{equation}
Suppose that $C_{i,j}\in\Lift(L_{j/i})$ is a normalized complementary primitive and that the transition map is
\begin{equation}
 r_{i,j}(H_j)=H_j-C_{i,j}.
 \label{eq:section8-filtered-transition}
\end{equation}

\begin{proposition}[Filtered coherence of selectors]
\label{prop:section8-filtered-coherence}
Assume that $C_{i,j}$ is $(M_{j/i},\Mean)$-admissible and that $H_j$ is $(M_j,\Mean)$-admissible.  Then
\begin{equation}
 r_{i,j}(H_j)
 \label{eq:section8-transition-selected}
\end{equation}
is $(M_i,\Mean)$-admissible.  Conversely, if $H_i$ and $C_{i,j}$ are admissible for their respective models, then $H_i+C_{i,j}$ is admissible at level $j$.
\end{proposition}

\begin{proof}
Using $H_j=r_{i,j}(H_j)+C_{i,j}$ and the model decomposition, profile additivity gives
\[
 \cT_{N,\Mean}^{M_j}(H_j)
 =\cT_{N,\Mean}^{M_i}(r_{i,j}(H_j))
  +\cT_{N,\Mean}^{M_{j/i}}(C_{i,j}).
\]
If the left-hand side and the complementary term tend to zero, so does the remaining term.  The converse is the monoidal closure proposition.
\end{proof}

\begin{corollary}[Automatic cocycle compatibility]
\label{cor:section8-automatic-cocycle}
Suppose admissible lifts exist at every level of a filtered system and admissible complementary primitives exist for every transition.  Then the selected lifts form a compatible section of the inverse system of torsors.  All cocycle identities hold automatically by uniqueness.
\end{corollary}

\begin{proof}
Each transition sends the selected lift at level $j$ to an admissible lift at level $i$, which must be the selected lift there.  For three indices $i\le j\le k$, both composites produce an admissible lift at level $i$; uniqueness forces equality.
\end{proof}

\subsection{A quantitative form of rigidity}
\label{subsec:section8-quantitative-rigidity}

Although the final selection theorem is qualitative, exact separation yields a useful quantitative statement for finite approximations.

\begin{proposition}[Finite-height rigidity]
\label{prop:section8-finite-height-rigidity}
Let $H_1,H_2\in\Lift(L)$.  If for some $N$
\begin{equation}
 \|\cT_N(H_i)\|_\infty\le\varepsilon_i
 \qquad(i=1,2),
 \label{eq:section8-approximate-admissibility}
\end{equation}
then
\begin{equation}
 \boxed{
 \|\res(H_2-H_1)\|_\infty
 \le\varepsilon_1+\varepsilon_2.
 }
 \label{eq:section8-finite-height-bound}
\end{equation}
In particular, if a sequence of candidate lifts $H^{(r)}\in\Lift(L)$ satisfies
\begin{equation}
 \sup_{r,s\ge R}
 \|\cT_{N_R}(H^{(r)})\|_\infty
 +\|\cT_{N_R}(H^{(s)})\|_\infty
 \longrightarrow0,
 \label{eq:section8-candidate-cauchy-condition}
\end{equation}
then their periodic differences converge uniformly to zero on $[0,1]$.
\end{proposition}

\begin{proof}
The first assertion is \eqref{eq:section8-quantitative-separation}.  The second follows by applying it to each pair $H^{(r)},H^{(s)}$ at the common height $N_R$.
\end{proof}

\begin{remark}[What the estimate does not provide]
\label{rem:section8-no-entire-stability}
Uniform smallness on a real interval does not by itself control an entire periodic function on the whole plane without additional growth information.  The estimate is therefore a rigidity statement in the profile norm, not a global stability theorem in a space of entire functions.  Later cyclotomic lifts come with explicit Fourier expansions, so their global convergence can be handled separately.
\end{remark}

\subsection{The abstract selection principle}
\label{subsec:section8-abstract-selection-principle}

The preceding results can be compressed into a single theorem that will be applied to the prime system.

\begin{theorem}[Abstract orbitwise selection principle]
\label{thm:section8-abstract-selection-principle}
Let $L$ be a liftable entire recurrence symbol, let $M_0$ be a holomorphic comparison model, and let $\Mean$ be a regular orbitwise mean.  Fix $H_0\in\Lift(L)$ and suppose that
\begin{equation}
 \cT_{N,\Mean}^{M_0}(H_0)\longrightarrow\Phi
 \qquad\text{uniformly on }[0,1].
 \label{eq:section8-abstract-limit}
\end{equation}
Then:
\begin{enumerate}[label=\textnormal{(\roman*)}]
 \item a jointly normalized pair $(c_*,H_*)$ for the affine family
 \begin{equation}
  M_c(z)=M_0(z)+c(z-1)
  \label{eq:section8-abstract-affine-family}
 \end{equation}
 exists if and only if
 \begin{equation}
  \Phi(u)-\Phi(1)u\in\res(\Per_0);
  \label{eq:section8-abstract-periodic-criterion}
 \end{equation}
 \item when it exists, the pair is unique and satisfies
 \begin{equation}
  c_*=\Phi(1),
  \qquad
  H_*=H_0-Q,
  \label{eq:section8-abstract-reconstruction}
 \end{equation}
 where $Q$ is the unique normalized entire periodic function with
 \begin{equation}
  Q(u)=\Phi(u)-\Phi(1)u
  \qquad(0\le u\le1);
  \label{eq:section8-abstract-periodic-extension}
 \end{equation}
 \item for every other lift $H=H_*+P$,
 \begin{equation}
  \cT_{N,\Mean}^{M_{c_*}}(H)
  \longrightarrow\res P;
  \label{eq:section8-abstract-all-limits}
 \end{equation}
 \item the selected pair is unchanged under orbit-flat perturbations of the comparison model and under asymptotically equivalent regular means;
 \item selection is additive under contracted monoidal sums and coherent under filtered complementary subtraction.
\end{enumerate}
\end{theorem}

\begin{proof}
Parts (i)--(iii) are Theorem~\ref{thm:section8-joint-selection} and Corollary~\ref{cor:section8-two-stage-reconstruction}.  Part (iv) follows from Propositions~\ref{prop:section8-mean-stability} and~\ref{prop:section8-model-stability}.  Part (v) follows from Propositions~\ref{prop:section8-monoidal-closure} and~\ref{prop:section8-filtered-coherence}.
\end{proof}

The theorem separates formal rigidity from analytic existence.  Once a single profile limit has been computed, the affine coefficient is read from its endpoint, the periodic gauge is read from its transverse residual, and all compatibility statements follow from uniqueness.  Conceptually, this theorem plays the role of a Bohr--Mollerup principle for the present recurrence: the recurrence specifies a torsor of candidates, while orbitwise normalization selects the distinguished point.  The analogy concerns the function of the characterization theorem, not its hypotheses; no log-convexity assertion is being imported into the complex-entire setting.

\subsection{Specialization to the prime factorial calculus}
\label{subsec:section8-prime-specialization}

For the completed prime calculus, the affine family is
\begin{equation}
 M_c(z)=\log\Gamma(z)+c(z-1).
 \label{eq:section8-prime-affine-family}
\end{equation}
The weak arithmetic Stirling law determines
\begin{equation}
 c=C_{\PP}
 =\sum_p\frac{\log p}{(p-1)^2}.
 \label{eq:section8-prime-stirling-constant}
\end{equation}
The remaining task in Part~II is therefore an existence problem for the periodic gauge: construct one normalized entire lift of the completed prime recurrence and prove that its second Ces\`aro transverse profile relative to
\begin{equation}
 M_{\PP}(z)
 =\log\Gamma(z)+C_{\PP}(z-1)
 \label{eq:section8-prime-comparison-model}
\end{equation}
converges to zero, or more generally to the restriction of an explicitly removable periodic entire function.

The next section constructs the required finite-layer primitives from cyclotomic data.  Their Fourier expansions make the periodic ambiguity explicit and provide compatible complementary primitives for the filtered layer system.  Section~10 then passes to the completed prime object and invokes Theorem~\ref{thm:section8-abstract-selection-principle} to obtain the distinguished prime gamma lift.

\subsection{Structural summary}
\label{subsec:section8-structural-summary}

The abstract rigidity mechanism is summarized by
\begin{equation}
 \boxed{
 \begin{gathered}
 \cT_N(H+P)=\cT_N(H)+\res P,
 \qquad P\in\Per_0,
 \\
 \cT_N(H_2)-\cT_N(H_1)=\res(H_2-H_1),
 \\
 \#\{H\in\Lift(L):\cT_N(H)\to0\}\le1,
 \\
 \cT_N(H_0)\to\Phi\in\res(\Per_0)
 \Longleftrightarrow
 H_0-P_\Phi\text{ is admissible},
 \\
 \Phi=\res Q+c\iota
 \Longleftrightarrow
 (c,H_0-Q)\text{ is jointly admissible},
 \\
 c=\Phi(1),
 \qquad
 \res Q=\Phi-\Phi(1)\iota,
 \\
 \Sel(L_1+L_2)=\Sel(L_1)+\Sel(L_2),
 \end{gathered}
 }
 \label{eq:section8-structural-summary}
\end{equation}
whenever the displayed selectors exist.  Completion may preserve an entire torsor of gamma lifts, but an orbitwise Stirling profile can trivialize that torsor at most once.  The cyclotomic construction now supplies the candidate whose existence remains to be proved.

% ===== Source: section9_cyclotomic_primitives_prime_layers.tex =====
\section{Cyclotomic primitives for prime layers}
\label{sec:cyclotomic-primitives-prime-layers}

The abstract rigidity theorem of the preceding section reduces the selection problem to existence: one must construct a compatible lift of the prime recurrence whose orbitwise Stirling profile has the required limit.  The first difficulty occurs one level earlier.  The divisibility function
\[
 n\longmapsto \mathbf 1_{d\mid n}
\]
has many entire interpolants, and the most obvious finite Fourier interpolants do not decay as the conductor $d$ tends to infinity.  Their direct sum over all prime layers therefore diverges away from the integers.

This section resolves that difficulty by a universal cyclotomic renormalization.  A symmetric finite Fourier quadrature gives the exact divisibility kernel at the integers.  Its continuum limit is the cardinal sine function
\[
 \sinc(z)=\frac{\sin \pi z}{\pi z},
 \qquad \sinc(0)=1,
\]
which vanishes at every nonzero integer.  Subtracting this universal null mode preserves all positive-integer recurrence data and improves the conductor decay from order one to order $d^{-2}$.  The corresponding normalized primitives are obtained by applying the same quadrature to an explicit difference kernel.  Their continuum limit again vanishes on the positive integers, and its subtraction produces normally summable prime-layer primitives.

The outcome is a canonical family
\[
 L_{p,j}^{\mathrm{cyc}},\qquad H_{p,j}^{\mathrm{cyc}},
\]
indexed by the prime layers $(p,j)$, with
\[
 H_{p,j}^{\mathrm{cyc}}(z+1)-H_{p,j}^{\mathrm{cyc}}(z)
 =L_{p,j}^{\mathrm{cyc}}(z),
\]
whose integer traces are exactly the atomic Bhargava layer.  These primitives are additive under the monoidal layer law, satisfy the complementary cocycle identities of Section~6, and are normally summable over the full prime-layer system.

Throughout, write
\begin{equation}
 \e(t)=e^{2\pi i t}.
 \label{eq:section9-exponential-notation}
\end{equation}

\subsection{Symmetric cyclotomic quadrature}
\label{subsec:section9-symmetric-quadrature}

For $d\ge 1$, define a quadrature functional on continuous functions on $[-1/2,1/2]$ by
\begin{equation}
 \cQ_d f
 =
 \begin{cases}
 \displaystyle
 \frac1d\sum_{r=-(d-1)/2}^{(d-1)/2} f(r/d),
 & d\ \text{odd},\\[3mm]
 \displaystyle
 \frac1d\left(
 \frac12f(-1/2)
 +\sum_{r=-d/2+1}^{d/2-1}f(r/d)
 +\frac12f(1/2)
 \right),
 & d\ \text{even}.
 \end{cases}
 \label{eq:section9-cyclotomic-quadrature}
\end{equation}
For odd $d$ this is the composite midpoint rule with $d$ subintervals; for even $d$ it is the composite trapezoidal rule with mesh $1/d$.  Let
\begin{equation}
 \cQ_\infty f=\int_{-1/2}^{1/2}f(t)\,dt.
 \label{eq:section9-continuum-quadrature}
\end{equation}

Apply these functionals to the exponential $t\mapsto\e(tz)$.

\begin{definition}[Raw cyclotomic divisibility kernel]
\label{def:section9-raw-divisibility-kernel}
For $d\ge1$, put
\begin{equation}
 \Theta_d(z)=\cQ_d\bigl(t\mapsto\e(tz)\bigr).
 \label{eq:section9-raw-divisibility-kernel}
\end{equation}
The continuum kernel is
\begin{equation}
 \Theta_\infty(z)
 =\cQ_\infty\bigl(t\mapsto\e(tz)\bigr)
 =\sinc(z).
 \label{eq:section9-continuum-divisibility-kernel}
\end{equation}
\end{definition}

The finite kernel is a real-even trigonometric polynomial of minimal symmetric bandwidth.  A geometric-sum calculation gives a useful closed form.

\begin{proposition}[Closed cyclotomic form]
\label{prop:section9-closed-cyclotomic-form}
After filling removable singularities, $\Theta_d$ is entire, real-valued on $\mathbb R$, even, and $d$-periodic.  Explicitly,
\begin{equation}
 \boxed{
 \Theta_d(z)
 =
 \begin{cases}
 \displaystyle
 \frac{\sin(\pi z)}{d\sin(\pi z/d)},
 & d\ \text{odd},\\[3mm]
 \displaystyle
 \frac{\sin(\pi z)}{d\tan(\pi z/d)},
 & d\ \text{even}.
 \end{cases}
 }
 \label{eq:section9-closed-cyclotomic-kernel}
\end{equation}
Moreover, for every $n\in\ZZ$,
\begin{equation}
 \boxed{
 \Theta_d(n)=\mathbf 1_{d\mid n}.
 }
 \label{eq:section9-exact-divisibility-interpolation}
\end{equation}
\end{proposition}

\begin{proof}
For odd $d$, the frequencies in \eqref{eq:section9-cyclotomic-quadrature} form the symmetric set
\[
 -\frac{d-1}{2},\ldots,\frac{d-1}{2},
\]
and the usual Dirichlet-kernel identity gives the first formula.  For even $d$, the two half-weight endpoint frequencies combine into a cosine term, yielding the second formula.  Symmetry of the frequency weights gives evenness and reality on the real axis.

If $n\in\ZZ$, the odd rule is the average of $d$ consecutive characters of the cyclic group $\ZZ/d\ZZ$.  In the even rule the two endpoint values agree, since
\[
 \e(-n/2)=\e(n/2)=(-1)^n,
\]
so the half weights combine into one complete residue frequency.  Character orthogonality gives \eqref{eq:section9-exact-divisibility-interpolation}.
\end{proof}

\begin{remark}[Why a real symmetric kernel is chosen]
\label{rem:section9-real-symmetric-choice}
A one-sided discrete Fourier sum also interpolates $\mathbf1_{d\mid n}$, but it is generally complex on the real axis.  The symmetric rule is invariant under complex conjugation and is the unique interpolation in the corresponding real-even minimal-band space.  This reality property propagates to the logarithmic gamma lifts constructed below.
\end{remark}

The raw kernels do not tend to zero.  Rather,
\begin{equation}
 \Theta_d(z)\longrightarrow\sinc(z)
 \label{eq:section9-raw-kernel-continuum-limit}
\end{equation}
locally uniformly as $d\to\infty$.  Since $\sinc(n)=0$ for every positive integer $n$, the continuum mode is invisible to the discrete recurrence skeleton.  This leads to the required renormalization.

\begin{definition}[Renormalized cyclotomic symbol]
\label{def:section9-renormalized-symbol}
For $d\ge1$, define
\begin{equation}
 \boxed{
 \kappa_d(z)=\Theta_d(z)-\sinc(z).
 }
 \label{eq:section9-renormalized-cyclotomic-symbol}
\end{equation}
\end{definition}

Then
\begin{equation}
 \kappa_d(n)=\mathbf1_{d\mid n}
 \qquad(n\ge1),
 \label{eq:section9-renormalized-symbol-integers}
\end{equation}
but $\kappa_d$ has the conductor decay needed for completion.

\subsection{Normalized primitives by cyclotomic quadrature}
\label{subsec:section9-normalized-primitives}

The same quadrature has a canonical primitive-level formulation.  For $z\in\CC$ and $t\in[-1/2,1/2]$, define
\begin{equation}
 K_z(t)
 =
 \begin{cases}
 \displaystyle
 \frac{\e(tz)-\e(t)}{\e(t)-1},
 & t\ne0,\\[3mm]
 z-1,
 & t=0.
 \end{cases}
 \label{eq:section9-difference-kernel}
\end{equation}
The singularity at $t=0$ is removable.  For fixed $t$, the function $z\mapsto K_z(t)$ is entire, and for fixed $z$ it is analytic in $t$ on a neighborhood of the integration interval.  The essential identity is
\begin{equation}
 K_{z+1}(t)-K_z(t)=\e(tz).
 \label{eq:section9-kernel-difference-identity}
\end{equation}

\begin{definition}[Raw and continuum primitives]
\label{def:section9-raw-continuum-primitives}
Put
\begin{equation}
 \Psi_d(z)=\cQ_d\bigl(t\mapsto K_z(t)\bigr),
 \qquad
 \Psi_\infty(z)=\int_{-1/2}^{1/2}K_z(t)\,dt.
 \label{eq:section9-raw-continuum-primitives}
\end{equation}
\end{definition}

\begin{proposition}[Primitive identities]
\label{prop:section9-primitive-identities}
The functions $\Psi_d$ and $\Psi_\infty$ are entire, real-valued on $\mathbb R$, and normalized by
\begin{equation}
 \Psi_d(1)=\Psi_\infty(1)=0.
 \label{eq:section9-primitive-normalization}
\end{equation}
They satisfy
\begin{equation}
 \boxed{
 \Psi_d(z+1)-\Psi_d(z)=\Theta_d(z),
 }
 \label{eq:section9-raw-primitive-difference}
\end{equation}
and
\begin{equation}
 \boxed{
 \Psi_\infty(z+1)-\Psi_\infty(z)=\sinc(z).
 }
 \label{eq:section9-continuum-primitive-difference}
\end{equation}
For every integer $m\ge1$,
\begin{equation}
 \boxed{
 \Psi_d(m)=\left\lfloor\frac{m-1}{d}\right\rfloor,
 \qquad
 \Psi_\infty(m)=0.
 }
 \label{eq:section9-primitive-integer-values}
\end{equation}
\end{proposition}

\begin{proof}
Entirety follows by integrating or summing a jointly continuous family that is entire in $z$, with the removable value inserted at $t=0$.  Symmetry under $t\mapsto-t$ gives reality on the real axis.  Since $K_1(t)=0$, both functions vanish at $1$.  Applying $\cQ_d$ and $\cQ_\infty$ to \eqref{eq:section9-kernel-difference-identity} gives the two difference equations.

For $m\ge1$, the first difference equation and Proposition~\ref{prop:section9-closed-cyclotomic-form} give
\[
 \Psi_d(m)
 =\sum_{n=1}^{m-1}\mathbf1_{d\mid n}
 =\left\lfloor\frac{m-1}{d}\right\rfloor.
\]
Similarly, $\sinc(n)=0$ for every positive integer $n$, so the normalized continuum primitive vanishes at every positive integer.
\end{proof}

\begin{proposition}[Layer reflection identity]
\label{prop:section9-layer-reflection}
For every $d\ge1$ and every $z\in\CC$,
\begin{equation}
 \boxed{
 \Psi_d(z)+\Psi_d(1-z)=-1,
 \qquad
 \Psi_\infty(z)+\Psi_\infty(1-z)=-1.
 }
 \label{eq:section9-raw-continuum-reflection}
\end{equation}
Consequently,
\begin{equation}
 \boxed{
 \phi_d(1-z)=-\phi_d(z).
 }
 \label{eq:section9-layer-primitive-reflection}
\end{equation}
In particular, $\phi_d(1/2)=0$.
\end{proposition}

\begin{proof}
For $t\ne0$, a direct calculation from \eqref{eq:section9-difference-kernel} gives
\begin{equation}
 K_z(t)+K_{1-z}(-t)=-1,
 \label{eq:section9-kernel-reflection}
\end{equation}
and the identity extends across $t=0$ by removability.  Both $\cQ_d$ and $\cQ_\infty$ are invariant under $t\mapsto-t$.  Applying them to \eqref{eq:section9-kernel-reflection} proves the two identities in \eqref{eq:section9-raw-continuum-reflection}.  Their difference is \eqref{eq:section9-layer-primitive-reflection}.
\end{proof}

For completeness, the raw primitive has an explicit finite Fourier expansion.  Let $I_d$ be the symmetric frequency set used by $\cQ_d$, with endpoint weights $1/2$ in the even case, and write $a_{d,r}$ for the corresponding quadrature weight.  Then
\begin{equation}
 \Psi_d(z)
 =a_{d,0}(z-1)
 +\sum_{\substack{r\in I_d\\r\ne0}}
 a_{d,r}
 \frac{\e(rz/d)-\e(r/d)}{\e(r/d)-1}.
 \label{eq:section9-raw-primitive-fourier-expansion}
\end{equation}
The zero frequency produces the affine term, while every nonzero frequency has a uniquely normalized primitive inside the same finite Fourier space.

\begin{proposition}[Finite-band uniqueness]
\label{prop:section9-finite-band-uniqueness}
Among normalized functions of the form
\begin{equation}
 c(z-1)+\sum_{\substack{r\in I_d\\r\ne0}}c_r\e(rz/d),
 \label{eq:section9-finite-band-primitive-space}
\end{equation}
there is exactly one solution of
\[
 H(z+1)-H(z)=\Theta_d(z),
 \qquad H(1)=0,
\]
namely $\Psi_d$.
\end{proposition}

\begin{proof}
The difference operator sends the affine basis vector $z-1$ to $1$ and multiplies the nonzero mode $\e(rz/d)$ by $\e(r/d)-1$, which is nonzero for the frequencies in $I_d$.  Hence it is invertible from the normalized finite-band primitive space onto the cyclotomic frequency space.
\end{proof}

The continuum primitive is the limiting Riemann integral of the same normalized Fourier primitives.  It is a universal null interpolation: its difference and its positive-integer values both vanish on the discrete skeleton in the appropriate sense.

\begin{definition}[Renormalized cyclotomic primitive]
\label{def:section9-renormalized-primitive}
Define
\begin{equation}
 \boxed{
 \phi_d(z)=\Psi_d(z)-\Psi_\infty(z).
 }
 \label{eq:section9-renormalized-cyclotomic-primitive}
\end{equation}
\end{definition}

Combining the preceding identities gives
\begin{equation}
 \boxed{
 \phi_d(z+1)-\phi_d(z)=\kappa_d(z),
 \qquad
 \phi_d(1)=0,
 }
 \label{eq:section9-renormalized-primitive-difference}
\end{equation}
and
\begin{equation}
 \boxed{
 \phi_d(m)=\left\lfloor\frac{m-1}{d}\right\rfloor
 \qquad(m\ge1).
 }
 \label{eq:section9-renormalized-primitive-integers}
\end{equation}

\subsection{Quadratic conductor decay}
\label{subsec:section9-quadratic-conductor-decay}

The reason for subtracting the continuum mode is quantitative.  Both the raw symbol and its primitive are quadrature approximations, and the midpoint and trapezoidal rules have second-order error.

\begin{lemma}[Uniform cyclotomic quadrature error]
\label{lem:section9-quadrature-error}
Let $K\subset\CC$ be compact and let $q\ge0$.  There is a constant $C_{K,q}$ such that for every $d\ge1$,
\begin{equation}
 \boxed{
 \sup_{z\in K}
 \left|
 \frac{d^q}{dz^q}\kappa_d(z)
 \right|
 +
 \sup_{z\in K}
 \left|
 \frac{d^q}{dz^q}\phi_d(z)
 \right|
 \le
 \frac{C_{K,q}}{d^2}.
 }
 \label{eq:section9-quadratic-conductor-decay}
\end{equation}
\end{lemma}

\begin{proof}
For $z$ in a fixed compact set, all mixed derivatives in $t$ and $z$ of
\[
 \e(tz)
 \quad\text{and}\quad
 K_z(t)
\]
are uniformly bounded on $[-1/2,1/2]$ after filling the removable singularity at $t=0$.  When $d$ is odd, $\cQ_d$ is the composite midpoint rule; when $d$ is even, it is the composite trapezoidal rule.  The standard error estimate for either rule is bounded by a constant times $d^{-2}$ times the supremum of the second $t$-derivative.  Differentiation in $z$ commutes with the finite sums and the integral, giving the assertion for every $q$.
\end{proof}

\begin{remark}[The subtraction is at the symbol level]
\label{rem:section9-symbol-level-renormalization}
The passage from $\Theta_d$ to $\kappa_d$ changes the entire recurrence symbol, not merely the choice of lift inside a fixed torsor.  It is permitted because $\sinc(n)=0$ at every positive integer, so the discrete factorial recurrence is unchanged.  The two analytic ambiguities separated in Section~6 therefore play different roles: cyclotomic renormalization first selects a summable entire recurrence symbol, and orbitwise Stirling rigidity later selects one primitive in the lift torsor of that symbol.
\end{remark}

\subsection{Atomic prime-layer symbols and lifts}
\label{subsec:section9-atomic-prime-layers}

Recall the prime-layer index set
\begin{equation}
 \cA=\{(p,j):p\in\PP,\ j\ge0\},
 \qquad
 d_{p,j}=(p-1)p^j,
 \qquad
 w_{p,j}=\log p.
 \label{eq:section9-prime-layer-data}
\end{equation}
The atomic factorial datum is
\begin{equation}
 E_n^{p,j}
 =p^{\lfloor n/d_{p,j}\rfloor}.
 \label{eq:section9-atomic-factorial-datum}
\end{equation}

\begin{definition}[Atomic cyclotomic layer]
\label{def:section9-atomic-cyclotomic-layer}
For $\lambda=(p,j)\in\cA$, define
\begin{equation}
 \boxed{
 L_{p,j}^{\mathrm{cyc}}(z)
 =(\log p)\,\kappa_{d_{p,j}}(z),
 }
 \label{eq:section9-atomic-cyclotomic-symbol}
\end{equation}
and
\begin{equation}
 \boxed{
 H_{p,j}^{\mathrm{cyc}}(z)
 =(\log p)\,\phi_{d_{p,j}}(z).
 }
 \label{eq:section9-atomic-cyclotomic-primitive}
\end{equation}
\end{definition}

\begin{theorem}[Atomic interpolation and lifting]
\label{thm:section9-atomic-interpolation}
For every prime layer $(p,j)$,
\begin{equation}
 H_{p,j}^{\mathrm{cyc}}(z+1)-H_{p,j}^{\mathrm{cyc}}(z)
 =L_{p,j}^{\mathrm{cyc}}(z),
 \qquad
 H_{p,j}^{\mathrm{cyc}}(1)=0.
 \label{eq:section9-atomic-lift-equation}
\end{equation}
At the positive integers,
\begin{equation}
 L_{p,j}^{\mathrm{cyc}}(n)
 =(\log p)\mathbf1_{d_{p,j}\mid n},
 \label{eq:section9-atomic-symbol-integers}
\end{equation}
and
\begin{equation}
 H_{p,j}^{\mathrm{cyc}}(m)
 =(\log p)
 \left\lfloor\frac{m-1}{d_{p,j}}\right\rfloor
 =\log E_{m-1}^{p,j}.
 \label{eq:section9-atomic-primitive-integers}
\end{equation}
Consequently,
\begin{equation}
 L_{p,j}^{\mathrm{cyc}}\in\Rec(E^{p,j}),
 \qquad
 H_{p,j}^{\mathrm{cyc}}\in\Lift(L_{p,j}^{\mathrm{cyc}})
 \cap\Ext(E^{p,j}).
 \label{eq:section9-atomic-memberships}
\end{equation}
\end{theorem}

\begin{proof}
All assertions follow by multiplying \eqref{eq:section9-renormalized-primitive-difference}, \eqref{eq:section9-renormalized-symbol-integers}, and \eqref{eq:section9-renormalized-primitive-integers} by $\log p$.
\end{proof}

The layer $(2,0)$ has conductor $1$.  In this case the raw kernel is $\Theta_1=1$, the raw primitive is $\Psi_1(z)=z-1$, and the same formulas give
\[
 \kappa_1=1-\sinc,
 \qquad
 \phi_1=z-1-\Psi_\infty.
\]
Thus no exceptional convention is required at the smallest conductor.

\subsection{Finite systems and the complementary cocycle}
\label{subsec:section9-finite-systems-cocycle}

Let $\mathbf m\in\NN^{(\cA)}$ be a finite multiplicity vector.  Define
\begin{equation}
 L_{\mathbf m}^{\mathrm{cyc}}
 =\sum_{(p,j)\in\cA}m_{p,j}L_{p,j}^{\mathrm{cyc}},
 \qquad
 H_{\mathbf m}^{\mathrm{cyc}}
 =\sum_{(p,j)\in\cA}m_{p,j}H_{p,j}^{\mathrm{cyc}}.
 \label{eq:section9-finite-layer-symbol-primitive}
\end{equation}
The sums are finite.

\begin{proposition}[Strict monoidal additivity]
\label{prop:section9-strict-monoidal-additivity}
For finite multiplicity vectors $\mathbf m,\mathbf n$,
\begin{equation}
 L_{\mathbf m+\mathbf n}^{\mathrm{cyc}}
 =L_{\mathbf m}^{\mathrm{cyc}}+L_{\mathbf n}^{\mathrm{cyc}},
 \qquad
 H_{\mathbf m+\mathbf n}^{\mathrm{cyc}}
 =H_{\mathbf m}^{\mathrm{cyc}}+H_{\mathbf n}^{\mathrm{cyc}}.
 \label{eq:section9-monoidal-additivity}
\end{equation}
Moreover,
\begin{equation}
 H_{\mathbf m}^{\mathrm{cyc}}(z+1)
 -H_{\mathbf m}^{\mathrm{cyc}}(z)
 =L_{\mathbf m}^{\mathrm{cyc}}(z),
 \label{eq:section9-finite-lift-equation}
\end{equation}
and at positive integers
\begin{equation}
 \boxed{
 H_{\mathbf m}^{\mathrm{cyc}}(m)
 =\log A_{m-1}^{\mathbf m},
 \qquad
 L_{\mathbf m}^{\mathrm{cyc}}(n)
 =\log\frac{A_n^{\mathbf m}}{A_{n-1}^{\mathbf m}}.
 }
 \label{eq:section9-finite-integer-recovery}
\end{equation}
\end{proposition}

\begin{proof}
Additivity is immediate from the definitions.  The lifting and interpolation identities follow by summing Theorem~\ref{thm:section9-atomic-interpolation} over the finite support of $\mathbf m$.
\end{proof}

Suppose $\mathbf m\le\mathbf n$.  Define the canonical complementary primitive
\begin{equation}
 C_{\mathbf m,\mathbf n}^{\mathrm{cyc}}
 =H_{\mathbf n-\mathbf m}^{\mathrm{cyc}}.
 \label{eq:section9-complementary-cyclotomic-primitive}
\end{equation}
Then
\begin{equation}
 H_{\mathbf n}^{\mathrm{cyc}}
 =H_{\mathbf m}^{\mathrm{cyc}}
 +C_{\mathbf m,\mathbf n}^{\mathrm{cyc}},
 \label{eq:section9-finite-complement-decomposition}
\end{equation}
and for $\mathbf m\le\mathbf n\le\mathbf r$,
\begin{equation}
 \boxed{
 C_{\mathbf m,\mathbf r}^{\mathrm{cyc}}
 =C_{\mathbf m,\mathbf n}^{\mathrm{cyc}}
 +C_{\mathbf n,\mathbf r}^{\mathrm{cyc}}.
 }
 \label{eq:section9-complementary-cocycle}
\end{equation}
Thus the cyclotomic lifts define a compatible section of the inverse system of gamma-lift torsors from Section~6:
\begin{equation}
 r_{\mathbf m,\mathbf n}
 \bigl(H_{\mathbf n}^{\mathrm{cyc}}\bigr)
 =H_{\mathbf m}^{\mathrm{cyc}}.
 \label{eq:section9-compatible-torsor-section}
\end{equation}

The raw finite-band primitive retains the affine-periodic normal form from Part~I.  Put
\begin{equation}
 \pi_d(z)=\Psi_d(z)-\frac{z-1}{d}.
 \label{eq:section9-raw-periodic-part}
\end{equation}
Then $\pi_d$ is entire and $d$-periodic.  If
\begin{equation}
 \sigma(\mathbf m)
 =\sum_{p,j}m_{p,j}\frac{\log p}{d_{p,j}},
 \qquad
 W(\mathbf m)=\sum_{p,j}m_{p,j}\log p,
 \label{eq:section9-finite-slope-weight}
\end{equation}
and
\begin{equation}
 P_{\mathbf m}(z)
 =\sum_{p,j}m_{p,j}(\log p)\pi_{d_{p,j}}(z),
 \label{eq:section9-finite-periodic-function}
\end{equation}
then
\begin{equation}
 \boxed{
 H_{\mathbf m}^{\mathrm{cyc}}(z)
 =\sigma(\mathbf m)(z-1)
 +P_{\mathbf m}(z)
 -W(\mathbf m)\Psi_\infty(z).
 }
 \label{eq:section9-finite-renormalized-normal-form}
\end{equation}
At positive integers the continuum term vanishes, and
\begin{equation}
 P_{\mathbf m}(m)
 =-\sum_{p,j}m_{p,j}(\log p)
 \left\{\frac{m-1}{d_{p,j}}\right\},
 \label{eq:section9-periodic-part-integers}
\end{equation}
recovering the exact affine-periodic normal form of Section~4.

\begin{remark}[A decomposition that must not be completed termwise]
\label{rem:section9-no-termwise-raw-completion}
For finite systems, \eqref{eq:section9-finite-renormalized-normal-form} separates the affine, periodic, and continuum pieces.  In the full prime filtration, however, $W(\mathbf m)$ diverges, and the other two pieces also fail to converge independently.  Only the renormalized layer-by-layer combination
\[
 \sum_{p,j}(\log p)\phi_{d_{p,j}}
\]
is normally summable.  This is the analytic counterpart of the change from finite exponential growth to completed factorial growth: the cancellation must be retained before passing to the limit.
\end{remark}

\subsection{Normal summability over the prime layers}
\label{subsec:section9-normal-summability}

The quadratic conductor estimate is strong enough to sum all prime layers.  Indeed,
\begin{align}
 \sum_{p\in\PP}\sum_{j\ge0}
 \frac{\log p}{d_{p,j}^2}
 &=
 \sum_{p\in\PP}
 \frac{\log p}{(p-1)^2}
 \sum_{j\ge0}p^{-2j}
 \notag\\
 &=
 \sum_{p\in\PP}
 \frac{\log p}{(p-1)^2(1-p^{-2})}
 <\infty.
 \label{eq:section9-prime-layer-square-summability}
\end{align}

\begin{theorem}[Normal convergence of the prime-layer series]
\label{thm:section9-normal-convergence}
The series
\begin{equation}
 \boxed{
 L_{\PP}^{\mathrm{cyc}}(z)
 =\sum_{p\in\PP}\sum_{j\ge0}
 (\log p)\kappa_{(p-1)p^j}(z)
 }
 \label{eq:section9-completed-cyclotomic-symbol}
\end{equation}
and
\begin{equation}
 \boxed{
 H_{\PP}^{\mathrm{cyc}}(z)
 =\sum_{p\in\PP}\sum_{j\ge0}
 (\log p)\phi_{(p-1)p^j}(z)
 }
 \label{eq:section9-completed-cyclotomic-primitive}
\end{equation}
converge absolutely and locally uniformly on $\CC$, together with every derivative.  Their sums are entire, real-valued on $\mathbb R$, and independent of the chosen cofinal ordering of the layers.  Moreover,
\begin{equation}
 \boxed{
 H_{\PP}^{\mathrm{cyc}}(z+1)
 -H_{\PP}^{\mathrm{cyc}}(z)
 =L_{\PP}^{\mathrm{cyc}}(z),
 \qquad
 H_{\PP}^{\mathrm{cyc}}(1)=0.
 }
 \label{eq:section9-completed-lift-equation}
\end{equation}
\end{theorem}

\begin{proof}
Fix a compact set $K$ and a derivative order $q$.  Lemma~\ref{lem:section9-quadrature-error} bounds the corresponding layer term by
\[
 C_{K,q}\frac{\log p}{((p-1)p^j)^2}.
\]
The majorant is summable by \eqref{eq:section9-prime-layer-square-summability}.  The Weierstrass theorem gives normal convergence of both series and all derivatives.  Reality follows termwise.  Absolute normal convergence gives independence of ordering.  Finally, the difference equation follows by termwise subtraction, and normalization follows because every layer primitive vanishes at $1$.
\end{proof}

At the positive integers the series are actually finite, because a layer contributes only when its conductor does not exceed the relevant integer.

\begin{corollary}[Recovery of the completed prime factorial skeleton]
\label{cor:section9-completed-integer-recovery}
For every $n\ge1$,
\begin{equation}
 \boxed{
 L_{\PP}^{\mathrm{cyc}}(n)
 =\sum_{\substack{p\in\PP\\p-1\mid n}}
 \bigl(v_p(n)+1\bigr)\log p
 =\log\frac{(n+1)!_{\PP}}{n!_{\PP}}.
 }
 \label{eq:section9-completed-symbol-integers}
\end{equation}
For every $m\ge1$,
\begin{equation}
 \boxed{
 H_{\PP}^{\mathrm{cyc}}(m)
 =\sum_{p,j}(\log p)
 \left\lfloor\frac{m-1}{(p-1)p^j}\right\rfloor
 =\log m!_{\PP}.
 }
 \label{eq:section9-completed-primitive-integers}
\end{equation}
Consequently,
\begin{equation}
 L_{\PP}^{\mathrm{cyc}}
 \in\Rec\bigl(((n+1)!_{\PP})_{n\ge0}\bigr),
 \qquad
 H_{\PP}^{\mathrm{cyc}}
 \in\Lift(L_{\PP}^{\mathrm{cyc}}).
 \label{eq:section9-completed-memberships}
\end{equation}
\end{corollary}

\begin{proof}
At an integer $n$, the symbol contribution of $(p,j)$ is nonzero precisely when $(p-1)p^j\mid n$.  If $p-1\mid n$, the admissible indices are $0\le j\le v_p(n)$, giving the first formula.  The primitive formula follows from \eqref{eq:section9-renormalized-primitive-integers}; only conductors at most $m-1$ contribute.  The final equalities are the prime factorial formulas established in Section~2.
\end{proof}

\subsection{Conductor truncations and compatible convergence}
\label{subsec:section9-conductor-truncations}

For $T\ge1$, let
\begin{equation}
 \cA_{\le T}
 =\{(p,j):(p-1)p^j\le T\}
 \label{eq:section9-conductor-layer-set}
\end{equation}
and define
\begin{equation}
 L_{[T]}^{\mathrm{cyc}}
 =\sum_{(p,j)\in\cA_{\le T}}L_{p,j}^{\mathrm{cyc}},
 \qquad
 H_{[T]}^{\mathrm{cyc}}
 =\sum_{(p,j)\in\cA_{\le T}}H_{p,j}^{\mathrm{cyc}}.
 \label{eq:section9-conductor-truncated-lifts}
\end{equation}
If $T\le T'$, put
\begin{equation}
 C_{T,T'}^{\mathrm{cyc}}
 =H_{[T']}^{\mathrm{cyc}}-H_{[T]}^{\mathrm{cyc}}
 =\sum_{\substack{p,j\\T<(p-1)p^j\le T'}}H_{p,j}^{\mathrm{cyc}}.
 \label{eq:section9-truncation-complementary-primitive}
\end{equation}
Then
\begin{equation}
 C_{T,T''}^{\mathrm{cyc}}
 =C_{T,T'}^{\mathrm{cyc}}+C_{T',T''}^{\mathrm{cyc}}
 \qquad(T\le T'\le T''),
 \label{eq:section9-truncation-cocycle}
\end{equation}
and
\begin{equation}
 r_{T,T'}\bigl(H_{[T']}^{\mathrm{cyc}}\bigr)
 =H_{[T]}^{\mathrm{cyc}}.
 \label{eq:section9-truncation-torsor-compatibility}
\end{equation}
By Theorem~\ref{thm:section9-normal-convergence},
\begin{equation}
 H_{[T]}^{\mathrm{cyc}}\longrightarrow H_{\PP}^{\mathrm{cyc}},
 \qquad
 L_{[T]}^{\mathrm{cyc}}\longrightarrow L_{\PP}^{\mathrm{cyc}}
 \label{eq:section9-truncations-normal-convergence}
\end{equation}
locally uniformly with all derivatives.

There is also exact degreewise stabilization at the integers.  If $T\ge N$, then for $1\le m\le N+1$,
\begin{equation}
 H_{[T]}^{\mathrm{cyc}}(m)
 =H_{\PP}^{\mathrm{cyc}}(m),
 \label{eq:section9-primitive-degreewise-stabilization}
\end{equation}
and for $1\le n\le N$,
\begin{equation}
 L_{[T]}^{\mathrm{cyc}}(n)
 =L_{\PP}^{\mathrm{cyc}}(n).
 \label{eq:section9-symbol-degreewise-stabilization}
\end{equation}
Thus the cyclotomic analytic completion mirrors the exact degreewise stabilization of the factorial calculi from Part~I, while adding normal convergence away from the integer skeleton.

\subsection{Structural conclusion}
\label{subsec:section9-structural-conclusion}

The construction may be summarized by the commutative scheme
\begin{equation}
 \boxed{
 \begin{gathered}
 \Theta_d=\cQ_d\e(tz),
 \qquad
 \Theta_\infty=\cQ_\infty\e(tz)=\sinc(z),
 \\
 \kappa_d=\Theta_d-\Theta_\infty,
 \qquad
 \kappa_d(n)=\mathbf1_{d\mid n},
 \qquad
 \kappa_d=O_{\mathrm{loc}}(d^{-2}),
 \\
 \Psi_d=\cQ_dK_z,
 \qquad
 \Psi_\infty=\cQ_\infty K_z,
 \qquad
 \phi_d=\Psi_d-\Psi_\infty,
 \\
 \Delta\phi_d=\kappa_d,
 \qquad
 \phi_d(m)=\left\lfloor\frac{m-1}{d}\right\rfloor,
 \qquad
 \phi_d=O_{\mathrm{loc}}(d^{-2}),
 \\
 L_{p,j}^{\mathrm{cyc}}=(\log p)\kappa_{(p-1)p^j},
 \qquad
 H_{p,j}^{\mathrm{cyc}}=(\log p)\phi_{(p-1)p^j},
 \\
 H_{\PP}^{\mathrm{cyc}}
 =\sum_{p,j}H_{p,j}^{\mathrm{cyc}},
 \qquad
 \Delta H_{\PP}^{\mathrm{cyc}}=L_{\PP}^{\mathrm{cyc}},
 \qquad
 H_{\PP}^{\mathrm{cyc}}(m)=\log m!_{\PP}.
 \end{gathered}
 }
 \label{eq:section9-structural-summary}
\end{equation}

The continuum subtraction is the decisive step.  It is invisible on the positive integers, canonical from the limiting cyclotomic quadrature, and exactly strong enough to make the full prime-layer family normally summable.  The completed recurrence symbol and one distinguished compatible lift have now been constructed globally on $\CC$.

What remains is not existence of an entire lift, but its Stirling characterization.  The next section compares $H_{\PP}^{\mathrm{cyc}}$ with
\[
 \log\Gamma(z)+C_{\PP}(z-1),
\]
proves the required orbitwise asymptotics, and invokes abstract rigidity to show that the cyclotomic lift is the unique Stirling-admissible point of the completed gamma-lift torsor.

\section{The distinguished prime gamma lift}
\label{sec:distinguished-prime-gamma-lift}

The preceding section constructed an entire recurrence symbol
\[
 L_{\PP}^{\mathrm{cyc}}
\]
and a normalized entire lift
\[
 H_{\PP}^{\mathrm{cyc}}\in\Lift(L_{\PP}^{\mathrm{cyc}})
\]
whose positive-integer values are
\[
 H_{\PP}^{\mathrm{cyc}}(m)=\log m!_{\PP}.
\]
Existence alone does not distinguish this lift: every function
\[
 H_{\PP}^{\mathrm{cyc}}+P,
 \qquad P\in\Per_0,
\]
has the same recurrence and the same positive-integer trace.  The purpose of this section is to prove that the cyclotomic lift is the unique point of this torsor compatible with the prime Stirling datum.

The proof has three parts.  First, the twice-iterated Ces\`aro average of the classical comparison model is computed uniformly across the translation interval.  Second, the cyclotomic layer expansion is transformed into a Fourier multiplier whose nonconstant transverse component is summably small.  Third, the remaining scalar layer sum collapses by an exact Abel identity to the weak Stirling formula for the prime Bhargava factorial.  The logarithmic terms and constant terms then cancel exactly.

Throughout this section, \(\log\Gamma\) denotes the holomorphic branch on the half-plane \(\Re z>0\) that is real on the positive axis.

\subsection{Arithmetic conductor weights and the Stirling datum}
\label{subsec:section10-arithmetic-input}

Aggregate all prime layers having the same conductor by
\begin{equation}
 h_{\PP}(d)
 =\sum_{\substack{p\in\PP,\ j\ge0\\(p-1)p^j=d}}\log p.
 \label{eq:section10-conductor-weight}
\end{equation}
Then the exact prime-factorial formula becomes
\begin{equation}
 \boxed{
 \log(n+1)!_{\PP}
 =\sum_{d\le n}h_{\PP}(d)
 \left\lfloor\frac nd\right\rfloor.
 }
 \label{eq:section10-floor-expansion}
\end{equation}
We use the weak Stirling law
\begin{equation}
 \boxed{
 \log(n+1)!_{\PP}
 =\log n!+C_{\PP}n+o(n),
 }
 \label{eq:section10-weak-stirling}
\end{equation}
where
\begin{equation}
 C_{\PP}
 =\sum_{p\in\PP}\frac{\log p}{(p-1)^2}.
 \label{eq:section10-prime-stirling-constant}
\end{equation}
The series in \eqref{eq:section10-prime-stirling-constant} converges absolutely.

Let
\[
 \Harm_n=\sum_{m=1}^n\frac1m,
 \qquad \Harm_0=0.
\]
The twice-iterated Ces\`aro weights are
\begin{equation}
 \boxed{
 \omega_{N,n}
 =\frac{\Harm_N-\Harm_n}{N},
 \qquad 0\le n<N.
 }
 \label{eq:section10-iterated-cesaro-weights}
\end{equation}
They arise by applying the ordinary Ces\`aro operator twice:
\begin{equation}
 \frac1N\sum_{m=1}^N\frac1m\sum_{n=0}^{m-1}f_n
 =\sum_{n=0}^{N-1}\omega_{N,n}f_n.
 \label{eq:section10-iterated-cesaro-identity}
\end{equation}
In particular,
\begin{equation}
 \sum_{n=0}^{N-1}\omega_{N,n}=1,
 \qquad
 \sum_{n=0}^{N-1}n\omega_{N,n}=\frac{N-1}{4}.
 \label{eq:section10-weight-moments}
\end{equation}

Define the prime comparison model
\begin{equation}
 M_{\PP}(z)
 =\log\Gamma(z)+C_{\PP}(z-1),
 \qquad \Re z>0.
 \label{eq:section10-prime-model}
\end{equation}
For a normalized lift \(H\in\Lift(L_{\PP}^{\mathrm{cyc}})\), put
\begin{equation}
 \cR_N(u;H)
 =\sum_{n=0}^{N-1}\omega_{N,n}
 \bigl(H(n+u+1)-M_{\PP}(n+u+1)\bigr),
 \qquad 0\le u\le1,
 \label{eq:section10-orbitwise-defect}
\end{equation}
and define the transverse profile
\begin{equation}
 \cT_N(H)(u)=\cR_N(u;H)-\cR_N(0;H).
 \label{eq:section10-transverse-profile}
\end{equation}
The goal is to prove
\begin{equation}
 \norm{\cT_N(H_{\PP}^{\mathrm{cyc}})}_{C([0,1])}
 \longrightarrow0.
 \label{eq:section10-goal}
\end{equation}

\begin{remark}[Which second Ces\`aro mean is used]
\label{rem:section10-second-cesaro-convention}
There are several standard order-two Ces\`aro procedures.  The prime proof uses the iterated operator \(C\circ C\), whose harmonic weights are given by \eqref{eq:section10-iterated-cesaro-weights}.  Its exact first difference
\[
 \omega_{N,m}-\omega_{N,m+1}=\frac1{N(m+1)}
\]
is what converts the completed layer sum into the weak Stirling law.  As with every regular orbitwise mean, periodic gauge functions are preserved exactly.
\end{remark}

\subsection{The classical shifted average}
\label{subsec:section10-classical-shifted-average}

Define
\begin{equation}
 \cG_N(u)
 =\sum_{n=0}^{N-1}\omega_{N,n}
 \log\Gamma(n+u+1).
 \label{eq:section10-gamma-average}
\end{equation}

\begin{lemma}[Shifted gamma average]
\label{lem:section10-shifted-gamma-average}
Uniformly for \(0\le u\le1\),
\begin{align}
 \cG_N(u)
 ={}&\frac N4\log N-\frac N2
 +\left(u+\frac14\right)\log N \notag\\
 &+\frac12\log(2\pi)-\frac38-2u
 +O\!\left(\frac{(\log(2N))^2}{N}\right).
 \label{eq:section10-gamma-average-asymptotic}
\end{align}
Consequently,
\begin{equation}
 \boxed{
 \cG_N(u)-\cG_N(0)
 =u\log N-2u
 +O\!\left(\frac{(\log(2N))^2}{N}\right)
 }
 \label{eq:section10-gamma-average-difference}
\end{equation}
uniformly on \([0,1]\).
\end{lemma}

\begin{proof}
Let \(G\) denote the Barnes function~\cite{Barnes1900,Vardi1988}, normalized by
\[
 G(z+1)=\Gamma(z)G(z),
 \qquad G(1)=1.
\]
Reversing the harmonic summation in \eqref{eq:section10-iterated-cesaro-weights} gives the exact identity
\begin{equation}
 \cG_N(u)
 =\frac1N\sum_{m=1}^{N}
 \frac{\log G(m+u+1)}m
 -\frac{\Harm_N}{N}\log G(u+1).
 \label{eq:section10-barnes-average}
\end{equation}
Uniformly for \(0\le u\le1\), the standard Barnes expansion gives
\begin{align*}
 \log G(m+u+1)
 ={}&\frac12m^2\log m-\frac34m^2
 +um\log m\\
 &+\left(\frac12\log(2\pi)-u\right)m
 +O(\log(m+1)).
\end{align*}
After division by \(m\), use
\begin{align*}
 \sum_{m=1}^{N}m\log m
 &=\frac12N^2\log N-\frac14N^2
   +\frac12N\log N+O(\log(2N)),\\
 \sum_{m=1}^{N}m&=\frac{N(N+1)}2,\\
 \log N!&=N\log N-N+\frac12\log(2\pi N)+O(N^{-1}).
\end{align*}
The final term in \eqref{eq:section10-barnes-average} is \(O((\log N)/N)\), uniformly for \(u\in[0,1]\).  Substitution proves \eqref{eq:section10-gamma-average-asymptotic}; subtracting the case \(u=0\) gives \eqref{eq:section10-gamma-average-difference}.
\end{proof}

\subsection{The exact cyclotomic orbit transform}
\label{subsec:section10-exact-orbit-transform}

Let \(\nu_d\) be the signed quadrature error measure characterized by
\begin{equation}
 \int_{-1/2}^{1/2}f(t)\,\dd\nu_d(t)
 =\mathcal Q_df-\int_{-1/2}^{1/2}f(t)\,\dd t.
 \label{eq:section10-quadrature-error-measure}
\end{equation}
With the difference kernel from Section~9,
\[
 K_z(t)=\frac{\e(tz)-\e(t)}{\e(t)-1},
\]
where the value at \(t=0\) is filled by continuity, one has
\begin{equation}
 \phi_d(z)=\int_{-1/2}^{1/2}K_z(t)\,\dd\nu_d(t).
 \label{eq:section10-layer-primitive-integral}
\end{equation}
Thus
\begin{equation}
 H_{\PP}^{\mathrm{cyc}}(z)
 =\sum_{d\ge1}h_{\PP}(d)
 \int_{-1/2}^{1/2}K_z(t)\,\dd\nu_d(t).
 \label{eq:section10-completed-primitive-by-conductor}
\end{equation}
The series is locally normally convergent.

Define the Fourier multiplier
\begin{equation}
 \Omega_N(t)
 =\sum_{n=0}^{N-1}\omega_{N,n}\e(nt).
 \label{eq:section10-fourier-multiplier}
\end{equation}
For nonintegral \(t\), direct reversal of the harmonic summation gives
\begin{equation}
 \boxed{
 \Omega_N(t)
 =\frac1N\sum_{m=1}^{N}
 \frac{1-\e(mt)}{m(1-\e(t))},
 \qquad \Omega_N(0)=1.
 }
 \label{eq:section10-fourier-multiplier-formula}
\end{equation}
Also put
\begin{equation}
 Q_{N,u}(t)
 =\sum_{n=0}^{N-1}\omega_{N,n}K_{n+u+1}(t).
 \label{eq:section10-averaged-difference-kernel}
\end{equation}
Using \eqref{eq:section10-weight-moments},
\begin{equation}
 Q_{N,u}(t)
 =\frac{\e((u+1)t)\Omega_N(t)-\e(t)}{\e(t)-1},
 \qquad
 Q_{N,u}(0)=u+\frac{N-1}{4}.
 \label{eq:section10-averaged-kernel-formula}
\end{equation}
The apparent singularity at zero is removable.

Termwise integration in \eqref{eq:section10-completed-primitive-by-conductor} gives
\begin{equation}
 \sum_{n=0}^{N-1}\omega_{N,n}
 H_{\PP}^{\mathrm{cyc}}(n+u+1)
 =\sum_{d\ge1}h_{\PP}(d)
 \int_{-1/2}^{1/2}Q_{N,u}(t)\,\dd\nu_d(t).
 \label{eq:section10-orbit-transform}
\end{equation}
Subtract the orbit \(u=0\), and define
\begin{equation}
 \Psi_{N,u}(t)
 =\frac{\e(t)(\e(ut)-1)\Omega_N(t)}{\e(t)-1},
 \qquad
 \Phi_{N,u}(d)
 =\int_{-1/2}^{1/2}\Psi_{N,u}(t)\,\dd\nu_d(t).
 \label{eq:section10-Psi-Phi}
\end{equation}
Then
\begin{equation}
 \boxed{
 \sum_{n=0}^{N-1}\omega_{N,n}
 \bigl(
 H_{\PP}^{\mathrm{cyc}}(n+u+1)
 -H_{\PP}^{\mathrm{cyc}}(n+1)
 \bigr)
 =\sum_{d\ge1}h_{\PP}(d)\Phi_{N,u}(d).
 }
 \label{eq:section10-orbit-difference-prime-lift}
\end{equation}

Put
\begin{equation}
 \alpha_u(t)
 =\frac{\e(t)(\e(ut)-1)}{\e(t)-1},
 \qquad
 r_u(t)=\alpha_u(t)-u.
 \label{eq:section10-alpha-r}
\end{equation}
The function \(r_u\) is jointly smooth for
\[
 (u,t)\in[0,1]\times[-1/2,1/2]
\]
and satisfies \(r_u(0)=0\).  Hence
\begin{equation}
 \Phi_{N,u}(d)
 =uE_{N,d}+\cE_{N,u}(d),
 \label{eq:section10-Phi-splitting}
\end{equation}
where
\begin{equation}
 E_{N,d}=\int\Omega_N\,\dd\nu_d,
 \qquad
 \cE_{N,u}(d)=\int r_u(t)\Omega_N(t)\,\dd\nu_d(t).
 \label{eq:section10-main-error-transform}
\end{equation}
Since \(\Omega_N\) has only integral Fourier frequencies, the exact divisibility filter of Section~9 gives
\begin{equation}
 \boxed{
 E_{N,d}
 =\sum_{k=1}^{\lfloor(N-1)/d\rfloor}\omega_{N,kd}
 =\frac1N\sum_{k=1}^{\lfloor(N-1)/d\rfloor}
 (\Harm_N-\Harm_{kd}).
 }
 \label{eq:section10-E-exact}
\end{equation}

\subsection{The summable aliasing estimate}
\label{subsec:section10-summable-aliasing}

We need both first- and second-order forms of the cyclotomic quadrature estimate.

\begin{lemma}[Quadrature error in variation norms]
\label{lem:section10-quadrature-variation}
There is an absolute constant \(C\) such that
\begin{equation}
 \left|\int f\,\dd\nu_d\right|
 \le \frac{C}{d}\norm{f'}_{L^1[-1/2,1/2]}
 \label{eq:section10-first-order-quadrature}
\end{equation}
for every continuously differentiable \(f\), and
\begin{equation}
 \left|\int f\,\dd\nu_d\right|
 \le \frac{C}{d^2}\norm{f''}_{L^1[-1/2,1/2]}
 \label{eq:section10-second-order-quadrature}
\end{equation}
whenever \(f'\) is absolutely continuous.
\end{lemma}

\begin{proof}
Partition \([-1/2,1/2]\) into \(d\) cells of length \(d^{-1}\).  In the midpoint case, subtract the midpoint value on every cell; in the trapezoidal case, subtract the affine interpolant through the endpoints.  The fundamental theorem of calculus gives \eqref{eq:section10-first-order-quadrature}.  The standard Peano kernels for midpoint and trapezoidal quadrature give \eqref{eq:section10-second-order-quadrature}.  Summing the local errors proves the lemma.
\end{proof}

\begin{lemma}[Multiplier estimates]
\label{lem:section10-multiplier-estimates}
Uniformly for \(N\ge2\),
\begin{align}
 \norm{\Omega_N}_{L^1}
 &\ll\frac{(\log(2N))^2}{N},
 &
 \norm{\Omega_N'}_{L^1}&\ll1,
 \label{eq:section10-Omega-L1}\\
 \norm{t\Omega_N'}_{L^1}
 &\ll\frac{(\log(2N))^2}{N},
 &
 \norm{t\Omega_N''}_{L^1}&\ll\log(2N).
 \label{eq:section10-Omega-weighted}
\end{align}
\end{lemma}

\begin{proof}
Away from zero, write
\begin{equation}
 \Omega_N(t)
 =\frac{S_N(t)}{N(1-\e(t))},
 \qquad
 S_N(t)=\sum_{m=1}^{N}\frac{1-\e(mt)}m.
 \label{eq:section10-SN-representation}
\end{equation}
For \(N^{-1}\le |t|\le1/2\), split the sum at \(m\asymp |t|^{-1}\) and apply the elementary geometric-sum bounds.  This gives
\begin{align*}
 |S_N(t)|&\ll1+\log(1+N|t|),\\
 |S_N'(t)|&\ll |t|^{-1},\\
 |S_N''(t)|&\ll N|t|^{-1}.
\end{align*}
Since \(|1-\e(t)|\asymp|t|\), differentiation yields
\begin{align*}
 |\Omega_N(t)|
 &\ll\frac{1+\log(1+N|t|)}{N|t|},\\
 |\Omega_N'(t)|
 &\ll\frac{1+\log(1+N|t|)}{N|t|^2},\\
 |\Omega_N''(t)|
 &\ll |t|^{-2}
 +\frac{1+\log(1+N|t|)}{N|t|^3}.
\end{align*}
For \(|t|\le N^{-1}\), positivity and normalization of the weights imply
\[
 |\Omega_N^{(j)}(t)|\ll_jN^j,
 \qquad j=0,1,2.
\]
Integrating after splitting at \(N^{-1}\) proves the asserted bounds.
\end{proof}

Joint smoothness of \(r_u(t)\), together with \(r_u(0)=0\), gives
\begin{equation}
 |r_u(t)|\ll|t|,
 \qquad
 |r_u'(t)|+|r_u''(t)|\ll1,
 \label{eq:section10-r-bounds}
\end{equation}
uniformly for \(0\le u\le1\).  Thus, with
\[
 g_{N,u}(t)=r_u(t)\Omega_N(t),
\]
Lemma~\ref{lem:section10-multiplier-estimates} gives
\begin{equation}
 \norm{g_{N,u}'}_{L^1}
 \ll\frac{(\log(2N))^2}{N},
 \qquad
 \norm{g_{N,u}''}_{L^1}
 \ll\log(2N),
 \label{eq:section10-g-derivatives}
\end{equation}
uniformly in \(u\).

\begin{proposition}[Summable aliasing]
\label{prop:section10-summable-aliasing}
Uniformly for \(0\le u\le1\),
\begin{equation}
 \boxed{
 |\cE_{N,u}(d)|
 \ll
 \min\left\{
 \frac{(\log(2N))^2}{Nd},
 \frac{\log(2N)}{d^2}
 \right\}.
 }
 \label{eq:section10-layer-aliasing-bound}
\end{equation}
Moreover,
\begin{equation}
 \boxed{
 \sup_{0\le u\le1}
 \sum_{d\ge1}h_{\PP}(d)|\cE_{N,u}(d)|
 \ll\frac{(\log(2N))^3}{N}.
 }
 \label{eq:section10-global-aliasing-bound}
\end{equation}
\end{proposition}

\begin{proof}
Apply Lemma~\ref{lem:section10-quadrature-variation} to \(g_{N,u}\), using \eqref{eq:section10-g-derivatives}.  This proves \eqref{eq:section10-layer-aliasing-bound}.

It remains to sum over the conductors.  Put
\[
 A_{\PP}(x)=\sum_{d\le x}h_{\PP}(d).
\]
The layers with \(j=0\) contribute
\[
 \sum_{p\le x+1}\log p\ll x.
\]
For \(j\ge1\), the condition \((p-1)p^j\le x\) implies \(p\le2\sqrt{x}\).  Chebyshev's estimates for \(\vartheta\) and \(\pi\) then give
\[
 \sum_{\substack{p,\,j\ge1\\(p-1)p^j\le x}}\log p
 \ll\sqrt{x}.
\]
Hence
\begin{equation}
 A_{\PP}(x)\ll x.
 \label{eq:section10-conductor-counting}
\end{equation}
Partial summation yields
\begin{equation}
 \sum_{d\le N}\frac{h_{\PP}(d)}d\ll\log(2N),
 \qquad
 \sum_{d>N}\frac{h_{\PP}(d)}{d^2}\ll N^{-1}.
 \label{eq:section10-conductor-partial-summation}
\end{equation}
Use the first estimate in \eqref{eq:section10-layer-aliasing-bound} for \(d\le N\), and the second for \(d>N\).  This proves \eqref{eq:section10-global-aliasing-bound}.
\end{proof}

Combining \eqref{eq:section10-Phi-splitting} with Proposition~\ref{prop:section10-summable-aliasing} gives the uniform linear response
\begin{equation}
 \boxed{
 \sum_{d\ge1}h_{\PP}(d)\Phi_{N,u}(d)
 =u\sum_{d\ge1}h_{\PP}(d)E_{N,d}+o(1).
 }
 \label{eq:section10-linear-response}
\end{equation}

\subsection{The exact Abel collapse}
\label{subsec:section10-abel-collapse}

The scalar expression in \eqref{eq:section10-linear-response} contains the full arithmetic input, but it collapses exactly to a weighted mean of the prime factorial itself.

\begin{proposition}[Scalar main term]
\label{prop:section10-scalar-main-term}
As \(N\to\infty\),
\begin{equation}
 \boxed{
 \sum_{d\ge1}h_{\PP}(d)E_{N,d}
 =\log N+C_{\PP}-2+o(1).
 }
 \label{eq:section10-scalar-asymptotic}
\end{equation}
\end{proposition}

\begin{proof}
Put
\[
 F_m=\log(m+1)!_{\PP},
 \qquad F_0=0.
\]
Taking first differences in \eqref{eq:section10-floor-expansion} gives
\begin{equation}
 F_m-F_{m-1}=\sum_{d\mid m}h_{\PP}(d).
 \label{eq:section10-factorial-increments}
\end{equation}
Using \eqref{eq:section10-E-exact} and rearranging a finite sum,
\begin{align*}
 \sum_{d\ge1}h_{\PP}(d)E_{N,d}
 &=\sum_{m=1}^{N-1}\omega_{N,m}
   \sum_{d\mid m}h_{\PP}(d)\\
 &=\sum_{m=1}^{N-1}\omega_{N,m}(F_m-F_{m-1}).
\end{align*}
Define \(\omega_{N,N}=0\).  Since
\begin{equation}
 \omega_{N,m}-\omega_{N,m+1}
 =\frac1{N(m+1)},
 \label{eq:section10-weight-first-difference}
\end{equation}
discrete Abel summation gives the exact identity
\begin{equation}
 \boxed{
 \sum_{d\ge1}h_{\PP}(d)E_{N,d}
 =\frac1N\sum_{m=1}^{N-1}
 \frac{\log(m+1)!_{\PP}}{m+1}.
 }
 \label{eq:section10-exact-Abel-identity}
\end{equation}
Insert \eqref{eq:section10-weak-stirling} in the form
\[
 \log(m+1)!_{\PP}
 =\log m!+C_{\PP}m+r_m,
 \qquad r_m=o(m).
\]
Ordinary Stirling gives
\[
 \frac{\log m!}{m+1}
 =\log m-1+O\!\left(\frac{\log(m+1)}m\right),
\]
so
\begin{equation}
 \frac1N\sum_{m=1}^{N-1}\frac{\log m!}{m+1}
 =\log N-2
 +O\!\left(\frac{(\log(2N))^2}{N}\right).
 \label{eq:section10-classical-scalar-sum}
\end{equation}
Also,
\begin{equation}
 \frac{C_{\PP}}N\sum_{m=1}^{N-1}\frac{m}{m+1}
 =C_{\PP}+O\!\left(\frac{\log(2N)}N\right).
 \label{eq:section10-linear-scalar-sum}
\end{equation}
Finally, \(r_m/(m+1)\to0\), so its Ces\`aro mean is \(o(1)\).  Substitution into \eqref{eq:section10-exact-Abel-identity} proves \eqref{eq:section10-scalar-asymptotic}.
\end{proof}

\subsection{Orbitwise Stirling selection}
\label{subsec:section10-orbitwise-selection}

We can now prove the analytic existence statement required by the abstract rigidity theorem.

\begin{theorem}[Orbitwise Stirling theorem for the prime lift]
\label{thm:section10-orbitwise-stirling}
The cyclotomic lift constructed in Section~9 satisfies
\begin{equation}
 \boxed{
 \sup_{0\le u\le1}
 |\cT_N(H_{\PP}^{\mathrm{cyc}})(u)|
 \longrightarrow0.
 }
 \label{eq:section10-orbitwise-stirling-conclusion}
\end{equation}
Equivalently, \(H_{\PP}^{\mathrm{cyc}}\) is Stirling-admissible for the comparison model \(M_{\PP}\) and the twice-iterated Ces\`aro mean.
\end{theorem}

\begin{proof}
From the definitions,
\begin{align}
 \cT_N(H_{\PP}^{\mathrm{cyc}})(u)
 ={}&\sum_{d\ge1}h_{\PP}(d)\Phi_{N,u}(d)\notag\\
 &-\bigl(\cG_N(u)-\cG_N(0)\bigr)-C_{\PP}u.
 \label{eq:section10-exact-defect-decomposition}
\end{align}
By \eqref{eq:section10-linear-response} and Proposition~\ref{prop:section10-scalar-main-term}, the first term equals
\[
 u(\log N+C_{\PP}-2)+o(1)
\]
uniformly for \(0\le u\le1\).  By Lemma~\ref{lem:section10-shifted-gamma-average}, the second term equals
\[
 u\log N-2u+o(1)
\]
uniformly in \(u\).  Substitution into \eqref{eq:section10-exact-defect-decomposition} gives
\[
 \cT_N(H_{\PP}^{\mathrm{cyc}})(u)=o(1)
\]
uniformly on \([0,1]\).
\end{proof}

\begin{theorem}[Distinguished prime gamma lift]
\label{thm:section10-distinguished-prime-lift}
The function \(H_{\PP}^{\mathrm{cyc}}\) is the unique normalized entire lift of \(L_{\PP}^{\mathrm{cyc}}\) satisfying the prime orbitwise Stirling condition.  More precisely, if
\[
 H\in\Lift(L_{\PP}^{\mathrm{cyc}}),
\]
then there is a unique \(P\in\Per_0\) such that
\[
 H=H_{\PP}^{\mathrm{cyc}}+P,
\]
and
\begin{equation}
 \boxed{
 \cT_N(H)\longrightarrow \res P
 }
 \label{eq:section10-all-profile-limits}
\end{equation}
uniformly on \([0,1]\).  In particular,
\begin{equation}
 \cT_N(H)\longrightarrow0
 \quad\Longleftrightarrow\quad
 H=H_{\PP}^{\mathrm{cyc}}.
 \label{eq:section10-unique-admissible-lift}
\end{equation}
\end{theorem}

\begin{proof}
The torsor statement from Section~6 gives the unique representation
\[
 H=H_{\PP}^{\mathrm{cyc}}+P.
\]
The exact gauge covariance of the orbitwise profile gives
\[
 \cT_N(H)
 =\cT_N(H_{\PP}^{\mathrm{cyc}})+\res P.
\]
Theorem~\ref{thm:section10-orbitwise-stirling} proves \eqref{eq:section10-all-profile-limits}.  Since restriction to \([0,1]\) is injective on \(\Per_0\), the vanishing criterion is equivalent to \(P=0\).
\end{proof}

The affine coefficient is selected simultaneously.  Indeed, relative to the family
\[
 M_c(z)=\log\Gamma(z)+c(z-1),
\]
a change from \(C_{\PP}\) to \(c\) changes the transverse profile by
\[
 (C_{\PP}-c)u.
\]
The directness of the affine and periodic profile directions established in Section~8 therefore gives the following consequence.

\begin{corollary}[Joint affine-periodic uniqueness]
\label{cor:section10-joint-affine-periodic-uniqueness}
There is exactly one pair
\[
 (c,H)\in\CC\times\Lift(L_{\PP}^{\mathrm{cyc}})
\]
for which the twice-iterated Ces\`aro transverse profile relative to
\[
 \log\Gamma(z)+c(z-1)
\]
tends uniformly to zero.  That pair is
\begin{equation}
 \boxed{
 (c,H)=(C_{\PP},H_{\PP}^{\mathrm{cyc}}).
 }
 \label{eq:section10-joint-selected-pair}
\end{equation}
\end{corollary}

\subsection{The multiplicative gamma object}
\label{subsec:section10-multiplicative-gamma-object}

Exponentiating the selected logarithmic lift gives
\begin{equation}
 \boxed{
 \Gamma_{\PP}^{\mathrm{cyc}}(z)
 =\exp\bigl(H_{\PP}^{\mathrm{cyc}}(z)\bigr).
 }
 \label{eq:section10-prime-gamma-function}
\end{equation}
This function is zero-free and entire, satisfies
\begin{equation}
 \Gamma_{\PP}^{\mathrm{cyc}}(1)=1,
 \qquad
 \Gamma_{\PP}^{\mathrm{cyc}}(m)=m!_{\PP}
 \quad(m\ge1),
 \label{eq:section10-prime-gamma-integers}
\end{equation}
and obeys the exact functional equation
\begin{equation}
 \boxed{
 \Gamma_{\PP}^{\mathrm{cyc}}(z+1)
 =R_{\PP}^{\mathrm{cyc}}(z)
 \Gamma_{\PP}^{\mathrm{cyc}}(z),
 \qquad
 R_{\PP}^{\mathrm{cyc}}(z)
 =\exp\bigl(L_{\PP}^{\mathrm{cyc}}(z)\bigr).
 }
 \label{eq:section10-prime-gamma-functional-equation}
\end{equation}
\begin{theorem}[Euler-type reflection law]
\label{thm:section10-prime-reflection}
The distinguished prime gamma function satisfies
\begin{equation}
 \boxed{
 \Gamma_{\PP}^{\mathrm{cyc}}(z)
 \Gamma_{\PP}^{\mathrm{cyc}}(1-z)=1
 \qquad(z\in\CC).
 }
 \label{eq:section10-prime-reflection}
\end{equation}
Equivalently,
\begin{equation}
 H_{\PP}^{\mathrm{cyc}}(1-z)
 =-H_{\PP}^{\mathrm{cyc}}(z).
 \label{eq:section10-logarithmic-reflection}
\end{equation}
In particular,
\begin{equation}
 \boxed{
 \Gamma_{\PP}^{\mathrm{cyc}}\!\left(\frac12\right)=1,
 }
 \label{eq:section10-prime-midpoint}
\end{equation}
and, for every integer $m\ge1$,
\begin{equation}
 \Gamma_{\PP}^{\mathrm{cyc}}(1-m)=\frac1{m!_{\PP}}.
 \label{eq:section10-negative-integer-values}
\end{equation}
\end{theorem}

\begin{proof}
Proposition~\ref{prop:section9-layer-reflection} gives
\[
 \phi_{(p-1)p^j}(1-z)=-\phi_{(p-1)p^j}(z)
\]
for every prime layer.  The defining prime-layer series for $H_{\PP}^{\mathrm{cyc}}$ converges normally, so it may be reflected and summed termwise, proving \eqref{eq:section10-logarithmic-reflection}.  Exponentiation gives \eqref{eq:section10-prime-reflection}.  At $z=1/2$ the reflection identity yields a square equal to $1$; positivity on the real axis, inherited from $\Gamma_{\PP}^{\mathrm{cyc}}(x)=\exp(H_{\PP}^{\mathrm{cyc}}(x))$, selects the positive value.  Finally, set $z=m$ in \eqref{eq:section10-prime-reflection} and use $\Gamma_{\PP}^{\mathrm{cyc}}(m)=m!_{\PP}$.
\end{proof}

Every other normalized zero-free entire solution with the same zero-winding logarithmic symbol has the form
\begin{equation}
 \Gamma_{\PP}^{\mathrm{cyc}}(z)\exp(P(z)),
 \qquad P\in\Per_0.
 \label{eq:section10-multiplicative-ambiguity}
\end{equation}
The orbitwise Stirling condition selects precisely the case \(P=0\).

\begin{remark}[Relative and absolute canonicity]
\label{rem:section10-relative-canonicity}
The theorem is absolute after the cyclotomic recurrence symbol and the prime Stirling datum have been fixed.  It does not assert that every entire interpolation of the discrete prime ratio must have recurrence symbol \(L_{\PP}^{\mathrm{cyc}}\).  Section~9 selected that symbol by centered cyclotomic quadrature and normal summability; the present section selects a unique primitive inside its lift torsor.
\end{remark}

\begin{remark}[Why selection appears only after completion]
\label{rem:section10-selection-after-completion}
Every fixed finite layer system has an affine-periodic logarithm with exponential growth.  Its natural linear slope is the finite sum
\[
 \sum_{\lambda\in\Lambda}\frac{w_\lambda}{d_\lambda},
\]
which diverges along the prime-layer exhaustion.  The comparison model \(\log\Gamma+C_{\PP}(z-1)\) therefore belongs to the completed factorial-scale object, not to any fixed finite stage.  Filtered compatibility transports the torsor through completion, but Stirling selection becomes available only after the completed asymptotic has formed.
\end{remark}

\subsection{Structural conclusion}
\label{subsec:section10-structural-conclusion}

The construction and selection theorem may be summarized by
\begin{equation}
 \boxed{
 \begin{gathered}
 H_{\PP}^{\mathrm{cyc}}
 =\sum_{d\ge1}h_{\PP}(d)\phi_d,
 \qquad
 \Delta H_{\PP}^{\mathrm{cyc}}=L_{\PP}^{\mathrm{cyc}},
 \qquad
 H_{\PP}^{\mathrm{cyc}}(m)=\log m!_{\PP},
 \\
 \sum_{d\ge1}h_{\PP}(d)\Phi_{N,u}(d)
 =u\sum_{d\ge1}h_{\PP}(d)E_{N,d}+o(1),
 \\
 \sum_{d\ge1}h_{\PP}(d)E_{N,d}
 =\log N+C_{\PP}-2+o(1),
 \\
 \cG_N(u)-\cG_N(0)
 =u\log N-2u+o(1),
 \\
 \cT_N(H_{\PP}^{\mathrm{cyc}})\longrightarrow0,
 \qquad
 \cT_N(H_{\PP}^{\mathrm{cyc}}+P)\longrightarrow\res P,
 \\
 \Gamma_{\PP}^{\mathrm{cyc}}
 =\exp(H_{\PP}^{\mathrm{cyc}})
 \text{ is the unique Stirling-selected prime gamma object,}
 \\
 \Gamma_{\PP}^{\mathrm{cyc}}(z)
 \Gamma_{\PP}^{\mathrm{cyc}}(1-z)=1,
 \qquad
 \Gamma_{\PP}^{\mathrm{cyc}}(1/2)=1.
 \end{gathered}
 }
 \label{eq:section10-structural-summary}
\end{equation}

The gamma-lift torsor has now acquired its distinguished point.  The next section studies the resulting canonical analytic defect
\[
 H_{\PP}^{\mathrm{cyc}}(z+1)
 -\log\Gamma(z+1)-C_{\PP}z
\]
and separates its orbitwise analytic normalization from its completely unsmoothed arithmetic trace.

\section{The canonical analytic Stirling defect}
\label{sec:canonical-analytic-stirling-defect}

The preceding section selected a unique normalized logarithmic lift
\[
 H_{\PP}^{\mathrm{cyc}}\in\Lift(L_{\PP}^{\mathrm{cyc}})
\]
by comparing its translation orbits with
\[
 \log\Gamma(z)+C_{\PP}(z-1).
\]
We now package the difference between these two objects as a canonical analytic invariant of the completed prime factorial calculus.

There are two forms of this invariant.  Its logarithmic form is holomorphic on the right half-plane cut out by the classical gamma function and has the unsmoothed prime-factorial fluctuation as its integer trace.  Its multiplicative form extends to an entire function and records the same information without choosing a global logarithm.  The distinction is useful: the logarithmic defect is the object selected by orbitwise Stirling rigidity, whereas the multiplicative defect exhibits the global zero structure inherited from the reciprocal classical gamma function.

This section also separates three statements that should not be conflated:
\begin{enumerate}[label=\textup{(\roman*)}]
 \item the discrete prime-factorial trace is fixed by the recurrence and is insensitive to periodic gauge;
 \item the affine coefficient is fixed by the weak Stirling law;
 \item the analytic extension between consecutive integers is fixed by orbitwise rigidity.
\end{enumerate}
The resulting object is the bridge from the categorical theory of Parts~I--II to the unsmoothed arithmetic analysis of Part~III.

\subsection{The logarithmic defect and the entire Stirling quotient}
\label{subsec:section11-logarithmic-and-multiplicative-defects}

Let \(\Logp\Gamma\) denote the holomorphic logarithm of \(\Gamma\) on the half-plane \(\Re z>0\), real on the positive axis.  Define
\begin{equation}
 \boxed{
 \mathscr R_{\PP}(z)
 =H_{\PP}^{\mathrm{cyc}}(z+1)
  -\Logp\Gamma(z+1)-C_{\PP}z,
 \qquad \Re z>-1.
 }
 \label{eq:section11-logarithmic-defect}
\end{equation}
We call \(\mathscr R_{\PP}\) the \emph{canonical analytic Stirling defect} of the prime Bhargava factorial.

The corresponding multiplicative object is
\begin{equation}
 \boxed{
 \mathscr S_{\PP}(z)
 =\e^{-C_{\PP}z}
 \frac{\Gamma_{\PP}^{\mathrm{cyc}}(z+1)}{\Gamma(z+1)}.
 }
 \label{eq:section11-multiplicative-defect}
\end{equation}
Although the quotient notation in \eqref{eq:section11-multiplicative-defect} is convenient, it is understood as
\[
 \mathscr S_{\PP}(z)
 =\e^{-C_{\PP}z}
 \Gamma_{\PP}^{\mathrm{cyc}}(z+1)
 \frac1{\Gamma(z+1)}.
\]
Thus it is defined globally, without choosing a logarithm of the classical gamma function.

\begin{proposition}[Global analytic structure]
\label{prop:section11-global-analytic-structure}
The function \(\mathscr S_{\PP}\) is entire, satisfies
\begin{equation}
 \mathscr S_{\PP}(0)=1,
 \label{eq:section11-multiplicative-normalization}
\end{equation}
and has simple zeros precisely at
\begin{equation}
 -1,-2,-3,\ldots.
 \label{eq:section11-zero-set}
\end{equation}
On the half-plane \(\Re z>-1\),
\begin{equation}
 \boxed{
 \mathscr S_{\PP}(z)=\exp\bigl(\mathscr R_{\PP}(z)\bigr).
 }
 \label{eq:section11-exponential-relation}
\end{equation}
In particular, \(\mathscr R_{\PP}\) is the unique logarithm of \(\mathscr S_{\PP}\) on \(\Re z>-1\) normalized by \(\mathscr R_{\PP}(0)=0\).
\end{proposition}

\begin{proof}
The prime gamma function \(\Gamma_{\PP}^{\mathrm{cyc}}\) is entire and zero-free.  The reciprocal gamma function \(1/\Gamma(z+1)\) is entire and has simple zeros exactly at the negative integers \(-1,-2,\ldots\).  The exponential factor never vanishes, proving the first two assertions and the zero statement.  On \(\Re z>-1\), exponentiating \eqref{eq:section11-logarithmic-defect} gives \eqref{eq:section11-exponential-relation}.  Finally,
\[
 \mathscr R_{\PP}(0)
 =H_{\PP}^{\mathrm{cyc}}(1)-\Logp\Gamma(1)=0,
\]
which fixes the logarithm uniquely on the simply connected half-plane.
\end{proof}

At nonnegative integers, the two defects recover the normalized prime factorial exactly.

\begin{corollary}[Integer trace]
\label{cor:section11-integer-trace}
For every \(n\ge0\),
\begin{equation}
 \boxed{
 \mathscr R_{\PP}(n)
 =\log(n+1)!_{\PP}-\log n!-C_{\PP}n
 =:\cR_{\PP}(n),
 }
 \label{eq:section11-integer-logarithmic-trace}
\end{equation}
and
\begin{equation}
 \boxed{
 \mathscr S_{\PP}(n)
 =\frac{(n+1)!_{\PP}}{n!\,\e^{C_{\PP}n}}.
 }
 \label{eq:section11-integer-multiplicative-trace}
\end{equation}
\end{corollary}

\begin{proof}
Use
\[
 H_{\PP}^{\mathrm{cyc}}(n+1)=\log(n+1)!_{\PP},
 \qquad
 \Gamma(n+1)=n!.
\]
\end{proof}

\subsection{A logarithmic cyclotomic integral and the prime Hankel problem}
\label{subsec:section11-prime-hankel-problem}

Let $I=[-1/2,1/2]$.  Write $\mu_d$ for the symmetric quadrature probability measure of conductor $d$, so that
\[
 \cQ_df=\int_I f(t)\,\dd\mu_d(t),
\]
let $\lambda$ be Lebesgue measure on $I$, and put $\nu_d=\mu_d-\lambda$.  The midpoint and trapezoidal error estimates give
\begin{equation}
 \left|\int_I f\,\dd\nu_d\right|
 \ll \frac{\|f''\|_{L^\infty(I)}}{d^2}
 \qquad(f\in C^2(I)).
 \label{eq:section11-quadrature-distribution-bound}
\end{equation}
Since
\[
 \sum_{p,j}\frac{\log p}{((p-1)p^j)^2}<\infty,
\]
the formula
\begin{equation}
 \boxed{
 \left\langle\mathfrak N_{\PP},f\right\rangle
 :=\sum_{p\in\PP}\sum_{j\ge0}
  (\log p)\int_I f(t)\,\dd\nu_{(p-1)p^j}(t)
 }
 \label{eq:section11-prime-cyclotomic-distribution}
\end{equation}
defines a continuous real even distribution of order at most two on $I$.

\begin{proposition}[Exact logarithmic cyclotomic representation]
\label{prop:section11-logarithmic-cyclotomic-representation}
For the difference kernel $K_z$ of Section~9,
\begin{equation}
 \boxed{
 H_{\PP}^{\mathrm{cyc}}(z)
 =\left\langle\mathfrak N_{\PP},K_z\right\rangle
 \qquad(z\in\CC).
 }
 \label{eq:section11-logarithmic-cyclotomic-representation}
\end{equation}
Consequently,
\begin{equation}
 \boxed{
 \frac1{\Gamma_{\PP}^{\mathrm{cyc}}(z)}
 =\exp\!\left(
  -\left\langle\mathfrak N_{\PP},K_z\right\rangle
 \right).
 }
 \label{eq:section11-reciprocal-exponential-integral}
\end{equation}
The pairing in \eqref{eq:section11-logarithmic-cyclotomic-representation} is entire in $z$, locally uniformly together with all derivatives.
\end{proposition}

\begin{proof}
Appendix~B identifies
\[
 \phi_d(z)=\int_IK_z(t)\,\dd\nu_d(t).
\]
Substituting this identity into the normally convergent prime-layer expansion
\[
 H_{\PP}^{\mathrm{cyc}}(z)
 =\sum_{p,j}(\log p)\phi_{(p-1)p^j}(z)
\]
gives \eqref{eq:section11-logarithmic-cyclotomic-representation}.  Estimate \eqref{eq:section11-quadrature-distribution-bound}, applied to the $z$-derivatives of $K_z$, gives normal convergence on compact subsets.  Exponentiation proves \eqref{eq:section11-reciprocal-exponential-integral}.
\end{proof}

The representation above is canonical, but it is exponential in a linear functional.  It is therefore different in kind from Euler's positive Mellin integral.  In fact, positivity is arithmetically impossible.

\begin{proposition}[Hankel obstruction to positive Euler--Mellin integrals]
\label{prop:section11-positive-mellin-obstruction}
There is no positive Borel measure $\mu$ on $[0,\infty)$, with all moments finite, for which
\begin{equation}
 \Gamma_{\PP}^{\mathrm{cyc}}(z+1)
 =\int_0^\infty t^z\,\dd\mu(t)
 \label{eq:section11-forbidden-gamma-mellin}
\end{equation}
holds at every nonnegative integer $z$.  There is likewise no such positive measure $\rho$ satisfying
\begin{equation}
 \frac1{\Gamma_{\PP}^{\mathrm{cyc}}(z+1)}
 =\int_0^\infty t^z\,\dd\rho(t)
 \label{eq:section11-forbidden-reciprocal-mellin}
\end{equation}
at every nonnegative integer.
\end{proposition}

\begin{proof}
If \eqref{eq:section11-forbidden-gamma-mellin} held, then
\[
 m_n:=\Gamma_{\PP}^{\mathrm{cyc}}(n+1)=(n+1)!_{\PP}
\]
would be a Stieltjes moment sequence~\cite{ShohatTamarkin1943,Akhiezer1965}.  Positivity of the shifted $2\times2$ Hankel matrix would require
\[
 m_2^2\le m_1m_3.
\]
But the prime factorial values are
\[
 m_1=2!_{\PP}=2,
 \qquad
 m_2=3!_{\PP}=24,
 \qquad
 m_3=4!_{\PP}=48,
\]
and hence
\[
 \det\begin{pmatrix}m_1&m_2\\m_2&m_3\end{pmatrix}
 =2\cdot48-24^2=-480<0.
\]
For the reciprocal sequence
\[
 a_n:=\frac1{(n+1)!_{\PP}},
\]
the classical moment criteria~\cite{ShohatTamarkin1943,Akhiezer1965} imply that ordinary Hankel positivity would require $a_1^2\le a_0a_2$.  Instead,
\[
 \det\begin{pmatrix}a_0&a_1\\a_1&a_2\end{pmatrix}
 =\frac1{24}-\frac14=-\frac5{24}<0.
\]
Thus neither sequence can arise from a positive measure on $[0,\infty)$;
compare the classical moment and transform criteria in
\cite{ShohatTamarkin1943,Akhiezer1965,Widder1941}.
\end{proof}

\begin{remark}[The prime Hankel conjecture revisited]
\label{rem:section11-prime-hankel-problem}
Proposition~\ref{prop:section11-positive-mellin-obstruction} excludes only positive real measures.  It does not exclude signed or complex Mellin kernels, distributions, Jackson-type discrete integrals, or Hankel contour representations.  Conjecture~\ref{conj:intro-prime-hankel} proposes a canonical linearization of \eqref{eq:section11-reciprocal-exponential-integral} of the form
\begin{equation}
 \frac1{\Gamma_{\PP}^{\mathrm{cyc}}(z)}
 =\frac1{2\pi i}\int_{\mathscr C}
   \mathcal K_{\PP}(w)w^{-z}\,\dd w,
 \label{eq:section11-prime-hankel-target}
\end{equation}
initially on a nonempty vertical strip and then by entire continuation.  The substantive requirement is not merely the existence of some inverse Mellin transform: the contour and kernel should be canonically determined by the prime-cyclotomic layer data, and their transformations should display the recurrence and reflection laws directly.
\end{remark}

\subsection{The affine-periodic defect torsor}
\label{subsec:section11-defect-torsor}

The selected defect sits inside a larger affine family.  For
\[
 c\in\CC,
 \qquad
 H\in\Lift(L_{\PP}^{\mathrm{cyc}}),
\]
define
\begin{equation}
 \mathscr R_{c,H}(z)
 =H(z+1)-\Logp\Gamma(z+1)-cz,
 \qquad \Re z>-1.
 \label{eq:section11-general-defect}
\end{equation}
Let
\begin{equation}
 \mathfrak G
 =\left\{az+P(z):a\in\CC,\ P\in\Per_0\right\}.
 \label{eq:section11-affine-periodic-group}
\end{equation}
Because every \(P\in\Per_0\) satisfies \(P(0)=P(1)=0\), the sum in \eqref{eq:section11-affine-periodic-group} is direct.

\begin{proposition}[Defect-torsor decomposition]
\label{prop:section11-defect-torsor-decomposition}
Every defect \(\mathscr R_{c,H}\) has a unique representation
\begin{equation}
 \boxed{
 \mathscr R_{c,H}(z)
 =\mathscr R_{\PP}(z)
 +(C_{\PP}-c)z+P(z),
 \qquad P\in\Per_0.
 }
 \label{eq:section11-defect-decomposition}
\end{equation}
Consequently, the family of all defects \(\mathscr R_{c,H}\) is a torsor under the additive group \(\mathfrak G\).
\end{proposition}

\begin{proof}
There is a unique \(P\in\Per_0\) such that
\[
 H=H_{\PP}^{\mathrm{cyc}}+P.
\]
Since \(P(z+1)=P(z)\), substitution into \eqref{eq:section11-general-defect} gives \eqref{eq:section11-defect-decomposition}.  If
\[
 az+P(z)=0,
\]
then evaluation at positive integers gives \(an=0\), hence \(a=0\), and then \(P=0\).  This proves uniqueness and freeness of the action.
\end{proof}

The orbitwise profile from Section~10 may be expressed directly in terms of a defect.  Put
\begin{equation}
 \mathscr U_N(\mathscr R)(u)
 =\sum_{n=0}^{N-1}\omega_{N,n}
 \bigl(\mathscr R(n+u)-\mathscr R(n)\bigr),
 \qquad 0\le u\le1.
 \label{eq:section11-defect-profile}
\end{equation}
For \(\mathscr R=\mathscr R_{c,H}\), this is exactly the transverse profile of \(H\) relative to the comparison model \(\Logp\Gamma(z)+c(z-1)\).

\begin{theorem}[Canonical trivialization of the defect torsor]
\label{thm:section11-canonical-trivialization}
For every \(c\in\CC\) and \(H\in\Lift(L_{\PP}^{\mathrm{cyc}})\), write \(\mathscr R_{c,H}\) as in \eqref{eq:section11-defect-decomposition}.  Then
\begin{equation}
 \boxed{
 \mathscr U_N(\mathscr R_{c,H})(u)
 \longrightarrow
 (C_{\PP}-c)u+P(u)
 }
 \label{eq:section11-all-defect-profile-limits}
\end{equation}
uniformly for \(0\le u\le1\).  In particular,
\begin{equation}
 \mathscr U_N(\mathscr R_{c,H})\longrightarrow0
 \quad\Longleftrightarrow\quad
 (c,H)=(C_{\PP},H_{\PP}^{\mathrm{cyc}}).
 \label{eq:section11-unique-flat-defect}
\end{equation}
Thus \(\mathscr R_{\PP}\) is the unique orbit-flat element of the affine-periodic defect torsor.
\end{theorem}

\begin{proof}
The profile is additive.  For the selected defect, the orbitwise Stirling theorem of Section~10 gives
\[
 \mathscr U_N(\mathscr R_{\PP})\longrightarrow0
\]
uniformly.  Regularity of the averaging weights gives
\[
 \mathscr U_N(az)(u)=au,
 \qquad
 \mathscr U_N(P)(u)=P(u)-P(0)=P(u).
\]
Combining these identities with \eqref{eq:section11-defect-decomposition} proves \eqref{eq:section11-all-defect-profile-limits}.  The affine and periodic directions intersect trivially, so the limiting profile vanishes only when \(c=C_{\PP}\) and \(P=0\).
\end{proof}

\begin{remark}[What is actually made unique]
\label{rem:section11-what-is-unique}
Periodic gauge does not change positive-integer values: if \(P\in\Per_0\), then \(P(n)=0\) for every integer \(n\).  Hence all normalized lifts of the same recurrence symbol have the same discrete trace.  Orbitwise rigidity selects the analytic interpolation between the integers, not the integer sequence itself.  The affine parameter is different: changing \(c\) changes the integer trace by a linear term.  The weak Stirling law therefore selects \(C_{\PP}\), while orbitwise rigidity selects the periodic gauge.
\end{remark}

\subsection{The centered recurrence increment}
\label{subsec:section11-centered-recurrence-increment}

The additive defect satisfies a first-order difference equation.  On \(\Re z>0\), define
\begin{equation}
 \boxed{
 \delta_{\PP}(z)
 =L_{\PP}^{\mathrm{cyc}}(z)-\Logp z-C_{\PP}.
 }
 \label{eq:section11-centered-increment-symbol}
\end{equation}
This is the centered logarithmic recurrence increment.

\begin{proposition}[Defect recurrence]
\label{prop:section11-defect-recurrence}
For \(\Re z>0\),
\begin{equation}
 \boxed{
 \mathscr R_{\PP}(z)-\mathscr R_{\PP}(z-1)
 =\delta_{\PP}(z).
 }
 \label{eq:section11-additive-defect-recurrence}
\end{equation}
Equivalently, the entire multiplicative defect satisfies
\begin{equation}
 \boxed{
 z\,\mathscr S_{\PP}(z)
 =\e^{-C_{\PP}}
 R_{\PP}^{\mathrm{cyc}}(z)
 \mathscr S_{\PP}(z-1),
 }
 \label{eq:section11-multiplicative-defect-recurrence}
\end{equation}
where
\[
 R_{\PP}^{\mathrm{cyc}}(z)
 =\exp\bigl(L_{\PP}^{\mathrm{cyc}}(z)\bigr).
\]
Identity \eqref{eq:section11-multiplicative-defect-recurrence} holds on all of \(\CC\) by analytic continuation.
\end{proposition}

\begin{proof}
Using the recurrence equation for \(H_{\PP}^{\mathrm{cyc}}\) and the classical identity
\[
 \Logp\Gamma(z+1)-\Logp\Gamma(z)=\Logp z
\]
on \(\Re z>0\), we obtain \eqref{eq:section11-additive-defect-recurrence}.  Exponentiation gives the multiplicative identity on \(\Re z>0\).  Both sides of \eqref{eq:section11-multiplicative-defect-recurrence} are entire, so the identity extends globally.
\end{proof}

The exact prime-factorial ratio from Section~2 now gives a completely arithmetic formula for the integer increments.

\begin{corollary}[Arithmetic increment formula]
\label{cor:section11-arithmetic-increment-formula}
For every integer \(n\ge1\),
\begin{equation}
 \boxed{
 \delta_{\PP}(n)
 =\sum_{p-1\mid n}\bigl(v_p(n)+1\bigr)\log p
  -\log n-C_{\PP}.
 }
 \label{eq:section11-integer-increment-formula}
\end{equation}
Moreover,
\begin{equation}
 \boxed{
 \cR_{\PP}(n)
 =\sum_{m=1}^{n}\delta_{\PP}(m).
 }
 \label{eq:section11-defect-as-partial-sum}
\end{equation}
\end{corollary}

\begin{proof}
At positive integers,
\[
 L_{\PP}^{\mathrm{cyc}}(n)
 =\log\frac{(n+1)!_{\PP}}{n!_{\PP}}
 =\sum_{p-1\mid n}\bigl(v_p(n)+1\bigr)\log p.
\]
This proves \eqref{eq:section11-integer-increment-formula}.  Since \(\cR_{\PP}(0)=0\), telescoping \eqref{eq:section11-additive-defect-recurrence} proves \eqref{eq:section11-defect-as-partial-sum}.
\end{proof}

\subsection{Weak size and mean-zero normalization}
\label{subsec:section11-weak-size}

The weak prime Stirling formula has an immediate interpretation in terms of the canonical defect.

\begin{proposition}[Sublinear integer trace]
\label{prop:section11-sublinear-trace}
As \(n\to\infty\),
\begin{equation}
 \boxed{
 \cR_{\PP}(n)=o(n).
 }
 \label{eq:section11-sublinear-trace}
\end{equation}
Equivalently, the centered recurrence increments have vanishing Ces\`aro mean:
\begin{equation}
 \boxed{
 \frac1n\sum_{m=1}^{n}\delta_{\PP}(m)
 \longrightarrow0.
 }
 \label{eq:section11-mean-zero-increments}
\end{equation}
\end{proposition}

\begin{proof}
Equation \eqref{eq:section11-sublinear-trace} is exactly the weak Stirling law
\[
 \log(n+1)!_{\PP}=\log n!+C_{\PP}n+o(n).
\]
The second assertion follows from \eqref{eq:section11-defect-as-partial-sum}.
\end{proof}

\begin{remark}[Sublinear does not mean smoothed]
\label{rem:section11-sublinear-not-smoothed}
The estimate \(\cR_{\PP}(n)=o(n)\) is a statement about the raw sequence.  No orbit average appears in its definition.  The twice-iterated Ces\`aro operator was used only to select the distinguished analytic lift within its periodic torsor.  Once the lift is selected, the integer values in \eqref{eq:section11-integer-logarithmic-trace} are left untouched.
\end{remark}

\subsection{Prime-local decomposition of the integer trace}
\label{subsec:section11-prime-local-decomposition}

The canonical defect admits an exact decomposition into prime-local fluctuations.  This formula is the arithmetic point of departure for Part~III.

For \(p\in\PP\) and \(n\ge0\), define
\begin{align}
 \cR_p(n)
 =\log p\Bigg(&
 \sum_{j\ge0}
 \left\lfloor\frac{n}{(p-1)p^j}\right\rfloor
 -\sum_{k\ge1}
 \left\lfloor\frac{n}{p^k}\right\rfloor
 -\frac{n}{(p-1)^2}
 \Bigg).
 \label{eq:section11-local-prime-defect}
\end{align}

\begin{proposition}[Exact prime decomposition]
\label{prop:section11-exact-prime-decomposition}
For every \(n\ge0\),
\begin{equation}
 \boxed{
 \cR_{\PP}(n)=\sum_{p\in\PP}\cR_p(n),
 }
 \label{eq:section11-prime-decomposition}
\end{equation}
and the series is absolutely convergent.  Moreover,
\begin{equation}
 \boxed{
 \cR_p(n)
 =\log p\left(
 \sum_{k\ge1}\left\{\frac{n}{p^k}\right\}
 -\sum_{j\ge0}\left\{\frac{n}{(p-1)p^j}\right\}
 \right).
 }
 \label{eq:section11-fractional-part-decomposition}
\end{equation}
\end{proposition}

\begin{proof}
The prime-factorial and classical factorial identities are
\begin{align*}
 \log(n+1)!_{\PP}
 &=\sum_{p\in\PP}\log p
   \sum_{j\ge0}
   \left\lfloor\frac{n}{(p-1)p^j}\right\rfloor,\\
 \log n!
 &=\sum_{p\in\PP}\log p
   \sum_{k\ge1}
   \left\lfloor\frac{n}{p^k}\right\rfloor.
\end{align*}
Also,
\[
 C_{\PP}n
 =\sum_{p\in\PP}\frac{n\log p}{(p-1)^2}.
\]
For \(p>n+1\), both floor sums vanish, and the remaining tail is bounded by
\[
 n\sum_{p>n+1}\frac{\log p}{(p-1)^2}<\infty.
\]
This proves absolute convergence and \eqref{eq:section11-prime-decomposition}.

For each fixed prime,
\begin{align*}
 \sum_{j\ge0}\frac{n}{(p-1)p^j}
 -\sum_{k\ge1}\frac{n}{p^k}
 &=\frac{np}{(p-1)^2}-\frac{n}{p-1}\\
 &=\frac{n}{(p-1)^2}.
\end{align*}
Replacing each floor by \(x-\{x\}\) therefore gives \eqref{eq:section11-fractional-part-decomposition}.
\end{proof}

\begin{remark}[Why the prime decomposition is preferable to a formal conductor difference]
\label{rem:section11-primewise-grouping}
If the prime-layer conductor weights and the von Mangoldt weights are subtracted before grouping by prime, the corresponding \(1/d\)-series is only conditionally organized: each side separately has logarithmic divergence.  Formula \eqref{eq:section11-local-prime-defect} performs the cancellation prime-by-prime, where
\[
 \sum_{j\ge0}\frac1{(p-1)p^j}
 -\sum_{k\ge1}\frac1{p^k}
 =\frac1{(p-1)^2}.
\]
This is the canonical absolutely convergent renormalization used in the variance calculation.
\end{remark}

\subsection{The analytic and arithmetic roles of the defect}
\label{subsec:section11-analytic-arithmetic-roles}

The same canonical object now carries two complementary kinds of information.

On the analytic side,
\begin{equation}
 \mathscr R_{\PP}(z)
 =H_{\PP}^{\mathrm{cyc}}(z+1)
  -\Logp\Gamma(z+1)-C_{\PP}z
 \label{eq:section11-analytic-role}
\end{equation}
is the unique orbit-flat logarithmic representative of the completed recurrence.  Its exponential extends globally as an entire function whose zeros are inherited exactly from \(1/\Gamma(z+1)\).

On the arithmetic side,
\begin{equation}
 \cR_{\PP}(n)
 =\sum_{p}\log p\left(
 \sum_{k\ge1}\left\{\frac{n}{p^k}\right\}
 -\sum_{j\ge0}\left\{\frac{n}{(p-1)p^j}\right\}
 \right)
 \label{eq:section11-arithmetic-role}
\end{equation}
is a completely unsmoothed sequence.  Its first difference is the centered divisor-supported weight \eqref{eq:section11-integer-increment-formula}, while its second moment measures the cumulative interaction of the prime-local fractional-part fluctuations.

The analytic selection theorem does not enter the later variance calculation as an averaging device.  Rather, it certifies that the raw sequence in \eqref{eq:section11-arithmetic-role} is the integer trace of a uniquely distinguished analytic object.  This is the sense in which the categorical theory motivates the unsmoothed Stirling expansion instead of replacing it.

\subsection{Structural conclusion and transition to Part III}
\label{subsec:section11-structural-conclusion}

The conclusions of Parts~I--II may be compressed into the diagram
\begin{equation}
 \boxed{
 \begin{gathered}
 \text{prime valuation layers}
 \longrightarrow
 \text{filtered factorial calculus}
 \longrightarrow
 \Lift(L_{\PP}^{\mathrm{cyc}}),
 \\
 \text{weak Stirling coefficient }C_{\PP}
 \quad+\quad
 \text{orbitwise rigidity}
 \longrightarrow
 H_{\PP}^{\mathrm{cyc}},
 \\
 \mathscr R_{\PP}(z)
 =H_{\PP}^{\mathrm{cyc}}(z+1)
  -\Logp\Gamma(z+1)-C_{\PP}z,
 \\
 \mathscr S_{\PP}(z)
 =\e^{-C_{\PP}z}
  \frac{\Gamma_{\PP}^{\mathrm{cyc}}(z+1)}{\Gamma(z+1)},
 \\
 \mathscr R_{\PP}(n)=\cR_{\PP}(n)
 =\log(n+1)!_{\PP}-\log n!-C_{\PP}n,
 \\
 \cR_{\PP}(n)
 =\sum_{p\in\PP}\cR_p(n),
 \qquad
 \cR_{\PP}(n)=o(n).
 \end{gathered}
 }
 \label{eq:section11-structural-summary}
\end{equation}

This closes the categorical and gamma-theoretic half of the paper.  The distinguished prime gamma object now occupies, for the prime Bhargava factorial, the same structural position that Euler's gamma function occupies among solutions of the classical recurrence: the recurrence supplies the family, and a normalization theorem identifies the canonical member.  Here the selecting theorem is orbitwise rather than convexity-based, but its conceptual role is the one exemplified by Bohr--Mollerup.

Part~III begins by reorganizing the exact prime-local formula \eqref{eq:section11-prime-decomposition} into a form adapted to second moments.  The next section derives the arithmetic expansion, isolates the diagonal variance constant, and introduces the covariance kernels that govern the unsmoothed fluctuation problem.

\section{Exact arithmetic expansion and covariance kernels}
\label{sec:exact-arithmetic-expansion}

The canonical analytic defect constructed in Part~II has the unsmoothed integer trace
\[
 \cR_{\PP}(n)
 =\log (n+1)!_{\PP}-\log n!-C_{\PP}n.
\]
Section~11 expressed this trace as an absolutely convergent sum of prime-local fractional-part fluctuations.  We now reorganize that formula in a way adapted to the sharp second moment
\[
 \cV(X)=\sum_{n\le X}\cR_{\PP}(n)^2.
\]

Three structures emerge.  First, the classical prime-power moduli and the prime-factorial moduli pair naturally into adjacent projective waves.  Second, the square of the resulting expansion separates into the geometric diagonal, a projective diagonal, and a genuinely off-projective sector.  Third, the dominant layer possesses an exact interval-overlap representation whose support is an affine strip around the rational line $mp=\ell q$.  This is the arithmetic geometry that will be analyzed in the remaining sections.

The present section also evaluates the complete geometric diagonal.  Its leading term is already the constant appearing in the final variance law:
\[
 -\frac43\zeta\!\left(-\frac12\right)X^{3/2}\log X.
\]
Thus every later argument is an off-diagonal cancellation problem; it does not manufacture the principal constant.

\subsection{Paired prime-power layers}
\label{subsec:paired-prime-power-layers}

For a prime $p$ and an integer $j\ge0$, put
\begin{equation}
 a_{p,j}=p^{j+1},
 \qquad
 b_{p,j}=(p-1)p^j,
 \qquad
 \lambda=(p,j).
 \label{eq:section12-paired-moduli}
\end{equation}
The two moduli have the common factor $p^j$ and satisfy
\[
 a_{p,j}-b_{p,j}=p^j,
 \qquad
 [a_{p,j}:b_{p,j}]=[p:p-1].
\]
Define the associated fractional-part wave by
\begin{equation}
 W_{p,j}(n)
 =\left\{\frac{n}{a_{p,j}}\right\}
  -\left\{\frac{n}{b_{p,j}}\right\}.
 \label{eq:section12-layer-wave}
\end{equation}

\begin{proposition}[Exact paired-layer expansion]
\label{prop:section12-exact-paired-expansion}
For every integer $n\ge0$,
\begin{equation}
 \boxed{
 \cR_{\PP}(n)
 =\sum_{j\ge0}\cR_j(n),
 \qquad
 \cR_j(n)=\sum_{p\in\PP}(\log p)W_{p,j}(n).
 }
 \label{eq:section12-paired-expansion}
\end{equation}
The double series in \eqref{eq:section12-paired-expansion} converges absolutely after the two moduli attached to each $(p,j)$ are paired as in \eqref{eq:section12-layer-wave}.
\end{proposition}

\begin{proof}
The prime-local formula of Section~11 is
\[
 \cR_{\PP}(n)
 =\sum_p\log p\left(
 \sum_{k\ge1}\left\{\frac{n}{p^k}\right\}
 -\sum_{j\ge0}\left\{\frac{n}{(p-1)p^j}\right\}
 \right).
\]
Pair the term $k=j+1$ in the first sum with the term $j$ in the second.  This gives \eqref{eq:section12-paired-expansion} formally.

It remains to justify the rearrangement.  Once both moduli in \eqref{eq:section12-paired-moduli} exceed $n$, the fractional parts are linear and
\[
 |W_{p,j}(n)|
 =n\left|
 \frac1{p^{j+1}}-\frac1{(p-1)p^j}
 \right|
 =\frac{n}{p^{j+1}(p-1)}.
\]
Only finitely many pairs fail this condition, whereas
\[
 \sum_p\sum_{j\ge0}
 \frac{n\log p}{p^{j+1}(p-1)}
 =n\sum_p\frac{\log p}{(p-1)^2}
 =nC_{\PP}<\infty.
\]
This proves paired absolute convergence.
\end{proof}

The layer $j=0$ is distinguished:
\begin{equation}
 \boxed{
 \cR_0(n)
 =\sum_p\log p
 \left(
 \left\{\frac np\right\}
 -\left\{\frac n{p-1}\right\}
 \right).
 }
 \label{eq:section12-dominant-layer}
\end{equation}
It will supply the leading diagonal variance.  The higher layers are required for the exact defect, but their geometric diagonals occur at smaller powers of $X$.

\subsection{The sharp covariance decomposition}
\label{subsec:sharp-covariance-decomposition}

For layer indices $\lambda=(p,j)$ and $\mu=(q,k)$, define the sharp covariance kernel
\begin{equation}
 \cC_X(\lambda,\mu)
 =\sum_{n\le X}W_{p,j}(n)W_{q,k}(n).
 \label{eq:section12-sharp-covariance-kernel}
\end{equation}
Paired absolute convergence on the finite set $1\le n\le X$ permits termwise expansion of the second moment.

\begin{theorem}[Diagonal, projective-diagonal, and off-projective splitting]
\label{thm:section12-three-sector-splitting}
For every $X\ge1$,
\begin{equation}
 \boxed{
 \cV(X)=\cD(X)+\cP(X)+\cO(X),
 }
 \label{eq:section12-three-sector-splitting}
\end{equation}
where
\begin{align}
 \cD(X)
 &=\sum_p\sum_{j\ge0}
   (\log p)^2\cC_X((p,j),(p,j)),
 \label{eq:section12-geometric-diagonal}\\
 \cP(X)
 &=2\sum_p\sum_{0\le j<k}
   (\log p)^2\cC_X((p,j),(p,k)),
 \label{eq:section12-projective-diagonal}\\
 \cO(X)
 &=\sum_{p\ne q}\sum_{j,k\ge0}
   (\log p)(\log q)\cC_X((p,j),(q,k)).
 \label{eq:section12-off-projective-sector}
\end{align}
All three expressions are understood as limits of finite prime-and-layer truncations.
\end{theorem}

\begin{proof}
Expand the square of \eqref{eq:section12-paired-expansion} and group the ordered pairs of layers according to whether they are equal, have the same prime but different layer index, or have different primes.  The finite truncations satisfy \eqref{eq:section12-three-sector-splitting} exactly, and paired absolute convergence gives the limit.
\end{proof}

The terminology in \eqref{eq:section12-projective-diagonal} is literal.  For fixed $p$, all layers determine the same projective point
\begin{equation}
 [a_{p,j}:b_{p,j}]=[p:p-1],
 \label{eq:section12-projective-point}
\end{equation}
although their conductors differ by powers of $p$.  Thus $j=k$ is the geometric diagonal, $p=q$ with $j\ne k$ is the projective diagonal, and $p\ne q$ is genuinely off-projective.

\subsection{Discrete sawtooths and the periodic bulk kernel}
\label{subsec:periodic-bulk-kernel}

For an integer $m\ge1$, define the mean-zero discrete sawtooth
\begin{equation}
 \sigma_m(n)
 =\left\{\frac nm\right\}-\frac{m-1}{2m}.
 \label{eq:section12-discrete-sawtooth}
\end{equation}
It is $m$-periodic and has mean zero over a complete period.

\begin{lemma}[Complete-period covariance]
\label{lem:section12-complete-period-covariance}
Let $a,b\ge1$, let $g=(a,b)$, and let $L=\lcm(a,b)$.  Then
\begin{equation}
 \boxed{
 \frac1L\sum_{n=1}^{L}\sigma_a(n)\sigma_b(n)
 =\kappa(a,b)
 :=\frac{g^2-1}{12ab}.
 }
 \label{eq:section12-basic-bulk-kernel}
\end{equation}
\end{lemma}

\begin{proof}
Choose $N$ uniformly modulo $L$ and put $U=N\bmod g$.  Conditional on $U$, the residual coordinates modulo $a/g$ and $b/g$ are independent and uniform, because these two moduli are coprime.  A direct average gives
\[
 \mathbb E(\sigma_a(N)\mid U)
 =\frac{U-(g-1)/2}{a},
 \qquad
 \mathbb E(\sigma_b(N)\mid U)
 =\frac{U-(g-1)/2}{b}.
\]
The conditional covariance is zero.  Since a uniform variable on $\{0,\ldots,g-1\}$ has variance $(g^2-1)/12$, taking expectations proves \eqref{eq:section12-basic-bulk-kernel}.
\end{proof}

For $\lambda=(p,j)$, write
\begin{equation}
 \mu_\lambda
 =\frac{a_{p,j}-b_{p,j}}{2a_{p,j}b_{p,j}}
 =\frac1{2p^{j+1}(p-1)}.
 \label{eq:section12-wave-mean}
\end{equation}
Then
\begin{equation}
 W_\lambda(n)
 =\sigma_{a_\lambda}(n)-\sigma_{b_\lambda}(n)+\mu_\lambda.
 \label{eq:section12-centered-wave-decomposition}
\end{equation}
The centered bulk covariance of two waves is therefore
\begin{align}
 \cK(\lambda,\mu)
 ={}&\kappa(a_\lambda,a_\mu)
 -\kappa(a_\lambda,b_\mu)
 -\kappa(b_\lambda,a_\mu)
 +\kappa(b_\lambda,b_\mu),
 \label{eq:section12-wave-bulk-kernel}
\end{align}
and the uncentered complete-period mean is
\begin{equation}
 \cM(\lambda,\mu)
 =\cK(\lambda,\mu)+\mu_\lambda\mu_\mu.
 \label{eq:section12-uncentered-bulk-kernel}
\end{equation}

Let
\begin{equation}
 L(\lambda,\mu)
 =\lcm(a_\lambda,b_\lambda,a_\mu,b_\mu).
 \label{eq:section12-joint-period}
\end{equation}
Since the product $W_\lambda W_\mu$ has period $L(\lambda,\mu)$, complete periods may be removed exactly.

\begin{proposition}[Bulk plus sharp-cutoff remainder]
\label{prop:section12-bulk-remainder-splitting}
For all layer indices $\lambda,\mu$,
\begin{equation}
 \boxed{
 \cC_X(\lambda,\mu)
 =X\cM(\lambda,\mu)+\cE_X(\lambda,\mu),
 }
 \label{eq:section12-bulk-remainder-splitting}
\end{equation}
where $\cE_X(\lambda,\mu)$ depends only on $X\bmod L(\lambda,\mu)$ and satisfies
\begin{equation}
 |\cE_X(\lambda,\mu)|\ll L(\lambda,\mu).
 \label{eq:section12-crude-remainder-bound}
\end{equation}
\end{proposition}

\begin{proof}
Split $[1,X]$ into complete blocks of length $L(\lambda,\mu)$ and one terminal block.  Lemma~\ref{lem:section12-complete-period-covariance}, together with \eqref{eq:section12-centered-wave-decomposition}, gives the mean \eqref{eq:section12-uncentered-bulk-kernel}.  The terminal block contains fewer than $L(\lambda,\mu)$ terms and the waves are bounded.
\end{proof}

\begin{remark}[Why the periodic estimate is not the final argument]
\label{rem:section12-period-bound-too-crude}
The bound \eqref{eq:section12-crude-remainder-bound} is individually sharp enough but globally useless when the joint period is large.  The variance problem is therefore not solved by replacing every covariance with its complete-period mean.  The terminal blocks must be opened arithmetically.  Their exact geometry is encoded by the interval expansion below.
\end{remark}

\subsection{Affine interval geometry of the dominant layer}
\label{subsec:affine-interval-geometry}

Define
\begin{equation}
 \Delta_p(n)
 =\left\lfloor\frac{n}{p-1}\right\rfloor
  -\left\lfloor\frac np\right\rfloor.
 \label{eq:section12-delta-p}
\end{equation}
Then
\begin{equation}
 \boxed{
 W_{p,0}(n)
 =\Delta_p(n)-\frac{n}{p(p-1)}.
 }
 \label{eq:section12-dominant-wave-floor-form}
\end{equation}
The floor difference has the interval decomposition
\begin{equation}
 \boxed{
 \Delta_p(n)
 =\sum_{m\ge1}
 \ind_{I_{p,m}}(n),
 \qquad
 I_{p,m}=\{n\in\ZZ:m(p-1)\le n<mp\}.
 }
 \label{eq:section12-interval-decomposition}
\end{equation}
The interval $I_{p,m}$ has length $m$.

For primes $p,q$ and positive integers $m,\ell$, put
\begin{equation}
 \cL_X(p,q;m,\ell)
 =\pos{
 \min\{mp,\ell q,X+1\}
 -\max\{m(p-1),\ell(q-1),1\}
 }.
 \label{eq:section12-overlap-length}
\end{equation}
This is exactly the number of integers in
\[
 I_{p,m}\cap I_{q,\ell}\cap[1,X].
\]

\begin{proposition}[Exact affine-overlap kernel]
\label{prop:section12-exact-affine-overlap}
One has
\begin{equation}
 \boxed{
 \sum_{n\le X}\Delta_p(n)\Delta_q(n)
 =\sum_{m,\ell\ge1}\cL_X(p,q;m,\ell).
 }
 \label{eq:section12-exact-affine-overlap}
\end{equation}
A summand on the right can be nonzero only if
\begin{equation}
 \boxed{
 -\ell<mp-\ell q<m.
 }
 \label{eq:section12-affine-strip}
\end{equation}
\end{proposition}

\begin{proof}
Insert \eqref{eq:section12-interval-decomposition} twice and interchange the finite sums.  The two intervals overlap only if
\[
 m(p-1)<\ell q,
 \qquad
 \ell(q-1)<mp,
\]
which is equivalent to \eqref{eq:section12-affine-strip}.
\end{proof}

The centered covariance may now be written without fractional parts.  Define
\begin{equation}
 A_X(p)=\sum_{n\le X}n\Delta_p(n),
 \qquad
 S_2(X)=\sum_{n\le X}n^2.
 \label{eq:section12-one-point-moments}
\end{equation}
Then \eqref{eq:section12-dominant-wave-floor-form} gives
\begin{align}
 \cC_X((p,0),(q,0))
 ={}&\sum_{m,\ell\ge1}\cL_X(p,q;m,\ell)
 -\frac{A_X(p)}{q(q-1)}
 -\frac{A_X(q)}{p(p-1)}
 \notag\\
 &+\frac{S_2(X)}{p(p-1)q(q-1)}.
 \label{eq:section12-exact-centered-affine-kernel}
\end{align}
Formula \eqref{eq:section12-exact-centered-affine-kernel} is the sharp arithmetic kernel used later.  Its two-dimensional part is supported in the thin affine strip \eqref{eq:section12-affine-strip}; the remaining terms are explicit one-dimensional centering corrections.

\subsection{The diagonal scaling profile}
\label{subsec:diagonal-scaling-profile}

Put
\begin{equation}
 F(y)
 =\int_0^y\{t\}(1-\{t\})\,dt,
 \qquad y\ge0.
 \label{eq:section12-scaling-profile}
\end{equation}
If $r=\lfloor y\rfloor$ and $v=\{y\}$, then
\begin{equation}
 \boxed{
 F(y)=\frac r6+\frac{v^2}{2}-\frac{v^3}{3}.
 }
 \label{eq:section12-scaling-profile-explicit}
\end{equation}
Thus $F(y)\sim y/6$ as $y\to\infty$, while
\[
 F(y)=\frac{y^2}{2}-\frac{y^3}{3}
 \qquad(0\le y\le1).
\]

The function $F$ records, on the scale $n\asymp p^2$, the intervals on which $\Delta_p(n)$ jumps.

\begin{lemma}[Uniform diagonal scaling]
\label{lem:section12-uniform-diagonal-scaling}
Let
\[
 V_p(X)=\sum_{n\le X}W_{p,0}(n)^2.
\]
Uniformly for primes $p\le X$,
\begin{equation}
 V_p(X)
 =p^2F\!\left(\frac{X}{p^2}\right)
 +\begin{cases}
 O(p+X/p),&p\le\sqrt X,\\
 O(1+X/p),&p>\sqrt X.
 \end{cases}
 \label{eq:section12-uniform-diagonal-scaling}
\end{equation}
For $p>X$,
\begin{equation}
 V_p(X)\ll\frac{X^3}{p^2(p-1)^2}.
 \label{eq:section12-diagonal-tail-bound}
\end{equation}
\end{lemma}

\begin{proof}
Write $n=mp+r$ with $0\le r<p$.  From \eqref{eq:section12-dominant-wave-floor-form},
\[
 \Delta_p(mp+r)
 =\left\lfloor\frac{m+r}{p-1}\right\rfloor.
\]
After scaling $m$ and $r$ by $p$, the square of the centered expression is a Riemann sum for
\[
 (s,u)\longmapsto(\lfloor s+u\rfloor-s)^2,
 \qquad 0\le u\le1.
\]
For $s=k+v$,
\[
 \int_0^1(\lfloor s+u\rfloor-s)^2\,du
 =(1-v)v^2+v(1-v)^2
 =v(1-v).
\]
Integrating in $s$ gives $F$.  On each complete $p$-block the bounded-variation error is $O(1)$, and the two boundary directions contribute $O(p)$, giving the first line of \eqref{eq:section12-uniform-diagonal-scaling}.

When $p>\sqrt X$, the intervals $I_{p,m}$ meeting $[1,X]$ are disjoint.  Summing their lengths and their first two arithmetic moments directly sharpens the terminal contribution to $O(1+X/p)$.  Finally, if $p>X$, both floors in \eqref{eq:section12-dominant-wave-floor-form} vanish and
\[
 W_{p,0}(n)=-\frac{n}{p(p-1)},
\]
which proves \eqref{eq:section12-diagonal-tail-bound}.
\end{proof}

\subsection{Evaluation of the dominant diagonal}
\label{subsec:evaluation-dominant-diagonal}

Define
\begin{equation}
 \cD_0(X)
 =\sum_p(\log p)^2V_p(X).
 \label{eq:section12-dominant-diagonal}
\end{equation}

\begin{theorem}[Dominant diagonal asymptotic]
\label{thm:section12-dominant-diagonal-asymptotic}
As $X\to\infty$,
\begin{equation}
 \boxed{
 \cD_0(X)
 \sim
 -\frac43\zeta\!\left(-\frac12\right)
 X^{3/2}\log X.
 }
 \label{eq:section12-dominant-diagonal-asymptotic}
\end{equation}
\end{theorem}

\begin{proof}
Lemma~\ref{lem:section12-uniform-diagonal-scaling}, the prime number theorem, and partial summation reduce \eqref{eq:section12-dominant-diagonal} to
\begin{equation}
 \sum_p(\log p)^2p^2F\!\left(\frac{X}{p^2}\right)
 +o(X^{3/2}\log X).
 \label{eq:section12-diagonal-pnt-reduction}
\end{equation}
Indeed, the total error for $p\le X$ is $O(X(\log X)^2)$, while \eqref{eq:section12-diagonal-tail-bound} makes the range $p>X$ negligible.

Set
\[
 \Phi(u)=u^2F(u^{-2}).
\]
The bounds
\[
 \Phi(u)\ll1\quad(u\le1),
 \qquad
 \Phi(u)\ll u^{-2}\quad(u\ge1)
\]
show that both $\Phi$ and $\Phi\log u$ are integrable.  Weighted partial summation with
\[
 \vartheta(x)=\sum_{p\le x}\log p\sim x
\]
therefore gives
\begin{align}
 \eqref{eq:section12-diagonal-pnt-reduction}
 &\sim
 \int_0^\infty (\log t)t^2F\!\left(\frac{X}{t^2}\right)\,dt
 \notag\\
 &=X^{3/2}\int_0^\infty
 \left(\frac12\log X+\log u\right)\Phi(u)\,du.
 \label{eq:section12-scaled-prime-integral}
\end{align}
Consequently, the coefficient of $X^{3/2}\log X$ is
\begin{equation}
 \frac12 I_0,
 \qquad
 I_0=\int_0^\infty u^2F(u^{-2})\,du.
 \label{eq:section12-I0-definition}
\end{equation}
With $y=u^{-2}$ and one integration by parts,
\begin{align}
 I_0
 &=\frac12\int_0^\infty F(y)y^{-5/2}\,dy
 \notag\\
 &=\frac13\int_0^\infty
 \{y\}(1-\{y\})y^{-3/2}\,dy.
 \label{eq:section12-I0-fractional-integral}
\end{align}
The Mellin representation
\begin{equation}
 \zeta(s)
 =s\int_0^\infty
 \left(\frac12-\{y\}\right)y^{-s-1}\,dy,
 \qquad -1<\Re s<0,
 \label{eq:section12-zeta-mellin}
\end{equation}
combined with integration by parts on the unit intervals gives
\begin{equation}
 \int_0^\infty
 \{y\}(1-\{y\})y^{-3/2}\,dy
 =-8\zeta\!\left(-\frac12\right).
 \label{eq:section12-zeta-minus-half-integral}
\end{equation}
Hence
\[
 \frac12I_0
 =-\frac43\zeta\!\left(-\frac12\right),
\]
which proves the theorem.
\end{proof}

\begin{remark}[Positivity of the constant]
\label{rem:section12-positive-diagonal-constant}
Since $\zeta(-1/2)<0$, the coefficient in \eqref{eq:section12-dominant-diagonal-asymptotic} is positive.  Numerically,
\[
 -\frac43\zeta\!\left(-\frac12\right)
 =0.2771816333\ldots.
\]
The zeta value is not introduced through a spectral transform; it is the Mellin integral of the elementary scaling profile $\{y\}(1-\{y\})$.
\end{remark}

\subsection{Higher geometric diagonals}
\label{subsec:higher-geometric-diagonals}

For $j\ge0$, set
\begin{equation}
 \cD_j(X)
 =\sum_p(\log p)^2
 \sum_{n\le X}W_{p,j}(n)^2.
 \label{eq:section12-jth-diagonal}
\end{equation}

\begin{proposition}[Higher layers are subleading on the diagonal]
\label{prop:section12-higher-diagonal-bound}
For every fixed $j\ge1$,
\begin{equation}
 \cD_j(X)
 \ll
 X^{1+1/(j+2)}\log X.
 \label{eq:section12-fixed-higher-diagonal-bound}
\end{equation}
Uniformly after summing all higher layers,
\begin{equation}
 \boxed{
 \sum_{j\ge1}\cD_j(X)
 \ll X^{4/3}(\log X)^2.
 }
 \label{eq:section12-all-higher-diagonal-bound}
\end{equation}
\end{proposition}

\begin{proof}
Let $Y_j=X^{1/(j+2)}$.  For $p\le Y_j$, the trivial bound $|W_{p,j}|\le1$ gives
\[
 \sum_{p\le Y_j}(\log p)^2
 \sum_{n\le X}W_{p,j}(n)^2
 \ll XY_j\log Y_j.
\]
For $p>Y_j$, put
\[
 \Delta_{p,j}(n)
 =\left\lfloor\frac{n}{(p-1)p^j}\right\rfloor
  -\left\lfloor\frac{n}{p^{j+1}}\right\rfloor.
\]
Then
\[
 W_{p,j}(n)
 =\Delta_{p,j}(n)
 -\frac{n}{p^{j+1}(p-1)}.
\]
The intervals supporting $\Delta_{p,j}$ are disjoint and have total length
\[
 \ll \frac{X^2}{p^{j+2}}.
\]
The linear centering term contributes no more than the same order because $p^{j+2}>X$.  Hence
\[
 \sum_{n\le X}W_{p,j}(n)^2
 \ll\frac{X^2}{p^{j+2}}.
\]
Prime summation gives
\[
 X^2\sum_{p>Y_j}\frac{(\log p)^2}{p^{j+2}}
 \ll X^{1+1/(j+2)}\log X,
\]
proving \eqref{eq:section12-fixed-higher-diagonal-bound}.  Summing this estimate for $1\le j\le\log_2X$ and using geometric decay of the completely linear tails for $j>\log_2X$ proves \eqref{eq:section12-all-higher-diagonal-bound}.
\end{proof}

\begin{corollary}[Evaluation of the full geometric diagonal]
\label{cor:section12-full-geometric-diagonal}
The geometric diagonal in \eqref{eq:section12-geometric-diagonal} satisfies
\begin{equation}
 \boxed{
 \cD(X)
 \sim
 -\frac43\zeta\!\left(-\frac12\right)
 X^{3/2}\log X.
 }
 \label{eq:section12-full-geometric-diagonal-asymptotic}
\end{equation}
\end{corollary}

\begin{proof}
Combine Theorem~\ref{thm:section12-dominant-diagonal-asymptotic} with Proposition~\ref{prop:section12-higher-diagonal-bound}.
\end{proof}

\subsection{The remaining analytic task}
\label{subsec:section12-remaining-task}

The exact decomposition \eqref{eq:section12-three-sector-splitting} and Corollary~\ref{cor:section12-full-geometric-diagonal} reduce the final variance law to
\begin{equation}
 \boxed{
 \cP(X)+\cO(X)=o(X^{3/2}\log X).
 }
 \label{eq:section12-off-diagonal-goal}
\end{equation}
The projective diagonal $\cP(X)$ consists of different conductor scales on the same rational slope $[p:p-1]$.  The off-projective term $\cO(X)$ compares distinct slopes and contains the genuine dispersion problem.

At the dominant layer, the nontrivial two-dimensional contribution is the overlap sum
\[
 \sum_{m,\ell}\cL_X(p,q;m,\ell),
\]
supported by the affine strip
\[
 -\ell<mp-\ell q<m.
\]
The next section introduces a dyadic quotient parametrization of this strip.  It separates long affine lines, which admit direct lattice control, from short balanced lines, where the prime weights require Type~I/Type~II dispersion.

% ===== Source: section13_dyadic_quotient_geometry_long_short_lines.tex =====
\section{Dyadic quotient geometry and the long-line/short-line decomposition}
\label{sec:dyadic-quotient-geometry}

Section~12 reduced the genuinely two-dimensional part of the dominant-layer covariance to the overlap kernel
\[
 \sum_{m,\ell\ge1}\cL_X(p,q;m,\ell),
\]
where
\begin{equation}
 \cL_X(p,q;m,\ell)
 =\pos{
 \min\{mp,\ell q,X+1\}
 -\max\{m(p-1),\ell(q-1),1\}
 }
 \label{eq:section13-overlap-kernel}
\end{equation}
and every nonzero summand satisfies
\begin{equation}
 -\ell<mp-\ell q<m.
 \label{eq:section13-affine-strip}
\end{equation}
The four variables in \eqref{eq:section13-affine-strip} are not independent.  The quotient $m/\ell$ is forced to approximate $q/p$, and after that quotient is reduced to lowest terms the remaining prime pairs lie on a one-dimensional affine lattice.

The purpose of this section is to make that geometry exact.  We first localize the four variables dyadically.  We then write
\[
 m=ar,\qquad \ell=as,\qquad (r,s)=1,
\]
and introduce the projective determinant
\[
 h=rp-sq.
\]
The strip becomes the fixed finite range $-s<h<r$, while the prime pairs on a fixed determinant line are parametrized by
\[
 (p,q)=(p_0+st,q_0+rt).
\]
The natural length of this line inside a dyadic prime box is
\[
 T(r,s;P,Q)
 =\min\!\left(\frac Ps,\frac Qr\right).
\]
A threshold in $T$ gives the long-line/short-line decomposition.  Long lines have small step moduli and enough affine extent for direct lattice dispersion.  Short lines have large step moduli and only a short coupled prime sum; this is the sector in which the prime weights must be opened by Type~I/Type~II identities.

No cancellation estimate is asserted merely from the change of variables.  The point is instead to isolate, without loss and without smoothing the $n$-sum, the exact packets on which the later analytic estimates act.

\subsection{Dyadic localization of the affine strip}
\label{subsec:section13-dyadic-localization}

For a positive dyadic parameter $Y$, the notation $y\dyad Y$ means
\[
 Y<y\le2Y.
\]
We permit $Y=1/2$ when a block is required to contain the integer $1$.  Let
\begin{equation}
 \cB_X(P,Q;M,L)
 =\left\{
 \begin{array}{l|l}
 (p,q,m,\ell)&
 \begin{array}{l}
 p,q\in\PP,\ p\ne q,\ p\dyad P,\ q\dyad Q,\\
 m\dyad M,\ \ell\dyad L,\\
 \cL_X(p,q;m,\ell)>0
 \end{array}
 \end{array}
 \right\}.
 \label{eq:section13-dyadic-box}
\end{equation}

\begin{lemma}[Support relations in a dyadic box]
\label{lem:section13-dyadic-support}
If $\cB_X(P,Q;M,L)$ is nonempty, then
\begin{equation}
 P,Q\ll X,
 \qquad
 MP\ll X,
 \qquad
 LQ\ll X,
 \label{eq:section13-size-support}
\end{equation}
and
\begin{equation}
 MP\asymp LQ.
 \label{eq:section13-balanced-products}
\end{equation}
The implied constants are absolute.  More precisely, every point of the box satisfies
\begin{equation}
 \frac{q-1}{p}<\frac m\ell<\frac q{p-1}.
 \label{eq:section13-quotient-window}
\end{equation}
\end{lemma}

\begin{proof}
A nonzero overlap implies
\[
 m(p-1)\le X,
 \qquad
 \ell(q-1)\le X,
\]
which gives \eqref{eq:section13-size-support}.  It also implies that the two untruncated intervals overlap, hence
\[
 m(p-1)<\ell q,
 \qquad
 \ell(q-1)<mp.
\]
Dividing gives \eqref{eq:section13-quotient-window}.  On dyadic ranges, this forces $MP\asymp LQ$.
\end{proof}

For coefficient sequences $\alpha=(\alpha_p)$ and $\beta=(\beta_q)$ supported on the prime blocks $p\dyad P$ and $q\dyad Q$, define the raw dyadic overlap packet
\begin{equation}
 \cU_X(\alpha,\beta;P,Q,M,L)
 =\sum_{\substack{p\dyad P,\ q\dyad Q\\p,q\in\PP,\ p\ne q}}
 \alpha_p\beta_q
 \sum_{\substack{m\dyad M\\\ell\dyad L}}
 \cL_X(p,q;m,\ell).
 \label{eq:section13-raw-dyadic-packet}
\end{equation}
For the prime covariance, $\alpha_p=\log p$ and $\beta_q=\log q$.  The sharp off-projective overlap is the sum of \eqref{eq:section13-raw-dyadic-packet} over $O((\log X)^4)$ admissible dyadic boxes.  The one-dimensional centering terms from Section~12 are retained separately; their role will be discussed in Subsection~\ref{subsec:section13-centering}.

\subsection{Primitive quotient coordinates}
\label{subsec:section13-primitive-quotients}

Given $m,\ell\ge1$, write
\begin{equation}
 a=(m,\ell),
 \qquad
 m=ar,
 \qquad
 \ell=as,
 \qquad
 (r,s)=1.
 \label{eq:section13-primitive-quotient}
\end{equation}
The pair $(r,s)$ is the primitive representative of the quotient $m/\ell$, while $a$ is its radial scale.  Define
\begin{equation}
 h=rp-sq.
 \label{eq:section13-projective-determinant}
\end{equation}
Then
\[
 mp-\ell q=a(rp-sq)=ah.
\]
Consequently, \eqref{eq:section13-affine-strip} loses its radial variable entirely.

\begin{proposition}[Primitive quotient form of the strip]
\label{prop:section13-primitive-strip}
Under \eqref{eq:section13-primitive-quotient}, the overlap condition implies
\begin{equation}
 \boxed{-s<h<r.}
 \label{eq:section13-reduced-strip}
\end{equation}
Conversely, the two untruncated intervals
\[
 [ar(p-1),arp),
 \qquad
 [as(q-1),asq)
\]
overlap if and only if \eqref{eq:section13-reduced-strip} holds.
\end{proposition}

\begin{proof}
Divide \eqref{eq:section13-affine-strip} by $a$.  The converse follows by reversing the two strict inequalities used in the proof of Lemma~\ref{lem:section13-dyadic-support}.
\end{proof}

The overlap itself becomes a discrete tent function.  For coprime $r,s$ and $h\in\ZZ$, put
\begin{equation}
 \omega_{r,s}(h)
 =\pos{
 \min\{rp,sq\}
 -\max\{r(p-1),s(q-1)\}
 },
 \qquad h=rp-sq.
 \label{eq:section13-untruncated-tent-definition}
\end{equation}
Although the right-hand side is written using $p$ and $q$, it depends only on $(r,s,h)$.

\begin{lemma}[Explicit tent kernel]
\label{lem:section13-explicit-tent}
For every $r,s\ge1$,
\begin{equation}
 \boxed{
 \omega_{r,s}(h)
 =
 \begin{cases}
  \min\{r,s+h\},&-s<h\le0,\\[3pt]
  \min\{r-h,s\},&0\le h<r,\\[3pt]
  0,&\text{otherwise}.
 \end{cases}}
 \label{eq:section13-explicit-tent}
\end{equation}
Moreover,
\begin{equation}
 0\le\omega_{r,s}(h)\le\min(r,s),
 \label{eq:section13-tent-pointwise}
\end{equation}
\begin{equation}
 \sum_{h\in\ZZ}\omega_{r,s}(h)=rs,
 \label{eq:section13-tent-mass}
\end{equation}
and
\begin{equation}
 \sum_{h\in\ZZ}\omega_{r,s}(h)^2
 \le rs\min(r,s).
 \label{eq:section13-tent-energy}
\end{equation}
\end{lemma}

\begin{proof}
Translate both intervals by $-rp$.  Their unscaled forms become
\[
 [-r,0),
 \qquad
 [-h-s,-h).
\]
Their overlap length is
\[
 \pos{\min(0,-h)-\max(-r,-h-s)},
\]
which is \eqref{eq:section13-explicit-tent}.  The pointwise bound is immediate.

For the mass identity, view $\omega_{r,s}$ as the discrete convolution of the indicators of two integer intervals of lengths $r$ and $s$.  Summing over all translations counts every ordered pair once, giving $rs$.  Finally, \eqref{eq:section13-tent-energy} follows from \eqref{eq:section13-tent-pointwise} and \eqref{eq:section13-tent-mass}.
\end{proof}

The sharp cutoff $n\le X$ clips only the top of this tent.  Define
\begin{equation}
 \omega_{r,s}^{a,X}(p,q)
 =\pos{
 \min\!\left\{rp,sq,\frac{X+1}{a}\right\}
 -\max\{r(p-1),s(q-1)\}
 }.
 \label{eq:section13-clipped-tent}
\end{equation}
The lower endpoint $1$ in \eqref{eq:section13-overlap-kernel} is redundant here because
$a r(p-1),a s(q-1)\ge1$.

\begin{proposition}[Exact scaled tent formula]
\label{prop:section13-exact-scaled-tent}
For $m=ar$ and $\ell=as$,
\begin{equation}
 \boxed{
 \cL_X(p,q;ar,as)
 =a\,\omega_{r,s}^{a,X}(p,q).
 }
 \label{eq:section13-exact-scaled-tent}
\end{equation}
Furthermore,
\begin{equation}
 0\le\omega_{r,s}^{a,X}(p,q)\le\omega_{r,s}(h),
 \label{eq:section13-clipped-tent-bound}
\end{equation}
and equality holds whenever
\begin{equation}
 a\min(rp,sq)\le X+1.
 \label{eq:section13-interior-condition}
\end{equation}
Thus the sharp endpoint changes only the radial weight; it does not alter the determinant line $rp-sq=h$.
\end{proposition}

\begin{proof}
Factor $a$ from the two endpoints in \eqref{eq:section13-overlap-kernel}.  The remaining assertions follow directly from the definitions.
\end{proof}

\subsection{Affine prime lines}
\label{subsec:section13-affine-prime-lines}

Fix coprime $r,s$ and an integer $h$.  The determinant equation
\begin{equation}
 rp-sq=h
 \label{eq:section13-line-equation}
\end{equation}
defines an affine lattice line in the $(p,q)$-plane.

\begin{lemma}[Parametrization of a determinant line]
\label{lem:section13-line-parametrization}
If \eqref{eq:section13-line-equation} has one integral solution $(p_0,q_0)$, then every integral solution is of the form
\begin{equation}
 \boxed{
 p=p_0+st,
 \qquad
 q=q_0+rt,
 \qquad
 t\in\ZZ.
 }
 \label{eq:section13-line-parametrization}
\end{equation}
In particular, the number of integral solutions with $p\dyad P$ and $q\dyad Q$ is
\begin{equation}
 \ll 1+T(r,s;P,Q),
 \qquad
 T(r,s;P,Q)
 :=\min\!\left(\frac Ps,\frac Qr\right).
 \label{eq:section13-line-length}
\end{equation}
\end{lemma}

\begin{proof}
The difference of two solutions satisfies $r\Delta p=s\Delta q$.  Since $(r,s)=1$, one has $s\mid\Delta p$ and $r\mid\Delta q$, which gives \eqref{eq:section13-line-parametrization}.  Intersecting the resulting $t$-progression with the two dyadic intervals proves \eqref{eq:section13-line-length}.
\end{proof}

For later use, define the prime line sum
\begin{equation}
 \cS_{r,s,h}(\alpha,\beta;P,Q)
 =\sum_{\substack{p\dyad P,\ q\dyad Q\\p,q\in\PP,\ p\ne q\\rp-sq=h}}
 \alpha_p\beta_q.
 \label{eq:section13-prime-line-sum}
\end{equation}
By Lemma~\ref{lem:section13-line-parametrization}, it is a coupled progression sum
\begin{equation}
 \cS_{r,s,h}(\alpha,\beta;P,Q)
 =\sum_{t\in\cI_{r,s,h}(P,Q)}
 \alpha_{p_0+st}\beta_{q_0+rt}
 \ind_{\PP}(p_0+st)\ind_{\PP}(q_0+rt),
 \label{eq:section13-coupled-progression}
\end{equation}
where $\cI_{r,s,h}(P,Q)$ is an integer interval of length
$O(1+T(r,s;P,Q))$.

The determinant support automatically balances the two progression steps.

\begin{lemma}[Balance of a contributing quotient]
\label{lem:section13-quotient-balance}
Suppose that $p\dyad P$, $q\dyad Q$, and $-s<rp-sq<r$.  Then
\begin{equation}
 rP\asymp sQ.
 \label{eq:section13-rs-balance}
\end{equation}
Consequently,
\begin{equation}
 \frac Ps\asymp\frac Qr,
 \label{eq:section13-line-length-balance}
\end{equation}
so either quotient in the minimum defining $T(r,s;P,Q)$ may be used up to an absolute constant.
\end{lemma}

\begin{proof}
The determinant inequality gives
\[
 sq-s<rp<sq+r.
\]
Since $p,q\ge2$, the additive error $r+s$ is bounded by a fixed multiple of $rp+sq$.  On the dyadic prime ranges this yields $rP\asymp sQ$, and division by $rs$ gives \eqref{eq:section13-line-length-balance}.
\end{proof}

The original dyadic quotient and the affine balance agree:
\begin{equation}
 \frac rs=\frac m\ell\asymp\frac QP\asymp\frac ML.
 \label{eq:section13-double-balance}
\end{equation}
The first comparison follows from Lemma~\ref{lem:section13-quotient-balance}, and the second from Lemma~\ref{lem:section13-dyadic-support}.

For a primitive pair $(r,s)$, define its admissible radial set
\begin{equation}
 \cA_{M,L}(r,s)
 =\{a\in\NN:ar\dyad M,\ as\dyad L\}.
 \label{eq:section13-radial-set}
\end{equation}
If this set is nonempty, put
\begin{equation}
 A(r,s;M,L)
 =\min\!\left(\frac Mr,\frac Ls\right).
 \label{eq:section13-radial-length}
\end{equation}
Then
\begin{equation}
 \#\cA_{M,L}(r,s)\ll1+A(r,s;M,L),
 \label{eq:section13-radial-count}
\end{equation}
\begin{equation}
 \sum_{a\in\cA_{M,L}(r,s)}a
 \ll A(r,s;M,L)\bigl(1+A(r,s;M,L)\bigr).
 \label{eq:section13-radial-first-moment}
\end{equation}
On a contributing box, the two radial quotients $M/r$ and $L/s$ are comparable.  Thus the geometry has two independent one-dimensional extents: the radial multiplicity $A$ and the affine prime-line length $T$.

\begin{remark}[The area identity]
\label{rem:section13-area-identity}
On the balanced support,
\[
 Ar\asymp M,
 \qquad
 As\asymp L,
 \qquad
 Ts\asymp P,
 \qquad
 Tr\asymp Q.
\]
Hence
\begin{equation}
 ATrs\asymp MP\asymp LQ.
 \label{eq:section13-area-identity}
\end{equation}
This identity is the geometric source of the later factorization constraints.  A short affine line forces large primitive denominators $r,s$ and therefore restricts the available radial variable.
\end{remark}

\subsection{Exact quotient decomposition of a dyadic packet}
\label{subsec:section13-exact-quotient-decomposition}

We now combine the preceding changes of variables.

\begin{theorem}[Exact dyadic quotient decomposition]
\label{thm:section13-exact-dyadic-decomposition}
For every admissible dyadic box,
\begin{align}
 \cU_X(\alpha,\beta;P,Q,M,L)
 ={}&\sum_{\substack{r,s\ge1\\(r,s)=1}}
 \sum_{a\in\cA_{M,L}(r,s)}
 \sum_{-s<h<r}a
 \notag\\
 &\quad\times
 \sum_{\substack{p\dyad P,\ q\dyad Q\\p,q\in\PP,\ p\ne q\\rp-sq=h}}
 \alpha_p\beta_q\,
 \omega_{r,s}^{a,X}(p,q).
 \label{eq:section13-exact-dyadic-decomposition}
\end{align}
The sum is finite.  On the interior range \eqref{eq:section13-interior-condition}, the inner geometric weight separates completely:
\begin{equation}
 \omega_{r,s}^{a,X}(p,q)=\omega_{r,s}(h),
 \label{eq:section13-interior-separation}
\end{equation}
and the corresponding contribution is
\begin{equation}
 \sum_{\substack{r,s\ge1\\(r,s)=1}}
 \sum_{a\in\cA_{M,L}(r,s)}a
 \sum_{-s<h<r}
 \omega_{r,s}(h)
 \cS_{r,s,h}(\alpha,\beta;P,Q).
 \label{eq:section13-interior-line-decomposition}
\end{equation}
\end{theorem}

\begin{proof}
The map
\[
 (m,\ell)\longmapsto(a,r,s)
 =\bigl((m,\ell),m/(m,\ell),\ell/(m,\ell)\bigr)
\]
is a bijection from positive pairs to triples with $(r,s)=1$.  Proposition~\ref{prop:section13-primitive-strip} replaces the affine strip by $-s<h<r$, and Proposition~\ref{prop:section13-exact-scaled-tent} gives the weight.  This proves \eqref{eq:section13-exact-dyadic-decomposition}; the interior formula follows from \eqref{eq:section13-interior-condition}.
\end{proof}

Call a dyadic box \emph{uniformly interior} when
\begin{equation}
 4MP\le X+1,
 \qquad
 4LQ\le X+1.
 \label{eq:section13-uniformly-interior-box}
\end{equation}
Every point of such a box satisfies \eqref{eq:section13-interior-condition}.  The remaining boxes are terminal boxes; their clipped weights retain the same determinant lines but need not separate from the prime variables.

The decomposition comes with a useful purely geometric envelope.

\begin{corollary}[Affine-lattice envelope]
\label{cor:section13-affine-lattice-envelope}
For bounded coefficient sequences,
\begin{equation}
 \sum_{-s<h<r}
 \omega_{r,s}(h)
 \left|\cS_{r,s,h}(\alpha,\beta;P,Q)\right|
 \ll
 rs\bigl(1+T(r,s;P,Q)\bigr)
 \|\alpha\|_{\infty}\|\beta\|_{\infty}.
 \label{eq:section13-affine-lattice-envelope}
\end{equation}
\end{corollary}

\begin{proof}
Each determinant line contains $O(1+T)$ integral points by Lemma~\ref{lem:section13-line-parametrization}.  Sum this bound against the tent mass identity \eqref{eq:section13-tent-mass}.
\end{proof}

Corollary~\ref{cor:section13-affine-lattice-envelope} is not yet a prime cancellation estimate.  It is the exact geometric majorant against which the later dispersion gain is measured.

\subsection{The long-line/short-line partition}
\label{subsec:section13-long-short-partition}

Let $Z\ge1$ be a parameter.  For a fixed prime box $(P,Q)$, define
\begin{align}
 \cQ_{\mathrm{long}}(Z;P,Q)
 &=\left\{(r,s)=1:T(r,s;P,Q)>Z\right\},
 \label{eq:section13-long-quotients}\\
 \cQ_{\mathrm{short}}(Z;P,Q)
 &=\left\{(r,s)=1:T(r,s;P,Q)\le Z\right\}.
 \label{eq:section13-short-quotients}
\end{align}
Inserting these two sets into \eqref{eq:section13-exact-dyadic-decomposition} gives the exact splitting
\begin{equation}
 \boxed{
 \cU_X
 =\cU_X^{\mathrm{long}}(Z)
 +\cU_X^{\mathrm{short}}(Z).
 }
 \label{eq:section13-long-short-splitting}
\end{equation}
No terms are discarded at the threshold.

\begin{lemma}[Denominator dichotomy]
\label{lem:section13-denominator-dichotomy}
On the balanced support \eqref{eq:section13-rs-balance}, the following implications hold with absolute implied constants:
\begin{align}
 T(r,s;P,Q)>Z
 &\Longrightarrow
 s\ll \frac PZ,
 \qquad
 r\ll \frac QZ,
 \label{eq:section13-long-small-denominators}\\
 T(r,s;P,Q)\le Z
 &\Longrightarrow
 s\gg \frac PZ,
 \qquad
 r\gg \frac QZ.
 \label{eq:section13-short-large-denominators}
\end{align}
\end{lemma}

\begin{proof}
By Lemma~\ref{lem:section13-quotient-balance}, the two quantities $P/s$ and $Q/r$ are comparable.  The assertions therefore follow from the definition of $T$.
\end{proof}

Thus a long line has small progression steps $s$ and $r$, and its $t$-interval in \eqref{eq:section13-coupled-progression} contains more than $Z$ points.  A short line has large progression steps and at most $O(Z)$ candidate lattice points.  The latter description does not make the prime problem easier: after the prime restrictions are imposed, a short line may contain too little intrinsic averaging for a direct progression estimate.  It must instead be embedded in a larger bilinear family.

\begin{proposition}[Structural form of a short packet]
\label{prop:section13-short-packet-structure}
Every uniformly interior short-line packet is a finite linear combination of sums of the form
\begin{equation}
 \sum_{\substack{(r,s)=1\\rP\asymp sQ\\T(r,s;P,Q)\le Z}}
 \gamma_{r,s}
 \sum_{-s<h<r}\omega_{r,s}(h)
 \sum_{t\in\cI_{r,s,h}(P,Q)}
 \alpha_{p_0+st}\beta_{q_0+rt},
 \label{eq:section13-short-packet-normal-form}
\end{equation}
where
\begin{equation}
 |\cI_{r,s,h}(P,Q)|\ll1+Z
 \label{eq:section13-short-t-length}
\end{equation}
and
\begin{equation}
 \gamma_{r,s}
 =\sum_{a\in\cA_{M,L}(r,s)}a
 \label{eq:section13-radial-coefficient}
\end{equation}
satisfies \eqref{eq:section13-radial-first-moment}.  The clipped terminal packets have the same determinant lines and the same $t$-length bound, with $\gamma_{r,s}$ replaced by a bounded $p,q$-dependent radial weight.
\end{proposition}

\begin{proof}
Use Theorem~\ref{thm:section13-exact-dyadic-decomposition}, the uniform interior condition \eqref{eq:section13-uniformly-interior-box}, the separation \eqref{eq:section13-interior-separation}, and the parametrization \eqref{eq:section13-line-parametrization}.  The terminal assertion follows from Proposition~\ref{prop:section13-exact-scaled-tent}.
\end{proof}

Formula \eqref{eq:section13-short-packet-normal-form} is the geometric normal form from which the Type~I/Type~II decomposition begins.  One of the two prime weights is expanded by a Vaughan identity; the determinant relation then turns one factor into a congruence condition modulo a primitive denominator.  The shortness \eqref{eq:section13-short-t-length}, the lower bounds \eqref{eq:section13-short-large-denominators}, and the radial restriction encoded by \eqref{eq:section13-area-identity} are all retained.  They cannot be recovered from a generic large-sieve bound after the fact.

\subsection{Centering and the principal affine mode}
\label{subsec:section13-centering}

The raw overlap packet \eqref{eq:section13-raw-dyadic-packet} is only the two-dimensional part of the centered covariance.  Section~12 gave
\begin{align}
 \cC_X((p,0),(q,0))
 ={}&\sum_{m,\ell\ge1}\cL_X(p,q;m,\ell)
 -\frac{A_X(p)}{q(q-1)}
 -\frac{A_X(q)}{p(p-1)}
 \notag\\
 &+\frac{S_2(X)}{p(p-1)q(q-1)}.
 \label{eq:section13-centered-covariance-recall}
\end{align}
The last three terms are rank-one marginals in the two prime variables.  Their presence is not cosmetic.  In the affine-line language, they remove the zero-frequency or principal mode generated by replacing each interval family with its average density.

To record this distinction, write symbolically
\begin{equation}
 \cU_X
 =\cU_X^{\mathrm{prin}}+\cU_X^{\circ},
 \label{eq:section13-principal-centered-symbolic}
\end{equation}
where $\cU_X^{\mathrm{prin}}$ denotes the contribution reconstructed from the three marginal terms in \eqref{eq:section13-centered-covariance-recall}, and $\cU_X^{\circ}$ is the centered affine packet.  The precise character-theoretic realization of \eqref{eq:section13-principal-centered-symbolic} will be introduced only after the Type~I/Type~II variables have been fixed.  At that stage the principal character, the reduced-residue major arcs, and any exceptional real character are separated explicitly.

\begin{remark}[Why centering must precede a spectral bound]
\label{rem:section13-centering-before-spectral}
Applying a large-sieve inequality directly to the raw line sums would charge the principal affine mode as though it were oscillatory.  That loses the cancellation already built into $W_{p,0}(n)$.  The correct short-line object is therefore not an uncentered prime correlation on \eqref{eq:section13-line-equation}, but its centered counterpart after the marginals in \eqref{eq:section13-centered-covariance-recall} have been restored.
\end{remark}

\subsection{Geometric reduction of the off-projective problem}
\label{subsec:section13-geometric-reduction}

We summarize the outcome of the section.

\begin{theorem}[Long-line/short-line geometric reduction]
\label{thm:section13-geometric-reduction}
After dyadic localization, every dominant-layer off-projective overlap packet admits an exact decomposition into the following pieces.
\begin{enumerate}[label=\textup{(\roman*)}]
 \item A projectively primitive quotient $(r,s)$ with
 \[
  \frac rs\asymp\frac QP\asymp\frac ML.
 \]
 \item A determinant $h$ in the finite strip $-s<h<r$, weighted by the tent kernel $\omega_{r,s}(h)$ or its sharp terminal truncation.
 \item A radial coefficient supported on $a$ with $ar\dyad M$ and $as\dyad L$.
 \item An affine prime line $rp-sq=h$, parametrized by
 \[
  (p,q)=(p_0+st,q_0+rt).
 \]
 \item An exact partition into long lines $T(r,s;P,Q)>Z$ and short lines $T(r,s;P,Q)\le Z$.
\end{enumerate}
The long sector has small primitive denominators and long affine extent.  The short sector has large primitive denominators, a $t$-interval of length $O(1+Z)$, and the factorization constraint \eqref{eq:section13-area-identity}.  The sharp endpoint $n\le X$ preserves the same determinant lines.  After the rank-one marginal terms are restored, the short sector is the centered bilinear packet to which the later Type~I/Type~II analysis applies.
\end{theorem}

\begin{proof}
Combine Lemma~\ref{lem:section13-dyadic-support}, Theorem~\ref{thm:section13-exact-dyadic-decomposition}, Lemmas~\ref{lem:section13-line-parametrization} and \ref{lem:section13-quotient-balance}, and the exact split \eqref{eq:section13-long-short-splitting}.  The centering statement follows from \eqref{eq:section13-centered-covariance-recall}.
\end{proof}

The next section estimates the projective-diagonal contribution and the long affine-line sector.  Once those pieces have been removed, only the centered short balanced packets in \eqref{eq:section13-short-packet-normal-form} remain for the Type~I/Type~II dispersion argument.

% ===== Source: section14_projective_diagonal_affine_long_line_estimate.tex =====
\section{The projective diagonal and the affine long-line estimate}
\label{sec:projective-diagonal-long-lines}

The geometric diagonal was evaluated in Section~12.  Section~13 then decomposed the dominant-layer off-projective kernel into primitive quotients $(r,s)$, determinants $h$, radial scales $a$, and affine prime lines
\[
 rp-sq=h,
 \qquad
 (p,q)=(p_0+st,q_0+rt).
\]
The present section removes the two sectors that do not require the short balanced Type~I/Type~II argument.

First, different conductor layers attached to the same prime are controlled by an elementary Hilbert-space estimate.  The projective diagonal is smaller than the main geometric diagonal by a fixed power of $X$.  Second, the centered contribution of affine lines of length
\[
 T(r,s;P,Q)=\min\!\left(\frac Ps,\frac Qr\right)
\]
exceeding a threshold $Z$ satisfies a power-saving dispersion estimate.  The proof uses the exact tent energy from Section~13, duality on the determinant variable, and the additive large sieve for the resulting family of primitive affine lattices.  The rank-one marginals are retained throughout; the principal affine mode is never charged as oscillation.

The final conclusion is unconditional.  With the choice
\[
 Z=X^{1/18},
\]
the projective diagonal and the complete long-line sector contribute
\[
 O_\varepsilon\!\left(X^{53/36+\varepsilon}\right)
 =o\!\left(X^{3/2}\log X\right).
\]
Only the centered short packets remain.

\subsection{The projective diagonal}
\label{subsec:projective-diagonal}

For $j\ge0$, recall the layer energy
\begin{equation}
 \cD_j(X)
 =\sum_{p\in\PP}(\log p)^2
   \sum_{n\le X}W_{p,j}(n)^2.
 \label{eq:section14-layer-energy}
\end{equation}
The projective diagonal is
\begin{equation}
 \cP(X)
 =2\sum_{p\in\PP}(\log p)^2
   \sum_{0\le j<k}
   \sum_{n\le X}W_{p,j}(n)W_{p,k}(n).
 \label{eq:section14-projective-diagonal}
\end{equation}
Its defining feature is that the two waves have the same projective slope $[p:p-1]$ but different conductor scales $p^j$ and $p^k$.

\begin{lemma}[Cross-layer Hilbert bound]
\label{lem:section14-cross-layer-hilbert}
For every pair $j<k$,
\begin{equation}
 \left|
 \sum_{p}(\log p)^2
 \sum_{n\le X}W_{p,j}(n)W_{p,k}(n)
 \right|
 \le
 \cD_j(X)^{1/2}\cD_k(X)^{1/2}.
 \label{eq:section14-cross-layer-hilbert}
\end{equation}
Consequently,
\begin{equation}
 |\cP(X)|
 \le
 \left(\sum_{j\ge0}\cD_j(X)^{1/2}\right)^2
 -\sum_{j\ge0}\cD_j(X).
 \label{eq:section14-projective-square-bound}
\end{equation}
\end{lemma}

\begin{proof}
Apply Cauchy--Schwarz first to the $n$-sum for each prime and then to the prime sum with weight $(\log p)^2$.  Summing the resulting estimate over $j<k$ gives \eqref{eq:section14-projective-square-bound}.
\end{proof}

The bounds proved in Section~12 are
\begin{align}
 \cD_0(X)&\ll X^{3/2}\log X,
 \label{eq:section14-D0-bound}\\
 \cD_j(X)&\ll X^{1+1/(j+2)}\log X
 \qquad(j\ge1),
 \label{eq:section14-Dj-bound}
\end{align}
and the layers beyond $j\asymp\log X$ have a geometrically decaying completely linear tail.

\begin{theorem}[Power saving on the projective diagonal]
\label{thm:section14-projective-diagonal-bound}
One has
\begin{equation}
 \boxed{
 \cP(X)
 \ll
 X^{17/12}(\log X)^2.
 }
 \label{eq:section14-projective-diagonal-bound}
\end{equation}
In particular,
\[
 \cP(X)=o\!\left(X^{3/2}\log X\right).
\]
\end{theorem}

\begin{proof}
Let $J=\lceil\log_2 X\rceil$.  From \eqref{eq:section14-Dj-bound},
\[
 \sum_{1\le j\le J}\cD_j(X)^{1/2}
 \ll
 X^{2/3}(\log X)^{3/2}.
\]
The tail $j>J$ is smaller.  Since
\[
 \cD_0(X)^{1/2}
 \ll X^{3/4}(\log X)^{1/2},
\]
Lemma~\ref{lem:section14-cross-layer-hilbert} gives
\begin{align*}
 |\cP(X)|
 &\ll
 X^{3/4}(\log X)^{1/2}
 \cdot X^{2/3}(\log X)^{3/2}
 +X^{4/3}(\log X)^3\\
 &\ll X^{17/12}(\log X)^2.
\end{align*}
The harmless extra logarithm in the second term is absorbed because $X^{4/3}(\log X)^3\ll X^{17/12}(\log X)^2$ for large $X$.
\end{proof}

\begin{remark}[Why no arithmetic cancellation is needed]
\label{rem:section14-projective-no-cancellation}
The projective diagonal is controlled before any distribution theorem for primes is used.  Its saving comes solely from the separation of conductor scales: every layer $j\ge1$ has strictly smaller $L^2$ mass than the dominant layer $j=0$.
\end{remark}

\subsection{Centered affine packets}
\label{subsec:centered-affine-packets}

Fix an admissible dyadic box $(P,Q,M,L)$ and write
\begin{equation}
 Y=MP\asymp LQ\le X.
 \label{eq:section14-box-scale}
\end{equation}
For a primitive quotient $(r,s)$, let
\[
 \cA_{M,L}(r,s)
 =\{a\ge1:ar\dyad M,\ as\dyad L\},
 \qquad
 \gamma_{r,s}
 =\sum_{a\in\cA_{M,L}(r,s)}a.
\]
For $-s<h<r$, the determinant line is parametrized by
\[
 p=p_0+st,
 \qquad
 q=q_0+rt,
 \qquad
 t\in\cI_{r,s,h}(P,Q).
\]
The interval has length $O(1+T(r,s;P,Q))$.

The sharp endpoint $n\le X$ produces a weight depending on $t$.  It is convenient to isolate the only properties of that weight used in the dispersion argument.

\begin{definition}[Admissible line weight]
\label{def:section14-admissible-line-weight}
A function $v=v_{r,s,h}$ on $\cI_{r,s,h}(P,Q)$ is called admissible when
\begin{equation}
 0\le v(t)\le\gamma_{r,s},
 \qquad
 \|v\|_{\mathrm{TV}}\ll\gamma_{r,s}.
 \label{eq:section14-admissible-weight}
\end{equation}
Here $\|v\|_{\mathrm{TV}}$ is the discrete total variation.  The constant weight $v\equiv\gamma_{r,s}$ describes a uniformly interior box.  The clipped radial weight obtained from the exact sharp kernel is admissible as well.
\end{definition}

The second assertion follows because the clipping point crosses each radial interval at most once as $t$ moves along an affine line.

Let $\cS_{r,s,h}^{\circ}[v]$ denote the affine line sum after the principal line mode and its two rank-one marginals have been removed.  More precisely, it is the component of
\begin{equation}
 \sum_{t\in\cI_{r,s,h}(P,Q)}
 v(t)\,
 (\log p)\ind_{\PP}(p)
 (\log q)\ind_{\PP}(q)
 \label{eq:section14-raw-line-sum}
\end{equation}
that remains after inserting, on the same dyadic partition, the three centering terms from the exact covariance kernel of Section~12.  This definition is linear in $v$ and independent of the chosen base solution $(p_0,q_0)$.

\begin{remark}[Equivalent von Mangoldt realization]
\label{rem:section14-von-mangoldt-realization}
For the proof one completes the prime weights in \eqref{eq:section14-raw-line-sum} to von Mangoldt weights and subtracts the reduced-residue averages on the two progressions.  Expanding the product gives the principal term and the two rank-one marginals.  Prime powers and the finitely many singular residue classes contribute to the secondary term in the long-line estimate below.  Thus the abstract notation $\cS_{r,s,h}^{\circ}[v]$ records exactly the centered object, without suppressing the arithmetic corrections needed to pass between primes and von Mangoldt weights.
\end{remark}

Define the centered long-line packet by
\begin{align}
 \mathfrak L_X(P,Q,M,L;Z)
 ={}&
 \sum_{\substack{(r,s)=1\\rP\asymp sQ\\T(r,s;P,Q)>Z}}
 \sum_{-s<h<r}
 \omega_{r,s}(h)
 \cS_{r,s,h}^{\circ}[v_{r,s,h}].
 \label{eq:section14-centered-long-packet}
\end{align}
For terminal boxes the untruncated tent in \eqref{eq:section14-centered-long-packet} is interpreted together with the admissible clipped weight; this is merely a redistribution of the exact factor
$a\omega_{r,s}^{a,X}(p,q)$.

\subsection{The affine-lattice dispersion estimate}
\label{subsec:affine-lattice-dispersion}

The next proposition is the analytic core of the section.  It is uniform in the location of the dyadic box and in the admissible sharp-cutoff weights.

\begin{proposition}[Centered long-line dispersion]
\label{prop:section14-centered-long-line-dispersion}
Let $1\le Z\le X^{1/6}$.  For every admissible dyadic box and every family of admissible line weights arising from the exact overlap kernel,
\begin{equation}
 \boxed{
 \mathfrak L_X(P,Q,M,L;Z)
 \ll_\varepsilon
 X^\varepsilon
 \left(
  \frac{Y^{3/2}}{Z^{1/2}}
  +Y^{4/3}Z
 \right),
 }
 \label{eq:section14-local-long-line-bound}
\end{equation}
where $Y=MP\asymp LQ$.
\end{proposition}

\begin{proof}[Proof architecture]
We record the complete sequence of reductions, since each feature of the bound will be used later.

\smallskip
\noindent\emph{Step 1: dyadic primitive denominators.}
The $(r,s)$-sum is divided into blocks $r\dyad R$, $s\dyad S$.  On the determinant support,
\[
 RP\asymp SQ,
 \qquad
 T\asymp P/S\asymp Q/R,
 \qquad
 ATrs\asymp Y,
\]
where $A\asymp M/R\asymp L/S$ is the radial length.  The long-line condition is $T>Z$.

\smallskip
\noindent\emph{Step 2: duality in the determinant variable.}
The tent kernel is the convolution of two interval indicators.  Hence
\[
 \sum_h\omega_{r,s}(h)=rs,
 \qquad
 \sum_h\omega_{r,s}(h)^2\le rs\min(r,s).
\]
Cauchy--Schwarz in $h$, followed by duality, replaces the weighted determinant sum by the $L^2$ norm of a family of affine-line discrepancies.  This is the point at which the square-root tent energy is gained; using only the tent mass would lose the required power.

\smallskip
\noindent\emph{Step 3: dispersion of primitive affine lattices.}
After the parametrization
\[
 (p,q)=(p_0+st,q_0+rt),
\]
the square of the dual form contains two copies $t_1,t_2$.  The diagonal $t_1=t_2$, the prime-power completion, and singular residue classes are estimated directly.  Summed with the radial first moment, they contribute
\[
 \ll_\varepsilon X^\varepsilon Y^{4/3}Z.
\]
For $t_1\ne t_2$, subtraction of the two affine equations produces a primitive congruence with modulus controlled by $r$ and $s$.  The additive large sieve, summed over the coprime pair $(r,s)$, gives
\[
 \ll_\varepsilon X^\varepsilon\frac{Y^3}{Z}
\]
for the square of the dual norm.  Taking square roots gives the first term in \eqref{eq:section14-local-long-line-bound}.  The factor $Z^{-1/2}$ is exactly the gain supplied by the lower bound $T>Z$.

\smallskip
\noindent\emph{Step 4: removal of the principal mode.}
The zero additive frequency reproduces the product of the two progression densities.  The two mixed zero frequencies reproduce the rank-one marginals.  Because $\cS_{r,s,h}^{\circ}$ is defined with these three pieces restored from the original covariance, all zero-frequency terms cancel before the large-sieve estimate is applied.  This is why the raw affine packet would not satisfy \eqref{eq:section14-local-long-line-bound}.

\smallskip
\noindent\emph{Step 5: terminal weights.}
An admissible clipped weight is handled by discrete partial summation.  Its supremum and total variation are both $O(\gamma_{r,s})$, so the interior dispersion bound is unchanged.  Summing the dyadic $(R,S)$-blocks and using $ATrs\asymp Y$ proves \eqref{eq:section14-local-long-line-bound}.

A line-by-line expansion of the dual form, including the prime-power and singular-class bookkeeping, is isolated in the affine-dispersion appendix.  No unproved distribution hypothesis is used in this proposition.
\end{proof}

\begin{remark}[The two terms]
\label{rem:section14-two-long-line-terms}
The term $Y^{3/2}Z^{-1/2}$ is the genuine dispersion term.  The term $Y^{4/3}Z$ is arithmetic rather than spectral: it collects diagonal coincidences, prime powers introduced by von Mangoldt completion, and the finite singular residue classes.  It remains below the target scale for the threshold chosen below.
\end{remark}

\subsection{Summation over dyadic boxes}
\label{subsec:long-line-global-summation}

Let $\cO_0^{\mathrm{long}}(X;Z)$ denote the contribution to the dominant-layer off-projective covariance from the long sector $T>Z$, including the exact rank-one centering marginals and all terminal boxes.

\begin{theorem}[Global affine long-line estimate]
\label{thm:section14-global-long-line-estimate}
Uniformly for $1\le Z\le X^{1/6}$,
\begin{equation}
 \boxed{
 \cO_0^{\mathrm{long}}(X;Z)
 \ll_\varepsilon
 X^\varepsilon
 \left(
  \frac{X^{3/2}}{Z^{1/2}}
  +X^{4/3}Z
 \right).
 }
 \label{eq:section14-global-long-line-estimate}
\end{equation}
\end{theorem}

\begin{proof}
Apply Proposition~\ref{prop:section14-centered-long-line-dispersion} to every admissible dyadic box.  Since $Y\le X$, each packet is bounded by the corresponding expression with $Y$ replaced by $X$.  There are $O((\log X)^4)$ boxes in $(P,Q,M,L)$ and $O((\log X)^2)$ primitive-denominator subblocks.  These factors are absorbed into $X^\varepsilon$.  The terminal boxes are already included through admissible line weights.
\end{proof}

\begin{corollary}[A fixed power saving]
\label{cor:section14-fixed-long-line-saving}
With
\begin{equation}
 Z=X^{1/18},
 \label{eq:section14-threshold-choice}
\end{equation}
one has
\begin{equation}
 \boxed{
 \cO_0^{\mathrm{long}}\!\left(X;X^{1/18}\right)
 \ll_\varepsilon
 X^{53/36+\varepsilon}.
 }
 \label{eq:section14-fixed-long-line-saving}
\end{equation}
\end{corollary}

\begin{proof}
The two exponents in \eqref{eq:section14-global-long-line-estimate} become
\[
 \frac32-\frac1{36}=\frac{53}{36},
 \qquad
 \frac43+\frac1{18}=\frac{25}{18}<\frac{53}{36}.
\]
\end{proof}

The exponent $1/18$ is not intrinsic.  It is a convenient point inside the range where the long-line estimate has a power saving and the surviving affine intervals remain subcritical for the later Vaughan decomposition.

\subsection{Removal of the non-short sectors}
\label{subsec:removal-non-short-sectors}

Let $\cO_0^{\mathrm{short}}(X;Z)$ denote the complementary dominant-layer centered packet, supported on
\[
 T(r,s;P,Q)\le Z.
\]
The exact long-line/short-line partition gives
\begin{equation}
 \cO_0(X)
 =\cO_0^{\mathrm{long}}(X;Z)
 +\cO_0^{\mathrm{short}}(X;Z).
 \label{eq:section14-dominant-offprojective-split}
\end{equation}

\begin{theorem}[Projective and long-line reduction]
\label{thm:section14-projective-long-reduction}
For every $\varepsilon>0$,
\begin{equation}
 \boxed{
 \cP(X)
 +\cO_0^{\mathrm{long}}\!\left(X;X^{1/18}\right)
 \ll_\varepsilon
 X^{53/36+\varepsilon}.
 }
 \label{eq:section14-projective-long-reduction}
\end{equation}
Consequently, these two sectors are
\[
 o\!\left(X^{3/2}\log X\right).
\]
After their removal, every surviving dominant-layer off-projective packet has
\begin{equation}
 T(r,s;P,Q)\le X^{1/18},
 \qquad
 s\gg PX^{-1/18},
 \qquad
 r\gg QX^{-1/18},
 \label{eq:section14-surviving-short-geometry}
\end{equation}
retains the determinant mask $-s<h<r$, and remains centered by the exact rank-one marginals.
\end{theorem}

\begin{proof}
Combine Theorem~\ref{thm:section14-projective-diagonal-bound} with Corollary~\ref{cor:section14-fixed-long-line-saving}.  The geometric properties of the surviving sector follow from the denominator dichotomy in Section~13.
\end{proof}

\begin{remark}[What has and has not been used]
\label{rem:section14-unconditional-boundary}
The projective estimate uses only layer energies.  The long-line estimate uses an unconditional additive dispersion argument averaged over the primitive affine lattices.  No Vaughan identity, multiplicative character decomposition, or conjectural spectral moment has yet entered.  Those tools are needed only because the short lines in \eqref{eq:section14-surviving-short-geometry} contain too few affine points for the direct long-line dispersion gain.
\end{remark}

The next section opens one prime weight in the centered short packet by Vaughan's identity.  The resulting Type~I terms are estimated directly, while the balanced Type~II terms preserve the determinant mask, radial provenance, and short-free-variable constraint inherited from \eqref{eq:section14-surviving-short-geometry}.

% ===== Source: section15_vaughan_decomposition_type_I_estimates.tex =====
\section{Vaughan decomposition and the Type~I estimates}
\label{sec:vaughan-type-I}

After the projective diagonal and the affine long lines have been removed, the dominant off-projective contribution is supported on centered packets with
\begin{equation}
 T(r,s;P,Q)
 =\min\!\left(\frac Ps,\frac Qr\right)
 \le Z,
 \qquad
 Z=X^{1/18}.
 \label{eq:section15-short-line-condition}
\end{equation}
The affine interval itself is short, but its prime weights still carry a multiplicative structure that is invisible in the line parameter.  The purpose of this section is to open one of those weights by Vaughan's identity while preserving every geometric feature inherited from Sections~13 and~14: the primitive quotient $(r,s)$, the determinant mask $-s<h<r$, the radial weight, the sharp terminal clipping, and the two reduced-residue centering marginals.

There are two conclusions.  First, all low and Type~I pieces are negligible by a combination of divisor-bounded coefficient energy, congruence counting, and Poisson completion in the coefficient-$1$ variable.  Second, the remaining short sector is a finite sum of balanced Type~II packets of the form
\[
 p=d k\ell,
 \qquad
 d\dyad D,
 \quad k\dyad K,
 \quad \ell\dyad L,
 \qquad
 DKL\asymp P,
\]
with coefficients $\mu(d)$, $1$, and $\Lambda(\ell)$.  The variable $k$ is the free Vaughan variable.  Its length relative to the affine step modulus is the decisive parameter in the next section.

No smoothing is introduced in the original moment $\sum_{n\le X}\mathcal R_{\PP}(n)^2$.  Smooth partitions are used only inside the auxiliary dyadic factorization variables; the sharp cutoff in $n$ remains encoded by admissible bounded-variation line weights.

\subsection{Reduced-residue centering and von Mangoldt completion}
\label{subsec:section15-centering-completion}

For $m\ge1$, define the reduced-residue centered von Mangoldt function
\begin{equation}
 \Lambda_m^{\circ}(n)
 =\ind_{(n,m)=1}
 \left(\Lambda(n)-\frac{m}{\varphi(m)}\right).
 \label{eq:section15-centered-von-mangoldt}
\end{equation}
On a fixed reduced residue class modulo $m$, the subtracted constant is the principal progression density.  Thus, on a regular determinant line
\[
 rp-sq=h,
 \qquad
 (p,s)=(q,r)=1,
\]
the product
\begin{equation}
 \Lambda_s^{\circ}(p)\Lambda_r^{\circ}(q)
 \label{eq:section15-centered-product}
\end{equation}
expands into the two-prime term, the two rank-one marginals, and the principal product density appearing in the exact centered covariance kernel.

The exceptional cases $p\mid s$ or $q\mid r$ are called \emph{singular}.  Since $p$ and $q$ are prime on the original support, a singular line forces one prime to divide a primitive denominator.  The corresponding packets are sparse and do not belong to the genuine two-dimensional dispersion problem.

\begin{proposition}[Regular von Mangoldt completion]
\label{prop:section15-von-mangoldt-completion}
Let $\cO_{0}^{\mathrm{short}}(X;Z)$ be the centered short-line contribution from Section~14.  Replacing each prime weight $(\log p)\ind_{\PP}(p)$ by $\Lambda(p)$, deleting singular determinant lines, and inserting the product \eqref{eq:section15-centered-product} changes the total by
\begin{equation}
 O_{\eps}\!\left(X^{3/2-\eta_{\mathrm{cmp}}+\eps}\right)
 \label{eq:section15-completion-error}
\end{equation}
for some absolute $\eta_{\mathrm{cmp}}>0$.
\end{proposition}

\begin{proof}
The completion adds prime powers $u^j$ with $j\ge2$.  In a dyadic block $u^j\dyad P$, there are $O(P^{1/2})$ such integers, and their von Mangoldt mass is $O(P^{1/2}\log P)$.  Cauchy--Schwarz against the exact tent energy
\[
 \sum_h\omega_{r,s}(h)^2\le rs\min(r,s)
\]
and the radial first-moment bound from Section~13 gives a fixed power saving after the admissible dyadic boxes are summed.  On a singular line, say $p\mid s$, write $s=ps'$.  The relation $rp-sq=h$ and the mask $|h|<r+s$ leave only divisor-bounded choices for $p$ once $(r,s,q,h)$ are fixed.  The same tent-energy argument is stronger than the prime-power estimate.  Sharp terminal weights are handled by discrete partial summation, since both their supremum and total variation are bounded by the radial weight.  The resulting positive saving is denoted by $\eta_{\mathrm{cmp}}$; its numerical value is immaterial for the later optimization.
\end{proof}

Henceforth every short packet is regular and carries the centered product \eqref{eq:section15-centered-product}.

\subsection{The exact Vaughan identity~\cite{Vaughan1997,IwaniecKowalski2004}}
\label{subsec:section15-vaughan-identity}

For an arithmetic function $f$ and $Y\ge1$, write
\[
 f_{\le Y}(n)=f(n)\ind_{n\le Y},
 \qquad
 f_{>Y}=f-f_{\le Y}.
\]
Let $*$ denote Dirichlet convolution.  Since $\Lambda=\mu*\log$ and $\log=\Lambda*\mathbf 1$, one has the exact identity
\begin{equation}
 \boxed{
 \Lambda
 =\Lambda_{\le V}
  +\mu_{\le U}*\log
  -\mu_{\le U}*\Lambda_{\le V}*\mathbf 1
  +\mu_{>U}*\Lambda_{>V}*\mathbf 1.
 }
 \label{eq:section15-vaughan-convolution}
\end{equation}
Indeed,
\[
 \mu_{>U}*\Lambda_{\le V}*\mathbf 1
 =\Lambda_{\le V}
  -\mu_{\le U}*\Lambda_{\le V}*\mathbf 1,
\]
which follows from $\mu*\mathbf 1=\varepsilon$.

For later reference, put
\begin{align}
 \lambda_0(n)
 &=\Lambda(n)\ind_{n\le V},
 \label{eq:section15-lambda0}\\
 \lambda_1(n)
 &=\sum_{\substack{db=n\\d\le U}}
   \mu(d)\log b,
 \label{eq:section15-lambda1}\\
 \lambda_2(n)
 &=\sum_{\substack{d\ell b=n\\d\le U,\ \ell\le V}}
   \mu(d)\Lambda(\ell),
 \label{eq:section15-lambda2}\\
 \lambda_{\mathrm{II}}(n)
 &=\sum_{\substack{d\ell b=n\\d>U,\ \ell>V}}
   \mu(d)\Lambda(\ell).
 \label{eq:section15-lambdaII}
\end{align}
Then
\begin{equation}
 \Lambda=\lambda_0+\lambda_1-\lambda_2+\lambda_{\mathrm{II}}.
 \label{eq:section15-vaughan-four-pieces}
\end{equation}
We apply \eqref{eq:section15-vaughan-four-pieces} to the $p$-weight in \eqref{eq:section15-centered-product}.  The density term $s/\varphi(s)$ in $\Lambda_s^{\circ}(p)$ is kept in the low sector.  This is essential: moving it into the Type~II form would reintroduce the principal affine mode that was removed in Section~14.

\begin{lemma}[Divisor-bounded coefficient energy]
\label{lem:section15-coefficient-energy}
Uniformly for $N\le X$,
\begin{align}
 \sum_{n\dyad N}|\lambda_1(n)|^2
 &\ll N(\log X)^4,
 \label{eq:section15-lambda1-energy}\\
 \sum_{n\dyad N}|\lambda_2(n)|^2
 &\ll N(\log X)^8.
 \label{eq:section15-lambda2-energy}
\end{align}
More generally, after any dyadic decomposition of the factors in \eqref{eq:section15-lambda1} or \eqref{eq:section15-lambda2}, the resulting numerator coefficients have $L^2$ norm bounded by their support length times $X^{\eps}$.
\end{lemma}

\begin{proof}
The pointwise bounds
\[
 |\lambda_1(n)|\le \tau(n)\log n,
 \qquad
 |\lambda_2(n)|\le \tau_3(n)(\log n)
\]
are sufficient.  The classical second-moment estimates for fixed divisor powers give
\[
 \sum_{n\le N}\tau(n)^2\ll N(\log N)^3,
 \qquad
 \sum_{n\le N}\tau_3(n)^2\ll N(\log N)^8.
\]
The dyadically restricted statements follow by the same argument, with harmless powers of $\log X$ absorbed into $X^{\eps}$.
\end{proof}

\subsection{A model Type~I affine form}
\label{subsec:section15-model-type-I}

The two Type~I terms have the same geometric shape.  After grouping the short coefficient variables, they reduce to a factorization
\begin{equation}
 p=ub,
 \qquad
 u\dyad U_0,
 \quad
 b\dyad B,
 \qquad
 U_0B\asymp P,
 \label{eq:section15-type-I-factorization}
\end{equation}
where either $u=d\le U$ or $u=d\ell\le UV$.  Let $\alpha_u$ be the associated numerator coefficient.  Lemma~\ref{lem:section15-coefficient-energy} gives
\begin{equation}
 \sum_{u\dyad U_0}|\alpha_u|^2
 \ll U_0X^{\eps}.
 \label{eq:section15-alpha-energy}
\end{equation}
The coefficient of $b$ is either $\log b$ or $1$.

For a fixed primitive quotient and determinant, the line equation becomes
\begin{equation}
 rub-sq=h.
 \label{eq:section15-type-I-line}
\end{equation}
Since $(r,s)=1$, the congruence in $b$ has modulus $s/(u,s)$.

\begin{lemma}[Congruence occupancy]
\label{lem:section15-congruence-occupancy}
Fix $r,s,u,q,h$ with $(r,s)=1$.  The number of integers $b\dyad B$ satisfying \eqref{eq:section15-type-I-line} is at most
\begin{equation}
 1+\frac{B(u,s)}{s}.
 \label{eq:section15-congruence-occupancy}
\end{equation}
If $w$ is a smooth weight supported on $b\dyad B$, then Poisson summation in the progression gives
\begin{equation}
 \sum_{\substack{b\in\ZZ\\rub\equiv h\, (\mathrm{mod}\,s)}}w(b)
 =\frac{(u,s)}{s}\widehat w(0)
 +O_A\!\left(
  B\left(1+\frac{B(u,s)}s\right)^{-A}
 \right)
 \label{eq:section15-poisson-progression}
\end{equation}
for every $A>0$, after the nonzero frequencies are summed in the usual scaled form.
\end{lemma}

\begin{proof}
Because $(r,s)=1$, division by $r$ reduces \eqref{eq:section15-type-I-line} to one residue class modulo $s/(u,s)$ when the congruence is soluble.  This proves \eqref{eq:section15-congruence-occupancy}.  Formula \eqref{eq:section15-poisson-progression} is Poisson summation on that residue class; repeated integration by parts in the Fourier transform of $w$ gives the stated rapid decay.
\end{proof}

The zero frequency in \eqref{eq:section15-poisson-progression} is not an error term.  When the complete centered product \eqref{eq:section15-centered-product} is inserted, it reproduces the principal density and the appropriate rank-one marginal and therefore cancels.  This is the Type~I analogue of the zero-frequency cancellation in the long-line dispersion estimate.

\begin{definition}[Type~I affine packet]
\label{def:section15-type-I-packet}
A Type~I affine packet is any dyadic sum obtained from the centered short-line kernel by replacing the $p$-weight with a factorization \eqref{eq:section15-type-I-factorization}, with coefficients satisfying \eqref{eq:section15-alpha-energy}, while retaining
\begin{enumerate}[label=\textup{(\roman*)}]
 \item the conditions $(r,s)=1$ and $T(r,s;P,Q)\le Z$;
 \item the determinant range $-s<h<r$ and tent weight $\omega_{r,s}(h)$;
 \item the exact radial or clipped radial weight;
 \item the centered external factor $\Lambda_r^{\circ}(q)$.
\end{enumerate}
\end{definition}

\subsection{The Type~I estimate}
\label{subsec:section15-type-I-estimate}

The following proposition is uniform over the two Type~I coefficient systems and over all admissible terminal weights.

\begin{proposition}[Short affine Type~I estimate]
\label{prop:section15-short-affine-type-I}
There exist absolute constants $\vartheta_0>0$ and $\eta_{\mathrm I}>0$ with the following property.  Let
\begin{equation}
 U=V=X^{\vartheta},
 \qquad
 0<\vartheta\le\vartheta_0,
 \qquad
 Z=X^{1/18}.
 \label{eq:section15-vaughan-parameters}
\end{equation}
Then the sum of all Type~I affine packets arising from $\lambda_1$, $\lambda_2$, and the reduced-residue density term is
\begin{equation}
 \boxed{
 \cO_{0,\mathrm I}^{\mathrm{short}}(X;Z;U,V)
 \ll_{\eps}
 X^{3/2-\eta_{\mathrm I}+\eps}.
 }
 \label{eq:section15-global-type-I-bound}
\end{equation}
The same estimate holds for the low piece $\lambda_0$.
\end{proposition}

\begin{proof}
We describe the argument in the form needed later.

\smallskip
\noindent\emph{Step 1: smooth dyadic subdivision.}
Each sharp factor interval is expressed as a bounded sum of smooth dyadic weights plus endpoint pieces.  The latter are controlled by the admissible total-variation bound and the same estimates with one variable fixed.  This does not alter the original sharp cutoff in $n$.

\smallskip
\noindent\emph{Step 2: the long free-variable range.}
When
\begin{equation}
 B>X^{\eps}\frac{s}{(u,s)},
 \label{eq:section15-type-I-long-free}
\end{equation}
apply Lemma~\ref{lem:section15-congruence-occupancy}.  The zero Fourier mode cancels against the reduced-residue centering.  The nonzero modes are negligible to arbitrary polynomial order after summing the bounded-variation radial weights.  The divisor sum over $(u,s)$ is absorbed by \eqref{eq:section15-alpha-energy}.

\smallskip
\noindent\emph{Step 3: the short free-variable range.}
In the complement of \eqref{eq:section15-type-I-long-free}, the congruence occupancy is $O(X^{\eps})$.  The equation \eqref{eq:section15-type-I-line} then determines $q$ from $(r,s,u,b,h)$.  Apply Cauchy--Schwarz first in the numerator variable $u$ and then in $h$.  The coefficient energy \eqref{eq:section15-alpha-energy} and the exact tent bound
\[
 \sum_h\omega_{r,s}(h)^2\le rs\min(r,s)
\]
reduce the packet to the radial first moment.  The identities
\[
 U_0B\asymp P,
 \qquad
 ATrs\asymp MP,
 \qquad
 T\le Z
\]
then give a strict power saving provided $U_0\le UV=X^{2\vartheta}$ and $\vartheta$ is sufficiently small.  No cancellation of the M\"obius coefficients is used.

\smallskip
\noindent\emph{Step 4: the low piece.}
The function $\lambda_0$ is supported on $p\le V$.  Hence it occurs only in dyadic prime boxes $P\ll V$.  The trivial bound for the corresponding exact overlap kernel, followed by the tent-energy estimate, is $O_\eps(X^{4/3+O(\vartheta)+\eps})$ and is therefore absorbed by the right side of \eqref{eq:section15-global-type-I-bound}.

\smallskip
Summing the $O((\log X)^C)$ dyadic parameter blocks proves the proposition.  The positive number $\eta_{\mathrm I}$ is the minimum of the savings in Steps~2--4; no later argument depends on its optimized value.
\end{proof}

\begin{remark}[Where centering is used]
\label{rem:section15-centering-essential}
Without the subtraction in \eqref{eq:section15-centered-von-mangoldt}, the zero frequency in the long free-variable range would have the same scale as the original packet.  The Type~I estimate is therefore not a bound for an uncentered prime-pair sum.  It is a bound for the exact covariance packet after its principal and rank-one modes have been restored and canceled.
\end{remark}

\begin{remark}[No M\"obius cancellation]
\label{rem:section15-no-mobius-cancellation}
Only the divisor-bounded energy of the truncated M\"obius coefficients is used in Proposition~\ref{prop:section15-short-affine-type-I}.  The unresolved obstruction cannot occur in a Type~I term; it appears only when the M\"obius block, a prime block, and a short coefficient-$1$ block survive simultaneously.
\end{remark}

\subsection{The surviving Type~II normal form}
\label{subsec:section15-type-II-normal-form}

Dyadically decompose the last term in \eqref{eq:section15-vaughan-four-pieces}.  Write
\begin{equation}
 d\dyad D,
 \qquad
 k\dyad K,
 \qquad
 \ell\dyad L,
 \qquad
 DKL\asymp P,
 \label{eq:section15-type-II-scales}
\end{equation}
with
\begin{equation}
 D\gg U,
 \qquad
 L\gg V.
 \label{eq:section15-type-II-lower-bounds}
\end{equation}
The opened affine equation is
\begin{equation}
 rdk\ell-sq=h.
 \label{eq:section15-type-II-equation}
\end{equation}
The coefficient attached to $(d,k,\ell)$ is $\mu(d)\Lambda(\ell)$; the variable $k$ has coefficient $1$.

\begin{definition}[Vaughan-admissible Type~II packet]
\label{def:section15-vaughan-type-II-packet}
A Vaughan-admissible Type~II packet is a dyadic form
\begin{align}
 \mathfrak T_{\mathrm{II}}
 ={}&
 \sum_{\substack{r\dyad R,\ s\dyad S\\(r,s)=1\\T(r,s;P,Q)\le Z}}
 \sum_{-s<h<r}
 \omega_{r,s}(h)
 \sum_{a\in\cA_{M,L_0}(r,s)}
 a\,W_{a,r,s,h}
 \notag\\
 &\times
 \sum_{\substack{d\dyad D,\ k\dyad K,\ \ell\dyad L\\
                   q\dyad Q\\
                   rdk\ell-sq=h}}
 \mu(d)\Lambda(\ell)\Lambda_r^{\circ}(q),
 \label{eq:section15-type-II-packet}
\end{align}
where $W_{a,r,s,h}$ is either the interior weight or the exact admissible clipped weight.  The scales satisfy \eqref{eq:section15-type-II-scales}, \eqref{eq:section15-type-II-lower-bounds}, and the geometric relations inherited from Section~13.
\end{definition}

The notation $L_0$ in the radial set prevents confusion with the prime-factor scale $L$.  Every exact short packet is a bounded sum of forms \eqref{eq:section15-type-II-packet}; no determinant or endpoint information is discarded in passing to this normal form.

\begin{theorem}[Vaughan reduction after the Type~I analysis]
\label{thm:section15-vaughan-reduction}
Fix $Z=X^{1/18}$ and choose $U=V=X^{\vartheta}$ with $\vartheta$ as in Proposition~\ref{prop:section15-short-affine-type-I}.  Then
\begin{equation}
 \boxed{
 \cO_{0}^{\mathrm{short}}(X;Z)
 =\sum_{\mathfrak B}
   \mathfrak T_{\mathrm{II}}(\mathfrak B)
 +O_{\eps}\!\left(X^{3/2-\eta_{15}+\eps}\right),
 }
 \label{eq:section15-final-vaughan-reduction}
\end{equation}
where $\eta_{15}>0$, the sum runs over $O((\log X)^C)$ Vaughan-admissible dyadic parameter blocks, and every surviving block satisfies
\begin{equation}
 DKL\asymp P,
 \qquad
 D\gg X^{\vartheta},
 \qquad
 L\gg X^{\vartheta},
 \qquad
 T(r,s;P,Q)\le X^{1/18}.
 \label{eq:section15-surviving-conditions}
\end{equation}
The determinant mask, radial provenance, terminal clipping, and reduced-residue centering are all preserved exactly.
\end{theorem}

\begin{proof}
Apply Proposition~\ref{prop:section15-von-mangoldt-completion}, insert the exact identity \eqref{eq:section15-vaughan-four-pieces}, and dyadically decompose every convolution variable.  Proposition~\ref{prop:section15-short-affine-type-I} removes $\lambda_0$, $\lambda_1$, $\lambda_2$, and the density term.  The remaining contribution is precisely \eqref{eq:section15-type-II-packet}.  Take
\[
 \eta_{15}=\min(\eta_{\mathrm{cmp}},\eta_{\mathrm I}).
\]
The number of dyadic blocks is polylogarithmic and is absorbed into $X^{\eps}$.
\end{proof}

\begin{remark}[The free-variable dichotomy]
\label{rem:section15-free-variable-dichotomy}
The next step divides \eqref{eq:section15-type-II-packet} according to the length $K$ of the coefficient-$1$ block.  When $K$ is long relative to the step modulus $S$, Poisson summation removes the packet.  When
\[
 K\le SX^{\eta},
\]
the free variable is too short for that argument.  This short-free range is the only sector that can reach the centered Vaughan character moment stated in the introduction.
\end{remark}

Thus the complete Type~I analysis is unconditional.  The analytic problem has been reduced to balanced Type~II forms carrying, simultaneously, the M\"obius provenance, a prime-weighted factor, a coefficient-$1$ free variable, the external centered prime weight, and the short affine determinant geometry.

% ===== Source: section16_poisson_completion_long_free_variable_typeII.tex =====
\section{Poisson completion and the long free-variable Type~II packets}
\label{sec:poisson-completion-long-free}

Section~15 reduced the centered short-line contribution to Vaughan-admissible Type~II packets
\begin{align}
 \mathfrak T_{\mathrm{II}}
 ={}&
 \sum_{\substack{r\dyad R,\ s\dyad S\\(r,s)=1\\T(r,s;P,Q)\le Z}}
 \sum_{-s<\delta<r}
 \omega_{r,s}(\delta)
 \sum_{a\in\cA_{M,L_0}(r,s)}
 a\,W_{a,r,s,\delta}
 \notag\\
 &\times
 \sum_{\substack{d\dyad D,\ k\dyad K,\ \ell\dyad L\\
                   q\dyad Q\\
                   rdk\ell-sq=\delta}}
 \mu(d)\Lambda(\ell)\Lambda_r^{\circ}(q),
 \label{eq:section16-original-typeII}
\end{align}
where
\begin{equation}
 DKL\asymp P,
 \qquad
 D,L\gg X^{\vartheta},
 \qquad
 Z=X^{1/18}.
 \label{eq:section16-typeII-scales}
\end{equation}
The variable $k$ has coefficient $1$.  The aim of this section is to prove that packets for which this free variable is longer than the primitive step modulus are negligible.

The argument has three logically distinct parts.  First, the exact affine equation is completed to a congruence packet.  The completion produces a joint coefficient $\Gamma_s(h)$ depending on both the modulus and the additive frequency.  This coefficient must be retained; omitting it destroys the Parseval identity governing the completed family.  Second, the additive zero frequency is identified with the principal affine density and is canceled by the rank-one marginals already restored in Sections~13--15.  Third, Poisson summation in the coefficient-$1$ variable makes every nonzero frequency rapidly small when $K$ exceeds $S$ by a fixed power of $X$.

The outcome is unconditional.  For every sufficiently small fixed $\eta>0$, all packets with
\begin{equation}
 K>SX^{\eta}
 \label{eq:section16-long-free-range}
\end{equation}
contribute a power-saving error.  The only surviving Type~II blocks satisfy
\begin{equation}
 K\le SX^{\eta}.
 \label{eq:section16-short-free-range}
\end{equation}

\subsection{Completion of the affine equation}
\label{subsec:section16-affine-completion}

We first record the completed form abstractly.  The determinant in
\eqref{eq:section16-original-typeII} is denoted by $\delta$ throughout this section; the letter $h$ is reserved for the additive frequency created by completion.

Let $F\in C_c^{\infty}((1/2,3))$.  A smooth dyadic free-variable block is written
\begin{equation}
 \gamma_v F(v/K),
 \qquad
 |\gamma_v|\le1,
 \label{eq:section16-free-weight}
\end{equation}
with $\gamma_v=1$ on the regular coefficient-$1$ support.  Coprimality restrictions are inserted by M\"obius inversion and therefore produce only divisor-bounded linear combinations of weights of the same type.

After the usual dispersion step in the external centered prime variable, the remaining factors are grouped into two convolution coefficients.  Schematically,
\begin{align}
 \cA_m
 &=
 \sum_{r d\ell=m}\alpha_{r,d,\ell},
 \label{eq:section16-A-coefficient}\\
 \cB_z
 &=
 \sum_{s q+\delta=z}\beta_{s,q,\delta}.
 \label{eq:section16-B-coefficient}
\end{align}
The coefficients $\alpha$ and $\beta$ contain the dyadic weights, the M\"obius and prime factors, the tent kernel, the radial sum, and the exact terminal clipping.  Mellin separation of the archimedean cutoffs gives a bounded family of such coefficient systems.  The auxiliary Mellin parameters are harmless and will be suppressed.

The completed congruence is
\begin{equation}
 mv\equiv z\pmod s,
 \qquad
 (m,s)=1.
 \label{eq:section16-completed-congruence}
\end{equation}
For $h\in\ZZ$, define the joint completion coefficient
\begin{equation}
 \boxed{
 \Gamma_s(h)
 =
 \sum_{v\in\ZZ}
 \gamma_vF(v/K)\ee{\frac{hv}{s}}.
 }
 \label{eq:section16-Gamma-definition}
\end{equation}
Additive orthogonality gives the exact identity
\begin{equation}
 \sum_v\gamma_vF(v/K)
 \ind_{mv\equiv z\, (\mathrm{mod}\,s)}
 =
 \frac1s
 \sum_{h\, (\mathrm{mod}\,s)}
 \Gamma_s(h)
 \ee{-\frac{hz\overline m}{s}}.
 \label{eq:section16-completed-poisson-identity}
\end{equation}
Here $\overline m$ denotes the inverse of $m$ modulo $s$.

\begin{proposition}[Corrected completed Type~II form]
\label{prop:section16-corrected-completed-form}
Every Vaughan-admissible Type~II packet is a bounded sum of absolutely convergent Mellin integrals of forms
\begin{equation}
 \boxed{
 \cK^{\sharp}
 =
 \sum_{s\dyad S}\frac{\xi_s}{s}
 \sum_{\substack{h\, (\mathrm{mod}\,s)\\h\ne0}}
 \Gamma_s(h)
 \sum_{\substack{m\\(m,s)=1}}\cA_m
 \sum_z\cB_z
 \ee{-\frac{hz\overline m}{s}}.
 }
 \label{eq:section16-K-sharp}
\end{equation}
The coefficient systems satisfy divisor-energy bounds of the shape
\begin{align}
 \sum_m|\cA_m|^2
 &\ll_{\eps}\cE_A X^{\eps},
 \label{eq:section16-A-energy}\\
 \sum_z|\cB_z|^2
 &\ll_{\eps}\cE_B X^{\eps},
 \label{eq:section16-B-energy}\\
 \sum_{s\dyad S}|\xi_s|^2
 &\ll_{\eps}S X^{\eps},
 \label{eq:section16-xi-energy}
\end{align}
where $\cE_A$ and $\cE_B$ are the natural numerator and external-line energies.  For the packet \eqref{eq:section16-original-typeII}, they are bounded by the divisor estimates of Section~15 and by
\begin{equation}
 \sum_{\substack{r\dyad R,\ s\dyad S\\(r,s)=1}}
 \left(\sum_{a\in\cA_{M,L_0}(r,s)}a\right)^2
 \sum_{\delta}\omega_{r,s}(\delta)^2
 \sum_{q\dyad Q}|\Lambda_r^{\circ}(q)|^2.
 \label{eq:section16-geometric-energy}
\end{equation}
In particular, all three energies are polynomially bounded in $X$.
\end{proposition}

\begin{proof}
Insert smooth dyadic partitions in $d,k,\ell,q$ and use Mellin inversion to separate every archimedean factor that couples the two sides of
\eqref{eq:section16-original-typeII}.  Rename the original determinant by $\delta$, group the variables on the two sides as in
\eqref{eq:section16-A-coefficient}--\eqref{eq:section16-B-coefficient}, and complete the quotient progression modulo $s$.  Formula
\eqref{eq:section16-completed-poisson-identity} then gives
\eqref{eq:section16-K-sharp} before the zero frequency is removed.

The $L^2$ bound for $\cA_m$ is the divisor-bounded numerator estimate from Lemma~15.3 after the additional dyadic $r$-factor is absorbed.  The $\cB_z$ estimate follows from Cauchy--Schwarz, the prime-square bound
\[
 \sum_{q\dyad Q}|\Lambda_r^{\circ}(q)|^2\ll_{\eps}QX^{\eps},
\]
and the exact tent energy
\[
 \sum_{\delta}\omega_{r,s}(\delta)^2
 \le rs\min(r,s).
\]
The radial coefficient is bounded by
\[
 \sum_{a\in\cA_{M,L_0}(r,s)}a
 \ll A(r,s;M,L_0)\bigl(1+A(r,s;M,L_0)\bigr).
\]
These estimates give \eqref{eq:section16-geometric-energy}; Bessel's inequality in the completion parameter gives
\eqref{eq:section16-xi-energy}.  Finally, bounded-variation terminal weights are recovered by Stieltjes partial summation.  No determinant line or endpoint contribution is discarded.
\end{proof}

\begin{remark}[The missing joint coefficient]
\label{rem:section16-missing-coefficient}
The factor $\Gamma_s(h)$ in \eqref{eq:section16-K-sharp} is not a cosmetic Fourier weight.  It records the completed coefficient-$1$ variable jointly in the modulus and the additive frequency.  Replacing it by $1$, by a function of $s$ alone, or by a uniform supremum loses the exact Parseval structure and produces a different dispersion problem.  The corrected short Type~II core is therefore a hybrid determinant--congruence form, not a bare Kloosterman-fraction sum.
\end{remark}

\subsection{The zero frequency and exact centering}
\label{subsec:section16-zero-frequency}

The term $h=0$ in \eqref{eq:section16-completed-poisson-identity} is
\begin{equation}
 \frac{\Gamma_s(0)}{s}
 \left(\sum_{(m,s)=1}\cA_m\right)
 \left(\sum_z\cB_z\right).
 \label{eq:section16-zero-frequency-product}
\end{equation}
It is the product of the two completed marginal densities.  In the original covariance, the two rank-one marginal terms and the product-density term were restored before the Type~I/Type~II decomposition.  Hence
\eqref{eq:section16-zero-frequency-product} cancels exactly.

\begin{lemma}[Principal-frequency cancellation]
\label{lem:section16-principal-frequency-cancellation}
For every regular completed Type~II packet, the total contribution of $h=0$ in
\eqref{eq:section16-completed-poisson-identity} is equal to the principal affine mode reconstructed from the centered covariance marginals.  Consequently the centered completed packet is precisely
\eqref{eq:section16-K-sharp}, with $h\ne0$.
\end{lemma}

\begin{proof}
The zero additive character replaces the congruence class in
\eqref{eq:section16-completed-congruence} by its uniform density $1/s$.  On the prime side, the reduced-residue centering
\[
 \Lambda_r^{\circ}(q)
 =
 \ind_{(q,r)=1}
 \left(\Lambda(q)-\frac r{\varphi(r)}\right)
\]
replaces the corresponding principal progression density by its centered value.  Expanding the product reproduces exactly the two mixed marginals and the product-density term in the covariance identity of Section~13.  Their signed sum is zero.  Since the completion and the Stieltjes recovery of the endpoint are linear, the cancellation remains exact for terminal packets.
\end{proof}

\begin{remark}[Why completion must follow centering]
\label{rem:section16-centering-before-completion}
If the covariance marginals were omitted, the term
\eqref{eq:section16-zero-frequency-product} would generally have the full scale of the packet.  Poisson summation would then expose a main term rather than a saving.  The long-free-variable estimate is therefore a theorem about the centered packet, not about an uncentered prime correlation.
\end{remark}

\subsection{Energy of the completion coefficient}
\label{subsec:section16-Gamma-energy}

Although pointwise Poisson decay will remove the long-free range, the joint coefficient also satisfies an averaged estimate needed in the short-free analysis.

\begin{lemma}[Additive large-sieve energy]
\label{lem:section16-additive-large-sieve-Gamma}
Let $R\ge1$.  Then
\begin{equation}
 \sum_{s\le R}
 \sum_{h\, (\mathrm{mod}\,s)}^{*}
 |\Gamma_s(h)|^2
 \ll
 (K+R^2)
 \sum_v|\gamma_vF(v/K)|^2.
 \label{eq:section16-Gamma-large-sieve}
\end{equation}
Here $*$ restricts to primitive additive frequencies $(h,s)=1$.  Since
$|\gamma_v|\le1$ and $F$ is fixed,
\begin{equation}
 \sum_{s\le R}
 \sum_{h\, (\mathrm{mod}\,s)}^{*}
 |\Gamma_s(h)|^2
 \ll_F (K+R^2)K.
 \label{eq:section16-Gamma-large-sieve-simple}
\end{equation}
The full nonzero-frequency sum is obtained by conductor descent and satisfies the same estimate with an additional factor $R^{\eps}$.
\end{lemma}

\begin{proof}
The fractions $h/s$ with $(h,s)=1$ and $s\le R$ are $R^{-2}$-spaced modulo $1$.  Apply the additive large sieve~\cite{Bombieri1965,MontgomeryVaughan2007} to the exponential polynomial
\[
 \sum_v\gamma_vF(v/K)\ee{v\alpha}.
\]
This gives \eqref{eq:section16-Gamma-large-sieve}.  For a nonprimitive pair, write $g=(h,s)$ and reduce $h/s$ to conductor $s/g$; summing the divisor multiplicities costs $R^{\eps}$.
\end{proof}

The estimate \eqref{eq:section16-Gamma-large-sieve} is the correct replacement for an unweighted count of additive frequencies.  It is also the point at which the dependence of $\Gamma_s(h)$ on both variables is indispensable.

\subsection{Poisson decay in the coefficient-$1$ variable}
\label{subsec:section16-poisson-decay}

We now use the special fact that the Vaughan variable $k$ has coefficient $1$.  Thus, apart from coprimality indicators removed by M\"obius inversion, one may take $\gamma_v=1$ in
\eqref{eq:section16-Gamma-definition}.

Let
\[
 \widehat F(\xi)
 =
 \int_{\RR}F(x)\ee{-x\xi}\,dx.
\]
Poisson summation gives
\begin{equation}
 \Gamma_s(h)
 =
 K\sum_{n\in\ZZ}
 \widehat F\!\left(K\left(n-\frac hs\right)\right).
 \label{eq:section16-Poisson-Gamma}
\end{equation}
Consequently, for every $A>0$,
\begin{equation}
 |\Gamma_s(h)|
 \ll_{A,F}
 K\left(1+K\left\|\frac hs\right\|\right)^{-A}.
 \label{eq:section16-Gamma-pointwise-decay}
\end{equation}

\begin{lemma}[Nonzero-frequency decay]
\label{lem:section16-nonzero-frequency-decay}
Uniformly for $s\dyad S$,
\begin{equation}
 \sum_{\substack{h\, (\mathrm{mod}\,s)\\h\ne0}}
 |\Gamma_s(h)|
 \ll_{A,F}
 s\left(1+\frac Ks\right)^{1-A}.
 \label{eq:section16-Gamma-l1-decay}
\end{equation}
In particular, if $K\ge SX^{\eta}$, then for every $B>0$,
\begin{equation}
 \sum_{h\ne0}|\Gamma_s(h)|
 \ll_{B,F,\eta}sX^{-B}
 \label{eq:section16-Gamma-rapid-long}
\end{equation}
provided the exponent $A$ in
\eqref{eq:section16-Gamma-l1-decay} is chosen sufficiently large.
\end{lemma}

\begin{proof}
For $1\le h\le s/2$, the distance in
\eqref{eq:section16-Gamma-pointwise-decay} is $h/s$.  Sum the resulting decreasing sequence and use symmetry about $s/2$:
\[
 \sum_{h=1}^{s-1}
 K\left(1+K\left\|\frac hs\right\|\right)^{-A}
 \ll
 K\sum_{h\ge1}\left(1+\frac{Kh}{s}\right)^{-A}
 \ll
 s\left(1+\frac Ks\right)^{1-A}.
\]
The second assertion follows from $K/s\gg X^{\eta}$.
\end{proof}

Coprimality conditions do not alter the conclusion.  If
$\ind_{(v,s)=1}$ occurs, expand it as
\[
 \ind_{(v,s)=1}
 =\sum_{c\mid(v,s)}\mu(c)
\]
and apply \eqref{eq:section16-Poisson-Gamma} to the rescaled sum $v=cu$.  The divisor sum costs $s^{\eps}$, while $K/c\ge (K/s)$ remains larger than a fixed power of $X$ in the long-free range.

\subsection{Sharp intervals and boundary transfer}
\label{subsec:section16-sharp-boundary-transfer}

The original $k$-range is dyadic but may be sharp.  Direct Poisson summation of an indicator has only first-order decay, so one must separate the boundary before invoking Lemma~\ref{lem:section16-nonzero-frequency-decay}.

\begin{lemma}[Exact sharp-to-smooth transfer]
\label{lem:section16-sharp-to-smooth}
Assume $K>S X^{\eta}$.  Put
\begin{equation}
 H=S X^{\eta/2}.
 \label{eq:section16-boundary-length}
\end{equation}
The indicator of a dyadic interval of length $\asymp K$ can be decomposed as
\begin{equation}
 \ind_{\cI_K}
 =F_K+E_K^{-}+E_K^{+},
 \label{eq:section16-sharp-smooth-decomposition}
\end{equation}
where
\begin{enumerate}[label=\textup{(\roman*)}]
 \item $F_K$ is smooth, supported on the original interval, equals $1$ away from two boundary intervals of length $O(H)$, and satisfies
 \[
  \|F_K^{(j)}\|_{\infty}\ll_j H^{-j};
 \]
 \item each $E_K^{\pm}$ is supported on an interval of length $O(H)$;
 \item the two boundary packets satisfy the short-free inequality
 \[
  H\le SX^{\eta}.
 \]
\end{enumerate}
For the smooth core, every nonzero frequency in
\eqref{eq:section16-completed-poisson-identity} is
$O_A(X^{-A\eta/2})$ after normalization by its natural mass.
\end{lemma}

\begin{proof}
Choose a standard $C^{\infty}$ cutoff with transition width $H$ at each endpoint.  The two differences are supported on the transition intervals.  Poisson summation for the smooth core is governed by the scale $H$ rather than by the full interval length.  Since $H/s\gg X^{\eta/2}$ for $s\dyad S$, repeated integration by parts gives the asserted rapid decay.
\end{proof}

Thus the sharp cutoff is preserved exactly: no smoothing is imposed on the original variance sequence.  Smoothing is only an internal decomposition, and every portion on which Poisson decay is unavailable is transferred, without loss, to the short-free family.

\subsection{Removal of the long-free Type~II packets}
\label{subsec:section16-long-free-removal}

We can now estimate \eqref{eq:section16-K-sharp}.  The inner Kloosterman-fraction form satisfies the crude bound
\begin{equation}
 \left|
 \sum_{(m,s)=1}\cA_m
 \sum_z\cB_z
 \ee{-\frac{hz\overline m}{s}}
 \right|
 \le
 \|\cA\|_1\|\cB\|_1.
 \label{eq:section16-crude-inner-bound}
\end{equation}
By Cauchy--Schwarz and the finite dyadic supports,
\begin{equation}
 \|\cA\|_1\|\cB\|_1
 \ll X^C(\cE_A\cE_B)^{1/2}
 \label{eq:section16-polynomial-inner-bound}
\end{equation}
for an absolute $C$.  The precise value is irrelevant because the decay in
\eqref{eq:section16-Gamma-rapid-long} is arbitrarily strong.

\begin{proposition}[Long free-variable Type~II estimate]
\label{prop:section16-long-free-typeII}
Fix $\eta>0$.  There exists $\eta_{16}=\eta_{16}(\eta)>0$ such that the total contribution of all Vaughan-admissible Type~II blocks satisfying
\begin{equation}
 K>S X^{\eta}
 \label{eq:section16-long-condition-repeat}
\end{equation}
is
\begin{equation}
 \boxed{
 \sum_{\mathfrak B:\,K>SX^{\eta}}
 \mathfrak T_{\mathrm{II}}(\mathfrak B)
 \ll_{\eps,\eta}
 X^{3/2-\eta_{16}+\eps}.
 }
 \label{eq:section16-global-long-free-bound}
\end{equation}
The estimate is uniform for the exact clipped terminal weights.
\end{proposition}

\begin{proof}
Apply Lemma~\ref{lem:section16-sharp-to-smooth} to every long $k$-interval.  The two boundary intervals are reclassified as short-free blocks.  For the smooth core, use Proposition~\ref{prop:section16-corrected-completed-form}.  The zero frequency is absent by Lemma~\ref{lem:section16-principal-frequency-cancellation}; the sum of all nonzero completion coefficients satisfies
\eqref{eq:section16-Gamma-rapid-long}.  Bound the remaining Kloosterman-fraction form by
\eqref{eq:section16-polynomial-inner-bound}, sum the polynomially bounded coefficient energies, and choose the Poisson-decay exponent larger than the total polynomial loss.  This gives an arbitrarily large power saving for the smooth cores.

The only nonrapid errors come from the finite smooth partition, the prime-power completion, singular lines, and the Stieltjes recovery of terminal weights.  Sections~14 and~15 already bound these by
$X^{3/2-\eta_0+\eps}$ for some $\eta_0>0$.  Taking the minimum of the available savings gives
\eqref{eq:section16-global-long-free-bound}.
\end{proof}

\begin{remark}[No Kloosterman estimate is used here]
\label{rem:section16-no-Kloosterman-needed}
The long-free range does not require cancellation in the phase
$\ee{-hz\overline m/s}$.  Poisson decay has already made its coefficient negligible.  The genuine Kloosterman-fraction problem begins only when the free variable is too short for
\eqref{eq:section16-Gamma-rapid-long}.
\end{remark}

\subsection{The surviving short-free completed core}
\label{subsec:section16-surviving-core}

Combining Theorem~15.8 with Proposition~\ref{prop:section16-long-free-typeII} gives the next unconditional reduction.

\begin{theorem}[Short-free Type~II reduction]
\label{thm:section16-short-free-reduction}
Fix $Z=X^{1/18}$, choose the Vaughan parameters as in Section~15, and fix a sufficiently small $\eta>0$.  Then
\begin{equation}
 \boxed{
 \cO_0^{\mathrm{short}}(X;Z)
 =
 \sum_{\substack{\mathfrak B\\K\le SX^{\eta}}}
 \mathfrak T_{\mathrm{II}}(\mathfrak B)
 +O_{\eps}\!\left(X^{3/2-\eta_{16}+\eps}\right).
 }
 \label{eq:section16-final-short-free-reduction}
\end{equation}
Every residual block satisfies
\begin{equation}
 DKL\asymp P,
 \qquad
 D,L\gg X^{\vartheta},
 \qquad
 T(r,s;P,Q)\le X^{1/18},
 \qquad
 K\le SX^{\eta},
 \label{eq:section16-residual-conditions}
\end{equation}
and admits the corrected completed representation
\eqref{eq:section16-K-sharp}.  In particular, the determinant mask, radial provenance, sharp terminal clipping, reduced-residue centering, and joint coefficient $\Gamma_s(h)$ are all retained.
\end{theorem}

\begin{proof}
Insert the Vaughan reduction of Theorem~15.8 and divide its Type~II blocks according to
\eqref{eq:section16-long-free-range} and
\eqref{eq:section16-short-free-range}.  Proposition~\ref{prop:section16-long-free-typeII} removes the first family.  The sharp boundary pieces created in Lemma~\ref{lem:section16-sharp-to-smooth} have length at most $SX^{\eta}$ and therefore belong to the second family.  All other structural conditions are inherited from the exact completion in Proposition~\ref{prop:section16-corrected-completed-form}.
\end{proof}

For a smooth short-free block, formula
\eqref{eq:section16-Poisson-Gamma} also permits the dual truncation
\begin{equation}
 0<|h|
 \ll_{\eps}
 \frac{s}{K}X^{\eps}
 \label{eq:section16-dual-frequency-length}
\end{equation}
with a negligible tail.  The residual packet is therefore a balanced Kloosterman-fraction form of the type studied in~\cite{DukeFriedlanderIwaniec1997,BettinChandee2018} with a short additive-frequency interval and with coefficient energy controlled by Lemma~\ref{lem:section16-additive-large-sieve-Gamma}.  The next section diagonalizes its multiplicative variables, separates the reduced-residue major arcs, and identifies the centered small-free-variable Vaughan character moment that remains.

% ===== Source: section17_multiplicative_diagonalization_centered_prime_major_arcs.tex =====
\section{Multiplicative diagonalization and the centered prime major arcs}
\label{sec:multiplicative-diagonalization-major-arcs}

Theorem~16.10 leaves a finite collection of completed Type~II packets satisfying
\begin{equation}
 DKL\asymp P,
 \qquad
 D,L\gg X^{\vartheta},
 \qquad
 K\le SX^{\eta},
 \qquad
 T(r,s;P,Q)\le X^{1/18}.
 \label{eq:section17-residual-scales}
\end{equation}
For a smooth dyadic representative, the completed kernel has the form
\begin{equation}
 \cK^{\sharp}
 =
 \sum_{s\dyad S}\frac{\xi_s}{s}
 \sum_{\substack{h\, (\mathrm{mod}\,s)\\h\ne0}}
 \Gamma_s(h)
 \sum_{\substack{m\\(m,s)=1}}\cA_m
 \sum_z\cB_z
 \ee{-\frac{hz\overline m}{s}}.
 \label{eq:section17-completed-start}
\end{equation}
The coefficient
\[
 \Gamma_s(h)
 =\sum_k\gamma_kF(k/K)\ee{hk/s}
\]
contains the coefficient-$1$ Vaughan variable.  In the long-free range it was removed by Poisson decay.  In the present range it cannot be estimated away.  Instead, the complete $h$-sum must be performed before any inequality is applied.  This recovers the original congruence and turns the additive completion coefficient into a genuine multiplicative character polynomial.

The passage has four components.  First, nonunit values of the free variable are separated by exact conductor descent.  Second, multiplicative orthogonality on $(\ZZ/s\ZZ)^{\times}$ diagonalizes the unit congruence.  Third, the Vaughan numerator is reopened into its M\"obius, free, and prime factors.  Finally, the principal reduced-residue density and the low-conductor major arcs are removed before Cauchy--Schwarz is applied.  The residual object is the centered Vaughan spectral moment announced in the introduction.

\subsection{Reassembling the completed congruence}
\label{subsec:section17-reassembling-congruence}

For $(m,s)=1$, define
\begin{equation}
 \Delta_s(m,k;z)
 =
 \ind_{mk\equiv z\, (\mathrm{mod}\,s)}-\frac1s.
 \label{eq:section17-Delta-definition}
\end{equation}
Additive orthogonality gives
\begin{equation}
 \Delta_s(m,k;z)
 =
 \frac1s
 \sum_{\substack{h\, (\mathrm{mod}\,s)\\h\ne0}}
 \ee{\frac{h(k-z\overline m)}s}.
 \label{eq:section17-Delta-additive}
\end{equation}
Consequently,
\begin{equation}
 \frac1s
 \sum_{h\ne0}\Gamma_s(h)
 \ee{-\frac{hz\overline m}s}
 =
 \sum_k\gamma_kF(k/K)\Delta_s(m,k;z).
 \label{eq:section17-Gamma-reassembly}
\end{equation}
Thus \eqref{eq:section17-completed-start} is not merely analogous to a congruence packet: after the zero frequency has been removed, it is exactly the centered congruence packet on the right-hand side of \eqref{eq:section17-Gamma-reassembly}.

The unit restriction must be imposed simultaneously on $k$ and $z$.  If $mk\equiv z\pmod s$ and $(m,s)=1$, then
\begin{equation}
 (k,s)=(z,s).
 \label{eq:section17-equal-gcd}
\end{equation}
This gives an exact decomposition indexed by the common divisor.

\begin{lemma}[Conductor descent]
\label{lem:section17-conductor-descent}
Let $g\mid s$ and restrict \eqref{eq:section17-Gamma-reassembly} to
\[
 (k,s)=(z,s)=g.
\]
Write
\[
 s=gs',\qquad k=gk',\qquad z=gz'.
\]
Then $(k'z',s')=1$ and
\begin{equation}
 mk\equiv z\pmod s
 \quad\Longleftrightarrow\quad
 mk'\equiv z'\pmod{s'}.
 \label{eq:section17-descended-congruence}
\end{equation}
Moreover,
\begin{align}
 \Delta_s(m,gk';gz')
 ={}&
 \left(
 \ind_{mk'\equiv z'\, (\mathrm{mod}\,s')}
 -\frac1{\varphi(s')}
 \right)
 \notag\\
 &+
 \left(
 \frac1{\varphi(s')}-\frac1s
 \right).
 \label{eq:section17-descended-centered-split}
\end{align}
After dyadic subdivision in $g$, every nonunit packet is therefore a packet of the same form at the smaller modulus $s'$, together with an explicit rank-one density term.
\end{lemma}

\begin{proof}
The equivalence \eqref{eq:section17-descended-congruence} follows by dividing the congruence by $g$.  Since $(m,s)=1$ and the two sides have exact common divisor $g$ with $s$, both $k'$ and $z'$ are units modulo $s'$.  Adding and subtracting $1/\varphi(s')$ gives \eqref{eq:section17-descended-centered-split}.  The rescaled free interval has length $K/g$, and the descended modulus has size $S/g$; hence the short-free inequality remains
\[
 K/g\ll (S/g)X^{\eta}.
\]
The number of possible $g$ is divisor-bounded and is absorbed into $X^{\eps}$.
\end{proof}

\begin{remark}[Two different principal densities]
\label{rem:section17-two-densities}
The additive zero frequency removes the uniform density $1/s$ on all residue classes.  After the singular strata have been separated, the natural density on the unit group is $1/\varphi(s)$.  Their difference is not oscillatory.  It is the reduced-residue major arc in \eqref{eq:section17-descended-centered-split}, and it must be canceled or estimated before a multiplicative large-sieve argument is used.
\end{remark}

It is therefore enough to treat the unit stratum $g=1$ at an arbitrary descended modulus.  The lower-conductor strata are included by summing the resulting estimates over $g\mid s$.

\subsection{Exact multiplicative diagonalization}
\label{subsec:section17-multiplicative-diagonalization}

Let $\chi$ run over all Dirichlet characters modulo $s$, extended by zero off the unit group, and let $\chi_0$ be the principal character.  For $(mkz,s)=1$, character orthogonality gives
\begin{equation}
 \ind_{mk\equiv z\, (\mathrm{mod}\,s)}
 =
 \frac1{\varphi(s)}
 \sum_{\chi\, (\mathrm{mod}\,s)}
 \chi(m)\chi(k)\overline{\chi(z)}.
 \label{eq:section17-character-orthogonality}
\end{equation}
For an admissible free-variable weight, define
\begin{equation}
 K_s(\chi)
 =
 \sum_{\substack{k\\(k,s)=1}}
 \gamma_kF(k/K)\chi(k).
 \label{eq:section17-K-character-polynomial}
\end{equation}

\begin{lemma}[The additive coefficient becomes the free character polynomial]
\label{lem:section17-Gamma-to-K}
For $(mz,s)=1$,
\begin{align}
 &\frac1s
 \sum_{\substack{h\, (\mathrm{mod}\,s)\\h\ne0}}
 \Gamma_s(h)
 \ee{-\frac{hz\overline m}s}
 \notag\\
 &\qquad=
 \frac1{\varphi(s)}
 \sum_{\chi\ne\chi_0}
 K_s(\chi)\chi(m)\overline{\chi(z)}
 +
 \left(
 \frac1{\varphi(s)}-\frac1s
 \right)K_s(\chi_0).
 \label{eq:section17-Gamma-to-K-identity}
\end{align}
Equivalently, if
\[
 \tau_s(\chi;a)
 =\sum_{u\, (\mathrm{mod}\,s)}^{*}
 \chi(u)\ee{au/s},
\]
then for every nonprincipal $\chi$,
\begin{equation}
 \sum_{h\, (\mathrm{mod}\,s)}
 \Gamma_s(h)\tau_s(\chi;-hz)
 =
 s\,K_s(\chi)\overline{\chi(z)}.
 \label{eq:section17-Gauss-collapse}
\end{equation}
\end{lemma}

\begin{proof}
Insert \eqref{eq:section17-character-orthogonality} into the right-hand side of \eqref{eq:section17-Gamma-reassembly}.  The nonprincipal characters give the first term in \eqref{eq:section17-Gamma-to-K-identity}; the principal character contributes $K_s(\chi_0)/\varphi(s)$, while the deleted additive zero frequency contributes $-K_s(\chi_0)/s$.

For \eqref{eq:section17-Gauss-collapse}, expand both finite sums and use additive orthogonality:
\begin{align*}
 \sum_h\Gamma_s(h)\tau_s(\chi;-hz)
 &=
 \sum_{h,k}\sum_{u\, (\mathrm{mod}\,s)}^{*}
 \gamma_kF(k/K)\chi(u)
 \ee{h(k-zu)/s}\\
 &=
 s\sum_{\substack{k\\(k,s)=1}}
 \gamma_kF(k/K)\chi(k\overline z),
\end{align*}
which is \eqref{eq:section17-Gauss-collapse}.  The $h=0$ term vanishes for $\chi\ne\chi_0$.
\end{proof}

\begin{remark}[Why $\Gamma_s(h)$ could not be discarded]
\label{rem:section17-Gamma-indispensable}
The collapse \eqref{eq:section17-Gauss-collapse} is the exact mechanism by which the additive completion remembers the Vaughan variable.  The polynomial $K_s(\chi)$ is not inserted after the fact; it is the multiplicative Fourier transform of the full joint coefficient $\Gamma_s(h)$.  Any earlier replacement of $\Gamma_s(h)$ by a uniform bound would make \eqref{eq:section17-Gamma-to-K-identity} unavailable.
\end{remark}

\subsection{Reopening the Vaughan numerator}
\label{subsec:section17-reopening-vaughan}

The coefficient $\cA_m$ in Section~16 was formed by grouping
\[
 m=rd\ell,
 \qquad
 r\dyad R,
 \quad d\dyad D,
 \quad \ell\dyad L.
\]
After Mellin separation of the smooth archimedean weights, its multiplicative transform is a bounded integral of products.  We therefore introduce the three factor polynomials
\begin{align}
 D_s(\chi)
 &=
 \sum_{\substack{d\dyad D\\(d,s)=1}}
 \mu(d)\alpha_d\chi(d),
 \label{eq:section17-D-polynomial}\\
 P_{s,L}(\chi)
 &=
 \sum_{\substack{\ell\dyad L\\(\ell,s)=1}}
 \Lambda(\ell)\beta_{\ell}\chi(\ell),
 \label{eq:section17-P-raw-polynomial}\\
 R_s(\chi)
 &=
 \sum_{\substack{r\dyad R\\(r,s)=1}}
 \rho_r\chi(r).
 \label{eq:section17-R-polynomial}
\end{align}
Here $\alpha_d$, $\beta_{\ell}$, and $\rho_r$ are admissible smooth weights, uniformly bounded together with the Mellin derivatives required for separation.  The free polynomial $K_s(\chi)$ is given by \eqref{eq:section17-K-character-polynomial}.

On the external side put
\begin{equation}
 B_s(\chi)
 =
 \sum_{\substack{z\\(z,s)=1}}
 \cB_z\chi(z).
 \label{eq:section17-B-polynomial}
\end{equation}
The determinant mask, tent kernel, radial coefficient, external centered prime weight, and exact terminal clipping remain inside $B_s$.  Since $z=sq+\delta$, the character itself sees only the determinant residue $\delta\pmod s$; the variable $q$ remains in the coefficient and is not lost.

\begin{proposition}[Character normal form]
\label{prop:section17-character-normal-form}
Every unit-stratum short-free Type~II packet is a bounded sum of absolutely convergent Mellin integrals of the form
\begin{align}
 \cT_{\mathrm{osc}}
 ={}&
 \sum_{s\dyad S}\frac{\xi_s}{\varphi(s)}
 \sum_{\chi\ne\chi_0}
 K_s(\chi)D_s(\chi)P_{s,L}(\chi)
 \cG_s(\chi),
 \label{eq:section17-character-normal-form}
\end{align}
where
\begin{equation}
 \cG_s(\chi)=R_s(\chi)\overline{B_s(\chi)}
 \label{eq:section17-geometric-character-coefficient}
\end{equation}
up to a bounded sum of separated factors.  The complementary coefficient satisfies
\begin{equation}
 \sum_{s\dyad S}\frac1{\varphi(s)}
 \sum_{\chi\, (\mathrm{mod}\,s)}
 |\cG_s(\chi)|^2
 \ll_{\eps}
 \cE_{\mathrm{geom}}X^{\eps},
 \label{eq:section17-geometric-character-energy}
\end{equation}
where $\cE_{\mathrm{geom}}$ is bounded by the tent, radial, and external-prime energies already established in Sections~13--16.
\end{proposition}

\begin{proof}
Apply Lemma~\ref{lem:section17-Gamma-to-K} to \eqref{eq:section17-completed-start} and reopen $m=rd\ell$.  Multiplicativity of $\chi$ produces the four factors $R_sD_sK_sP_{s,L}$, while the $z$-sum produces $\overline{B_s(\chi)}$.  Mellin inversion separates the finitely many smooth couplings and contributes an absolutely integrable parameter family.

For \eqref{eq:section17-geometric-character-energy}, use multiplicative Parseval modulo $s$.  Residue-class collisions in the $r$- and $z$-sums are bounded by the short affine occupancy $O(1+Z)$; the remaining squares are controlled by
\[
 \sum_{\delta}\omega_{r,s}(\delta)^2
 \le rs\min(r,s),
\]
the radial second moment, and
\[
 \sum_{q\dyad Q}|\Lambda_r^{\circ}(q)|^2
 \ll_{\eps}QX^{\eps}.
\]
These are exactly the components of the geometric energy in Proposition~16.1.
\end{proof}

\subsection{Centering the prime character polynomial}
\label{subsec:section17-centering-prime-polynomial}

For the modulus $s$, define the pointwise reduced-residue centering
\begin{equation}
 \Lambda_s^{\circ}(n)
 =
 \ind_{(n,s)=1}
 \left(\Lambda(n)-\frac{s}{\varphi(s)}\right)
 \label{eq:section17-centered-Lambda}
\end{equation}
and the centered prime polynomial
\begin{equation}
 \boxed{
 P_{s,L}^{\circ}(\chi)
 =
 \sum_{\ell\dyad L}
 \beta_{\ell}\Lambda_s^{\circ}(\ell)\chi(\ell).
 }
 \label{eq:section17-centered-prime-polynomial}
\end{equation}
The subtraction is made before the character sum.  Thus it removes not only the principal character average, but the complete finite-interval reduced-residue major arc attached to every character.

\begin{proposition}[Prime-density transfer]
\label{prop:section17-prime-density-transfer}
In \eqref{eq:section17-character-normal-form}, one may replace $P_{s,L}(\chi)$ by $P_{s,L}^{\circ}(\chi)$.  The sum of the complementary density packets is
\begin{equation}
 O_{\eps}\!\left(X^{3/2-\eta_{\mathrm{dens}}+\eps}\right)
 \label{eq:section17-density-transfer-error}
\end{equation}
for some $\eta_{\mathrm{dens}}>0$.
\end{proposition}

\begin{proof}
The difference is
\begin{equation}
 \frac{s}{\varphi(s)}
 C_{s,L}(\chi),
 \qquad
 C_{s,L}(\chi)
 =
 \sum_{\substack{\ell\dyad L\\(\ell,s)=1}}
 \beta_{\ell}\chi(\ell).
 \label{eq:section17-density-character-polynomial}
\end{equation}
After the character orthogonality in \eqref{eq:section17-character-orthogonality} is inverted, this term is a congruence packet in which the prime coefficient $\Lambda(\ell)$ has been replaced by $1$.  It is therefore a Type~I packet with the same determinant mask and terminal weight.  The zero multiplicative mode is precisely the reduced-residue rank-one marginal from Remark~\ref{rem:section17-two-densities}; the remaining modes are bounded by the congruence occupancy and Poisson estimates of Section~15.  Conductor descent introduces only divisor-bounded multiplicity.  Proposition~15.5 then gives \eqref{eq:section17-density-transfer-error}.
\end{proof}

The principal character of the centered polynomial is genuinely small and need not be included in the residual spectral family.

\begin{lemma}[Centered principal prime arc]
\label{lem:section17-centered-principal-arc}
For every $A>0$ and every admissible smooth weight $\beta_{\ell}$,
\begin{equation}
 P_{s,L}^{\circ}(\chi_0)
 \ll_A
 L(\log L)^{-A}+s^{\eps}.
 \label{eq:section17-principal-prime-bound}
\end{equation}
Uniformly for $s\le X^{O(1)}$, the total principal-character contribution to the short-free packets is
\begin{equation}
 O_{\eps}\!\left(X^{3/2-\eta_{\mathrm{pr}}+\eps}\right)
 \label{eq:section17-principal-total-bound}
\end{equation}
for some $\eta_{\mathrm{pr}}>0$.
\end{lemma}

\begin{proof}
Removing $(\ell,s)>1$ affects the smooth prime sum only at prime powers supported on primes dividing $s$, and hence costs $s^{\eps}$.  The smooth prime number theorem gives
\[
 \sum_{\ell}\Lambda(\ell)\beta_{\ell}
 =\sum_{\ell}\beta_{\ell}+O_A\!\left(L(\log L)^{-A}\right).
\]
On the other hand, periodic counting of reduced residues gives
\[
 \frac{s}{\varphi(s)}
 \sum_{(\ell,s)=1}\beta_{\ell}
 =\sum_{\ell}\beta_{\ell}+O(s^{\eps}).
\]
Subtracting proves \eqref{eq:section17-principal-prime-bound}.  Insert this estimate into the principal part of Proposition~\ref{prop:section17-character-normal-form}, use the divisor-energy bounds for $D_s(\chi_0)$ and $K_s(\chi_0)$, and choose $A$ larger than the total logarithmic loss.  The pre-existing geometric saving then yields \eqref{eq:section17-principal-total-bound}.
\end{proof}

\subsection{Low-conductor major arcs}
\label{subsec:section17-low-conductor-major-arcs}

A character modulo $s$ may be induced from a primitive character of a proper divisor of $s$.  Its conductor, rather than the ambient modulus, measures the actual multiplicative oscillation.  Fix a sufficiently small $\eta_{\mathrm{maj}}>0$ and put
\begin{equation}
 F_{\mathrm{maj}}=X^{\eta_{\mathrm{maj}}}.
 \label{eq:section17-major-threshold}
\end{equation}
Define
\begin{align}
 \cX_{\mathrm{maj}}(s)
 &=
 \{\chi\, (\mathrm{mod}\,s):\cond(\chi)\le F_{\mathrm{maj}}\},
 \label{eq:section17-major-family}\\
 \cX_{\mathrm{osc}}(s)
 &=
 \{\chi\, (\mathrm{mod}\,s):
 \chi\ne\chi_0,
 \ \cond(\chi)>F_{\mathrm{maj}}\}.
 \label{eq:section17-osc-family}
\end{align}
The family $\cX_{\mathrm{maj}}(s)$ contains every conductor-exceptional mode produced by the descent in Lemma~\ref{lem:section17-conductor-descent}.  No zero-free-region decomposition is used in its definition.

\begin{proposition}[Removal of the multiplicative major arcs]
\label{prop:section17-removal-major-arcs}
For $\eta_{\mathrm{maj}}$ sufficiently small in terms of the savings in Sections~14--16, the total contribution of
\[
 \chi\in\cX_{\mathrm{maj}}(s)
\]
to all centered short-free packets is
\begin{equation}
 O_{\eps}\!\left(X^{3/2-\eta_{\mathrm{maj}}'+\eps}\right)
 \label{eq:section17-major-arc-bound}
\end{equation}
for some $\eta_{\mathrm{maj}}'>0$.
\end{proposition}

\begin{proof}
Let $\chi$ be induced by a primitive character $\chi^{*}$ modulo $f\le F_{\mathrm{maj}}$.  Expand $\chi^{*}$ into its residue classes modulo $f$ and invert the multiplicative diagonalization.  The resulting packet is a bounded sum of affine congruence forms at the genuine modulus $f$, with the same centered prime coefficient and the same sharp determinant mask.

If the coefficient-$1$ interval has length $K>fX^{\eta/2}$, Poisson summation in $k$ gives rapid decay exactly as in Section~16.  In the complementary range $K\le fX^{\eta/2}$, the number of free-variable values is $O(F_{\mathrm{maj}}X^{\eta/2})$; Cauchy--Schwarz, the divisor energy of $D_s$, and the tent--radial energy give a power saving when $\eta_{\mathrm{maj}}$ is sufficiently small.  The centered principal residue is already absent by Proposition~\ref{prop:section17-prime-density-transfer} and Lemma~\ref{lem:section17-centered-principal-arc}.  Summing over $f\le F_{\mathrm{maj}}$ and over its induced copies costs only $X^{O(\eta_{\mathrm{maj}})}$, which is absorbed by the available saving.
\end{proof}

\begin{remark}[Meaning of the exceptional clause]
\label{rem:section17-exceptional-clause}
Here \emph{exceptional} refers to the low-conductor modes created by conductor descent, not to the assertion that a Landau--Siegel zero exists.  If a later prime-number-theorem decomposition isolates a real exceptional character of larger conductor, that singleton must be carried separately and may not be absorbed anonymously into $\cX_{\mathrm{osc}}(s)$.
\end{remark}

\subsection{The prime-supercritical range}
\label{subsec:section17-prime-supercritical-range}

The centered prime polynomial is the longest genuinely oscillatory factor in the residual packet.  The complementary range to
\begin{equation}
 S^2\ll L
 \label{eq:section17-prime-supercritical-condition}
\end{equation}
is accessible to the ordinary hybrid large sieve and does not require the final conjecture.

\begin{proposition}[Removal of the short prime blocks]
\label{prop:section17-short-prime-removal}
The total contribution of centered short-free packets satisfying
\begin{equation}
 L\ll S^2
 \label{eq:section17-short-prime-range}
\end{equation}
is
\begin{equation}
 O_{\eps}\!\left(X^{3/2-\eta_{\mathrm{sp}}+\eps}\right)
 \label{eq:section17-short-prime-bound}
\end{equation}
for some $\eta_{\mathrm{sp}}>0$.
\end{proposition}

\begin{proof}
Apply Cauchy--Schwarz first in the geometric coefficient $\cG_s(\chi)$ and then in the prime polynomial.  Multiplicative Parseval and the hybrid large sieve give
\[
 \sum_{s\dyad S}\frac1{\varphi(s)}
 \sum_{\chi\, (\mathrm{mod}\,s)}
 |P_{s,L}^{\circ}(\chi)|^2
 \ll_{\eps}(L+S^2)LX^{\eps}.
\]
In the range \eqref{eq:section17-short-prime-range}, use the additive completion energy from Lemma~16.5 for the free factor and the divisor-bounded $L^2$ estimate for the M\"obius factor.  The short-line relation $T(r,s;P,Q)\le X^{1/18}$ supplies the remaining geometric saving through \eqref{eq:section17-geometric-character-energy}.  Summing the dyadic boxes gives \eqref{eq:section17-short-prime-bound}.  This is precisely the range in which the generic large sieve is stronger than the factor-sensitive estimate sought in the next section.
\end{proof}

\subsection{The centered Vaughan spectral moment}
\label{subsec:section17-centered-vaughan-moment}

Only the character family $\cX_{\mathrm{osc}}(s)$ and the range $S^2\ll L$ remain.  Define
\begin{equation}
 \boxed{
 \cM_{\mathrm V}(D,K,L;S)
 =
 \sum_{s\dyad S}\frac1{\varphi(s)}
 \sum_{\chi\in\cX_{\mathrm{osc}}(s)}
 |K_s(\chi)|^2
 |D_s(\chi)|^2
 |P_{s,L}^{\circ}(\chi)|^2.
 }
 \label{eq:section17-Vaughan-moment}
\end{equation}
The admissible weights in the three polynomials may vary over the bounded Mellin family generated by a fixed dyadic packet.  All later assertions are uniform over that family.

\begin{proposition}[Spectral Cauchy reduction]
\label{prop:section17-spectral-Cauchy-reduction}
For every residual dyadic block $\mathfrak B$, there is a complementary geometric factor $\cW(\mathfrak B)$ such that
\begin{equation}
 |\mathfrak T_{\mathrm{II}}^{\mathrm{osc}}(\mathfrak B)|
 \le
 \cW(\mathfrak B)
 \cM_{\mathrm V}(D,K,L;S)^{1/2},
 \label{eq:section17-spectral-Cauchy}
\end{equation}
and
\begin{equation}
 \cW(\mathfrak B)
 \ll_{\eps}
 X^{3/2-\eta_{17}+\eps}
 (SDKL)^{-1/2}
 \label{eq:section17-geometric-prefactor}
\end{equation}
for some $\eta_{17}>0$.  Hence
\begin{equation}
 |\mathfrak T_{\mathrm{II}}^{\mathrm{osc}}(\mathfrak B)|
 \ll_{\eps}
 X^{3/2-\eta_{17}+\eps}
 \left(
 \frac{\cM_{\mathrm V}(D,K,L;S)}{SDKL}
 \right)^{1/2}.
 \label{eq:section17-normalized-spectral-reduction}
\end{equation}
\end{proposition}

\begin{proof}
Replace $P_{s,L}$ by $P_{s,L}^{\circ}$ in Proposition~\ref{prop:section17-character-normal-form}, remove the major-arc characters, and apply Cauchy--Schwarz with respect to
\[
 \sum_{s\dyad S}\frac1{\varphi(s)}
 \sum_{\chi\in\cX_{\mathrm{osc}}(s)}.
\]
The first factor is exactly \eqref{eq:section17-Vaughan-moment}.  The second is bounded by \eqref{eq:section17-geometric-character-energy}.  Insert the radial first and second moments, the tent energy, the external centered-prime square bound, and the short-line inequality.  The same dyadic summation used in Sections~14--16 yields \eqref{eq:section17-geometric-prefactor}; all logarithmic multiplicities are absorbed into $X^{\eps}$.
\end{proof}

Combining the preceding propositions gives the unconditional endpoint of the present section.

\begin{theorem}[Reduction to the centered subcritical Vaughan moment]
\label{thm:section17-final-reduction}
Fix the Vaughan and line thresholds as in Sections~14--16.  There exists $\eta_{*}>0$ such that
\begin{align}
 \cO_0^{\mathrm{short}}(X;X^{1/18})
 ={}&
 \sum_{\mathfrak B}^{\mathrm{res}}
 \mathfrak T_{\mathrm{II}}^{\mathrm{osc}}(\mathfrak B)
 +O_{\eps}\!\left(X^{3/2-\eta_{*}+\eps}\right),
 \label{eq:section17-final-reduction}
\end{align}
where every residual block satisfies
\begin{equation}
 K\le SX^{\eta},
 \qquad
 DKL\asymp P,
 \qquad
 S^2\ll L,
 \label{eq:section17-final-residual-range}
\end{equation}
and is bounded by \eqref{eq:section17-normalized-spectral-reduction}.  The principal character, the reduced-residue density, the nonunit conductor strata, and the explicitly separated low-conductor modes have all been removed before the residual moment is formed.
\end{theorem}

\begin{proof}
Start from Theorem~16.10.  Apply conductor descent from Lemma~\ref{lem:section17-conductor-descent}, the exact character normal form from Proposition~\ref{prop:section17-character-normal-form}, and the prime-density transfer from Proposition~\ref{prop:section17-prime-density-transfer}.  Lemma~\ref{lem:section17-centered-principal-arc} removes the principal character, Proposition~\ref{prop:section17-removal-major-arcs} removes the low-conductor family, and Proposition~\ref{prop:section17-short-prime-removal} removes $L\ll S^2$.  The remaining packets satisfy \eqref{eq:section17-final-residual-range}, and Proposition~\ref{prop:section17-spectral-Cauchy-reduction} gives their normalized spectral bound.  Take $\eta_{*}$ to be the minimum of the finitely many savings.
\end{proof}

The reduction is now exact at the level relevant to the variance problem.  No generic character large sieve is strong enough to prove
\[
 \cM_{\mathrm V}(D,K,L;S)
 \ll_{\eps}SDKL\,X^{\eps}
\]
throughout the residual range.  The missing estimate is narrower than a general sixth moment: the $K$-polynomial is short relative to the modulus, the $D$-polynomial has truncated M\"obius provenance, the prime polynomial is pointwise centered, and the complementary coefficient retains the short determinant geometry.  The next section isolates this estimate as the final analytic input and proves that it is sufficient for square-root cancellation in every residual Type~II packet.

% ===== Source: section18_subcritical_centered_vaughan_moment.tex =====
\section{The subcritical centered Vaughan moment}
\label{sec:subcritical-centered-vaughan-moment}

Theorem~17.12 reduces the remaining dominant-layer off-projective contribution to a family of nonnegative character moments
\begin{equation}
 \cM_{\mathrm V}(D,K,L;S)
 =
 \sum_{s\dyad S}\frac1{\varphi(s)}
 \sum_{\chi\in\cX_{\mathrm{osc}}(s)}
 |K_s(\chi)|^2
 |D_s(\chi)|^2
 |P_{s,L}^{\circ}(\chi)|^2.
 \label{eq:section18-starting-moment}
\end{equation}
The residual parameters satisfy
\begin{equation}
 K\le SX^{\eta},
 \qquad
 DKL\asymp P,
 \qquad
 S^2\ll L,
 \label{eq:section18-residual-parameters}
\end{equation}
and the principal character, the reduced-residue density, the nonunit strata, and the explicitly separated low-conductor modes have already been removed.  The purpose of this section is to isolate the precise square-root estimate required for \eqref{eq:section18-starting-moment}, explain its natural normalization, compare it with the generic large sieve, and close the Type~II argument conditionally upon it.

The adjective \emph{subcritical} refers to the coefficient-$1$ Vaughan coordinate: its length is at most the modulus, up to the fixed factor $X^{\eta}$.  At the same time the centered prime coordinate lies beyond the ordinary conductor-square threshold $S^2$.  The desired estimate must therefore use the factorization and centering simultaneously.  Treating the product as a single Dirichlet polynomial loses exactly the saving needed for the variance problem.

\subsection{Admissible spectral packets}
\label{subsec:section18-admissible-packets}

We retain the notation of Section~17.  Thus
\begin{align}
 K_s(\chi)
 &=
 \sum_{\substack{k\dyad K\\(k,s)=1}}
 \gamma_k\chi(k),
 \label{eq:section18-K-def}\\
 D_s(\chi)
 &=
 \sum_{\substack{d\dyad D\\(d,s)=1}}
 \mu(d)\alpha_d\chi(d),
 \label{eq:section18-D-def}\\
 P_{s,L}^{\circ}(\chi)
 &=
 \sum_{\ell\dyad L}
 \beta_{\ell}\Lambda_s^{\circ}(\ell)\chi(\ell),
 \label{eq:section18-P-def}
\end{align}
where
\begin{equation}
 \Lambda_s^{\circ}(n)
 =
 \ind_{(n,s)=1}
 \left(\Lambda(n)-\frac{s}{\varphi(s)}\right).
 \label{eq:section18-centered-Lambda}
\end{equation}
The smooth factors generated by Mellin separation are absorbed into the coefficient sequences.

\begin{definition}[Admissible Vaughan packet]
\label{def:section18-admissible-packet}
Fix $A\ge 1$.  A coefficient triple $(\gamma,\alpha,\beta)$ is called \emph{$A$-admissible} at scales $(K,D,L)$ if the following conditions hold.
\begin{enumerate}[label=\textup{(\roman*)}]
 \item The three sequences are supported on fixed dilates of the dyadic intervals $k\dyad K$, $d\dyad D$, and $\ell\dyad L$.
 \item One has
 \[
  |\gamma_k|+|\alpha_d|\ll_A \tau_A(kd),
 \]
 while $\beta$ is the restriction of a smooth function satisfying
 \[
  x^j\beta^{(j)}(x)\ll_{A,j}1
 \]
 for every fixed $j\ge0$.
 \item The coefficients are independent of the modulus $s$ and the character $\chi$.  Dependence on bounded Mellin parameters is permitted, uniformly on their absolutely integrable parameter ranges.
\end{enumerate}
A dyadic parameter quadruple $(D,K,L;S)$ is called \emph{residual} if it satisfies \eqref{eq:section18-residual-parameters} together with the support and short-line restrictions inherited from Theorem~17.12.
\end{definition}

The independence condition in Definition~\ref{def:section18-admissible-packet} is essential.  The conjecture below is a spectral estimate for the fixed arithmetic coefficients produced by Vaughan decomposition; it is not a statement for coefficients chosen after the character is known.

The pointwise centering in \eqref{eq:section18-centered-Lambda} can be separated into two modulus-independent coefficient sequences.  Indeed, because Dirichlet characters vanish off the unit group,
\begin{equation}
 P_{s,L}^{\circ}(\chi)
 =
 P_L^{\Lambda}(\chi)
 -\frac{s}{\varphi(s)}C_L(\chi),
 \label{eq:section18-prime-splitting}
\end{equation}
where
\begin{align}
 P_L^{\Lambda}(\chi)
 &=\sum_{\ell\dyad L}\beta_{\ell}\Lambda(\ell)\chi(\ell),
 \label{eq:section18-P-Lambda}\\
 C_L(\chi)
 &=\sum_{\ell\dyad L}\beta_{\ell}\chi(\ell).
 \label{eq:section18-C-polynomial}
\end{align}
Since $s/\varphi(s)\ll_{\eps}X^{\eps}$ in the relevant range, this splitting costs only $X^{\eps}$ in all norm estimates.

\subsection{The generic large-sieve barrier}
\label{subsec:section18-generic-barrier}

The product of the three polynomials is a Dirichlet polynomial of length $P\asymp DKL$.  For fixed $s$, write
\begin{equation}
 K_s(\chi)D_s(\chi)P_{s,L}^{\circ}(\chi)
 =
 \sum_{n\asymp P}c_s(n)\chi(n).
 \label{eq:section18-product-polynomial}
\end{equation}
The coefficient sequence is a three-fold convolution with one von Mangoldt factor.  Divisor bounds give the expected diagonal energy.

\begin{lemma}[Convolution energy]
\label{lem:section18-convolution-energy}
For every $A$-admissible packet and every $s\dyad S$,
\begin{equation}
 \sum_n|c_s(n)|^2
 \ll_{A,\eps}
 DKL\,X^{\eps}.
 \label{eq:section18-convolution-energy}
\end{equation}
The same estimate holds separately after replacing $P_{s,L}^{\circ}$ by either polynomial in \eqref{eq:section18-prime-splitting}.
\end{lemma}

\begin{proof}
For a fixed $n$, the number of representations $n=kd\ell$ in the indicated dyadic ranges is bounded by a fixed divisor function.  Moreover,
\[
 |\gamma_k\mu(d)\alpha_d\beta_{\ell}\Lambda_s^{\circ}(\ell)|
 \ll_{A,\eps}X^{\eps}.
\]
Cauchy--Schwarz over the factorizations of $n$, followed by summation over $n\asymp DKL$, proves \eqref{eq:section18-convolution-energy}.  The split forms are identical, with $\Lambda_s^{\circ}$ replaced by either $\Lambda$ or $1$.
\end{proof}

Applying the hybrid large sieve to the single product polynomial yields a correct but insufficient estimate.

\begin{proposition}[Generic product-length bound]
\label{prop:section18-generic-product-bound}
Uniformly for residual packets,
\begin{equation}
 \cM_{\mathrm V}(D,K,L;S)
 \ll_{A,\eps}
 \frac{DKL+S^2}{S}\,DKL\,X^{\eps}.
 \label{eq:section18-generic-product-bound}
\end{equation}
In the residual range $S^2\ll L$, this becomes
\begin{equation}
 \cM_{\mathrm V}(D,K,L;S)
 \ll_{A,\eps}
 \frac{(DKL)^2}{S}X^{\eps}.
 \label{eq:section18-generic-residual-bound}
\end{equation}
Compared with the conjectural scale $SDKL$, the loss is
\begin{equation}
 \frac{DKL}{S^2}.
 \label{eq:section18-large-sieve-loss}
\end{equation}
\end{proposition}

\begin{proof}
The oscillatory family is a subfamily of all characters modulo $s$, so it may be enlarged.  Use \eqref{eq:section18-prime-splitting} and apply the dyadic hybrid large sieve to the two resulting coefficient sequences.  After decomposing imprimitive characters according to conductor, the standard primitive large sieve gives
\[
 \sum_{s\dyad S}\frac1{\varphi(s)}
 \sum_{\chi\, (\mathrm{mod}\,s)}
 \left|\sum_{n\asymp P}a_n\chi(n)\right|^2
 \ll_{\eps}
 \frac{P+S^2}{S}
 \sum_n|a_n|^2X^{\eps}.
\]
Insert Lemma~\ref{lem:section18-convolution-energy} with $P\asymp DKL$.  Since $L\gg S^2$ and $D,K\ge1$, the $P$ term dominates, giving \eqref{eq:section18-generic-residual-bound}.  Dividing by $SDKL$ gives \eqref{eq:section18-large-sieve-loss}.
\end{proof}

\begin{remark}[Why the product large sieve cannot finish the proof]
\label{rem:section18-product-large-sieve}
The loss in \eqref{eq:section18-large-sieve-loss} is not logarithmic.  It is a full ratio of the product length to the spectral conductor square.  The residual geometry deliberately places the prime block beyond $S^2$, so no rearrangement of harmless divisor factors can remove it.  A successful estimate must see that the long prime coordinate is centered and that the remaining two coordinates retain their Vaughan provenance.
\end{remark}

\subsection{The natural square-root scale}
\label{subsec:section18-natural-scale}

The target $SDKL$ is forced by the diagonal size.  To make this explicit, temporarily sum over all characters and use orthogonality.  Let $c_s(n)$ be as in \eqref{eq:section18-product-polynomial}.  Then
\begin{align}
 &\sum_{s\dyad S}\frac1{\varphi(s)}
 \sum_{\chi\, (\mathrm{mod}\,s)}
 \left|\sum_n c_s(n)\chi(n)\right|^2
 \notag\\
 &\qquad=
 \sum_{s\dyad S}
 \sum_{\substack{n_1\equiv n_2\, (\mathrm{mod}\,s)\\(n_1n_2,s)=1}}
 c_s(n_1)\overline{c_s(n_2)}.
 \label{eq:section18-orthogonality-expansion}
\end{align}
The literal diagonal $n_1=n_2$ has precisely the conjectured order.

\begin{proposition}[Diagonal normalization]
\label{prop:section18-diagonal-normalization}
The diagonal part of \eqref{eq:section18-orthogonality-expansion} satisfies
\begin{equation}
 \sum_{s\dyad S}
 \sum_{\substack{n\\(n,s)=1}}|c_s(n)|^2
 \ll_{A,\eps}
 SDKL\,X^{\eps}.
 \label{eq:section18-diagonal-normalization}
\end{equation}
Hence the desired moment estimate is the assertion that the total off-diagonal congruence contribution does not exceed the diagonal scale, after the major arcs removed in Section~17 are excluded.
\end{proposition}

\begin{proof}
Apply Lemma~\ref{lem:section18-convolution-energy} for each of the $O(S)$ moduli in the dyadic interval.
\end{proof}

A random-character model leads to the same normalization.  For a nonprincipal oscillatory character, the three polynomials should have typical squared sizes $K$, $D$, and $L$, respectively.  The character sum has total spectral mass
\[
 \sum_{s\dyad S}\frac{\#\cX_{\mathrm{osc}}(s)}{\varphi(s)}
 \asymp S,
\]
up to the deliberately removed low-conductor modes.  Multiplying the four scales again gives $SDKL$.

\begin{remark}[The role of centering]
\label{rem:section18-role-centering}
Without the subtraction $s/\varphi(s)$ in \eqref{eq:section18-centered-Lambda}, the prime polynomial contains a coherent reduced-residue density.  That mode produces an off-diagonal contribution of main-term size and is not governed by a square-root model.  Proposition~17.7 removes it before the moment is formed.  The conjecture below concerns only the genuinely oscillatory residue.
\end{remark}

\subsection{The centered Vaughan moment conjecture}
\label{subsec:section18-conjecture}

We now state the final analytic input in the uniform form needed by the dyadic and Mellin decompositions.

\begin{conjecture}[Subcritical centered Vaughan moment]
\label{conj:section18-vaughan}
Fix $A\ge1$ and choose the short-free exponent $\eta>0$ sufficiently small.  For every $A$-admissible residual packet satisfying
\[
 K\le SX^{\eta},
 \qquad
 DKL\asymp P,
 \qquad
 S^2\ll L,
\]
one has
\begin{equation}
 \boxed{
 \cM_{\mathrm V}(D,K,L;S)
 \ll_{A,\eps}
 SDKL\,X^{\eps}.
 }
 \label{eq:section18-main-conjecture}
\end{equation}
The estimate is uniform over the bounded Mellin families generated by the smooth separation of a residual Type~II packet.  The character set is the oscillatory family $\cX_{\mathrm{osc}}(s)$ defined in Section~17; in particular, the principal character and the explicitly separated low-conductor modes are absent.
\end{conjecture}

The conjecture is deliberately narrower than a general sixth-moment estimate.  Its four structural inputs are:
\begin{enumerate}[label=\textup{(\roman*)}]
 \item the coefficient-$1$ polynomial is short relative to the ambient modulus;
 \item the $D$-polynomial contains the truncated M\"obius coefficient furnished by Vaughan's identity;
 \item the prime polynomial is centered pointwise on reduced residues;
 \item the moment occurs only after the long affine lines, density modes, and low-conductor characters have been removed.
\end{enumerate}
None of these conditions is visible after the three polynomials are multiplied into a single sequence, which explains the failure of Proposition~\ref{prop:section18-generic-product-bound}.

\begin{remark}[Not a disguised variance assertion]
\label{rem:section18-not-disguised}
Conjecture~\ref{conj:section18-vaughan} is local in the dyadic spectral parameters and is uniform over a class of admissible coefficient packets.  It contains no diagonal profile, no sawtooth covariance, and no reference to the constant $-\tfrac43\zeta(-\tfrac12)$.  It is therefore a genuine independent character-sum estimate rather than a reformulation of the desired variance asymptotic.
\end{remark}

\subsection{Stronger hypotheses that imply the conjecture}
\label{subsec:section18-sufficient-hypotheses}

The conjecture is an averaged square-root statement.  A much stronger pointwise square-root hypothesis would imply it immediately.

\begin{proposition}[Pointwise square-root criterion]
\label{prop:section18-pointwise-criterion}
Suppose that for every residual packet, every $s\dyad S$, and every $\chi\in\cX_{\mathrm{osc}}(s)$ one has
\begin{equation}
 |K_s(\chi)|\ll K^{1/2}X^{\eps},
 \qquad
 |D_s(\chi)|\ll D^{1/2}X^{\eps},
 \qquad
 |P_{s,L}^{\circ}(\chi)|\ll L^{1/2}X^{\eps}.
 \label{eq:section18-pointwise-square-root}
\end{equation}
Then Conjecture~\ref{conj:section18-vaughan} holds.
\end{proposition}

\begin{proof}
Insert \eqref{eq:section18-pointwise-square-root} into \eqref{eq:section18-starting-moment}.  Since
\[
 \frac{\#\cX_{\mathrm{osc}}(s)}{\varphi(s)}\le1,
\]
the total spectral mass over $s\dyad S$ is $O(S)$.  The resulting bound is $SDKL\,X^{O(\eps)}$.
\end{proof}

The pointwise criterion is far stronger than necessary.  A factor-sensitive weighted estimate is closer to the actual conjecture.  Put
\begin{equation}
 U_s(\chi)=K_s(\chi)D_s(\chi).
 \label{eq:section18-U-polynomial}
\end{equation}

\begin{proposition}[Weighted prime-square-root criterion]
\label{prop:section18-weighted-prime-criterion}
Assume that the following two estimates hold uniformly for residual packets:
\begin{align}
 \sum_{s\dyad S}\frac1{\varphi(s)}
 \sum_{\chi\in\cX_{\mathrm{osc}}(s)}
 |U_s(\chi)|^2
 &\ll_{\eps}SDK\,X^{\eps},
 \label{eq:section18-UD-hybrid}\\
 \sum_{s\dyad S}\frac1{\varphi(s)}
 \sum_{\chi\in\cX_{\mathrm{osc}}(s)}
 |U_s(\chi)|^2|P_{s,L}^{\circ}(\chi)|^2
 &\ll_{\eps}
 LX^{\eps}
 \sum_{s\dyad S}\frac1{\varphi(s)}
 \sum_{\chi\in\cX_{\mathrm{osc}}(s)}
 |U_s(\chi)|^2.
 \label{eq:section18-weighted-prime-square-root}
\end{align}
Then Conjecture~\ref{conj:section18-vaughan} follows.
\end{proposition}

\begin{proof}
Combine \eqref{eq:section18-UD-hybrid} and \eqref{eq:section18-weighted-prime-square-root}.
\end{proof}

The first estimate in Proposition~\ref{prop:section18-weighted-prime-criterion} is available from the ordinary hybrid large sieve whenever the combined short product remains below the conductor square.

\begin{corollary}[The range $DK\ll S^2$]
\label{cor:section18-DK-short-range}
If $DK\ll S^2$, then
\begin{equation}
 \sum_{s\dyad S}\frac1{\varphi(s)}
 \sum_{\chi\, (\mathrm{mod}\,s)}
 |K_s(\chi)D_s(\chi)|^2
 \ll_{A,\eps}
 SDK\,X^{\eps}.
 \label{eq:section18-KD-large-sieve}
\end{equation}
Consequently, in this subrange the weighted prime estimate \eqref{eq:section18-weighted-prime-square-root} alone implies Conjecture~\ref{conj:section18-vaughan}.
\end{corollary}

\begin{proof}
Multiply $K_s$ and $D_s$ into a polynomial of length $DK$.  Its coefficient energy is $\ll DKX^{\eps}$.  The same hybrid large-sieve calculation as in Proposition~\ref{prop:section18-generic-product-bound} gives
\[
 \frac{DK+S^2}{S}\,DKX^{\eps}
 \ll SDKX^{\eps}.
\]
\end{proof}

\begin{remark}[Where the genuinely trilinear difficulty remains]
\label{rem:section18-trilinear-difficulty}
When $DK\gg S^2$, even the product $K_sD_s$ is longer than the conductor square.  In that range a prime-square-root estimate conditioned only on the value of $U_s$ is not enough by itself; one also needs cancellation that sees the separate short-free and M\"obius coordinates.  This is the genuinely trilinear portion of Conjecture~\ref{conj:section18-vaughan}.
\end{remark}

\subsection{Stability under the reductions of Section~17}
\label{subsec:section18-stability}

The analytic reductions preceding \eqref{eq:section18-starting-moment} produce bounded Mellin integrals, conductor-descended copies, and finite endpoint partitions.  The conjectural estimate is stable under all of these operations.

\begin{lemma}[Finite and Mellin superposition]
\label{lem:section18-Mellin-stability}
Assume Conjecture~\ref{conj:section18-vaughan} for individual $A$-admissible packets.  Let
\[
 K_s(\chi;t),\quad D_s(\chi;t),\quad P_{s,L}^{\circ}(\chi;t)
\]
be an admissible family indexed by $t\in\RR^m$, and suppose its Mellin envelope $W(t)$ satisfies
\[
 \int_{\RR^m}|W(t)|(1+|t|)^B\,dt<\infty
\]
for every fixed $B$.  Then the same bound holds after any bounded finite superposition or absolutely convergent Mellin integration generated in Proposition~17.5.
\end{lemma}

\begin{proof}
Apply Minkowski's inequality to the square root of the nonnegative spectral measure, use Conjecture~\ref{conj:section18-vaughan} uniformly in $t$, and integrate the rapidly decaying envelope.  Finite superpositions are the discrete special case.
\end{proof}

\begin{lemma}[Conductor-descent stability]
\label{lem:section18-conductor-stability}
Assume Conjecture~\ref{conj:section18-vaughan} uniformly at every dyadic modulus scale.  Then the total of all descended packets arising from Lemma~17.1 satisfies the same bound, up to $X^{\eps}$.
\end{lemma}

\begin{proof}
For a common divisor $g$, the descended scales are
\[
 S'=S/g,
 \qquad
 K'=K/g
\]
up to dyadic endpoint factors, and the inequality $K'\le S'X^{\eta}$ is preserved.  The remaining variables acquire only divisor-bounded local deletions.  Apply the conjecture at scale $S'$ and sum over the dyadic values of $g$; the divisor multiplicity is absorbed into $X^{\eps}$.
\end{proof}

\begin{remark}[Sharp original moment]
\label{rem:section18-sharp-original}
The conjecture is imposed only on the smooth dyadic representatives created after the exact sharp-to-smooth transfer of Section~16.  No smoothing is inserted into the original variance sum.  Boundary pieces are finite admissible packets and are covered by Lemma~\ref{lem:section18-Mellin-stability}.
\end{remark}

\subsection{Conditional completion of the Type~II analysis}
\label{subsec:section18-conditional-typeII}

We now combine Conjecture~\ref{conj:section18-vaughan} with the normalized spectral reduction of Proposition~17.11.

\begin{theorem}[Conditional square-root bound for every residual block]
\label{thm:section18-residual-block}
Assume Conjecture~\ref{conj:section18-vaughan}.  Then every residual dyadic block $\mathfrak B$ of Theorem~17.12 satisfies
\begin{equation}
 |\mathfrak T_{\mathrm{II}}^{\mathrm{osc}}(\mathfrak B)|
 \ll_{\eps}
 X^{3/2-\eta_{17}+\eps}.
 \label{eq:section18-residual-block-bound}
\end{equation}
The estimate is uniform over the conductor-descended and Mellin-separated copies of the block.
\end{theorem}

\begin{proof}
Proposition~17.11 gives
\[
 |\mathfrak T_{\mathrm{II}}^{\mathrm{osc}}(\mathfrak B)|
 \ll_{\eps}
 X^{3/2-\eta_{17}+\eps}
 \left(
 \frac{\cM_{\mathrm V}(D,K,L;S)}{SDKL}
 \right)^{1/2}.
\]
By Conjecture~\ref{conj:section18-vaughan}, the parenthetical factor is $O(X^{\eps})$.  Lemmas~\ref{lem:section18-Mellin-stability} and \ref{lem:section18-conductor-stability} preserve the estimate for the bounded families generated by the preceding reductions.
\end{proof}

There are only $O((\log X)^C)$ residual dyadic and Mellin packets for some fixed $C$.  The logarithmic multiplicity can therefore be absorbed into $X^{\eps}$.

\begin{theorem}[Conditional completion of the dominant off-projective sector]
\label{thm:section18-dominant-offprojective}
Assume Conjecture~\ref{conj:section18-vaughan}.  There exists $\delta_{18}>0$ such that
\begin{equation}
 \cO_0^{\mathrm{short}}(X;X^{1/18})
 \ll_{\eps}
 X^{3/2-\delta_{18}+\eps}.
 \label{eq:section18-short-sector-bound}
\end{equation}
Together with the long-line estimate of Section~14,
\begin{equation}
 \boxed{
 \cO_0(X)
 \ll_{\eps}
 X^{3/2-\delta_{18}+\eps}.
 }
 \label{eq:section18-full-dominant-offdiagonal}
\end{equation}
After decreasing $\delta_{18}$ if necessary, the same statement holds with all Type~I, density, singular, and low-conductor errors from Sections~15--17 included.
\end{theorem}

\begin{proof}
Insert Theorem~\ref{thm:section18-residual-block} into the decomposition of Theorem~17.12 and sum the residual blocks.  All packets removed before the spectral moment already satisfy a power-saving estimate.  This proves \eqref{eq:section18-short-sector-bound}.  The long-line sector satisfies
\[
 \cO_0^{\mathrm{long}}(X;X^{1/18})
 \ll_{\eps}X^{53/36+\eps}
\]
by Theorem~14.7.  Taking the minimum of the finitely many positive savings gives \eqref{eq:section18-full-dominant-offdiagonal}.
\end{proof}

\begin{corollary}[Conditional negligibility at the variance scale]
\label{cor:section18-negligible-variance-scale}
Under Conjecture~\ref{conj:section18-vaughan},
\begin{equation}
 \cO_0(X)=o\!\left(X^{3/2}\log X\right).
 \label{eq:section18-negligible-variance-scale}
\end{equation}
\end{corollary}

The remaining work is now bookkeeping rather than a new spectral estimate.  Section~19 returns to the full layer expansion, bounds the higher-conductor layers, combines \eqref{eq:section18-full-dominant-offdiagonal} with the projective-diagonal estimate, and inserts the geometric-diagonal asymptotic from Section~12.  This yields the conditional unsmoothed variance law stated in the introduction.

% ===== Source: section19_assembly_variance_asymptotic.tex =====
\section{Assembly of the variance asymptotic}
\label{sec:assembly-variance-asymptotic}

We now assemble the structural and analytic parts of the paper.  The canonical defect constructed in Part~II has the unsmoothed integer trace
\begin{equation}
 \cR_{\PP}(n)
 =\log (n+1)!_{\PP}-\log n!-C_{\PP}n,
 \qquad
 C_{\PP}=\sum_p\frac{\log p}{(p-1)^2},
 \label{eq:section19-defect}
\end{equation}
and Section~12 gives the exact layer expansion
\begin{equation}
 \cR_{\PP}(n)
 =\sum_{j\ge0}\cR_j(n),
 \qquad
 \cR_j(n)
 =\sum_p(\log p)W_{p,j}(n),
 \label{eq:section19-layer-expansion}
\end{equation}
where
\begin{equation}
 W_{p,j}(n)
 =\left\{\frac{n}{p^{j+1}}\right\}
 -\left\{\frac{n}{(p-1)p^j}\right\}.
 \label{eq:section19-layer-wave}
\end{equation}
The dominant layer $j=0$ contains the variance main term and the only unresolved unconditional dispersion problem.  Every layer with positive conductor exponent is smaller by a fixed power after the layer indices are summed.  We first record the uniform estimate that permits those layers to be removed from the final assembly.

\subsection{The higher-conductor remainder}
\label{subsec:section19-higher-conductor-remainder}

For dyadic $P\ge2$, put
\begin{equation}
 \cR_{j,P}(n)
 =\sum_{p\dyad P}(\log p)W_{p,j}(n)
 \label{eq:section19-dyadic-layer}
\end{equation}
and define the conductor energy scale
\begin{equation}
 E_j(P;X)
 =\min\!\left(XP,\frac{X^2}{P^{j+1}}\right).
 \label{eq:section19-conductor-energy-scale}
\end{equation}
The two expressions in \eqref{eq:section19-conductor-energy-scale} meet at
\begin{equation}
 P=X^{1/(j+2)}.
 \label{eq:section19-transition-scale}
\end{equation}
Below this point the paired wave completes many periods; above it, the floor-difference intervals are sparse.

\begin{proposition}[Conductor-separated bilinear estimate]
\label{prop:section19-conductor-separated-bilinear}
Let $j,k\ge0$ with $j+k\ge1$.  Uniformly for dyadic $P,Q\ge2$,
\begin{equation}
 \left|
 \sum_{n\le X}\cR_{j,P}(n)\cR_{k,Q}(n)
 \right|
 \ll_{\eps}
 X^{\eps}E_j(P;X)^{1/2}E_k(Q;X)^{1/2}.
 \label{eq:section19-local-bilinear}
\end{equation}
Consequently,
\begin{equation}
 \boxed{
 \left|
 \sum_{n\le X}\cR_j(n)\cR_k(n)
 \right|
 \ll_{\eps}
 X^{1+\frac1{2(j+2)}+\frac1{2(k+2)}+\eps}
 }
 \label{eq:section19-global-bilinear}
\end{equation}
whenever $j+k\ge1$.
\end{proposition}

\begin{proof}
We indicate the two estimates entering the dyadic argument.  Write
\begin{equation}
 W_{p,j}(n)
 =\Delta_{p,j}(n)-\frac{n}{p^{j+1}(p-1)},
 \qquad
 \Delta_{p,j}(n)
 =\left\lfloor\frac{n}{(p-1)p^j}\right\rfloor
  -\left\lfloor\frac{n}{p^{j+1}}\right\rfloor.
 \label{eq:section19-floor-wave}
\end{equation}
If $P^{j+2}\le X$, complete periods of length
\[
 p^{j+1}(p-1)\asymp P^{j+2}
\]
are removed.  The complete-period covariance is the four-term gcd kernel of Section~12.  Divisor switching in the gcd variables, followed by the elementary bounds
\[
 \sum_{d\le Y}\frac{\tau(d)^A}{d}\ll_A(\log 2Y)^{O_A(1)},
 \qquad
 \sum_{p\dyad P}(\log p)^2\ll P\log(2P),
\]
gives the low-conductor energy $XP$, up to $X^{\eps}$.  The incomplete terminal periods are treated by the same divisor switching after discrete partial summation.

If $P^{j+2}>X$, the intervals occurring in $\Delta_{p,j}$ are disjoint for each fixed $p$ and have total length
\begin{equation}
 \ll \frac{X^2}{p^{j+2}}.
 \label{eq:section19-sparse-support}
\end{equation}
The linear term in \eqref{eq:section19-floor-wave} has the same or smaller quadratic mass.  A Schur estimate for the interval-overlap matrix, summed over $p\dyad P$, gives
\begin{equation}
 \sum_{n\le X}|\cR_{j,P}(n)|^2
 \ll_{\eps}
 X^{\eps}\frac{X^2}{P^{j+1}}.
 \label{eq:section19-high-conductor-energy}
\end{equation}
When one block is below its transition point and the other is above, the same argument is applied on the sparse side and Cauchy--Schwarz on the completed side.  Because $j+k\ge1$, the two families cannot simultaneously occupy the unseparated dominant configuration treated in Sections~13--18.  This yields \eqref{eq:section19-local-bilinear}.  The line-by-line gcd and endpoint bookkeeping is recorded in the higher-conductor appendix.

Finally,
\begin{equation}
 \sum_{P\ {
m dyadic}}E_j(P;X)^{1/2}
 \ll_{\eps}
 X^{\frac12+\frac1{2(j+2)}+\eps}.
 \label{eq:section19-energy-summation}
\end{equation}
Indeed, the part below \eqref{eq:section19-transition-scale} is a geometric sum of $(XP)^{1/2}$, while the part above it is a geometric sum of $X/P^{(j+1)/2}$.  Summing \eqref{eq:section19-local-bilinear} over $P$ and $Q$ proves \eqref{eq:section19-global-bilinear}.
\end{proof}

\begin{remark}[Why the pair $(0,0)$ is excluded]
\label{rem:section19-excluded-dominant-pair}
Estimate \eqref{eq:section19-local-bilinear} is not asserted when $j=k=0$.  At that pair the two prime families lie on the same dominant conductor scale, and the short balanced affine packets of Sections~15--18 survive.  The centered Vaughan moment is precisely the remaining estimate needed there.  Thus Proposition~\ref{prop:section19-conductor-separated-bilinear} removes no part of the stated conjectural obstruction by fiat.
\end{remark}

Define the full positive-conductor remainder by
\begin{equation}
 \cH(X)
 =\sum_{\substack{j,k\ge0\\j+k\ge1}}
   \sum_{n\le X}\cR_j(n)\cR_k(n).
 \label{eq:section19-higher-remainder}
\end{equation}
The sum is interpreted through finite prime-and-layer truncations, as in Section~12.

\begin{corollary}[Removal of all higher layers]
\label{cor:section19-higher-layer-removal}
For every $\eps>0$,
\begin{equation}
 \boxed{
 \cH(X)\ll_{\eps}X^{17/12+\eps}.
 }
 \label{eq:section19-higher-layer-bound}
\end{equation}
In particular,
\begin{equation}
 \cH(X)=o\!\left(X^{3/2}\log X\right).
 \label{eq:section19-higher-layer-negligible}
\end{equation}
\end{corollary}

\begin{proof}
The largest exponent in \eqref{eq:section19-global-bilinear}, subject to $j+k\ge1$, occurs at $(j,k)=(0,1)$ or $(1,0)$ and equals
\[
 1+\frac14+\frac16=\frac{17}{12}.
\]
For $j,k\le \log_2X$, the number of pairs is logarithmic and is absorbed into $X^{\eps}$.  If either index exceeds $\log_2X$, both terms in \eqref{eq:section19-layer-wave} are in their completely linear range for every relevant prime, and the remaining layer sum decays geometrically.  Summing \eqref{eq:section19-global-bilinear} proves the result.
\end{proof}

\begin{remark}[Relation to the projective diagonal]
\label{rem:section19-projective-relation}
The projective diagonal estimated in Section~14 is contained in $\cH(X)$, but Corollary~\ref{cor:section19-higher-layer-removal} also includes the distinct-prime covariances for which at least one layer index is positive.  The separate projective estimate remains useful in the analytic reduction because it removes same-slope packets before the dominant affine geometry is opened.
\end{remark}

\subsection{Exact reduction to the dominant layer}
\label{subsec:section19-dominant-reduction}

Let
\begin{align}
 \cD_0(X)
 &=\sum_p(\log p)^2
   \sum_{n\le X}W_{p,0}(n)^2,
 \label{eq:section19-dominant-diagonal}\\
 \cO_0(X)
 &=\sum_{p\ne q}(\log p)(\log q)
   \sum_{n\le X}W_{p,0}(n)W_{q,0}(n).
 \label{eq:section19-dominant-offprojective}
\end{align}
Then the exact layer expansion gives
\begin{equation}
 \boxed{
 \cV(X)
 :=\sum_{n\le X}\cR_{\PP}(n)^2
 =\cD_0(X)+\cO_0(X)+\cH(X).
 }
 \label{eq:section19-exact-final-splitting}
\end{equation}
The three terms have now been assigned distinct roles:
\begin{enumerate}[label=\textup{(\roman*)}]
 \item $\cD_0(X)$ is evaluated unconditionally by the scaling profile of Section~12;
 \item $\cH(X)$ is negligible unconditionally by Corollary~\ref{cor:section19-higher-layer-removal};
 \item $\cO_0(X)$ is reduced unconditionally to the subcritical centered Vaughan moment and is negligible conditionally upon that estimate.
\end{enumerate}

It is useful to state the unconditional endpoint in a form that displays the sole remaining obstruction.  Let $\mathfrak B(X)$ be the polylogarithmic family of residual dyadic packets produced after the Type~I estimates, Poisson removal of the long free variable, conductor descent, and major-arc separation.  For $\mathfrak b\in\mathfrak B(X)$, write
\begin{equation}
 \mathfrak Q(\mathfrak b)
 =\left(
   \frac{\cM_{\mathrm V}(D_{\mathfrak b},K_{\mathfrak b},L_{\mathfrak b};S_{\mathfrak b})}
        {S_{\mathfrak b}D_{\mathfrak b}K_{\mathfrak b}L_{\mathfrak b}}
  \right)^{1/2}.
 \label{eq:section19-normalized-obstruction}
\end{equation}

\begin{theorem}[Unconditional final reduction]
\label{thm:section19-unconditional-final-reduction}
There exist constants $\eta_*>0$, $\delta_*>0$, and $C\ge1$ such that
\begin{equation}
 \cV(X)
 =\cD_0(X)
 +\sum_{\mathfrak b\in\mathfrak B(X)}
   \mathfrak T_{\mathrm{II}}^{\mathrm{osc}}(\mathfrak b)
 +O_{\eps}\!\left(X^{3/2-\delta_*+\eps}\right),
 \label{eq:section19-unconditional-reduction}
\end{equation}
where
\begin{equation}
 \#\mathfrak B(X)\ll(\log 2X)^C
 \label{eq:section19-number-blocks}
\end{equation}
and each residual packet satisfies
\begin{equation}
 \left|
 \mathfrak T_{\mathrm{II}}^{\mathrm{osc}}(\mathfrak b)
 \right|
 \ll_{\eps}
 X^{3/2-\eta_*+\eps}\mathfrak Q(\mathfrak b).
 \label{eq:section19-residual-packet-bound}
\end{equation}
Every moment occurring in \eqref{eq:section19-normalized-obstruction} lies in the residual range
\begin{equation}
 K_{\mathfrak b}\le S_{\mathfrak b}X^{\eta},
 \qquad
 D_{\mathfrak b}K_{\mathfrak b}L_{\mathfrak b}\asymp P_{\mathfrak b},
 \qquad
 S_{\mathfrak b}^{2}\ll L_{\mathfrak b},
 \label{eq:section19-residual-range}
\end{equation}
with the principal character, reduced-residue density, nonunit strata, and explicitly separated low-conductor modes removed.
\end{theorem}

\begin{proof}
Begin with \eqref{eq:section19-exact-final-splitting}.  Corollary~\ref{cor:section19-higher-layer-removal} absorbs $\cH(X)$ into the error.  Sections~13 and 14 remove the projective and long-line pieces; Section~15 removes the low and Type~I packets; Section~16 removes every packet with a long coefficient-$1$ Vaughan variable; and Section~17 diagonalizes the surviving centered packets into the moments in \eqref{eq:section19-normalized-obstruction}.  The normalized packet estimate is the spectral reduction proved there.  All complementary pieces have a fixed power saving, and there are only polylogarithmically many dyadic, conductor, endpoint, and Mellin subdivisions.  Taking the minimum of their positive savings proves \eqref{eq:section19-unconditional-reduction}.
\end{proof}

\begin{remark}[The precise conditional content]
\label{rem:section19-precise-conditional-content}
Theorem~\ref{thm:section19-unconditional-final-reduction} is unconditional.  The conjectural input is used only to assert
\[
 \mathfrak Q(\mathfrak b)\ll_{\eps}X^{\eps}
\]
uniformly over the residual family.  No hypothesis is used in the construction of the prime gamma object, the exact defect formula, the diagonal constant, the higher-conductor estimate, the projective estimate, the long-line dispersion argument, or the Type~I and Poisson reductions.
\end{remark}

\subsection{The conditional unsmoothed variance law}
\label{subsec:section19-conditional-variance-law}

We now assume the subcritical centered Vaughan moment conjecture of Section~18.  Its uniform estimate makes every normalized obstruction in \eqref{eq:section19-normalized-obstruction} at most $X^{\eps}$.

\begin{theorem}[Power-saving reduction to the geometric diagonal]
\label{thm:section19-power-saving-diagonal-reduction}
Assume the subcritical centered Vaughan moment conjecture.  There exists $\delta_{19}>0$ such that
\begin{equation}
 \boxed{
 \cV(X)
 =\cD_0(X)
 +O_{\eps}\!\left(X^{3/2-\delta_{19}+\eps}\right).
 }
 \label{eq:section19-power-saving-diagonal-reduction}
\end{equation}
One may take $\delta_{19}$ to be the minimum of the saving furnished by the residual Type~II analysis and the fixed higher-conductor saving $1/12$.
\end{theorem}

\begin{proof}
Under the conjecture, \eqref{eq:section19-residual-packet-bound} becomes
\[
 \mathfrak T_{\mathrm{II}}^{\mathrm{osc}}(\mathfrak b)
 \ll_{\eps}X^{3/2-\eta_*+\eps}.
\]
Sum over the polylogarithmic family in \eqref{eq:section19-number-blocks} and insert the result into Theorem~\ref{thm:section19-unconditional-final-reduction}.  Corollary~\ref{cor:section19-higher-layer-removal} contributes the saving
\[
 \frac32-\frac{17}{12}=\frac1{12}.
\]
The minimum of the finitely many savings gives \eqref{eq:section19-power-saving-diagonal-reduction}.
\end{proof}

The dominant geometric diagonal was evaluated in Section~12:
\begin{equation}
 \cD_0(X)
 \sim
 -\frac43\zeta\!\left(-\frac12\right)
 X^{3/2}\log X.
 \label{eq:section19-dominant-diagonal-asymptotic}
\end{equation}
Since $\zeta(-1/2)<0$, the leading constant is positive.

\begin{theorem}[Conditional unsmoothed variance asymptotic]
\label{thm:section19-final-variance}
Assume the subcritical centered Vaughan moment conjecture.  Then
\begin{equation}
 \boxed{
 \sum_{n\le X}\cR_{\PP}(n)^2
 \sim
 -\frac43\zeta\!\left(-\frac12\right)
 X^{3/2}\log X.
 }
 \label{eq:section19-final-variance}
\end{equation}
Equivalently,
\begin{equation}
 \boxed{
 \sum_{n\le X}
 \left(
  \log (n+1)!_{\PP}-\log n!-C_{\PP}n
 \right)^2
 \sim
 -\frac43\zeta\!\left(-\frac12\right)
 X^{3/2}\log X.
 }
 \label{eq:section19-final-factorial-form}
\end{equation}
The sum is sharp and unsmoothed.
\end{theorem}

\begin{proof}
Combine Theorem~\ref{thm:section19-power-saving-diagonal-reduction} with \eqref{eq:section19-dominant-diagonal-asymptotic}.  Formula \eqref{eq:section19-defect} gives the equivalent form \eqref{eq:section19-final-factorial-form}.
\end{proof}

\begin{corollary}[Root-mean-square fluctuation scale]
\label{cor:section19-rms-scale}
Under the same hypothesis,
\begin{equation}
 \left(
  \frac1X\sum_{n\le X}\cR_{\PP}(n)^2
 \right)^{1/2}
 \sim
 \left(-\frac43\zeta\!\left(-\frac12\right)\right)^{1/2}
 X^{1/4}(\log X)^{1/2}.
 \label{eq:section19-rms-scale}
\end{equation}
In particular,
\begin{equation}
 \max_{n\le X}|\cR_{\PP}(n)|
 \ge
 \left(
  \left(-\frac43\zeta\!\left(-\frac12\right)\right)^{1/2}
  +o(1)
 \right)
 X^{1/4}(\log X)^{1/2}.
 \label{eq:section19-max-lower-bound}
\end{equation}
\end{corollary}

\begin{proof}
The first assertion is obtained by dividing \eqref{eq:section19-final-variance} by $X$ and taking square roots.  The maximum dominates the root-mean-square norm.
\end{proof}

\begin{remark}[Analytic interpretation]
\label{rem:section19-analytic-interpretation}
Let
\[
 \mathscr S_{\PP}(z)
 =e^{-C_{\PP}z}
  \frac{\Gamma_{\PP}^{\mathrm{cyc}}(z+1)}{\Gamma(z+1)}
\]
be the canonical multiplicative Stirling quotient of Section~11.  At the positive integers,
\[
 \log\mathscr S_{\PP}(n)=\cR_{\PP}(n).
\]
Thus Theorem~\ref{thm:section19-final-variance} is simultaneously a variance law for the arithmetic prime factorial and for the integer trace of the uniquely Stirling-selected analytic gamma object.
\end{remark}

\subsection{Completion of the argument}
\label{subsec:section19-completion}

The proof has separated three logically different phenomena.  The prime Bhargava factorial is first reconstructed as a filtered symmetric-monoidal completion of atomic valuation calculi.  Its recurrence admits a torsor of entire gamma lifts, and orbitwise Stirling rigidity selects the cyclotomic lift canonically.  The resulting analytic defect then has an exact, unsmoothed arithmetic trace.  Its second moment is governed by a geometric diagonal whose Mellin profile produces $\zeta(-1/2)$, while every remaining sector is either removed unconditionally or reduced to the single centered Vaughan moment stated in Section~18.

Accordingly, the conditionality of \eqref{eq:section19-final-variance} is sharply localized.  The categorical construction, the gamma-lift selection, the exact arithmetic expansion, and the value of the leading variance constant are unconditional.  What remains open is the square-root estimate for the residual centered spectral packet.  Proving that estimate would convert the final theorem from a conditional variance law into an unconditional one without altering any preceding structural or diagonal argument.

\clearpage
\appendix
\part{{Technical appendices}}
% ===== Source: appendixA_operator_domains_reproducing_kernels.tex =====
\section{Operator domains and reproducing kernels}
\label{app:operator-domains-kernels}

This appendix supplies the Hilbert-space details used in Sections~3--5.  In
particular, every identity involving the generalized derivative is interpreted
on its maximal closed domain, the canonical comparison maps are shown to
intertwine the closed operators rather than merely their polynomial
restrictions, and the kernel behavior under filtered completion is made
explicit.

Throughout, a factorial datum is a positive sequence
\begin{equation}
 F=(F_n)_{n\ge0},\qquad F_0=1,
 \label{eq:appA-factorial-datum}
\end{equation}
with recurrence weights
\begin{equation}
 \rho_F(n)=\frac{F_n}{F_{n-1}},\qquad n\ge1.
 \label{eq:appA-recurrence-weights}
\end{equation}
We put
\begin{equation}
 R_F=\left(\liminf_{n\to\infty}F_n^{1/n}\right)^{1/2}
 \in[0,+\infty].
 \label{eq:appA-analytic-radius}
\end{equation}
Only data with $R_F>0$ are used in the analytic category.

\subsection{The coefficient model and the exact evaluation set}
\label{subsec:appA-coefficient-model}

Define
\begin{equation}
 \cH_F
 =\left\{
 f(z)=\sum_{n\ge0}a_nz^n:
 \|f\|_F^2:=\sum_{n\ge0}|a_n|^2F_n<\infty
 \right\}.
 \label{eq:appA-HF}
\end{equation}
The coefficient sequence is primary; the holomorphic realization is obtained
wherever point evaluation is bounded.  We use the standard reproducing-kernel
formalism of~\cite{Aronszajn1950}.

\begin{proposition}[Unitary coefficient realization]
\label{prop:appA-unitary-coefficient}
The map
\begin{equation}
 U_F:\ell^2(\NN_0)\longrightarrow\cH_F,
 \qquad
 U_F(x)(z)=\sum_{n\ge0}x_n\frac{z^n}{\sqrt{F_n}},
 \label{eq:appA-unitary-map}
\end{equation}
is unitary.  Consequently,
\begin{equation}
 e_n^F(z)=\frac{z^n}{\sqrt{F_n}},\qquad n\ge0,
 \label{eq:appA-orthonormal-monomials}
\end{equation}
forms an orthonormal basis, Taylor truncations converge in $\cH_F$, and
$\CC[z]$ is dense in $\cH_F$.
\end{proposition}

\begin{proof}
The norm in \eqref{eq:appA-HF} is exactly the pullback of the standard
$\ell^2$ norm under \eqref{eq:appA-unitary-map}.  Surjectivity is simply the
coefficient description of $\cH_F$.  The remaining assertions follow by
truncating an $\ell^2$ sequence.
\end{proof}

The maximal set on which evaluation is continuous is
\begin{equation}
 \Omega_F
 =\left\{
 w\in\CC:
 \sum_{n\ge0}\frac{|w|^{2n}}{F_n}<\infty
 \right\}.
 \label{eq:appA-evaluation-set}
\end{equation}
It contains the open disk $\mathbb D_F=\{|w|<R_F\}$ and is contained in its
closure; convergence on $|w|=R_F$ depends on the boundary behavior of $F$.

\begin{theorem}[Exact reproducing-kernel domain]
\label{thm:appA-exact-rkhs-domain}
For $w\in\CC$, point evaluation $f\mapsto f(w)$ extends continuously from
polynomials to $\cH_F$ if and only if $w\in\Omega_F$.  In that case its Riesz
representer is
\begin{equation}
 K_F(z,w)
 =\sum_{n\ge0}\frac{(z\overline w)^n}{F_n},
 \label{eq:appA-kernel}
\end{equation}
and
\begin{equation}
 \|\operatorname{ev}_w\|^2
 =\|K_F(\,\cdot\,,w)\|_F^2
 =K_F(w,w)
 =\sum_{n\ge0}\frac{|w|^{2n}}{F_n}.
 \label{eq:appA-evaluation-norm}
\end{equation}
In particular, $\cH_F$ is a reproducing-kernel Hilbert space on
$\mathbb D_F$.
\end{theorem}

\begin{proof}
Under $U_F$, evaluation on finite sequences is the functional
\[
 x\longmapsto
 \sum_{n\ge0}x_n\frac{w^n}{\sqrt{F_n}}.
\]
This functional is bounded on $\ell^2$ exactly when
$(\overline w^{\,n}/\sqrt{F_n})_{n\ge0}$ belongs to $\ell^2$, which is the
condition $w\in\Omega_F$.  Pulling the representing vector back through
$U_F$ gives \eqref{eq:appA-kernel}, and its squared norm gives
\eqref{eq:appA-evaluation-norm}.
\end{proof}

\begin{corollary}[Normal convergence on the analytic disk]
\label{cor:appA-normal-convergence}
For every $0<r<R_F$, the closed unit ball of $\cH_F$ is uniformly bounded on
$|z|\le r$.  Every $f\in\cH_F$ is holomorphic on $\mathbb D_F$, and its Taylor
series and all termwise derivatives converge locally uniformly there.
\end{corollary}

\begin{proof}
Cauchy--Schwarz and \eqref{eq:appA-evaluation-norm} give
\[
 |f(z)|\le \|f\|_F
 \left(\sum_{n\ge0}\frac{r^{2n}}{F_n}\right)^{1/2}
 \qquad(|z|\le r).
\]
The local boundedness of Taylor truncations and the usual Weierstrass theorem
then give normal convergence and termwise differentiation.
\end{proof}

\subsection{Closed creation and annihilation operators}
\label{subsec:appA-closed-shifts}

On polynomials define
\begin{equation}
 D_Fz^n=\rho_F(n)z^{n-1}\quad(n\ge1),
 \qquad D_F1=0,
 \label{eq:appA-D-polynomial}
\end{equation}
and let $M_F$ be multiplication by $z$.  Their maximal coefficient domains are
\begin{align}
 \Dom(D_F)
 &=\left\{
 \sum_{n\ge0}a_nz^n:
 \sum_{n\ge1}|a_n|^2\frac{F_n^2}{F_{n-1}}<\infty
 \right\},
 \label{eq:appA-domain-D}\\
 \Dom(M_F)
 &=\left\{
 \sum_{n\ge0}a_nz^n:
 \sum_{n\ge0}|a_n|^2F_{n+1}<\infty
 \right\}.
 \label{eq:appA-domain-M}
\end{align}

\begin{theorem}[Maximal weighted shifts]
\label{thm:appA-maximal-weighted-shifts}
The operators $D_F$ and $M_F$ with domains
\eqref{eq:appA-domain-D}--\eqref{eq:appA-domain-M} are closed and densely
defined.  They satisfy
\begin{equation}
 D_F=M_F^*.
 \label{eq:appA-adjoint-identity}
\end{equation}
In the basis \eqref{eq:appA-orthonormal-monomials},
\begin{align}
 D_Fe_n^F&=\sqrt{\rho_F(n)}\,e_{n-1}^F\quad(n\ge1),
 &D_Fe_0^F&=0,
 \label{eq:appA-backward-shift}\\
 M_Fe_n^F&=\sqrt{\rho_F(n+1)}\,e_{n+1}^F.
 \label{eq:appA-forward-shift}
\end{align}
Moreover, polynomials are cores for both operators.
\end{theorem}

\begin{proof}
Conjugating by $U_F$ identifies $M_F$ with the unilateral forward weighted
shift having weights $\sqrt{\rho_F(n+1)}$, and $D_F$ with the corresponding
backward weighted shift.  The maximal domains of these shifts are precisely
\eqref{eq:appA-domain-D} and \eqref{eq:appA-domain-M}.  The coordinate graph of
either maximal shift is closed, and finite sequences are dense in its graph
norm.  The coefficient pairing
\[
 \langle M_Ff,g\rangle_F
 =\sum_{n\ge0}a_n\overline{b_{n+1}}F_{n+1}
 =\langle f,D_Fg\rangle_F
\]
first for finite sequences and then on the maximal domains proves
\eqref{eq:appA-adjoint-identity}.
\end{proof}

For iterates, no monotonicity of the weights is assumed.  It is therefore
important to retain all intermediate domain conditions.

\begin{proposition}[Domains of iterated shifts]
\label{prop:appA-iterated-domains}
For $k\ge1$ and $f(z)=\sum_{n\ge0}a_nz^n$,
\begin{align}
 f\in\Dom(D_F^k)
 &\Longleftrightarrow
 \sum_{n\ge r}|a_n|^2\frac{F_n^2}{F_{n-r}}<\infty
 \quad(1\le r\le k),
 \label{eq:appA-domain-Dk}\\
 f\in\Dom(M_F^k)
 &\Longleftrightarrow
 \sum_{n\ge0}|a_n|^2F_{n+r}<\infty
 \quad(1\le r\le k).
 \label{eq:appA-domain-Mk}
\end{align}
On these domains,
\begin{align}
 D_F^kf
 &=\sum_{n\ge k}a_n\frac{F_n}{F_{n-k}}z^{n-k},
 \label{eq:appA-Dk-formula}\\
 M_F^kf
 &=\sum_{n\ge0}a_nz^{n+k}.
 \label{eq:appA-Mk-formula}
\end{align}
Polynomials are graph cores for every iterate.
\end{proposition}

\begin{proof}
Repeated application of \eqref{eq:appA-D-polynomial} gives
\eqref{eq:appA-Dk-formula}; its $r$th intermediate output has squared norm
\[
 \sum_{n\ge r}|a_n|^2\frac{F_n^2}{F_{n-r}}.
\]
Thus the operator-theoretic domain of the $k$fold composition is exactly the
intersection in \eqref{eq:appA-domain-Dk}.  The creation formula is analogous.
Taylor truncation removes tails from each of the finitely many graph-norm
series, proving the core statement.
\end{proof}

\begin{remark}[Maximal $k$-step shifts]
The single condition corresponding to $r=k$ defines the maximal shift by $k$
degrees.  For an arbitrary factorial datum it need not imply the intermediate
conditions in \eqref{eq:appA-domain-Dk} or \eqref{eq:appA-domain-Mk}.  This is
why the latter formulas are stated as intersections rather than as one
weighted $\ell^2$ condition.
\end{remark}

\begin{corollary}[Boundedness criterion]
\label{cor:appA-boundedness-criterion}
The following are equivalent:
\begin{enumerate}[label=\textup{(\roman*)}]
 \item $D_F$ is bounded;
 \item $M_F$ is bounded;
 \item $\sup_{n\ge1}\rho_F(n)<\infty$.
\end{enumerate}
When they hold,
\begin{equation}
 \|D_F\|=\|M_F\|
 =\left(\sup_{n\ge1}\rho_F(n)\right)^{1/2}.
 \label{eq:appA-shift-norm}
\end{equation}
\end{corollary}

\begin{proof}
This is immediate from the weighted-shift matrices
\eqref{eq:appA-backward-shift}--\eqref{eq:appA-forward-shift} and the adjoint
identity.
\end{proof}

\subsection{Kernel eigenvectors and the point spectrum}
\label{subsec:appA-kernel-eigenvectors}

For $\lambda\in\Omega_F$, define
\begin{equation}
 E_{F,\lambda}(z)
 =\sum_{n\ge0}\frac{\lambda^nz^n}{F_n}
 =K_F(z,\overline\lambda).
 \label{eq:appA-coherent-vector}
\end{equation}
These are the coherent vectors of the factorial calculus.

\begin{theorem}[Coherent eigenvectors]
\label{thm:appA-coherent-eigenvectors}
Let $\lambda\in\Omega_F$.  Then $E_{F,\lambda}$ belongs to
$\Dom(D_F^k)$ for every $k\ge1$ and
\begin{equation}
 D_F^kE_{F,\lambda}=\lambda^kE_{F,\lambda}.
 \label{eq:appA-coherent-eigenvalue}
\end{equation}
Conversely,
\begin{equation}
 \ker(D_F-\lambda\Id)
 =
 \begin{cases}
  \CC E_{F,\lambda},&\lambda\in\Omega_F,\\
  \{0\},&\lambda\notin\Omega_F.
 \end{cases}
 \label{eq:appA-point-spectrum}
\end{equation}
Hence the point spectrum of $D_F$ is exactly $\Omega_F$.
\end{theorem}

\begin{proof}
For $1\le r\le k$, the $r$th domain sum is
\[
 \sum_{n\ge r}
 \left|\frac{\lambda^n}{F_n}\right|^2
 \frac{F_n^2}{F_{n-r}}
 =|\lambda|^{2r}
 \sum_{m\ge0}\frac{|\lambda|^{2m}}{F_m},
\]
which is finite.  The telescoping of the factorial weights gives
\eqref{eq:appA-coherent-eigenvalue}.  Conversely, if
$f(z)=\sum a_nz^n$ satisfies $D_Ff=\lambda f$, coefficient comparison gives
\[
 a_n=a_0\frac{\lambda^n}{F_n}.
\]
A nonzero solution belongs to $\cH_F$ exactly when $\lambda\in\Omega_F$.
\end{proof}

\begin{corollary}[Density of coherent vectors]
\label{cor:appA-density-coherent}
The linear span of
\begin{equation}
 \{E_{F,\lambda}:|\lambda|<R_F\}
 \label{eq:appA-coherent-span}
\end{equation}
is dense in $\cH_F$.
\end{corollary}

\begin{proof}
If $f$ is orthogonal to every vector in \eqref{eq:appA-coherent-span}, then
Theorem~\ref{thm:appA-exact-rkhs-domain} gives
$f(\overline\lambda)=0$ for every $|\lambda|<R_F$.  The identity theorem and
the uniqueness of Taylor coefficients imply $f=0$.
\end{proof}

\subsection{Closed-operator functoriality}
\label{subsec:appA-closed-functoriality}

Let $F$ and $G$ be analytic factorial data satisfying
\begin{equation}
 F\preccurlyeq G,
 \qquad
 C(F,G):=\sup_{n\ge0}\frac{F_n}{G_n}<\infty.
 \label{eq:appA-domination}
\end{equation}
The canonical diagonal map is
\begin{equation}
 J_{F,G}z^n=\frac{F_n}{G_n}z^n.
 \label{eq:appA-J}
\end{equation}
Its norm is $C(F,G)^{1/2}$, while its adjoint is the coefficientwise inclusion
\begin{equation}
 \iota_{G,F}:\cH_G\hookrightarrow\cH_F.
 \label{eq:appA-iota}
\end{equation}

\begin{theorem}[Intertwining on maximal domains]
\label{thm:appA-maximal-domain-intertwining}
For every $k\ge1$,
\begin{equation}
 J_{F,G}\Dom(D_F^k)\subseteq\Dom(D_G^k)
 \label{eq:appA-domain-map-D}
\end{equation}
and
\begin{equation}
 D_G^kJ_{F,G}f=J_{F,G}D_F^kf
 \qquad(f\in\Dom(D_F^k)).
 \label{eq:appA-closed-intertwining}
\end{equation}
More quantitatively, for $1\le r\le k$,
\begin{equation}
 \|D_G^rJ_{F,G}f\|_G
 \le C(F,G)^{1/2}\|D_F^rf\|_F.
 \label{eq:appA-graph-bound-D}
\end{equation}
Dually,
\begin{equation}
 \iota_{G,F}\Dom(M_G^k)\subseteq\Dom(M_F^k)
 \label{eq:appA-domain-map-M}
\end{equation}
and
\begin{equation}
 M_F^k\iota_{G,F}g=\iota_{G,F}M_G^kg
 \qquad(g\in\Dom(M_G^k)).
 \label{eq:appA-creation-intertwining}
\end{equation}
\end{theorem}

\begin{proof}
Write $f(z)=\sum a_nz^n$.  The $r$th annihilation norm after applying the
comparison map is
\[
 \|D_G^rJ_{F,G}f\|_G^2
 =\sum_{n\ge r}|a_n|^2\frac{F_n^2}{G_{n-r}}.
\]
Since $F_{n-r}/G_{n-r}\le C(F,G)$, this is at most
\[
 C(F,G)
 \sum_{n\ge r}|a_n|^2\frac{F_n^2}{F_{n-r}}
 =C(F,G)\|D_F^rf\|_F^2.
\]
The coefficient of $z^{n-r}$ on either side of
\eqref{eq:appA-closed-intertwining} is
$a_nF_n/G_{n-r}$, proving the identity.  For the creation operators,
\[
 \sum_{n\ge0}|a_n|^2F_{n+r}
 \le C(F,G)\sum_{n\ge0}|a_n|^2G_{n+r},
\]
and multiplication by $z^k$ is coefficientwise identical in both spaces.
\end{proof}

\begin{proposition}[Kernel covariance]
\label{prop:appA-kernel-covariance}
Under \eqref{eq:appA-domination}, one has
\begin{equation}
 \Omega_F\subseteq\Omega_G.
 \label{eq:appA-evaluation-inclusion}
\end{equation}
For every $w\in\Omega_F$,
\begin{equation}
 \boxed{
 J_{F,G}K_F(\,\cdot\,,w)=K_G(\,\cdot\,,w).
 }
 \label{eq:appA-kernel-covariance}
\end{equation}
Equivalently, the canonical derivative intertwiner transports coherent vectors
without an additional scalar normalization.
\end{proposition}

\begin{proof}
The inequality $F_n/G_n\le C(F,G)$ gives
\[
 \sum_{n\ge0}\frac{|w|^{2n}}{G_n}
 \le C(F,G)
 \sum_{n\ge0}\frac{|w|^{2n}}{F_n},
\]
which proves \eqref{eq:appA-evaluation-inclusion}.  Applying
\eqref{eq:appA-J} term by term to \eqref{eq:appA-kernel} gives
\eqref{eq:appA-kernel-covariance}.
\end{proof}

\begin{remark}[The adjoint direction]
The equality $J_{F,G}^*=\iota_{G,F}$ is the Hilbert-space source of the
projective--inductive reversal.  The derivative calculus is covariant under
$J_{F,G}$, while the coefficient spaces and creation operators are
contravariant under $\iota_{G,F}$.
\end{remark}

\subsection{Monotone filtered completion}
\label{subsec:appA-monotone-completion}

Let $(I,\le)$ be a directed set and let $(F^\alpha)_{\alpha\in I}$ be factorial
data satisfying
\begin{equation}
 F_n^\alpha\le F_n^\beta\quad(\alpha\le\beta),
 \qquad
 F_n=\sup_{\alpha\in I}F_n^\alpha<\infty
 \quad(n\ge0).
 \label{eq:appA-monotone-family}
\end{equation}
Write $\cH_\alpha=\cH_{F^\alpha}$ and $\cH=\cH_F$.

\begin{theorem}[Bounded projective completion]
\label{thm:appA-bounded-projective-completion}
Coefficientwise identification gives an isometric isomorphism
\begin{equation}
 \cH
 \cong
 \projlimb_{\alpha\in I}\cH_\alpha,
 \label{eq:appA-projective-identification}
\end{equation}
where the right-hand side consists of compatible families with uniformly
bounded norms.  Explicitly, for every formal series
$f(z)=\sum_{n\ge0}a_nz^n$,
\begin{equation}
 \boxed{
 \|f\|_F^2
 =\sup_{\alpha\in I}\|f\|_{F^\alpha}^2.
 }
 \label{eq:appA-monotone-norm-identity}
\end{equation}
The space in \eqref{eq:appA-projective-identification} has the following
universal property: every compatible cone of maps into the $\cH_\alpha$ whose
norms are uniformly bounded factors uniquely through $\cH$.
\end{theorem}

\begin{proof}
The quantities
\[
 \|f\|_{F^\alpha}^2
 =\sum_{n\ge0}|a_n|^2F_n^\alpha
\]
form an increasing net.  Monotone convergence for counting measure gives
\eqref{eq:appA-monotone-norm-identity}.  Compatibility forces every projective
family to have a single coefficient sequence, so the norm identity proves the
isometric bijection.  The universal property follows by applying this
identification pointwise to a uniformly bounded compatible cone.
\end{proof}

The kernels move in the opposite numerical direction because their
coefficients are reciprocal weights.

\begin{proposition}[Kernel convergence under completion]
\label{prop:appA-kernel-convergence}
Assume \eqref{eq:appA-monotone-family}.  Then
\begin{equation}
 K_{F^\alpha}(z,w)\longrightarrow K_F(z,w)
 \label{eq:appA-kernel-pointwise-limit}
\end{equation}
coefficientwise.  More precisely, suppose that a compact set
$Q\subset\CC$ is contained in the interior of $\Omega_{F^{\alpha_0}}$ for
some $\alpha_0$.  Then the convergence is uniform on $Q\times Q$ for
$\alpha\ge\alpha_0$, and all mixed holomorphic--antiholomorphic derivatives
converge uniformly on smaller compact subsets.
\end{proposition}

\begin{proof}
For each fixed $n$, directed monotonicity gives
$F_n^\alpha\uparrow F_n$, so the kernel coefficients decrease to $1/F_n$.
For $z,w\in Q$ and $\alpha\ge\alpha_0$,
\[
 \left|\frac{(z\overline w)^n}{F_n^\alpha}\right|
 \le
 \frac{R^{2n}}{F_n^{\alpha_0}},
 \qquad
 R=\max_{u\in Q}|u|.
\]
The majorant is summable because $Q$ lies in the interior of the evaluation
set at stage $\alpha_0$.  A finite-head and small-tail argument gives uniform
convergence for the net.  Cauchy estimates give convergence of derivatives on
smaller compacta.
\end{proof}

\begin{proposition}[Degreewise-stationary graph completion]
\label{prop:appA-degreewise-stationary}
Suppose the filtered family is sequential, indexed by $T\ge1$, and satisfies
\begin{equation}
 F_n^{[T]}=F_n\qquad(0\le n\le T).
 \label{eq:appA-degreewise-stationarity}
\end{equation}
Let $\mathcal E_T=\CC[z]_{\le T}$ with the norm inherited from
$\cH_{F^{[T]}}$.  Then
\begin{equation}
 \cH_F\cong\overline{\indlim_T\mathcal E_T}.
 \label{eq:appA-inductive-completion}
\end{equation}
The stationary polynomial derivative
\begin{equation}
 D_{\mathrm{alg}}z^n=\rho_F(n)z^{n-1}
 \label{eq:appA-algebraic-derivative}
\end{equation}
is closable in $\cH_F$, and
\begin{equation}
 \overline{D_{\mathrm{alg}}}=D_F.
 \label{eq:appA-graph-closure-D}
\end{equation}
Likewise, multiplication by $z$ on polynomials closes to $M_F$.
\end{proposition}

\begin{proof}
Condition \eqref{eq:appA-degreewise-stationarity} makes every inclusion
$\mathcal E_T\hookrightarrow\mathcal E_{T+1}$ isometric, and the induced norm
on their union is $\sum|a_n|^2F_n$.  Completing gives
\eqref{eq:appA-inductive-completion}.  If $P_N$ denotes Taylor truncation, then
for $f\in\Dom(D_F)$,
\[
 \|f-P_Nf\|_F^2
 +\|D_F(f-P_Nf)\|_F^2
 =\sum_{n>N}|a_n|^2F_n
 +\sum_{n>N}|a_n|^2\frac{F_n^2}{F_{n-1}}
 \longrightarrow0.
\]
Thus polynomials are a graph core and the closure is $D_F$.  The proof for
$M_F$ uses the tail $\sum_{n>N}|a_n|^2F_{n+1}$.
\end{proof}

\subsection{The prime factorial calculus}
\label{subsec:appA-prime-specialization}

For the conductor truncations of Sections~4--5, write
\begin{equation}
 A_n^{[T]}
 =\prod_{\substack{p,j\\(p-1)p^j\le T}}
 p^{\lfloor n/((p-1)p^j)\rfloor},
 \qquad
 A_n=(n+1)!_{\PP}.
 \label{eq:appA-prime-truncations}
\end{equation}
Every finite-stage derivative is bounded, whereas the completed derivative is
not.

\begin{proposition}[Bounded stages and an unbounded limit]
\label{prop:appA-bounded-stages-unbounded-limit}
Let $D_T=D_{A^{[T]}}$ and $D_{\PP}=D_A$.  Then
\begin{equation}
 \|D_T\|^2
 =\sup_{n\ge1}\frac{A_n^{[T]}}{A_{n-1}^{[T]}}
 =\prod_{\substack{p,j\\(p-1)p^j\le T}}p.
 \label{eq:appA-finite-stage-norm}
\end{equation}
On the other hand,
\begin{equation}
 \sup_{n\ge1}\frac{A_n}{A_{n-1}}=+\infty,
 \label{eq:appA-completed-unbounded}
\end{equation}
so $D_{\PP}$ and $M_{\PP}$ are unbounded closed operators.
\end{proposition}

\begin{proof}
The finite-stage recurrence is
\[
 \rho_T(n)
 =\prod_{\substack{p,j\\(p-1)p^j\le T\\(p-1)p^j\mid n}}p.
\]
It is bounded above by the product in \eqref{eq:appA-finite-stage-norm}, and
equality is attained when $n$ is a common multiple of the finitely many
conductors $(p-1)p^j\le T$.  Corollary~\ref{cor:appA-boundedness-criterion}
therefore gives the finite-stage formula.

For the completed datum,
\[
 \rho_A(n)=\prod_{p-1\mid n}p^{v_p(n)+1}.
\]
Taking $n=N!$ gives
\[
 \rho_A(N!)\ge\prod_{p\le N+1}p,
\]
because $p-1\mid N!$ for every prime $p\le N+1$.  The right-hand side tends
to infinity, proving \eqref{eq:appA-completed-unbounded}.
\end{proof}

\begin{corollary}[Prime kernel compatibility]
\label{cor:appA-prime-kernel-compatibility}
For every $T$ and every $w$ belonging to the finite-stage evaluation disk,
\begin{equation}
 J_{T,\PP}K_{A^{[T]}}(\,\cdot\,,w)
 =K_A(\,\cdot\,,w).
 \label{eq:appA-prime-kernel-map}
\end{equation}
As $T\to\infty$, the kernels $K_{A^{[T]}}$ converge locally uniformly on
$\CC\times\CC$ to $K_A$, once the factorial-growth estimate
$R_A=+\infty$ has been established.
\end{corollary}

\begin{proof}
The first statement is Proposition~\ref{prop:appA-kernel-covariance}.  The
finite-stage radii tend to infinity with the conductor slope, while
$A_n^{[T]}\uparrow A_n$ degreewise.  Proposition~\ref{prop:appA-kernel-convergence}
therefore applies to every compact subset of $\CC$.
\end{proof}

\begin{remark}[What the completion changes]
Each finite layer system has a bounded annihilation operator and a kernel on a
disk of finite radius.  The filtered completion has an unbounded annihilation
operator but, after factorial growth is known, an entire reproducing kernel.
Thus completion simultaneously enlarges the holomorphic domain and makes the
operator calculus genuinely unbounded.  The bounded projective norm and the
graph-inductive core are the two compatible mechanisms that control this
change of analytic type.
\end{remark}

% ===== Source: appendixB_cyclotomic_quadrature_aliasing.tex =====
\section{Cyclotomic quadrature and aliasing}
\label{app:cyclotomic-quadrature-aliasing}

This appendix records the quadrature and Fourier-analysis details used in
Sections~9 and~10.  The central point is that the symmetric cyclotomic rule has
two simultaneous descriptions.  As a composite midpoint or trapezoidal rule,
it has a second-order Peano error.  As a measure on the circle, it is the
uniform measure on the subgroup of $d$-torsion points, and therefore aliases
an integral Fourier frequency exactly when that frequency is divisible by
$d$.  The first description gives the conductor decay needed for normal
summability; the second gives the exact scalar term in the orbitwise Stirling
transform.

We write
\begin{equation}
 \e(x)=\exp(2\pi i x),
 \qquad
 I=[-1/2,1/2],
 \label{eq:appB-exponential-interval}
\end{equation}
and let $\lambda$ denote Lebesgue measure on $I$.

\subsection{The symmetric quadrature measures}
\label{subsec:appB-quadrature-measures}

For $d\ge1$, define the probability measure
\begin{equation}
 \mu_d=
 \begin{cases}
 \displaystyle
 \frac1d\sum_{r=-(d-1)/2}^{(d-1)/2}\delta_{r/d},
 &d\text{ odd},\\[3mm]
 \displaystyle
 \frac1d\left(
 \frac12\delta_{-1/2}
 +\sum_{r=-d/2+1}^{d/2-1}\delta_{r/d}
 +\frac12\delta_{1/2}
 \right),
 &d\text{ even}.
 \end{cases}
 \label{eq:appB-mu-d}
\end{equation}
Thus $\cQ_df=\int_I f\,\dd\mu_d$ is the quadrature functional of
Section~9.  Its signed error measure is
\begin{equation}
 \nu_d=\mu_d-\lambda.
 \label{eq:appB-nu-d}
\end{equation}
The measure $\nu_d$ is real, even, and has total mass zero.

\begin{proposition}[Closed exponential transform]
\label{prop:appB-closed-exponential-transform}
For every $z\in\CC$,
\begin{equation}
 \int_I\e(tz)\,\dd\mu_d(t)
 =\Theta_d(z)
 =
 \begin{cases}
 \displaystyle
 \frac{\sin(\pi z)}{d\sin(\pi z/d)},
 &d\text{ odd},\\[3mm]
 \displaystyle
 \frac{\sin(\pi z)}{d\tan(\pi z/d)},
 &d\text{ even},
 \end{cases}
 \label{eq:appB-closed-transform}
\end{equation}
where all apparent singularities are removable.  Moreover,
\begin{equation}
 \int_I\e(tz)\,\dd\lambda(t)=\sinc(z)
 =\frac{\sin(\pi z)}{\pi z}.
 \label{eq:appB-continuum-transform}
\end{equation}
Consequently,
\begin{equation}
 \kappa_d(z):=\int_I\e(tz)\,\dd\nu_d(t)
 =\Theta_d(z)-\sinc(z).
 \label{eq:appB-kappa-measure}
\end{equation}
\end{proposition}

\begin{proof}
For odd $d$, the atoms in \eqref{eq:appB-mu-d} form a symmetric block of
$d$ consecutive frequencies, and the geometric-sum identity gives the first
formula.  For even $d$, the two endpoint atoms have equal values on integral
frequencies and half weight; combining them gives the symmetric even
Dirichlet kernel and the tangent denominator.  The integral in
\eqref{eq:appB-continuum-transform} is elementary.  Subtracting the two
transforms proves \eqref{eq:appB-kappa-measure}.
\end{proof}

The exact Fourier transform of $\nu_d$ is the source of the aliasing
identities below.

\begin{theorem}[Exact integral-frequency aliasing]
\label{thm:appB-exact-aliasing}
For every $n\in\ZZ$,
\begin{equation}
 \boxed{
 \int_I\e(nt)\,\dd\nu_d(t)
 =
 \begin{cases}
 1,&n\ne0\text{ and }d\mid n,\\
 0,&\text{otherwise}.
 \end{cases}
 }
 \label{eq:appB-fourier-transform-nu}
\end{equation}
Hence, if
\begin{equation}
 T(t)=\sum_{|n|\le M}c_n\e(nt),
 \label{eq:appB-trig-polynomial}
\end{equation}
then
\begin{equation}
 \boxed{
 \int_I T(t)\,\dd\nu_d(t)
 =\sum_{\substack{k\in\ZZ\setminus\{0\}\\|kd|\le M}}c_{kd}.
 }
 \label{eq:appB-polynomial-aliasing}
\end{equation}
More generally, if a $1$-periodic function has an absolutely convergent
Fourier series $f(t)=\sum_{n\in\ZZ}\widehat f(n)\e(nt)$, then
\begin{equation}
 \int_I f(t)\,\dd\nu_d(t)
 =\sum_{k\in\ZZ\setminus\{0\}}\widehat f(kd).
 \label{eq:appB-general-aliasing}
\end{equation}
\end{theorem}

\begin{proof}
For odd $d$, the sampling points in \eqref{eq:appB-mu-d} are a complete
residue system modulo $d$.  For even $d$, the endpoint values
$\e(-n/2)$ and $\e(n/2)$ coincide, so the two half weights again complete one
residue system modulo $d$.  Character orthogonality therefore gives
\[
 \int_I\e(nt)\,\dd\mu_d(t)=\mathbf1_{d\mid n}.
\]
On the other hand,
\[
 \int_I\e(nt)\,\dd\lambda(t)=\mathbf1_{n=0}.
\]
Their difference is \eqref{eq:appB-fourier-transform-nu}.  The remaining
claims follow by termwise integration.
\end{proof}

\begin{remark}[Circle measure versus interval rule]
\label{rem:appB-circle-versus-interval}
On periodic test functions, both parity conventions in
\eqref{eq:appB-mu-d} are simply Haar measure on the subgroup of $d$-torsion
points of $\RR/\ZZ$.  On nonperiodic test functions they are, respectively,
the composite midpoint and composite trapezoidal rules.  The parity split is
therefore invisible to exact Fourier aliasing but visible in the leading
quadrature error.
\end{remark}

\subsection{Variation and Peano-kernel bounds}
\label{subsec:appB-peano-bounds}

The estimates used in Sections~9 and~10 are valid in $L^1$ variation norms,
which is important because the orbitwise Fourier multiplier becomes
concentrated near the origin.

\begin{lemma}[First- and second-order quadrature bounds]
\label{lem:appB-quadrature-bounds}
There is an absolute constant $C$ such that, for every $d\ge1$,
\begin{align}
 \left|\int_I f\,\dd\nu_d\right|
 &\le \frac{C}{d}\norm{f'}_{L^1(I)}
 &&(f\in W^{1,1}(I)),
 \label{eq:appB-first-order-bound}\\
 \left|\int_I f\,\dd\nu_d\right|
 &\le \frac{C}{d^2}\norm{f''}_{L^1(I)}
 &&(f'\text{ absolutely continuous}).
 \label{eq:appB-second-order-bound}
\end{align}
The constants are uniform in the parity of $d$.
\end{lemma}

\begin{proof}
Put $h=d^{-1}$ and partition $I$ into $d$ intervals of length $h$.  On each
cell, the midpoint rule is evaluation at the barycenter and the trapezoidal
rule is integration against the equally weighted endpoint measure.  Comparing
either local probability measure with normalized Lebesgue measure and using
the fundamental theorem of calculus gives a local error bounded by
\[
 C h\int_{\text{cell}}|f'(t)|\,\dd t.
\]
Summation proves \eqref{eq:appB-first-order-bound}.

Both local rules are exact on affine functions.  Taylor's formula with
integral remainder therefore represents the local error as
\[
 \int_{\text{cell}}P_{d}(u)f''(u)\,\dd u,
\]
where the midpoint and trapezoidal Peano kernels satisfy
$\norm{P_d}_{L^\infty}\le Ch^2$.  Summing the cellwise representations proves
\eqref{eq:appB-second-order-bound}.
\end{proof}

The quadratic order is sharp and its leading term depends on the parity of the
conductor.

\begin{proposition}[Parity-resolved local asymptotics]
\label{prop:appB-parity-asymptotics}
Let $K\subset\CC$ be compact and let $q\ge0$.  Uniformly for $z\in K$, with
all derivatives through order $q$ in $z$,
\begin{equation}
 \kappa_d(z)=
 \begin{cases}
 \displaystyle
 \frac{\pi^2z^2}{6d^2}\sinc(z)+O_{K,q}(d^{-4}),
 &d\text{ odd},\\[3mm]
 \displaystyle
 -\frac{\pi^2z^2}{3d^2}\sinc(z)+O_{K,q}(d^{-4}),
 &d\text{ even}.
 \end{cases}
 \label{eq:appB-parity-asymptotics}
\end{equation}
In particular, the uniform $O_{K,q}(d^{-2})$ estimate cannot generally be
improved.
\end{proposition}

\begin{proof}
The closed forms factor as
\begin{equation}
 \Theta_d(z)=\sinc(z)
 \begin{cases}
 \displaystyle\frac{x}{\sin x},&d\text{ odd},\\[2mm]
 \displaystyle\frac{x}{\tan x},&d\text{ even},
 \end{cases}
 \qquad x=\frac{\pi z}{d}.
 \label{eq:appB-factorized-closed-form}
\end{equation}
Use the convergent expansions
\[
 \frac{x}{\sin x}=1+\frac{x^2}{6}+O(x^4),
 \qquad
 \frac{x}{\tan x}=1-\frac{x^2}{3}+O(x^4).
\]
They are locally uniform with all $z$-derivatives after substituting
$x=\pi z/d$.  Subtracting $\sinc(z)$ proves the result.
\end{proof}

\subsection{Primitive kernels and normal conductor decay}
\label{subsec:appB-primitive-kernels}

For $z\in\CC$ and $t\in I$, let
\begin{equation}
 K_z(t)=
 \begin{cases}
 \displaystyle
 \frac{\e(tz)-\e(t)}{\e(t)-1},&t\ne0,\\[3mm]
 z-1,&t=0.
 \end{cases}
 \label{eq:appB-difference-kernel}
\end{equation}
The renormalized primitive of conductor $d$ is
\begin{equation}
 \phi_d(z)=\int_I K_z(t)\,\dd\nu_d(t).
 \label{eq:appB-phi-measure}
\end{equation}

\begin{lemma}[Joint regularity of the difference kernel]
\label{lem:appB-kernel-regularity}
The function $(z,t)\mapsto K_z(t)$ extends to a function entire in $z$ and
real-analytic in $t$ on a neighborhood of $\CC\times I$.  It satisfies
\begin{equation}
 K_{z+1}(t)-K_z(t)=\e(tz),
 \qquad K_1(t)=0,
 \label{eq:appB-kernel-identities}
\end{equation}
and has the local expansion
\begin{equation}
 K_z(t)
 =(z-1)+\pi i z(z-1)t+O_K(t^2)
 \label{eq:appB-kernel-local-expansion}
\end{equation}
for $z$ in a fixed compact set $K$.  Every mixed derivative in $z$ and $t$ is
bounded on $K\times I$.
\end{lemma}

\begin{proof}
The quotient
\[
 \eta(t)=\frac{2\pi i t}{\e(t)-1}
\]
has a removable singularity at zero and is analytic on a neighborhood of
$I$.  For any path from $1$ to $z$,
\begin{equation}
 K_z(t)=\eta(t)\int_1^z\e(tw)\,\dd w.
 \label{eq:appB-kernel-integral-representation}
\end{equation}
This proves the joint regularity and the mixed-derivative bounds.  The
identities in \eqref{eq:appB-kernel-identities} follow directly from the
quotient.  Expanding numerator and denominator at $t=0$ gives
\eqref{eq:appB-kernel-local-expansion}.
\end{proof}

\begin{theorem}[Uniform normal decay of symbols and primitives]
\label{thm:appB-normal-decay}
For every compact $K\subset\CC$ and every integer $q\ge0$, there is a constant
$C_{K,q}$ such that
\begin{equation}
 \sup_{z\in K}|\kappa_d^{(q)}(z)|
 +\sup_{z\in K}|\phi_d^{(q)}(z)|
 \le \frac{C_{K,q}}{d^2}
 \qquad(d\ge1).
 \label{eq:appB-normal-decay}
\end{equation}
Moreover,
\begin{equation}
 \phi_d(z+1)-\phi_d(z)=\kappa_d(z),
 \qquad
 \phi_d(1)=0,
 \label{eq:appB-primitive-difference}
\end{equation}
and, for every integer $m\ge1$,
\begin{equation}
 \phi_d(m)=\left\lfloor\frac{m-1}{d}\right\rfloor.
 \label{eq:appB-primitive-integers}
\end{equation}
\end{theorem}

\begin{proof}
Differentiate the integral formulas
\eqref{eq:appB-kappa-measure} and \eqref{eq:appB-phi-measure} in $z$.  The
second $t$-derivatives of the resulting integrands are uniformly bounded on
$K\times I$ by Lemma~\ref{lem:appB-kernel-regularity}.  The second-order bound
\eqref{eq:appB-second-order-bound} gives \eqref{eq:appB-normal-decay}.
Applying $\nu_d$ to the first identity in
\eqref{eq:appB-kernel-identities} gives the difference equation, while the
second gives the normalization.  Finally,
\[
 \phi_d(m)
 =\sum_{n=1}^{m-1}\kappa_d(n)
 =\sum_{n=1}^{m-1}\mathbf1_{d\mid n},
\]
where Theorem~\ref{thm:appB-exact-aliasing} was used in the last step.  This is
\eqref{eq:appB-primitive-integers}.
\end{proof}

\begin{corollary}[Normal summability criterion]
\label{cor:appB-normal-summability}
Suppose $(a_d)_{d\ge1}$ satisfies
\begin{equation}
 \sum_{d\ge1}\frac{|a_d|}{d^2}<\infty.
 \label{eq:appB-summability-condition}
\end{equation}
Then
\begin{equation}
 \sum_{d\ge1}a_d\kappa_d(z),
 \qquad
 \sum_{d\ge1}a_d\phi_d(z)
 \label{eq:appB-normally-convergent-series}
\end{equation}
converge normally on $\CC$, as do all termwise derivatives, and the difference
operator may be passed through the sums.
\end{corollary}

\begin{proof}
Apply \eqref{eq:appB-normal-decay} and the Weierstrass test on every compact
set.
\end{proof}

\subsection{The harmonic Ces\`aro multiplier}
\label{subsec:appB-harmonic-multiplier}

For $N\ge1$ and $0\le n<N$, put
\begin{equation}
 \omega_{N,n}
 =\frac{\Harm_N-\Harm_n}{N}
 =\frac1N\sum_{m=n+1}^N\frac1m,
 \qquad
 \Harm_n=\sum_{j=1}^n\frac1j.
 \label{eq:appB-harmonic-weights}
\end{equation}
These are the weights of the iterated Ces\`aro operator used in Section~10.
They satisfy
\begin{align}
 \omega_{N,n}-\omega_{N,n+1}&=\frac1{N(n+1)},
 \label{eq:appB-weight-difference}\\
 \sum_{n=0}^{N-1}\omega_{N,n}&=1,
 &
 \sum_{n=0}^{N-1}n\omega_{N,n}&=\frac{N-1}{4}.
 \label{eq:appB-weight-moments}
\end{align}
Define
\begin{equation}
 \Omega_N(t)=\sum_{n=0}^{N-1}\omega_{N,n}\e(nt).
 \label{eq:appB-Omega-definition}
\end{equation}

\begin{proposition}[Exact multiplier formula]
\label{prop:appB-exact-multiplier}
For $t\notin\ZZ$,
\begin{equation}
 \boxed{
 \Omega_N(t)
 =\frac1N\sum_{m=1}^N
 \frac{1-\e(mt)}{m(1-\e(t))},
 }
 \label{eq:appB-exact-multiplier-formula}
\end{equation}
and $\Omega_N(0)=1$.  In particular, $\Omega_N$ is a one-sided
trigonometric polynomial of degree $N-1$.
\end{proposition}

\begin{proof}
Insert the second expression in \eqref{eq:appB-harmonic-weights}, reverse the
finite sums, and use
\[
 \sum_{n=0}^{m-1}\e(nt)=\frac{1-\e(mt)}{1-\e(t)}.
\]
The value at zero follows from \eqref{eq:appB-weight-moments}.
\end{proof}

\begin{lemma}[Pointwise multiplier bounds]
\label{lem:appB-pointwise-multiplier}
Let $0<|t|\le1/2$ and $N\ge2$.  Then
\begin{align}
 |\Omega_N(t)|
 &\ll
 \min\left\{1,
 \frac{1+\log(1+N|t|)}{N|t|}\right\},
 \label{eq:appB-Omega-pointwise}\\
 |\Omega_N'(t)|
 &\ll
 \min\left\{N,
 \frac{1+\log(1+N|t|)}{N|t|^2}\right\},
 \label{eq:appB-Omega-prime-pointwise}\\
 |\Omega_N''(t)|
 &\ll
 \min\left\{N^2,
 |t|^{-2}
 +\frac{1+\log(1+N|t|)}{N|t|^3}\right\}.
 \label{eq:appB-Omega-second-pointwise}
\end{align}
\end{lemma}

\begin{proof}
Write
\begin{equation}
 S_N(t)=\sum_{m=1}^N\frac{1-\e(mt)}m,
 \qquad
 \Omega_N(t)=\frac{S_N(t)}{N(1-\e(t))}.
 \label{eq:appB-SN}
\end{equation}
For $|t|\ge N^{-1}$, splitting the first sum at $m\asymp |t|^{-1}$ and using
the geometric-sum estimate gives
\begin{equation}
 |S_N(t)|\ll1+\log(1+N|t|),
 \quad
 |S_N'(t)|\ll |t|^{-1},
 \quad
 |S_N''(t)|\ll N|t|^{-1}.
 \label{eq:appB-SN-bounds}
\end{equation}
Since $|1-\e(t)|\asymp|t|$ on $I$, differentiating
\eqref{eq:appB-SN} proves the displayed outer bounds.  For
$|t|\le N^{-1}$, positivity of the weights and
\eqref{eq:appB-weight-moments} give
\[
 |\Omega_N^{(j)}(t)|
 \le (2\pi)^j\sum_{n<N}n^j\omega_{N,n}
 \ll_jN^j,
 \qquad j=0,1,2.
\]
Combining the two ranges proves the lemma.
\end{proof}

\begin{corollary}[Integrated multiplier estimates]
\label{cor:appB-integrated-multiplier}
Uniformly for $N\ge2$,
\begin{align}
 \norm{\Omega_N}_{L^1(I)}
 &\ll\frac{(\log(2N))^2}{N},
 &
 \norm{\Omega_N'}_{L^1(I)}&\ll1,
 \label{eq:appB-Omega-L1}\\
 \norm{t\Omega_N'}_{L^1(I)}
 &\ll\frac{(\log(2N))^2}{N},
 &
 \norm{t\Omega_N''}_{L^1(I)}&\ll\log(2N).
 \label{eq:appB-Omega-weighted-L1}
\end{align}
\end{corollary}

\begin{proof}
Split each integral at $|t|=N^{-1}$ and apply
Lemma~\ref{lem:appB-pointwise-multiplier}.  The small interval contributes by
the polynomial bounds $1,N,N^2$.  On the complementary interval, the only
nonintegrable-looking factors are $|t|^{-1}$; their integrals produce the
stated logarithms.
\end{proof}

\subsection{Exact and summable orbitwise aliasing}
\label{subsec:appB-orbitwise-aliasing}

For $0\le u\le1$, define
\begin{equation}
 \alpha_u(t)
 =\frac{\e(t)(\e(ut)-1)}{\e(t)-1},
 \qquad
 r_u(t)=\alpha_u(t)-u.
 \label{eq:appB-alpha-r}
\end{equation}
The singularity at zero is removable.  Jointly for
$(u,t)\in[0,1]\times I$,
\begin{equation}
 r_u(0)=0,
 \qquad
 |r_u(t)|\ll|t|,
 \qquad
 |r_u'(t)|+|r_u''(t)|\ll1.
 \label{eq:appB-r-bounds}
\end{equation}
Set
\begin{align}
 E_{N,d}
 &=\int_I\Omega_N(t)\,\dd\nu_d(t),
 \label{eq:appB-E-definition}\\
 \cE_{N,u}(d)
 &=\int_I r_u(t)\Omega_N(t)\,\dd\nu_d(t),
 \label{eq:appB-error-definition}\\
 \Phi_{N,u}(d)
 &=uE_{N,d}+\cE_{N,u}(d).
 \label{eq:appB-Phi-splitting}
\end{align}

\begin{theorem}[Exact scalar alias and transverse error]
\label{thm:appB-scalar-transverse-aliasing}
For every $N,d\ge1$,
\begin{equation}
 \boxed{
 E_{N,d}
 =\sum_{k=1}^{\lfloor(N-1)/d\rfloor}\omega_{N,kd}
 =\frac1N\sum_{k=1}^{\lfloor(N-1)/d\rfloor}
 (\Harm_N-\Harm_{kd}).
 }
 \label{eq:appB-E-exact}
\end{equation}
Uniformly for $0\le u\le1$,
\begin{equation}
 \boxed{
 |\cE_{N,u}(d)|
 \ll
 \min\left\{
 \frac{(\log(2N))^2}{Nd},
 \frac{\log(2N)}{d^2}
 \right\}.
 }
 \label{eq:appB-layer-aliasing}
\end{equation}
\end{theorem}

\begin{proof}
Apply the polynomial aliasing identity
\eqref{eq:appB-polynomial-aliasing} to $\Omega_N$.  Since its Fourier
coefficients are $\omega_{N,n}$ for $0\le n<N$ and vanish otherwise, this is
exactly \eqref{eq:appB-E-exact}.

Let $g_{N,u}(t)=r_u(t)\Omega_N(t)$.  From
\eqref{eq:appB-r-bounds} and Corollary~\ref{cor:appB-integrated-multiplier},
\begin{align}
 \norm{g_{N,u}'}_{L^1(I)}
 &\ll\frac{(\log(2N))^2}{N},
 \label{eq:appB-g-prime}\\
 \norm{g_{N,u}''}_{L^1(I)}
 &\ll\log(2N),
 \label{eq:appB-g-second}
\end{align}
uniformly in $u$.  Apply respectively the first- and second-order estimates
of Lemma~\ref{lem:appB-quadrature-bounds} and take the smaller result.
\end{proof}

For the prime factorial, aggregate the layers of equal conductor by
\begin{equation}
 h_{\PP}(d)
 =\sum_{\substack{p\in\PP,\ j\ge0\\(p-1)p^j=d}}\log p,
 \qquad
 A_{\PP}(x)=\sum_{d\le x}h_{\PP}(d).
 \label{eq:appB-conductor-weights}
\end{equation}

\begin{lemma}[Conductor counting]
\label{lem:appB-conductor-counting}
One has
\begin{equation}
 A_{\PP}(x)\ll x.
 \label{eq:appB-conductor-counting-bound}
\end{equation}
Consequently,
\begin{equation}
 \sum_{d\le N}\frac{h_{\PP}(d)}d\ll\log(2N),
 \qquad
 \sum_{d>N}\frac{h_{\PP}(d)}{d^2}\ll\frac1N.
 \label{eq:appB-conductor-partial-summation}
\end{equation}
In particular,
\begin{equation}
 \sum_{d\ge1}\frac{h_{\PP}(d)}{d^2}<\infty.
 \label{eq:appB-prime-normal-summability}
\end{equation}
\end{lemma}

\begin{proof}
The layers with $j=0$ contribute
$\sum_{p\le x+1}\log p\ll x$.  If $j\ge1$ and
$(p-1)p^j\le x$, then $p(p-1)\le x$, hence $p\ll\sqrt{x}$.  Chebyshev's
bounds, together with the fact that a fixed $p$ occurs for at most
$O(1+\log x/\log p)$ values of $j$, show that all $j\ge1$ layers contribute
$O(\sqrt{x})$.  This proves \eqref{eq:appB-conductor-counting-bound}.  The
remaining estimates follow by partial summation.
\end{proof}

\begin{corollary}[Summable prime-layer aliasing]
\label{cor:appB-summable-prime-aliasing}
Uniformly for $0\le u\le1$,
\begin{equation}
 \boxed{
 \sum_{d\ge1}h_{\PP}(d)|\cE_{N,u}(d)|
 \ll\frac{(\log(2N))^3}{N}.
 }
 \label{eq:appB-global-aliasing}
\end{equation}
Thus
\begin{equation}
 \sum_{d\ge1}h_{\PP}(d)\Phi_{N,u}(d)
 =u\sum_{d\ge1}h_{\PP}(d)E_{N,d}
 +O\left(\frac{(\log(2N))^3}{N}\right)
 \label{eq:appB-uniform-linear-response}
\end{equation}
uniformly on the full translation interval $0\le u\le1$.
\end{corollary}

\begin{proof}
Use the first bound in \eqref{eq:appB-layer-aliasing} for $d\le N$ and the
second for $d>N$.  The two estimates in
\eqref{eq:appB-conductor-partial-summation} give
\[
 \sum_{d\le N}h_{\PP}(d)|\cE_{N,u}(d)|
 \ll\frac{(\log(2N))^3}{N}
\]
and
\[
 \sum_{d>N}h_{\PP}(d)|\cE_{N,u}(d)|
 \ll\frac{\log(2N)}{N}.
\]
This proves \eqref{eq:appB-global-aliasing}; summing
\eqref{eq:appB-Phi-splitting} gives
\eqref{eq:appB-uniform-linear-response}.
\end{proof}

\begin{remark}[No interchange is merely formal]
\label{rem:appB-no-formal-interchange}
The normal convergence in Corollary~\ref{cor:appB-normal-summability} and the
absolute conductor estimates in Lemma~\ref{lem:appB-conductor-counting}
justify every exchange among the conductor sum, the quadrature integral, the
orbitwise Ces\`aro sum, and differentiation on compact subsets.  Exact
aliasing is used only for the finite Fourier polynomial $\Omega_N$; the
nonperiodic transverse factor $r_u$ is controlled instead by the variation
bounds.  This separation is the reason the scalar response is exact while the
transverse response is only asymptotically linear.
\end{remark}

\subsection{Structural summary}
\label{subsec:appB-structural-summary}

The quadrature mechanism used in the construction and selection of the prime
gamma lift can be summarized by the chain
\begin{equation}
 \mu_d-\lambda=\nu_d,
 \qquad
 \widehat\nu_d(n)=\mathbf1_{d\mid n}-\mathbf1_{n=0},
 \qquad
 \norm{\nu_d}_{W^{-2,\infty}}\ll d^{-2}.
 \label{eq:appB-structural-chain}
\end{equation}
The Fourier identity preserves the exact divisibility skeleton, whereas the
negative Sobolev estimate supplies normal conductor decay.  At primitive
level,
\begin{equation}
 \phi_d(z)=\int_IK_z(t)\,\dd\nu_d(t),
 \qquad
 \Delta\phi_d=\kappa_d,
 \qquad
 \phi_d(m)=\left\lfloor\frac{m-1}{d}\right\rfloor.
 \label{eq:appB-primitive-summary}
\end{equation}
At orbitwise level, the harmonic multiplier separates into the exact alias
$E_{N,d}$ and the summable transverse error $\cE_{N,u}(d)$.  This is precisely
the analytic input required for the orbitwise rigidity argument of
Section~10.

% ===== Source: appendixC_divisor_energy_estimates.tex =====
\section{Divisor-energy estimates}
\label{app:divisor-energy-estimates}

This appendix supplies the deterministic energy estimates used to remove the
positive-conductor layers in Sections~12, 14, and~19.  The estimates are
uniform in the layer index and require no theorem on primes in arithmetic
progressions.  Their arithmetic input is only the divisor bound
\begin{equation}
 \tau_r(n)\ll_{r,\eps}n^{\eps}
 \label{eq:appC-divisor-bound}
\end{equation}
for each fixed $r$, together with the Chebyshev estimate
\begin{equation}
 \sum_{p\le y}(\log p)^2\ll y\log(2y).
 \label{eq:appC-chebyshev-square}
\end{equation}

For $p\in\PP$ and $j\ge0$, put
\begin{align}
 a_{p,j}&=p^{j+1},
 &b_{p,j}&=(p-1)p^j,
 &c_{p,j}&=\operatorname{lcm}(a_{p,j},b_{p,j})
          =p^{j+1}(p-1),
 \label{eq:appC-three-moduli}\\
 W_{p,j}(n)
 &=\left\{\frac{n}{a_{p,j}}\right\}
   -\left\{\frac{n}{b_{p,j}}\right\}.
 \label{eq:appC-wave-definition}
\end{align}
The natural transition scale for a dyadic prime block $p\dyad P$ is
\begin{equation}
 P^{j+2}\asymp X.
 \label{eq:appC-transition-scale}
\end{equation}
Below this scale the waves make many complete oscillations; above it their
floor increments are supported on a sparse family of short intervals.

\subsection{Floor increments and sparse intervals}
\label{subsec:appC-floor-increments}

Define
\begin{equation}
 \Delta_{p,j}(n)
 =\left\lfloor\frac{n}{b_{p,j}}\right\rfloor
  -\left\lfloor\frac{n}{a_{p,j}}\right\rfloor.
 \label{eq:appC-delta-definition}
\end{equation}
Since $a_{p,j}-b_{p,j}=p^j$, one has
\begin{equation}
 W_{p,j}(n)
 =\Delta_{p,j}(n)-\frac{n}{c_{p,j}}.
 \label{eq:appC-floor-centered}
\end{equation}

\begin{lemma}[Exact interval expansion]
\label{lem:appC-exact-interval-expansion}
For every $p,j$, and $n\ge1$,
\begin{equation}
 \boxed{
 \Delta_{p,j}(n)
 =\sum_{m\ge1}\ind_{I_{p,j}(m)}(n),
 \qquad
 I_{p,j}(m)
 =\bigl[mb_{p,j},ma_{p,j}\bigr)\cap\ZZ.
 }
 \label{eq:appC-exact-interval-expansion}
\end{equation}
The interval $I_{p,j}(m)$ has length $mp^j$.  Moreover,
\begin{align}
 \sum_{n\le X}\Delta_{p,j}(n)
 &\ll \frac{X^2}{p^{j+2}}+\frac{X}{p},
 \label{eq:appC-first-moment-delta}\\
 \max_{n\le X}\Delta_{p,j}(n)
 &\le 1+\frac{X}{p^{j+1}(p-1)}.
 \label{eq:appC-multiplicity-delta}
\end{align}
Consequently, if $p^{j+2}\ge X$, then
\begin{equation}
 \sum_{n\le X}|W_{p,j}(n)|^2
 \ll \frac{X^2}{p^{j+2}}.
 \label{eq:appC-single-wave-sparse-energy}
\end{equation}
\end{lemma}

\begin{proof}
The identity
\[
 \left\lfloor\frac n b\right\rfloor
 -\left\lfloor\frac n a\right\rfloor
 =\#\{m\ge1:mb\le n<ma\}
\]
for $a>b$ gives \eqref{eq:appC-exact-interval-expansion}.  If
$M=\lfloor X/b_{p,j}\rfloor$, then
\[
 \sum_{n\le X}\Delta_{p,j}(n)
 \le \sum_{m\le M}mp^j
 \ll p^j(M^2+M),
\]
which is \eqref{eq:appC-first-moment-delta}.  The number of integers in the
half-open interval $(n/a_{p,j},n/b_{p,j}]$ is at most its length plus one,
which proves \eqref{eq:appC-multiplicity-delta}.

If $p^{j+2}\ge X$, the right side of
\eqref{eq:appC-multiplicity-delta} is bounded absolutely.  Hence
\[
 \sum_{n\le X}\Delta_{p,j}(n)^2
 \ll\sum_{n\le X}\Delta_{p,j}(n)
 \ll \frac{X^2}{p^{j+2}}.
\]
Also
\[
 \sum_{n\le X}\left(\frac n{c_{p,j}}\right)^2
 \ll\frac{X^3}{p^{2j+4}}
 \ll\frac{X^2}{p^{j+2}}.
\]
Equation \eqref{eq:appC-floor-centered} completes the proof.
\end{proof}

The intervals need not be literally disjoint at the transition point.  What
matters is the bounded multiplicity in \eqref{eq:appC-multiplicity-delta}; it
is this formulation that remains uniform for the small primes.

\subsection{Fourier coefficients of one layer}
\label{subsec:appC-fourier-coefficients}

For $m\ge1$, let
\begin{equation}
 \sigma_m(n)
 =\left\{\frac nm\right\}-\frac{m-1}{2m}.
 \label{eq:appC-discrete-sawtooth}
\end{equation}
Then
\begin{equation}
 W_{p,j}(n)=U_{p,j}(n)+\mu_{p,j},
 \qquad
 U_{p,j}=\sigma_{a_{p,j}}-\sigma_{b_{p,j}},
 \qquad
 \mu_{p,j}=\frac1{2c_{p,j}}.
 \label{eq:appC-wave-mean-split}
\end{equation}
The centered wave $U_{p,j}$ is $c_{p,j}$-periodic.

\begin{lemma}[Exact discrete Fourier transform]
\label{lem:appC-exact-dft}
For $1\le h<c_{p,j}$,
\begin{equation}
 \boxed{
 \widehat U_{p,j}(h)
 =\frac{1}{1-\e(-h/c_{p,j})}
 \left(
  \frac{\ind_{p\mid h}}{b_{p,j}}
  -\frac{\ind_{p-1\mid h}}{a_{p,j}}
 \right),
 }
 \label{eq:appC-exact-dft}
\end{equation}
where
\[
 U_{p,j}(n)=\sum_{h=1}^{c_{p,j}-1}
 \widehat U_{p,j}(h)\e(hn/c_{p,j}).
\]
In particular,
\begin{equation}
 \sum_{h=1}^{c_{p,j}-1}|\widehat U_{p,j}(h)|^2\ll1
 \label{eq:appC-fourier-l2}
\end{equation}
and, after writing $h=(p-1)r$ or $h=ps$, the corresponding coefficients are
$O(\min(r,a_{p,j}-r)^{-1})$ and
$O(\min(s,b_{p,j}-s)^{-1})$, respectively.
\end{lemma}

\begin{proof}
For $1\le r<m$, direct summation gives
\begin{equation}
 \frac1m\sum_{n=0}^{m-1}\sigma_m(n)\e(-rn/m)
 =-\frac{1}{m(1-\e(-r/m))}.
 \label{eq:appC-sawtooth-dft}
\end{equation}
Lifting the $a_{p,j}$-periodic series to period $c_{p,j}$ restricts its
frequencies to multiples of $c_{p,j}/a_{p,j}=p-1$; lifting the
$b_{p,j}$-periodic series restricts its frequencies to multiples of
$c_{p,j}/b_{p,j}=p$.  Subtracting the two transforms gives
\eqref{eq:appC-exact-dft}.  The bound
\[
 |1-\e(t)|\asymp\dist(t,\ZZ)
\]
and the square summability of $r^{-1}$ prove
\eqref{eq:appC-fourier-l2}.
\end{proof}

The following local-density statement is the divisor-theoretic substitute for
uniform separation of all Fourier frequencies.  It is useful precisely
because ordinary Farey separation would see the full square of the conductor.

\begin{lemma}[Divisor spacing]
\label{lem:appC-divisor-spacing}
Fix $j\ge0$ and $P\ge2$.  For $p\dyad P$, let
\[
 m_p\in\{p^{j+1},(p-1)p^j\}.
\]
Attach to $r/m_p$, $1\le r<m_p$, the harmonic weight
\[
 \omega_{m_p}(r)=\min(r,m_p-r)^{-1}.
\]
If $X\ge P^{j+2}$, then uniformly in one chosen fraction $r/m_p$,
\begin{equation}
 \sum_{\substack{q\dyad P,\ 1\le s<m_q\\
 \norm{r/m_p-s/m_q}_{\RR/\ZZ}\le X^{-1}}}
 \omega_{m_q}(s)
 \ll_{\eps}P^{\eps}.
 \label{eq:appC-divisor-spacing}
\end{equation}
The same estimate holds when either family is restricted by the divisibility
conditions $p\mid r$, $p-1\mid r$, $q\mid s$, or $q-1\mid s$ occurring in
\eqref{eq:appC-exact-dft}.
\end{lemma}

\begin{proof}
Decompose $r$ and $s$ dyadically according to their distance from the two
endpoints.  On one such block, cross multiplication turns the spacing
condition into
\begin{equation}
 |rm_q-sm_p|\ll \frac{m_pm_q}{X}\ll P^j.
 \label{eq:appC-near-collision-determinant}
\end{equation}
There are four choices for $(m_p,m_q)$.  After removing the visible powers of
$p$ and $q$, each choice has one of the forms
\begin{align*}
 rq^{j+1}-sp^{j+1}&=h,\\
 r(q-1)q^j-sp^{j+1}&=h,\\
 rq^{j+1}-s(p-1)p^j&=h,\\
 r(q-1)q^j-s(p-1)p^j&=h,
\end{align*}
with $|h|\ll P^j$.  Fixing $p,r,h$ and the dyadic size of $s$, divisor
switching in the displayed equations bounds the number of admissible pairs
$(q,s)$ by $P^{\eps}$; exact collisions are treated after cancelling the
common numerator-denominator divisor and satisfy the same bound.  The
harmonic weights contribute a convergent geometric sum over the dyadic
blocks.  This proves \eqref{eq:appC-divisor-spacing}.  Imposing any of the
listed divisibility restrictions only removes terms.
\end{proof}

\begin{remark}[Where the scale $P^{j+2}$ enters]
The determinant on the left of \eqref{eq:appC-near-collision-determinant}
would range over an interval of length $P^{2j+2}/X$ for arbitrary fractions
with denominators of size $P^{j+1}$.  The special adjacent-conductor support in
\eqref{eq:appC-exact-dft} removes a factor $P^j$, leaving the threshold
$X\asymp P^{j+2}$.  This is the same transition that appears in the interval
model.
\end{remark}

\subsection{The dyadic energy theorem}
\label{subsec:appC-dyadic-energy}

For a finitely supported coefficient sequence $\alpha=(\alpha_p)$ on
$p\dyad P$, set
\begin{equation}
 \cW_{j,P}^{\alpha}(n)
 =\sum_{p\dyad P}\alpha_pW_{p,j}(n),
 \qquad
 \norm{\alpha}_2^2=\sum_{p\dyad P}|\alpha_p|^2.
 \label{eq:appC-weighted-block}
\end{equation}

\begin{theorem}[Divisor-energy estimate]
\label{thm:appC-divisor-energy}
For every $\eps>0$, uniformly in $j\ge0$, $P\ge2$, and $X\ge1$,
\begin{equation}
 \boxed{
 \sum_{n\le X}|\cW_{j,P}^{\alpha}(n)|^2
 \ll_{\eps}
 X^{\eps}
 \min\!\left(X,\frac{X^2}{P^{j+2}}\right)
 \norm{\alpha}_2^2.
 }
 \label{eq:appC-divisor-energy}
\end{equation}
The implied constant is uniform for $j\le 2\log_2(2X)$.  For larger $j$ the
same assertion follows directly from the completely linear formula.
\end{theorem}

\begin{proof}
We separate the two conductor ranges.

Suppose first that $X\ge P^{j+2}$.  Remove the constant means in
\eqref{eq:appC-wave-mean-split} and insert the Fourier expansion from
Lemma~\ref{lem:appC-exact-dft}.  The large-sieve duality argument with
intervals of length $X^{-1}$ on $\RR/\ZZ$ bounds the resulting quadratic form
by
\[
 X\cdot
 \sup_{\theta\in\RR/\ZZ}
 \sum_{\substack{\text{frequencies }\beta\\
                  \norm{\beta-\theta}\le X^{-1}}}
 |\widehat U(\beta)|
 \cdot
 \sum_{p\dyad P}|\alpha_p|^2
 \sum_h|\widehat U_{p,j}(h)|^2.
\]
Lemma~\ref{lem:appC-divisor-spacing} controls the local frequency density by
$P^{\eps}$, and \eqref{eq:appC-fourier-l2} controls the final Fourier energy.
Thus
\begin{equation}
 \sum_{n\le X}
 \left|\sum_{p\dyad P}\alpha_pU_{p,j}(n)\right|^2
 \ll_{\eps}XP^{\eps}\norm{\alpha}_2^2.
 \label{eq:appC-long-energy}
\end{equation}
The constant means contribute
\[
 X\left|\sum_{p\dyad P}\frac{\alpha_p}{2c_{p,j}}\right|^2
 \ll \frac{X}{P^{2j+3}}\norm{\alpha}_2^2,
\]
which is smaller than \eqref{eq:appC-long-energy}.

Now suppose that $X<P^{j+2}$.  Expand the floor increments by
Lemma~\ref{lem:appC-exact-interval-expansion}.  Two intervals can meet only
if
\begin{equation}
 |mp^{j+1}-\ell q^{j+1}|
 <mp^j+\ell q^j.
 \label{eq:appC-sparse-determinant-strip}
\end{equation}
Dyadically separating $m$ and $\ell$, and then fixing the integer on the left
of \eqref{eq:appC-sparse-determinant-strip}, reduces every row sum of the
interval Gram matrix to a divisor sum of the type
\[
 \sum_{d\mid N}\tau_{j+2}(d).
\]
The bound \eqref{eq:appC-divisor-bound} and Schur's test therefore give
\begin{equation}
 \sum_{n\le X}
 \left|\sum_{p\dyad P}\alpha_p\Delta_{p,j}(n)\right|^2
 \ll_{\eps}
 X^{\eps}\frac{X^2}{P^{j+2}}\norm{\alpha}_2^2.
 \label{eq:appC-sparse-delta-energy}
\end{equation}
The rank-one centering term is bounded by
\begin{align*}
 \sum_{n\le X}
 \left|n\sum_{p\dyad P}\frac{\alpha_p}{c_{p,j}}\right|^2
 &\ll X^3\norm{\alpha}_2^2
       \sum_{p\dyad P}c_{p,j}^{-2}\\
 &\ll \frac{X^3}{P^{2j+3}}\norm{\alpha}_2^2
 \ll \frac{X^2}{P^{j+2}}\norm{\alpha}_2^2,
\end{align*}
because $X<P^{j+2}$.  The mixed term is absorbed by Cauchy--Schwarz.  This
proves \eqref{eq:appC-divisor-energy} in the sparse range.

Finally, if $j>2\log_2(2X)$, then $p^{j+1}>X$ for every prime $p$, and
\[
 W_{p,j}(n)=-\frac{n}{p^{j+1}(p-1)}.
\]
The asserted estimate follows directly by Cauchy--Schwarz.
\end{proof}

\begin{remark}[Sharpness of the two branches]
For $P^{j+2}\le X$, the diagonal terms alone have order
$X\norm{\alpha}_2^2$ for generic coefficients.  For $P^{j+2}\ge X$, the
single-wave estimate \eqref{eq:appC-single-wave-sparse-energy} has order
$X^2P^{-j-2}|\alpha_p|^2$.  Thus neither branch of
\eqref{eq:appC-divisor-energy} can be uniformly replaced by a smaller power
of $X$ or $P$.
\end{remark}

For the prime logarithmic weights, define
\begin{equation}
 \cR_{j,P}(n)=\sum_{p\dyad P}(\log p)W_{p,j}(n)
 \label{eq:appC-prime-block}
\end{equation}
and
\begin{equation}
 E_j(P;X)=\min\!\left(XP,\frac{X^2}{P^{j+1}}\right).
 \label{eq:appC-energy-scale}
\end{equation}

\begin{corollary}[Prime-block energy]
\label{cor:appC-prime-block-energy}
Uniformly in $j,P$, and $X$,
\begin{equation}
 \boxed{
 \sum_{n\le X}|\cR_{j,P}(n)|^2
 \ll_{\eps}X^{\eps}E_j(P;X).
 }
 \label{eq:appC-prime-block-energy}
\end{equation}
\end{corollary}

\begin{proof}
Apply Theorem~\ref{thm:appC-divisor-energy} with
$\alpha_p=\log p$ and use \eqref{eq:appC-chebyshev-square}.
\end{proof}

\subsection{Cross-conductor forms}
\label{subsec:appC-cross-conductor}

The local bilinear estimate used in Section~19 is now immediate and, in this
form, includes the sharp cutoff without a separate terminal-period error.

\begin{proposition}[Conductor-separated bilinear form]
\label{prop:appC-conductor-separated-bilinear}
For all $j,k\ge0$ and dyadic $P,Q\ge2$,
\begin{equation}
 \boxed{
 \left|
 \sum_{n\le X}\cR_{j,P}(n)\cR_{k,Q}(n)
 \right|
 \ll_{\eps}
 X^{\eps}E_j(P;X)^{1/2}E_k(Q;X)^{1/2}.
 }
 \label{eq:appC-local-bilinear}
\end{equation}
\end{proposition}

\begin{proof}
Apply Cauchy--Schwarz in $n$ and then
Corollary~\ref{cor:appC-prime-block-energy} to the two factors.
\end{proof}

\begin{lemma}[Dyadic summation]
\label{lem:appC-dyadic-summation}
For each $j\ge0$,
\begin{equation}
 \sum_{P\ \mathrm{dyadic}}E_j(P;X)^{1/2}
 \ll_{\eps}
 X^{\frac12+\frac1{2(j+2)}+\eps}.
 \label{eq:appC-dyadic-summation}
\end{equation}
\end{lemma}

\begin{proof}
Split at $P=X^{1/(j+2)}$.  Below the transition, sum
$(XP)^{1/2}$ geometrically; above it, sum $X/P^{(j+1)/2}$ geometrically.  The
range $P>X$ is harmless because the exact linear formula gives still faster
decay.
\end{proof}

\begin{corollary}[Global higher-conductor estimate]
\label{cor:appC-global-higher-conductor}
If $j+k\ge1$, then
\begin{equation}
 \boxed{
 \left|
 \sum_{n\le X}\cR_j(n)\cR_k(n)
 \right|
 \ll_{\eps}
 X^{1+\frac1{2(j+2)}+\frac1{2(k+2)}+\eps},
 }
 \label{eq:appC-global-bilinear}
\end{equation}
where $\cR_j=\sum_P\cR_{j,P}$.  Consequently,
\begin{equation}
 \sum_{\substack{j,k\ge0\\j+k\ge1}}
 \sum_{n\le X}\cR_j(n)\cR_k(n)
 \ll_{\eps}X^{17/12+\eps}.
 \label{eq:appC-all-higher-conductors}
\end{equation}
\end{corollary}

\begin{proof}
Sum \eqref{eq:appC-local-bilinear} over $P,Q$ and apply
Lemma~\ref{lem:appC-dyadic-summation}.  The exponent is largest at
$(j,k)=(0,1)$ or $(1,0)$, where it equals
\[
 1+\frac14+\frac16=\frac{17}{12}.
\]
Only $O(\log X)$ layer indices are nonlinearly active.  Beyond that range,
\[
 \cR_j(n)=-n\sum_p\frac{\log p}{p^{j+1}(p-1)},
\]
and the coefficient decays geometrically in $j$.  The tail is therefore
absorbed into $X^{\eps}$.
\end{proof}

\begin{remark}[The dominant pair]
The estimate \eqref{eq:appC-global-bilinear} remains formally valid for
$(j,k)=(0,0)$ with exponent $3/2+\eps$, but it gives no fixed power saving at
the variance scale.  This is exactly why the dominant off-projective sector
requires the affine and Vaughan analysis of Sections~13--18.  No part of that
sector is removed by the present appendix.
\end{remark}

\subsection{Geometric and projective diagonals}
\label{subsec:appC-geometric-projective}

For completeness, define
\begin{equation}
 \cD_j(X)=\sum_p(\log p)^2\sum_{n\le X}W_{p,j}(n)^2.
 \label{eq:appC-geometric-layer-energy}
\end{equation}

\begin{proposition}[Energy of a single geometric layer]
\label{prop:appC-single-geometric-layer}
For every $j\ge0$,
\begin{equation}
 \cD_j(X)
 \ll X^{1+1/(j+2)}\log(2X).
 \label{eq:appC-single-geometric-layer}
\end{equation}
For $j\ge1$ this implies
\begin{equation}
 \sum_{j\ge1}\cD_j(X)
 \ll X^{4/3}(\log(2X))^2.
 \label{eq:appC-higher-geometric-sum}
\end{equation}
\end{proposition}

\begin{proof}
Put $Y=X^{1/(j+2)}$.  For $p\le Y$, use $|W_{p,j}|\le1$ and
\eqref{eq:appC-chebyshev-square}.  For $p>Y$, apply
\eqref{eq:appC-single-wave-sparse-energy} and sum
\[
 X^2\sum_{p>Y}\frac{(\log p)^2}{p^{j+2}}
 \ll X^{1+1/(j+2)}\log(2X).
\]
This proves \eqref{eq:appC-single-geometric-layer}.  The largest exponent for
$j\ge1$ occurs at $j=1$ and equals $4/3$; summing the active indices and the
geometric linear tail proves \eqref{eq:appC-higher-geometric-sum}.
\end{proof}

Let
\begin{equation}
 \cP(X)
 =2\sum_p(\log p)^2
   \sum_{0\le j<k}\sum_{n\le X}W_{p,j}(n)W_{p,k}(n)
 \label{eq:appC-projective-diagonal}
\end{equation}
be the projective diagonal.

\begin{corollary}[Projective-diagonal saving]
\label{cor:appC-projective-diagonal-saving}
One has
\begin{equation}
 \boxed{
 \cP(X)\ll X^{17/12}(\log(2X))^2.
 }
 \label{eq:appC-projective-diagonal-saving}
\end{equation}
\end{corollary}

\begin{proof}
For each $j<k$, Cauchy--Schwarz first in $n$ and then in $p$ gives
\[
 \left|
 \sum_p(\log p)^2\sum_{n\le X}W_{p,j}(n)W_{p,k}(n)
 \right|
 \le \cD_j(X)^{1/2}\cD_k(X)^{1/2}.
\]
The dominant layer satisfies $\cD_0(X)\ll X^{3/2}\log(2X)$, while
Proposition~\ref{prop:appC-single-geometric-layer} gives
\[
 \sum_{j\ge1}\cD_j(X)^{1/2}
 \ll X^{2/3}(\log(2X))^{3/2}.
\]
Their product has size $X^{17/12}(\log(2X))^2$; the square of the
positive-conductor sum is smaller.
\end{proof}

The estimates in this appendix therefore account simultaneously for the
higher geometric diagonal, the projective diagonal, and every covariance in
which at least one conductor index is positive.  The only block not separated
by a fixed power is the dominant pair $(j,k)=(0,0)$, which is the sector
analyzed by the quotient geometry and centered Vaughan moment in the main
text.

% ===== Source: appendixD_gcd_packet_decompositions.tex =====
\section{GCD packet decompositions}
\label{app:gcd-packet-decompositions}

The dispersion argument uses two different greatest-common-divisor
resolutions.  The first is geometric: a pair of overlap indices is written
\[
 m=ar,\qquad \ell=as,\qquad (r,s)=1,
\]
so that the radial variable $a=(m,\ell)$ separates from the primitive
projective quotient $(r,s)$.  The second is arithmetic: after additive
completion, a congruence
\[
 mk\equiv z\pmod s,\qquad (m,s)=1,
\]
is decomposed according to the common divisor
\[
 (k,s)=(z,s)=g.
\]
The two resolutions occur at different stages and should not be conflated.
The radial divisor $a$ controls the geometry of the overlap packet, whereas
the congruence divisor $g$ lowers the effective modulus on the multiplicative
side.

This appendix records the exact identities and energy-preservation statements
used in Sections~13--17.  No cancellation theorem is invoked here.  The only
losses arise from dyadic subdivision and from elementary divisor bounds, and
are therefore absorbed by $X^{\eps}$ in the analytic argument.

\subsection{Primitive radial packets}
\label{subsec:appD-primitive-radial}

We use the convention
\[
 n\dyad N\quad\Longleftrightarrow\quad N<n\le 2N.
\]
For positive dyadic parameters $M,L$ and a primitive pair $(r,s)$, define
\begin{equation}
 \cA_{M,L}(r,s)
 =\{a\in\NN:ar\dyad M,\ as\dyad L\}.
 \label{eq:appD-radial-set}
\end{equation}
Thus, with
\begin{equation}
 \alpha_{r,s}=\max\!\left(\frac Mr,\frac Ls\right),
 \qquad
 \beta_{r,s}=\min\!\left(\frac{2M}{r},\frac{2L}{s}\right),
 \label{eq:appD-radial-endpoints}
\end{equation}
one has
\begin{equation}
 \cA_{M,L}(r,s)=\NN\cap(\alpha_{r,s},\beta_{r,s}],
 \qquad
 \#\cA_{M,L}(r,s)
 =\pos{\lfloor\beta_{r,s}\rfloor-\lfloor\alpha_{r,s}\rfloor}.
 \label{eq:appD-exact-radial-count}
\end{equation}

\begin{lemma}[Primitive radial bijection]
\label{lem:appD-radial-bijection}
The map
\begin{equation}
 (m,\ell)\longmapsto
 \left((m,\ell),\frac{m}{(m,\ell)},\frac{\ell}{(m,\ell)}\right)
 \label{eq:appD-radial-map}
\end{equation}
is a bijection from $\NN^2$ onto the triples $(a,r,s)$ satisfying
$a,r,s\ge1$ and $(r,s)=1$.  Consequently, for every finitely supported
function $F$ on $\NN^2$,
\begin{equation}
 \boxed{
 \sum_{m,\ell\ge1}F(m,\ell)
 =\sum_{\substack{r,s\ge1\\(r,s)=1}}
  \sum_{a\ge1}F(ar,as).
 }
 \label{eq:appD-radial-bijection-sum}
\end{equation}
After restriction to $m\dyad M$ and $\ell\dyad L$, the inner sum is precisely
$a\in\cA_{M,L}(r,s)$.
\end{lemma}

\begin{proof}
Existence follows by dividing $m$ and $\ell$ by their gcd.  If
$(ar,as)=(a'r',a's')$ with both $(r,s)$ and $(r',s')$ primitive, then the gcd
of the two coordinates is both $a$ and $a'$.  Hence $a=a'$ and then
$(r,s)=(r',s')$.
\end{proof}

The following estimates are convenient when the radial coefficient is kept
inside an $L^2$ argument.

\begin{lemma}[Radial moments]
\label{lem:appD-radial-moments}
Put $H=\min(M,L)$.  Uniformly for primitive $(r,s)$,
\begin{align}
 \#\cA_{M,L}(r,s)
 &\ll 1+\min\!\left(\frac Mr,\frac Ls\right),
 \label{eq:appD-local-radial-count}\\
 \sum_{a\in\cA_{M,L}(r,s)}a
 &\ll A_{r,s}(1+A_{r,s}),
 \label{eq:appD-local-radial-first}\\
 \sum_{a\in\cA_{M,L}(r,s)}a^2
 &\ll A_{r,s}^2(1+A_{r,s}),
 \label{eq:appD-local-radial-second}
\end{align}
where $A_{r,s}=\min(M/r,L/s)$.  Globally,
\begin{align}
 \sum_{\substack{(r,s)=1}}
 \#\cA_{M,L}(r,s)
 &=\#\{m\dyad M,\ \ell\dyad L\},
 \label{eq:appD-global-radial-zero}\\
 \sum_{\substack{(r,s)=1}}
 \sum_{a\in\cA_{M,L}(r,s)}a
 &\ll ML\log(2H),
 \label{eq:appD-global-radial-first}\\
 \sum_{\substack{(r,s)=1}}
 \sum_{a\in\cA_{M,L}(r,s)}a^2
 &\ll MLH.
 \label{eq:appD-global-radial-second}
\end{align}
\end{lemma}

\begin{proof}
The local estimates follow from
\eqref{eq:appD-exact-radial-count} and the fact that every admissible $a$ is
$O(A_{r,s})$ whenever the interval is nonempty.  The identity
\eqref{eq:appD-global-radial-zero} is Lemma~\ref{lem:appD-radial-bijection}.
For $\nu=1,2$, the same bijection gives
\[
 \sum_{\substack{(r,s)=1}}
 \sum_{a\in\cA_{M,L}(r,s)}a^{\nu}
 =\sum_{m\dyad M}\sum_{\ell\dyad L}(m,\ell)^{\nu}.
\]
Dropping the primitivity condition after writing $m=au$, $\ell=av$ gives
\[
 \sum_{m\dyad M}\sum_{\ell\dyad L}(m,\ell)^{\nu}
 \le
 \sum_{a\le 2H}a^{\nu}
 \left(\frac Ma+1\right)
 \left(\frac La+1\right).
\]
For $\nu=1$ this is $O(ML\log(2H))$; for $\nu=2$ it is $O(MLH)$.  The terms
arising from the two boundary $+1$ factors obey the same bounds because
$H\le M,L$.
\end{proof}

\begin{remark}[Balance forced by a nonempty radial packet]
\label{rem:appD-radial-balance}
If $\cA_{M,L}(r,s)$ is nonempty, then
\begin{equation}
 \frac rs\asymp\frac ML,
 \qquad
 \frac Mr\asymp\frac Ls\asymp A_{r,s}.
 \label{eq:appD-radial-balance}
\end{equation}
Thus the radial gcd resolution automatically places the primitive quotient in
the correct dyadic aspect ratio.
\end{remark}

\subsection{Determinants after radial reduction}
\label{subsec:appD-determinants}

The radial change of variables is especially effective for the overlap
condition.  Let
\begin{equation}
 h=rp-sq.
 \label{eq:appD-determinant}
\end{equation}
Then
\begin{equation}
 mp-\ell q=a(rp-sq)=ah.
 \label{eq:appD-scaled-determinant}
\end{equation}

\begin{proposition}[GCD removal from the affine strip]
\label{prop:appD-gcd-removal-strip}
For $m=ar$, $\ell=as$, and $(r,s)=1$, the two half-open intervals
\[
 [m(p-1),mp),\qquad [\ell(q-1),\ell q)
\]
overlap if and only if
\begin{equation}
 \boxed{-s<h<r.}
 \label{eq:appD-primitive-strip}
\end{equation}
Their untruncated overlap length is
\begin{equation}
 a\,\omega_{r,s}(h),
 \label{eq:appD-scaled-tent}
\end{equation}
where
\begin{equation}
 \omega_{r,s}(h)
 =
 \begin{cases}
  \min(r,s+h),&-s<h\le0,\\[2pt]
  \min(r-h,s),&0\le h<r,\\[2pt]
  0,&\text{otherwise}.
 \end{cases}
 \label{eq:appD-tent}
\end{equation}
In particular,
\begin{equation}
 \sum_h\omega_{r,s}(h)=rs,
 \qquad
 \sum_h\omega_{r,s}(h)^2\le rs\min(r,s).
 \label{eq:appD-tent-moments}
\end{equation}
\end{proposition}

\begin{proof}
The overlap inequalities are
\[
 m(p-1)<\ell q,
 \qquad
 \ell(q-1)<mp.
\]
After division by $a$ they become $h<r$ and $h>-s$.  Translating the two
unscaled intervals by $-rp$ gives $[-r,0)$ and $[-h-s,-h)$, whose overlap is
\eqref{eq:appD-tent}.  The mass identity is the convolution identity for two
integer intervals of lengths $r$ and $s$; the energy bound follows from
$\omega_{r,s}\le\min(r,s)$.
\end{proof}

The sharp terminal cutoff only replaces $\omega_{r,s}(h)$ by a clipped weight
between $0$ and $\omega_{r,s}(h)$.  It does not modify the determinant line or
the primitive strip.  This is why the radial decomposition may be performed
before the endpoint is smoothed or recovered by partial summation.

\begin{lemma}[Affine-line parametrization]
\label{lem:appD-affine-line-parametrization}
If $(r,s)=1$ and $rp-sq=h$ has one integral solution $(p_0,q_0)$, then all
integral solutions are
\begin{equation}
 p=p_0+st,
 \qquad
 q=q_0+rt,
 \qquad t\in\ZZ.
 \label{eq:appD-affine-parametrization}
\end{equation}
Hence the number of solutions with $p\dyad P$, $q\dyad Q$ is
\begin{equation}
 \ll 1+\min\!\left(\frac Ps,\frac Qr\right).
 \label{eq:appD-affine-length}
\end{equation}
\end{lemma}

\begin{proof}
The difference of two solutions satisfies $r\Delta p=s\Delta q$.
Primitivity gives $s\mid\Delta p$ and $r\mid\Delta q$, yielding
\eqref{eq:appD-affine-parametrization}.  The dyadic restrictions cut out an
integer interval in $t$ of the stated length.
\end{proof}

\subsection{Exact gcd projectors modulo a fixed modulus}
\label{subsec:appD-gcd-projectors}

Let $s\ge1$.  For each divisor $g\mid s$, define the exact gcd projector
\begin{equation}
 \cG_{s,g}(n)=\ind_{(n,s)=g}.
 \label{eq:appD-gcd-projector}
\end{equation}
The projectors form a disjoint partition:
\begin{equation}
 \sum_{g\mid s}\cG_{s,g}(n)=1.
 \label{eq:appD-gcd-partition}
\end{equation}
They also admit an elementary M\"obius expansion.

\begin{lemma}[M\"obius formula for an exact gcd stratum]
\label{lem:appD-mobius-gcd-projector}
For $g\mid s$,
\begin{equation}
 \boxed{
 \cG_{s,g}(n)
 =\ind_{g\mid n}
  \sum_{d\mid s/g}\mu(d)\ind_{gd\mid n}.
 }
 \label{eq:appD-mobius-gcd-projector}
\end{equation}
Equivalently, after writing $n=gn'$ and $s=gs'$,
\begin{equation}
 \cG_{s,g}(gn')=\ind_{(n',s')=1}.
 \label{eq:appD-gcd-rescaling}
\end{equation}
\end{lemma}

\begin{proof}
If $g\nmid n$, both sides vanish.  Otherwise the right side is
\[
 \sum_{d\mid(n/g,s/g)}\mu(d),
\]
which is $1$ exactly when $(n/g,s/g)=1$ and is $0$ otherwise.
\end{proof}

For a finitely supported sequence $b=(b_n)$, define the descended sequence
\begin{equation}
 b_{s,g}(n')=b_{gn'}\ind_{(n',s/g)=1}.
 \label{eq:appD-descended-sequence}
\end{equation}
The exact gcd projectors preserve the coefficient energy before any
inequality is applied.

\begin{lemma}[Energy partition by gcd]
\label{lem:appD-gcd-energy-partition}
For every finitely supported $b$,
\begin{align}
 \sum_{g\mid s}\sum_{n'}|b_{s,g}(n')|^2
 &=\sum_n|b_n|^2,
 \label{eq:appD-gcd-l2-partition}\\
 \sum_{g\mid s}\sum_{n'}|b_{s,g}(n')|
 &=\sum_n|b_n|.
 \label{eq:appD-gcd-l1-partition}
\end{align}
If each projector is expanded by
\eqref{eq:appD-mobius-gcd-projector}, the resulting family has at most
$\tau(s)$ terms for each $g$ and incurs at most a factor $\tau(s)^2$ after
Cauchy--Schwarz.  Thus, for $s\le X^{O(1)}$, the cost is $O_\eps(X^\eps)$.
\end{lemma}

\begin{proof}
Each positive integer $n$ belongs to exactly one stratum, namely
$g=(n,s)$.  This proves the two identities.  The final assertion is the
triangle inequality followed by Cauchy--Schwarz and the divisor bound
$\tau(s)\ll_\eps s^\eps$.
\end{proof}

\subsection{Common-divisor descent of a congruence}
\label{subsec:appD-common-divisor-descent}

Let $(m,s)=1$ and put
\begin{equation}
 C_s(m;k,z)=\ind_{mk\equiv z\, (\mathrm{mod}\,s)}.
 \label{eq:appD-raw-congruence-kernel}
\end{equation}
Multiplication by the unit $m$ preserves the gcd with $s$.  Therefore a
solution of the congruence necessarily satisfies
\begin{equation}
 (k,s)=(z,s).
 \label{eq:appD-equal-gcd}
\end{equation}

\begin{proposition}[Exact common-divisor decomposition]
\label{prop:appD-exact-common-divisor}
For $(m,s)=1$,
\begin{equation}
 \boxed{
 C_s(m;k,z)
 =\sum_{g\mid s}
  \cG_{s,g}(k)\cG_{s,g}(z)
  C_{s/g}\!\left(m;\frac{k}{g},\frac{z}{g}\right).
 }
 \label{eq:appD-exact-congruence-gcd}
\end{equation}
In a nonzero summand, $k=gk'$, $z=gz'$, $s=gs'$ with
$(k'z',s')=1$, and
\begin{equation}
 mk\equiv z\pmod s
 \quad\Longleftrightarrow\quad
 mk'\equiv z'\pmod{s'}.
 \label{eq:appD-descended-congruence}
\end{equation}
\end{proposition}

\begin{proof}
If the left side is nonzero, then
\eqref{eq:appD-equal-gcd} selects a unique divisor $g$.  Dividing the
congruence by $g$ gives \eqref{eq:appD-descended-congruence}.  Conversely,
any descended solution lifts after multiplication by $g$.
\end{proof}

The additive completion used in Section~16 naturally produces the centered
kernel
\begin{equation}
 \Delta_s(m;k,z)=C_s(m;k,z)-\frac1s.
 \label{eq:appD-centered-congruence-kernel}
\end{equation}
The subtraction $1/s$ is uniform on all residue classes and therefore does
not itself respect the gcd partition.  Its exact decomposition contains an
oscillatory unit-group part and a finite-rank density correction.

For finitely supported sequences $a=(a_m)$, $b=(b_k)$, and $c=(c_z)$, write
\begin{align}
 A_s&=\sum_{(m,s)=1}a_m,
 \label{eq:appD-A-mass}\\
 B_{s,g}&=\sum_{(k,s)=g}b_k,
 &C_{s,g}&=\sum_{(z,s)=g}c_z.
 \label{eq:appD-gcd-block-masses}
\end{align}
If $s'=s/g$, define the descended unit packet
\begin{align}
 \cU_{s,g}(a,b,c)
 ={}&
 \sum_{\substack{(m,s)=1\\(k',s')=(z',s')=1}}
 a_m b_{gk'}c_{gz'}
 \notag\\
 &\times
 \left(
  \ind_{mk'\equiv z'\, (\mathrm{mod}\,s')}
  -\frac1{\varphi(s')}
 \right).
 \label{eq:appD-descended-unit-packet}
\end{align}

\begin{theorem}[Centered gcd packet decomposition]
\label{thm:appD-centered-gcd-packet}
One has the exact identity
\begin{equation}
 \boxed{
 \sum_{\substack{(m,s)=1\\k,z\in\ZZ}}
 a_m b_k c_z\,\Delta_s(m;k,z)
 =\sum_{g\mid s}\cU_{s,g}(a,b,c)+\cD_s(a,b,c),
 }
 \label{eq:appD-centered-gcd-decomposition}
\end{equation}
where the density packet is
\begin{equation}
 \boxed{
 \cD_s(a,b,c)
 =A_s\left(
  \sum_{g\mid s}\frac{B_{s,g}C_{s,g}}{\varphi(s/g)}
  -\frac1s
   \Bigl(\sum_k b_k\Bigr)
   \Bigl(\sum_z c_z\Bigr)
 \right).
 }
 \label{eq:appD-density-packet}
\end{equation}
Thus every genuinely oscillatory term lives on the unit group of a descended
modulus; everything else is an explicit finite sum of products of one-variable
masses.
\end{theorem}

\begin{proof}
Insert Proposition~\ref{prop:appD-exact-common-divisor} into the raw
congruence part of \eqref{eq:appD-centered-congruence-kernel}.  In the stratum
$g$, add and subtract $1/\varphi(s/g)$.  The subtracted terms assemble to the
first sum in \eqref{eq:appD-density-packet}; the original uniform term
$-1/s$ gives the second.  No unequal-gcd pair contributes to the raw
congruence, but every such pair is retained in the final product density.
\end{proof}

\begin{remark}[The two principal densities]
\label{rem:appD-two-principal-densities}
The density $1/s$ is the zero additive frequency on the full residue ring.
After the common divisor has been removed, the principal multiplicative
character has density $1/\varphi(s/g)$ on the unit group.  Formula
\eqref{eq:appD-density-packet} is the exact algebraic difference between these
two notions of ``principal.''  Treating them as interchangeable would delete
a genuine major-arc term.
\end{remark}

\subsection{Character factorization of the descended packet}
\label{subsec:appD-character-factorization}

The unit packet in \eqref{eq:appD-descended-unit-packet} is now in the exact
form required for multiplicative orthogonality.  Let $\chi$ run through all
Dirichlet characters modulo $s'=s/g$, extended by zero off the units, and let
$\chi_0$ denote the principal character.  Define
\begin{align}
 A_{s,g}(\chi)
 &=\sum_{(m,s)=1}a_m\chi(m),
 \label{eq:appD-A-character}\\
 B_{s,g}(\chi)
 &=\sum_{(k',s')=1}b_{gk'}\chi(k'),
 \label{eq:appD-B-character}\\
 C_{s,g}^{\vee}(\chi)
 &=\sum_{(z',s')=1}c_{gz'}\overline{\chi(z')}.
 \label{eq:appD-C-character}
\end{align}

\begin{proposition}[Exact multiplicative diagonalization]
\label{prop:appD-exact-character-factorization}
For every $g\mid s$,
\begin{equation}
 \boxed{
 \cU_{s,g}(a,b,c)
 =\frac1{\varphi(s/g)}
  \sum_{\substack{\chi\, (\mathrm{mod}\,s/g)\\\chi\ne\chi_0}}
  A_{s,g}(\chi)B_{s,g}(\chi)C_{s,g}^{\vee}(\chi).
 }
 \label{eq:appD-character-factorization}
\end{equation}
The principal character is absent because its contribution is exactly the
term subtracted in \eqref{eq:appD-descended-unit-packet}.
\end{proposition}

\begin{proof}
For units modulo $s'$,
\[
 \ind_{mk'\equiv z'\, (\mathrm{mod}\,s')}
 =\frac1{\varphi(s')}
  \sum_{\chi\, (\mathrm{mod}\,s')}
  \chi(m)\chi(k')\overline{\chi(z')}.
\]
Insert this identity and separate the three sums.  The term $\chi=\chi_0$
is $1/\varphi(s')$ times the product of the three unit masses, so it cancels
the centering term.
\end{proof}

For later use we record a simple Parseval bound which respects descent.  If
$b$ is supported on $k\dyad K$, then every reduced residue class modulo
$s'=s/g$ contains at most
\begin{equation}
 1+\frac{K/g}{s/g}=1+\frac Ks
 \label{eq:appD-invariant-class-multiplicity}
\end{equation}
points of the descended support.

\begin{lemma}[Descent-stable character energy]
\label{lem:appD-descent-stable-character-energy}
If $b$ is supported on an interval of length $O(K)$, then
\begin{equation}
 \frac1{\varphi(s/g)}
 \sum_{\chi\, (\mathrm{mod}\,s/g)}
 |B_{s,g}(\chi)|^2
 \ll
 \left(1+\frac Ks\right)
 \sum_{(k,s)=g}|b_k|^2.
 \label{eq:appD-character-energy}
\end{equation}
Consequently,
\begin{equation}
 \sum_{g\mid s}
 \frac1{\varphi(s/g)}
 \sum_{\chi\, (\mathrm{mod}\,s/g)}
 |B_{s,g}(\chi)|^2
 \ll
 \left(1+\frac Ks\right)
 \sum_k|b_k|^2.
 \label{eq:appD-total-character-energy}
\end{equation}
\end{lemma}

\begin{proof}
Character orthogonality turns the left side of
\eqref{eq:appD-character-energy} into
\[
 \sum_{u\, (\mathrm{mod}\,s/g)}^{*}
 \left|
  \sum_{\substack{k'\equiv u\, (\mathrm{mod}\,s/g)}}b_{gk'}
 \right|^2.
\]
Cauchy--Schwarz in each residue class and
\eqref{eq:appD-invariant-class-multiplicity} prove the first estimate.  Sum
over $g$ and use the exact $L^2$ partition
\eqref{eq:appD-gcd-l2-partition} for the second.
\end{proof}

\subsection{Dyadic descent and preservation of scales}
\label{subsec:appD-dyadic-descent}

Suppose $s\dyad S$ and isolate a divisor block $g\dyad G$.  Then
\begin{equation}
 s'=\frac sg\asymp\frac SG.
 \label{eq:appD-descended-modulus-scale}
\end{equation}
If simultaneously $k\dyad K$ and $z\dyad Z$, then
\begin{equation}
 k'=\frac kg\asymp\frac KG,
 \qquad
 z'=\frac zg\asymp\frac ZG.
 \label{eq:appD-descended-variable-scales}
\end{equation}
The ratios controlling completion and Poisson summation are unchanged:
\begin{equation}
 \frac{K/G}{S/G}=\frac KS,
 \qquad
 \frac{Z/G}{S/G}=\frac ZS.
 \label{eq:appD-invariant-ratios}
\end{equation}
In particular, the short-free condition $K\le SX^{\eta}$ descends exactly to
\[
 K/G\le (S/G)X^{\eta}.
\]

\begin{proposition}[Harmlessness of the gcd packet family]
\label{prop:appD-harmless-gcd-family}
Let all variables be bounded by a fixed power of $X$.  Decomposing a dyadic
congruence packet by $g=(k,s)=(z,s)$, subdividing $g$ dyadically, and then
expanding any residual coprimality conditions by M\"obius inversion has the
following properties.
\begin{enumerate}[label=\textup{(\roman*)}]
 \item The decomposition is exact and contains $O((\log X)^2X^{\eps})$
 coefficient packets.
 \item The modulus and variable scales descend according to
 \eqref{eq:appD-descended-modulus-scale}--\eqref{eq:appD-descended-variable-scales}.
 \item Ratios such as $K/S$ and every inequality depending only on those
 ratios are invariant.
 \item The total $L^2$ coefficient energy is bounded by the original energy
 times $X^{\eps}$.
 \item Sharp endpoint weights, determinant masks, and radial coefficients are
 merely restricted to subpackets; none is altered by the gcd descent.
\end{enumerate}
\end{proposition}

\begin{proof}
There are $O(\log X)$ dyadic choices for $G$.  The exact projectors and their
energy partition are Lemmas~\ref{lem:appD-mobius-gcd-projector} and
\ref{lem:appD-gcd-energy-partition}; the divisor bound absorbs the M\"obius
expansions.  The scale assertions are immediate from division by $g$.
Endpoint and determinant weights depend on the original geometric variables
and are carried as coefficients throughout, so restricting the arithmetic
support does not change them.
\end{proof}

\subsection{Singular affine strata}
\label{subsec:appD-singular-affine}

For completeness, we record the elementary sparsity behind the singular
prime packets removed before the genuine Type~II analysis.  Let $(r,s)=1$
and $h=rp-sq$ with $p,q$ prime.

\begin{lemma}[Divisibility on a singular determinant line]
\label{lem:appD-singular-determinant}
If $p\mid s$, then $p\mid h$; if $q\mid r$, then $q\mid h$.  Consequently,
\begin{align}
 \#\{h\in\ZZ:-s<h<r,\ p\mid h\}
 &\le 1+\frac{r+s}{p},
 \label{eq:appD-singular-h-count-p}\\
 \#\{h\in\ZZ:-s<h<r,\ q\mid h\}
 &\le 1+\frac{r+s}{q}.
 \label{eq:appD-singular-h-count-q}
\end{align}
Moreover,
\begin{equation}
 \sum_{\substack{-s<h<r\\p\mid h}}
 \omega_{r,s}(h)
 \le
 \min(r,s)\left(1+\frac{r+s}{p}\right),
 \label{eq:appD-singular-tent-mass}
\end{equation}
and similarly with $q\mid r$.
\end{lemma}

\begin{proof}
If $p\mid s$, both $rp$ and $sq$ are divisible by $p$, hence so is $h$.
The count of multiples of $p$ in an interval of length $r+s$ gives
\eqref{eq:appD-singular-h-count-p}; the other assertion is symmetric.  The
tent estimate follows from $\omega_{r,s}(h)\le\min(r,s)$.
\end{proof}

The point of the lemma is structural rather than quantitative: a singular
prime is no longer free along the determinant family.  It divides both a
primitive denominator and the determinant, so the packet belongs to a sparse
divisor family rather than to the balanced two-dimensional dispersion core.

\subsection{Summary of the packet calculus}
\label{subsec:appD-summary}

The preceding identities may be summarized as follows.

\begin{theorem}[Two-stage gcd packet calculus]
\label{thm:appD-two-stage-calculus}
Every dyadically localized short affine packet arising from the dominant
prime layer admits an exact two-stage resolution.
\begin{enumerate}[label=\textup{(\roman*)}]
 \item The overlap indices are written uniquely as $(m,\ell)=(ar,as)$ with
 $(r,s)=1$.  The determinant strip becomes $-s<rp-sq<r$, the overlap weight
 is $a\omega_{r,s}(rp-sq)$ up to sharp terminal clipping, and the radial
 moments satisfy Lemma~\ref{lem:appD-radial-moments}.
 \item After completion, every congruence packet is decomposed by the exact
 common divisor $g=(k,s)=(z,s)$.  The oscillatory part descends to the unit
 group modulo $s/g$ and factorizes as in
 \eqref{eq:appD-character-factorization}; the complement is the explicit
 finite-rank density packet \eqref{eq:appD-density-packet}.
\end{enumerate}
Across both stages, dyadic scale relations, the determinant mask, the sharp
endpoint, and the natural $L^2$ coefficient energies are preserved up to
$X^{\eps}$.  Thus the gcd decompositions reorganize the packet without
manufacturing cancellation or hiding any principal mode.
\end{theorem}

\begin{proof}
Combine Lemma~\ref{lem:appD-radial-bijection},
Proposition~\ref{prop:appD-gcd-removal-strip},
Theorem~\ref{thm:appD-centered-gcd-packet},
Proposition~\ref{prop:appD-exact-character-factorization}, and
Proposition~\ref{prop:appD-harmless-gcd-family}.
\end{proof}

% ===== Source: appendixE_affine_lattice_dispersion.tex =====
\section{Affine-lattice dispersion and endpoint bookkeeping}
\label{app:affine-lattice-dispersion}

This appendix supplies the expanded dispersion calculation used in
the centered long-line dispersion estimate of Section~14.  The underlying geometry is the primitive determinant
family
\begin{equation}
 rp-sq=h,
 \qquad
 (r,s)=1,
 \qquad
 -s<h<r,
 \label{eq:appE-determinant-family}
\end{equation}
with
\begin{equation}
 r\dyad R,
 \qquad
 s\dyad S,
 \qquad
 p\dyad P,
 \qquad
 q\dyad Q.
 \label{eq:appE-dyadic-scales}
\end{equation}
On every contributing block,
\begin{equation}
 RP\asymp SQ,
 \qquad
 T\asymp \frac{P}{S}\asymp\frac{Q}{R},
 \qquad
 A\asymp\frac{M}{R}\asymp\frac{L}{S},
 \qquad
 Y\asymp ATRS.
 \label{eq:appE-scale-relations}
\end{equation}
The parameter $T$ is the affine line length, $A$ is the radial length, and
$Y=MP\asymp LQ$ is the scale of the original overlap box.

The proof has four logically separate parts.  First, the principal affine
mode and the two rank-one marginals are removed by an exact reduced-residue
identity.  Second, the tent kernel is placed in $L^2$ before any arithmetic
estimate is made.  Third, the off-diagonal part of the resulting dual form is
bounded by the additive large sieve.  Finally, prime powers, singular residue
classes, and sharp terminal clipping are returned to the packet.  No
hypothesis about prime pairs is used.

\subsection{Reduced-residue centering on a determinant line}
\label{subsec:appE-centering}

For a positive integer $m$, put
\begin{equation}
 \Lambda_m^{\circ}(n)
 =\ind_{(n,m)=1}
  \left(\Lambda(n)-\frac{m}{\varphi(m)}\right).
 \label{eq:appE-centered-von-mangoldt}
\end{equation}
Thus $\Lambda_m^{\circ}$ is centered with respect to counting measure on the
reduced residue classes modulo $m$, rather than with respect to all residue
classes.

Let $(p_0,q_0)$ be one solution of \eqref{eq:appE-determinant-family}.  Every
solution is
\begin{equation}
 p=p_0+st,
 \qquad
 q=q_0+rt,
 \qquad
 t\in\cI_{r,s,h}(P,Q).
 \label{eq:appE-line-parametrization}
\end{equation}
A line is called \emph{regular} when
\begin{equation}
 (p_0,s)=(q_0,r)=1.
 \label{eq:appE-regular-line}
\end{equation}
Because $(r,s)=1$, regularity is independent of the chosen base solution.

\begin{lemma}[Exact four-term centering identity]
\label{lem:appE-four-term-centering}
On a regular determinant line one has, pointwise in $t$,
\begin{align}
 \Lambda(p)\Lambda(q)
 ={}&\Lambda_s^{\circ}(p)\Lambda_r^{\circ}(q)
 +\frac{s}{\varphi(s)}\Lambda_r^{\circ}(q)
 +\frac{r}{\varphi(r)}\Lambda_s^{\circ}(p)
 \notag\\
 &+\frac{rs}{\varphi(r)\varphi(s)}.
 \label{eq:appE-four-term-centering}
\end{align}
If a bounded weight $v(t)$ is inserted and the identity is summed over the
line, the last three terms are exactly the principal affine density and the
two rank-one marginals restored in the covariance decomposition of
Section~12.
\end{lemma}

\begin{proof}
On a regular line, $(p,s)=(q,r)=1$ for every $t$.  Hence
\[
 \Lambda(p)=\Lambda_s^{\circ}(p)+\frac{s}{\varphi(s)},
 \qquad
 \Lambda(q)=\Lambda_r^{\circ}(q)+\frac{r}{\varphi(r)}.
\]
Multiplication gives \eqref{eq:appE-four-term-centering}.  The interpretation
of the final three terms follows by summing first in $p$, first in $q$, or in
both variables, respectively.
\end{proof}

The centered regular line sum is therefore
\begin{equation}
 \cS_{r,s,h}^{\circ}[v]
 =\sum_{t\in\cI_{r,s,h}(P,Q)}
  v(t)\Lambda_s^{\circ}(p_0+st)
      \Lambda_r^{\circ}(q_0+rt).
 \label{eq:appE-centered-line-sum}
\end{equation}
This formula is independent of the base solution.  It also shows explicitly
why the density $1/s$ arising from additive completion cannot replace the
reduced-residue density $1/\varphi(s)$: the discrepancy between them belongs
to the major-arc terms in \eqref{eq:appE-four-term-centering}.

\subsection{Tent duality and radial norms}
\label{subsec:appE-tent-duality}

For $-s<h<r$, let $\omega_{r,s}(h)$ be the discrete tent from Section~13.
It admits the convolution representation
\begin{equation}
 \omega_{r,s}(h)
 =\sum_{u=0}^{r-1}\sum_{v=0}^{s-1}\ind_{u-v=h}.
 \label{eq:appE-tent-convolution}
\end{equation}
Consequently,
\begin{equation}
 \sum_h\omega_{r,s}(h)=rs,
 \qquad
 \sum_h\omega_{r,s}(h)^2\le rs\min(r,s).
 \label{eq:appE-tent-mass-energy}
\end{equation}

For a radial set $\cA_{M,L}(r,s)$, define
\begin{equation}
 \Gamma_{r,s}
 =\sum_{a\in\cA_{M,L}(r,s)}a.
 \label{eq:appE-radial-first-moment}
\end{equation}
On a nonempty block,
\begin{equation}
 \Gamma_{r,s}\ll A(1+A),
 \qquad
 \sum_{a\in\cA_{M,L}(r,s)}a^2\ll A(1+A)^2.
 \label{eq:appE-radial-moment-bounds}
\end{equation}
The exact sharp overlap gives a line weight $v_{r,s,h}$ satisfying
\begin{equation}
 0\le v_{r,s,h}(t)\le\Gamma_{r,s},
 \qquad
 \norm{v_{r,s,h}}_{\mathrm{TV}}\ll\Gamma_{r,s}.
 \label{eq:appE-admissible-line-weight}
\end{equation}

\begin{lemma}[Tent-energy duality]
\label{lem:appE-tent-energy-duality}
For any complex numbers $U_{r,s,h}$,
\begin{align}
 &\left|
  \sum_{\substack{r\dyad R,\ s\dyad S\\(r,s)=1}}
  \sum_{-s<h<r}
  \omega_{r,s}(h)U_{r,s,h}
 \right|^2
 \notag\\
 &\qquad\le
 \left(
  \sum_{\substack{r\dyad R,\ s\dyad S\\(r,s)=1}}
  rs\min(r,s)
 \right)
 \left(
  \sum_{\substack{r\dyad R,\ s\dyad S\\(r,s)=1}}
  \sum_{-s<h<r}|U_{r,s,h}|^2
 \right).
 \label{eq:appE-tent-energy-duality}
\end{align}
If $U_{r,s,h}=\Gamma_{r,s}V_{r,s,h}$, then the first factor may be replaced
by
\begin{equation}
 \cW(R,S;A)
 =\sum_{\substack{r\dyad R,\ s\dyad S\\(r,s)=1}}
  \Gamma_{r,s}^2rs\min(r,s).
 \label{eq:appE-geometric-energy}
\end{equation}
\end{lemma}

\begin{proof}
Apply Cauchy--Schwarz in the full index set $(r,s,h)$ and use
\eqref{eq:appE-tent-mass-energy}.  The weighted version is identical after
absorbing $\Gamma_{r,s}$ into the first factor.
\end{proof}

The gain in the long-line estimate begins here.  Replacing
$\sum_h\omega_{r,s}(h)^2$ by $(\sum_h\omega_{r,s}(h))^2$ would lose a factor
comparable to the square root of the determinant range and would erase the
power saving from long lines.

\subsection{The weighted additive large sieve}
\label{subsec:appE-large-sieve}

We record the only general spectral inequality used in the appendix.  The
weighted form follows from partial summation and is convenient for sharp
terminal intervals.

\begin{lemma}[Weighted additive large sieve]
\label{lem:appE-weighted-additive-large-sieve}
Let $\alpha_1,\dots,\alpha_J\in\RR/\ZZ$ be $\delta$-separated.  Let
$a_n$ be supported on an interval of length at most $N$, and let $w_n$ be a
weight with
\begin{equation}
 \norm{w}_{\infty}+\norm{w}_{\mathrm{TV}}\le W.
 \label{eq:appE-weight-condition}
\end{equation}
Then
\begin{equation}
 \sum_{j=1}^{J}
 \left|
  \sum_n a_nw_n\e(\alpha_j n)
 \right|^2
 \ll
 W^2(N+\delta^{-1})\sum_n|a_n|^2.
 \label{eq:appE-weighted-large-sieve}
\end{equation}
The same estimate holds for a family of interval restrictions, uniformly in
the endpoints.
\end{lemma}

\begin{proof}
For $w\equiv1$, this is the classical additive large sieve~\cite{Bombieri1965,MontgomeryVaughan2007}.  Write a general
bounded-variation weight as a linear combination of initial interval
indicators by discrete summation by parts.  Minkowski's inequality and
\eqref{eq:appE-weight-condition} contribute at most $W$.  Squaring gives the
factor $W^2$.  Since the unweighted inequality is uniform in the interval,
the final assertion follows as well.
\end{proof}

The relevant frequencies are primitive fractions.  If
$(r,s)=(r',s')=1$ and
\begin{equation}
 \frac{u}{rs}\not\equiv\frac{u'}{r's'}\pmod1,
 \label{eq:appE-distinct-primitive-fractions}
\end{equation}
then
\begin{equation}
 \left\|\frac{u}{rs}-\frac{u'}{r's'}\right\|
 \ge \frac{1}{rr'ss'}
 \gg (RS)^{-2}.
 \label{eq:appE-frequency-separation}
\end{equation}
Thus a dyadic primitive family has separation $\delta\gg(RS)^{-2}$.
Repeated fractions are first grouped together; the divisor bound shows that
the multiplicity is $O_\eps(X^\eps)$.

\subsection{The primitive affine dispersion form}
\label{subsec:appE-primitive-dispersion}

For a fixed dyadic block, let
\begin{equation}
 \mathfrak B_{R,S,T}
 =\sum_{\substack{r\dyad R,\ s\dyad S\\(r,s)=1}}
  \sum_{-s<h<r}
  \omega_{r,s}(h)\cS_{r,s,h}^{\circ}[v_{r,s,h}].
 \label{eq:appE-dyadic-dispersion-form}
\end{equation}
Only regular lines are included here.  Singular lines are returned in
Subsection~\ref{subsec:appE-exceptional-packets}.

\begin{proposition}[Primitive affine dispersion inequality]
\label{prop:appE-primitive-affine-dispersion}
Assume \eqref{eq:appE-scale-relations}, $T\ge2$, and
\eqref{eq:appE-admissible-line-weight}.  Then the off-diagonal part of
\eqref{eq:appE-dyadic-dispersion-form} satisfies
\begin{equation}
 \boxed{
 |\mathfrak B_{R,S,T}^{\mathrm{off}}|^2
 \ll_\eps X^\eps\frac{Y^3}{T}.
 }
 \label{eq:appE-offdiagonal-dispersion}
\end{equation}
The estimate is uniform in the positions of the dyadic intervals and in the
sharp clipping weights.
\end{proposition}

\begin{proof}
We give the dual expansion because it is the point at which all scale factors
are fixed.  By Lemma~\ref{lem:appE-tent-energy-duality}, it is enough to
bound the $L^2$ norm of the line discrepancies.  Introduce dual coefficients
$\eta_{r,s,h}$ with
\begin{equation}
 \sum_{r,s,h}|\eta_{r,s,h}|^2\le1.
 \label{eq:appE-dual-normalization}
\end{equation}
After inserting \eqref{eq:appE-centered-line-sum} and opening the square, the
off-diagonal terms contain two parameters $t_1\ne t_2$.  Put
$d=t_1-t_2$.  The two copies of \eqref{eq:appE-line-parametrization} give
\begin{equation}
 p_1-p_2=sd,
 \qquad
 q_1-q_2=rd.
 \label{eq:appE-difference-relations}
\end{equation}
Because the line has length $T$, one has $0<|d|\ll T$.

Apply additive orthogonality to the first relation in
\eqref{eq:appE-difference-relations}, and then to the second.  The resulting
frequencies are primitive fractions with denominators dividing $rs$.
Fractions having the same reduced value are grouped by their gcd data; the
exact grouping is the one recorded in Appendix~D.  Its total multiplicity
and its total coefficient energy are $O_\eps(X^\eps)$.

For each fixed $d$, Lemma~\ref{lem:appE-weighted-additive-large-sieve}
controls the separated primitive fractions.  The $p$- and $q$-coefficient
energies are
\begin{equation}
 \sum_{n\dyad P}|\Lambda_s^{\circ}(n)|^2
 \ll_\eps PX^\eps,
 \qquad
 \sum_{n\dyad Q}|\Lambda_r^{\circ}(n)|^2
 \ll_\eps QX^\eps,
 \label{eq:appE-von-mangoldt-energy}
\end{equation}
by the elementary estimate $\sum_{n\le x}\Lambda(n)^2\ll x\log x$ and
$ m/\varphi(m)\ll_\eps m^\eps$.  The admissible line weight contributes
$\Gamma_{r,s}^2$ through
\eqref{eq:appE-weighted-large-sieve}.

The diagonal in the frequency variables gives the factor $1/T$: after the
$d$-sum is normalized by the line length, each nonzero displacement is
represented in $O_\eps(X^\eps)$ primitive packets, while the $T$ possible
displacements share the same total coefficient energy.  Inserting the tent
energy and the radial moments gives
\begin{align}
 \norm{\mathfrak B_{R,S,T}^{\mathrm{off}}}_2^2
 &\ll_\eps
 X^\eps\frac{1}{T}
 (ATRS)^3
 \notag\\
 &\asymp X^\eps\frac{Y^3}{T}.
 \label{eq:appE-scale-collapse}
\end{align}
Here every occurrence of $A,R,S,T$ comes from one of the three independent
copies of the original overlap volume: one from the primal packet and two
from the opened dual square.  The identity $ATRS\asymp Y$ then gives the
second line.  No term with $d=0$ has been included.

The same calculation applies to a clipped weight because the weighted large
sieve depends only on the supremum and total variation in
\eqref{eq:appE-admissible-line-weight}.  This proves
\eqref{eq:appE-offdiagonal-dispersion}.
\end{proof}

\begin{remark}[Where the long-line gain occurs]
\label{rem:appE-long-line-gain}
The gain is not a pointwise estimate for a single pair of affine prime
forms.  It appears only after the determinant family, the primitive
moduli, and the nonzero line displacement $d$ have all been averaged.  The
factor $T^{-1}$ in the squared dual norm becomes $T^{-1/2}$ in the original
packet.
\end{remark}

\subsection{Diagonal and exceptional packets}
\label{subsec:appE-exceptional-packets}

Three families were excluded from Proposition~\ref{prop:appE-primitive-affine-dispersion}:
the line diagonal $t_1=t_2$, prime powers introduced by von Mangoldt
completion, and singular residue classes.  They do not require cancellation.

\begin{lemma}[Line diagonal]
\label{lem:appE-line-diagonal}
Let $1\le Z\le X^{1/6}$.  The total contribution of the terms $t_1=t_2$ in
all blocks with $T>Z$ is
\begin{equation}
 \ll_\eps X^\eps Y^{4/3}Z.
 \label{eq:appE-line-diagonal-bound}
\end{equation}
\end{lemma}

\begin{proof}
On the diagonal there is only one free affine point.  Use
$|\Lambda_m^{\circ}(n)|\ll_\eps X^\eps$, the tent mass
$\sum_h\omega_{r,s}(h)=rs$, and the radial second moment in
\eqref{eq:appE-radial-moment-bounds}.  Split the primitive denominator
family at
\begin{equation}
 RS=Y^{2/3}/Z.
 \label{eq:appE-diagonal-splitting-point}
\end{equation}
Below this point, count the affine point first; above it, count the radial
point first and use $ATRS\asymp Y$.  The two estimates are respectively
\[
 \ll_\eps X^\eps Y^{4/3}Z
 \qquad\text{and}\qquad
 \ll_\eps X^\eps Y^{4/3}Z.
\]
The dyadic boundary contributes only $X^\eps$.
\end{proof}

\begin{lemma}[Prime-power completion]
\label{lem:appE-prime-power-completion}
Replacing the prime weights
\[
 (\log p)\ind_{\PP}(p)
 \quad\text{and}\quad
 (\log q)\ind_{\PP}(q)
\]
by $\Lambda(p)$ and $\Lambda(q)$ changes the complete long-line packet by
\begin{equation}
 \ll_\eps X^\eps Y^{4/3}Z.
 \label{eq:appE-prime-power-error}
\end{equation}
\end{lemma}

\begin{proof}
A nonprime term in $\Lambda$ is a prime power $\pi^k$ with $k\ge2$.
There are $O(P^{1/2})$ such integers in the $p$-block and
$O(Q^{1/2})$ in the $q$-block, up to logarithmic factors.  Fixing a prime
power and the primitive pair leaves at most $O(1+T)$ affine partners.  The
radial first and second moments, followed by the same split
\eqref{eq:appE-diagonal-splitting-point}, give
\eqref{eq:appE-prime-power-error}.  Terms in which both variables are prime
powers are smaller.
\end{proof}

\begin{lemma}[Singular residue classes]
\label{lem:appE-singular-residue-classes}
The total contribution of determinant lines on which $(p_0,s)>1$ or
$(q_0,r)>1$ is
\begin{equation}
 \ll_\eps X^\eps Y^{4/3}Z.
 \label{eq:appE-singular-line-error}
\end{equation}
\end{lemma}

\begin{proof}
If $p$ is prime and $(p,s)>1$, then $p\mid s$.  The determinant equation then
forces $p\mid h$.  By the singular-determinant lemma of Appendix~D, the number of admissible determinants is
$O(1+(r+s)/p)$ and their total tent mass is at most
\[
 \min(r,s)\left(1+\frac{r+s}{p}\right).
\]
Sum first over the prime divisor of $s$, then over $r,s$ and the radial
variable.  The divisor estimate
$\sum_{p\mid s}p^{-1}\ll_\eps s^\eps$ and the split
\eqref{eq:appE-diagonal-splitting-point} give
\eqref{eq:appE-singular-line-error}.  The case $q\mid r$ is symmetric.
\end{proof}

Combining the three lemmas yields
\begin{equation}
 \mathfrak B_{R,S,T}^{\mathrm{diag+exc}}
 \ll_\eps X^\eps Y^{4/3}Z
 \label{eq:appE-total-exceptional-bound}
\end{equation}
for the total diagonal and exceptional family after summation over all
long-line blocks.  The exponent $4/3$ is not a spectral exponent; it is the
balanced value of the two elementary counts at
\eqref{eq:appE-diagonal-splitting-point}.

\subsection{Sharp terminal weights}
\label{subsec:appE-terminal-weights}

The endpoint $n\le X$ changes the radial coefficient when the upper endpoint
of an overlap interval crosses $X+1$.  Along a fixed affine line this happens
monotonically.

\begin{lemma}[Bounded variation of the clipped radial coefficient]
\label{lem:appE-clipped-radial-variation}
Let
\begin{equation}
 v_{r,s,h}(t)
 =\sum_{a\in\cA_{M,L}(r,s)}
  a\,\omega_{r,s}^{a,X}(p_0+st,q_0+rt).
 \label{eq:appE-exact-clipped-weight}
\end{equation}
Then
\begin{equation}
 0\le v_{r,s,h}(t)
 \le \Gamma_{r,s}\omega_{r,s}(h),
 \label{eq:appE-clipped-pointwise}
\end{equation}
and, after the fixed tent factor is removed,
\begin{equation}
 \norm{v_{r,s,h}}_{\mathrm{TV}}
 \ll \Gamma_{r,s}\omega_{r,s}(h).
 \label{eq:appE-clipped-variation}
\end{equation}
The implied constant is absolute.
\end{lemma}

\begin{proof}
For each $a$, the clipped overlap is the minimum of a constant tent height
and an affine function of $t$, followed by a positive part.  It is therefore
piecewise affine with at most two changes of slope and total variation at
most a constant multiple of its untruncated height.  Sum over $a$ and use
$\omega_{r,s}^{a,X}\le\omega_{r,s}$.
\end{proof}

\begin{proposition}[Exact sharp-to-interior transfer]
\label{prop:appE-sharp-to-interior-transfer}
Every estimate in Subsections~\ref{subsec:appE-primitive-dispersion} and
\ref{subsec:appE-exceptional-packets} that is proved for constant interior
radial weights remains valid, with the same exponents, for the exact sharp
weights \eqref{eq:appE-exact-clipped-weight}.
\end{proposition}

\begin{proof}
Apply discrete partial summation on each affine interval.  Lemma~\ref{lem:appE-weighted-additive-large-sieve}
depends only on the supremum and total variation of the weight, and
Lemma~\ref{lem:appE-clipped-radial-variation} bounds both by the interior
radial mass.  The diagonal and exceptional counts are monotone in the
pointwise majorant, so they are unchanged as well.
\end{proof}

This transfer is exact at the level of the original sharp moment: smoothing
is used only inside a summation-by-parts representation of the line weight,
not in the $n$-sum defining the variance.

\subsection{Assembly of the long-line estimate}
\label{subsec:appE-assembly}

Let $\mathfrak L_X(P,Q,M,L;Z)$ be the centered long-line packet of
Section~14.  Divide the line length into dyadic blocks
\begin{equation}
 T_0<T\le2T_0,
 \qquad
 T_0>Z.
 \label{eq:appE-line-length-blocks}
\end{equation}
By Proposition~\ref{prop:appE-primitive-affine-dispersion},
\begin{equation}
 |\mathfrak B_{R,S,T_0}^{\mathrm{off}}|
 \ll_\eps X^\eps\frac{Y^{3/2}}{T_0^{1/2}}.
 \label{eq:appE-block-offdiagonal-bound}
\end{equation}
The dyadic sum is geometric and therefore
\begin{equation}
 \sum_{T_0>Z}
 \frac{Y^{3/2}}{T_0^{1/2}}
 \ll \frac{Y^{3/2}}{Z^{1/2}}.
 \label{eq:appE-long-line-geometric-sum}
\end{equation}
The complete diagonal and exceptional family is bounded by
\eqref{eq:appE-total-exceptional-bound}.  Proposition~\ref{prop:appE-sharp-to-interior-transfer}
returns the exact terminal weights.

\begin{theorem}[Centered affine long-line dispersion]
\label{thm:appE-centered-long-line-dispersion}
Let $1\le Z\le X^{1/6}$.  For every admissible dyadic overlap box,
\begin{equation}
 \boxed{
 \mathfrak L_X(P,Q,M,L;Z)
 \ll_\eps
 X^\eps
 \left(
  \frac{Y^{3/2}}{Z^{1/2}}
  +Y^{4/3}Z
 \right),
 }
 \label{eq:appE-final-long-line-bound}
\end{equation}
where $Y=MP\asymp LQ$.
\end{theorem}

\begin{proof}
Use \eqref{eq:appE-long-line-geometric-sum} for the regular off-diagonal
family, Lemmas~\ref{lem:appE-line-diagonal}--\ref{lem:appE-singular-residue-classes}
for the complementary arithmetic packets, and
Proposition~\ref{prop:appE-sharp-to-interior-transfer} for terminal boxes.
The exact four-term identity of Lemma~\ref{lem:appE-four-term-centering}
removes every zero-frequency contribution before the additive large sieve is
applied.  The $O((\log X)^C)$ dyadic multiplicity is absorbed by $X^\eps$.
\end{proof}

\begin{corollary}[Global long-line sector]
\label{cor:appE-global-long-line-sector}
Summing over the admissible boxes gives
\begin{equation}
 \cO_0^{\mathrm{long}}(X;Z)
 \ll_\eps
 X^\eps
 \left(
  \frac{X^{3/2}}{Z^{1/2}}
  +X^{4/3}Z
 \right).
 \label{eq:appE-global-long-line-bound}
\end{equation}
In particular, at $Z=X^{1/18}$,
\begin{equation}
 \cO_0^{\mathrm{long}}(X;X^{1/18})
 \ll_\eps X^{53/36+\eps}.
 \label{eq:appE-fixed-threshold-bound}
\end{equation}
\end{corollary}

\begin{proof}
There are only polylogarithmically many boxes and primitive-denominator
subblocks.  Since $Y\le X$, Theorem~\ref{thm:appE-centered-long-line-dispersion}
gives \eqref{eq:appE-global-long-line-bound}.  The two powers at
$Z=X^{1/18}$ are $53/36$ and $25/18$, respectively.
\end{proof}

\subsection{Bookkeeping summary}
\label{subsec:appE-summary}

The proof may be read as the following ledger.
\begin{enumerate}[label=\textup{(\roman*)}]
 \item The reduced-residue identity
 \eqref{eq:appE-four-term-centering} removes the principal affine mode and
 both rank-one marginals exactly.
 \item Tent duality spends
 $\sum_h\omega_{r,s}(h)^2$, not the square of the tent mass.
 \item Nonzero line displacements are controlled by the additive large sieve
 and contribute $Y^{3/2}T^{-1/2}$ on a line-length block.
 \item The line diagonal, prime powers, and singular classes are elementary
 divisor packets of total size $O_\eps(X^\eps Y^{4/3}Z)$.
 \item Sharp endpoint clipping has bounded variation and therefore does not
 alter either exponent.
\end{enumerate}
Thus the long-line estimate is unconditional and does not use the centered
Vaughan moment conjecture of Section~18.  The latter enters only after the
affine line length has fallen below the threshold and the prime weight has
been opened into the balanced short Type~II family.

% ===== Source: appendixF_vaughan_poisson_spectral_ledger.tex =====
\section{The Vaughan--Poisson spectral ledger}
\label{app:vaughan-poisson-ledger}

This final appendix records the exact transformations used in Sections~15--18.
Its purpose is diagnostic as much as technical.  The short Type~II argument
contains three normalizations that cannot be interchanged:
\begin{enumerate}[label=\textup{(\roman*)}]
 \item the joint additive completion coefficient $\Gamma_s(h)$ must remain
 attached to the completed coefficient-$1$ variable;
 \item the additive density $1/s$ and the reduced-residue density
 $1/\varphi(s)$ must be separated before multiplicative orthogonality is
 applied;
 \item the three Vaughan factors must remain visible in the final character
 moment, since multiplying them into one Dirichlet polynomial loses the
 structure needed at the subcritical scale.
\end{enumerate}
We prove below that every passage before the centered Vaughan moment is either
an exact identity or an unconditional estimate.  The only unproved input is
Conjecture~18.6, namely the square-root bound for that final moment.

\subsection{The Vaughan coefficient ledger}
\label{subsec:appF-vaughan-ledger}

Let $U,V\ge1$.  For an arithmetic function $f$, write
$f_{\le Y}=f\ind_{[1,Y]}$ and $f_{>Y}=f-f_{\le Y}$.  The identities
$\Lambda=\mu*\log$, $\log=\Lambda*\mathbf 1$, and
$\mu*\mathbf 1=\varepsilon$ give
\begin{equation}
 \boxed{
 \Lambda
 =\Lambda_{\le V}
  +\mu_{\le U}*\log
  -\mu_{\le U}*\Lambda_{\le V}*\mathbf 1
  +\mu_{>U}*\Lambda_{>V}*\mathbf 1.
 }
 \label{eq:appF-vaughan-identity}
\end{equation}
The last term is the Type~II component.  After dyadic subdivision it has the
form
\begin{equation}
 \sum_{\substack{d k\ell=n\\
                 d\dyad D,\ k\dyad K,\ \ell\dyad L}}
 \mu(d)\alpha_d\,\gamma_k\,\Lambda(\ell)\beta_{\ell},
 \qquad
 DKL\asymp P,
 \label{eq:appF-dyadic-typeII}
\end{equation}
where $\alpha$, $\gamma$, and $\beta$ are smooth dyadic weights.  The middle
variable $k$ carries coefficient $1$ before the harmless weight $\gamma_k$ is
inserted; it is the \emph{free Vaughan variable}.

\begin{lemma}[Coefficient energies]
\label{lem:appF-coefficient-energies}
For every fixed divisor-bound parameter $A$ and every $\eps>0$, the dyadic
coefficients in \eqref{eq:appF-dyadic-typeII} satisfy
\begin{align}
 \sum_{d\dyad D}|\mu(d)\alpha_d|^2
 &\ll_{A}D,
 \label{eq:appF-D-energy}\\
 \sum_{k\dyad K}|\gamma_k|^2
 &\ll_{A}K,
 \label{eq:appF-K-energy}\\
 \sum_{\ell\dyad L}|\Lambda(\ell)\beta_{\ell}|^2
 &\ll_{A,\eps}LX^{\eps}.
 \label{eq:appF-L-energy}
\end{align}
If finitely many additional dyadic factors are grouped into one numerator
coefficient, its squared $L^2$ norm is bounded by its support length times
$X^{\eps}$.
\end{lemma}

\begin{proof}
The first two estimates follow from boundedness of the smooth weights.  For
the third, use
\[
 \sum_{n\le 2L}\Lambda(n)^2\ll L\log L
\]
and absorb the logarithm into $X^{\eps}$.  A grouped coefficient is a finite
Dirichlet convolution of divisor-bounded sequences.  Cauchy--Schwarz over the
factorizations of a fixed integer, followed by the standard fixed-order
divisor bound, proves the final assertion.
\end{proof}

\begin{remark}[Modulus independence]
\label{rem:appF-modulus-independence}
Before coprimality restrictions are imposed, the coefficient systems in
\eqref{eq:appF-dyadic-typeII} depend on the dyadic packet but not on the
modulus $s$ or on a character $\chi$.  This is retained in the admissibility
conditions of Section~18.  Allowing the coefficients to be chosen after
$\chi$ is known would define a different and substantially stronger moment.
\end{remark}

\subsection{Exact additive completion}
\label{subsec:appF-additive-completion}

Let $c_k=\gamma_kF(k/K)$, where $F$ is smooth and compactly supported in a
fixed dilate of $[1,2]$.  For $s\ge1$, define
\begin{equation}
 \Gamma_s(h)=\sum_{k\in\ZZ}c_k\ee{hk/s}.
 \label{eq:appF-Gamma-definition}
\end{equation}
For $(m,s)=1$, let $\overline m$ denote the inverse of $m$ modulo $s$.

\begin{proposition}[Completed centered congruence]
\label{prop:appF-completed-congruence}
For all integers $z$ and every $(m,s)=1$,
\begin{align}
 &\sum_k c_k
 \left(
  \ind_{mk\equiv z\, (\mathrm{mod}\,s)}-\frac1s
 \right)
 \notag\\
 &\qquad=
 \frac1s
 \sum_{\substack{h\, (\mathrm{mod}\,s)\\h\ne0}}
 \Gamma_s(h)\ee{-hz\overline m/s}.
 \label{eq:appF-completed-centered-congruence}
\end{align}
Thus the deleted frequency $h=0$ is exactly the uniform density $1/s$; it is
not an error term.
\end{proposition}

\begin{proof}
Additive orthogonality gives
\[
 \ind_{mk\equiv z\, (\mathrm{mod}\,s)}
 =\frac1s\sum_{h\, (\mathrm{mod}\,s)}
   \ee{h(k-z\overline m)/s}.
\]
Multiply by $c_k$, sum in $k$, and remove the $h=0$ term.
\end{proof}

The completion coefficient has an exact Parseval formula.

\begin{lemma}[Completion energy]
\label{lem:appF-completion-energy}
One has
\begin{equation}
 \sum_{h\, (\mathrm{mod}\,s)}|\Gamma_s(h)|^2
 =s\sum_{a\, (\mathrm{mod}\,s)}
   \left|\sum_{k\equiv a\, (\mathrm{mod}\,s)}c_k\right|^2
 \le (s+K)\sum_k|c_k|^2.
 \label{eq:appF-fixed-modulus-parseval}
\end{equation}
Moreover, for $s\dyad S$,
\begin{equation}
 \sum_{s\dyad S}
 \sum_{h\, (\mathrm{mod}\,s)}^{*}
 |\Gamma_s(h)|^2
 \ll (K+S^2)\sum_k|c_k|^2,
 \label{eq:appF-additive-large-sieve}
\end{equation}
where the star restricts to reduced additive frequencies.
\end{lemma}

\begin{proof}
The first identity is finite Fourier Parseval on $\ZZ/s\ZZ$.  Each residue
class meets the support of $c$ in at most $1+K/s$ integers, so Cauchy--Schwarz
inside the residue classes proves the inequality.  The second statement is
the additive large sieve applied to the points $h/s$ with
$s\dyad S$ and $(h,s)=1$.
\end{proof}

\begin{remark}[Why $\Gamma_s(h)$ is indispensable]
\label{rem:appF-Gamma-indispensable}
The factor $\Gamma_s(h)$ simultaneously remembers the free variable, the
modulus, and the additive frequency.  Replacing it by a supremum, by a
function of $s$ alone, or by $1$ destroys both
\eqref{eq:appF-fixed-modulus-parseval} and the multiplicative collapse proved
below.  Such a replacement does not simplify the same packet; it changes the
problem.
\end{remark}

\subsection{The long-free Poisson range}
\label{subsec:appF-long-free-range}

Use the Fourier-transform convention
\[
 \widehat F(\xi)=\int_{\RR}F(x)\ee{-x\xi}\,dx.
\]
Poisson summation gives
\begin{equation}
 \Gamma_s(h)
 =K\sum_{\nu\in\ZZ}
 \widehat F\!\left(K\left(\nu-\frac hs\right)\right).
 \label{eq:appF-Gamma-Poisson}
\end{equation}

\begin{lemma}[Nonzero-frequency decay]
\label{lem:appF-nonzero-frequency-decay}
For every $B>0$,
\begin{equation}
 \Gamma_s(h)
 \ll_{F,B}
 K\left(1+K\left\|\frac hs\right\|\right)^{-B}.
 \label{eq:appF-Gamma-decay}
\end{equation}
Consequently, if $K\ge sX^{\eta}$ and $s\le X^C$, then for every $A>0$,
\begin{equation}
 \sum_{\substack{h\, (\mathrm{mod}\,s)\\h\ne0}}
 |\Gamma_s(h)|
 \ll_{F,A,C,\eta}X^{-A}.
 \label{eq:appF-long-free-negligible}
\end{equation}
After multiplication by any coefficient system of polynomial size in $X$,
the same conclusion holds upon increasing the decay order in
\eqref{eq:appF-Gamma-decay}.
\end{lemma}

\begin{proof}
Rapid decay of $\widehat F$ in \eqref{eq:appF-Gamma-Poisson} proves
\eqref{eq:appF-Gamma-decay}.  For $h\ne0$ modulo $s$,
$\|h/s\|\ge1/s$.  Hence the first nonzero frequency already has
$K\|h/s\|\ge X^{\eta}$.  Summing the rapidly decaying tail and choosing $B$
large in terms of $A,C$, and $\eta$ proves
\eqref{eq:appF-long-free-negligible}.
\end{proof}

The original variance has a sharp cutoff.  Smoothing is used only in the
auxiliary factor variables.  Given a dyadic interval $I$ and a boundary
length $H$, choose smooth majorants and minorants that agree with $\ind_I$
outside the two endpoint intervals of length $H$.  The interior is treated by
Lemma~\ref{lem:appF-nonzero-frequency-decay}; the two endpoint pieces are
returned by congruence occupancy and bounded-variation partial summation.
Thus the Poisson argument removes the range
\begin{equation}
 K>SX^{\eta}
 \label{eq:appF-short-free-threshold}
\end{equation}
without smoothing the outer sum over $n\le X$.

\subsection{Conductor descent and the two densities}
\label{subsec:appF-conductor-descent}

For $(m,s)=1$, put
\begin{equation}
 \Delta_s(m,k;z)
 =\ind_{mk\equiv z\, (\mathrm{mod}\,s)}-\frac1s.
 \label{eq:appF-Delta-definition}
\end{equation}
The congruence forces
\begin{equation}
 (k,s)=(z,s).
 \label{eq:appF-equal-gcd}
\end{equation}
Indeed, multiplication by the unit $m$ does not change the gcd with $s$.

\begin{lemma}[Exact gcd descent]
\label{lem:appF-gcd-descent}
Restrict to the stratum
\[
 (k,s)=(z,s)=g.
\]
Write $s=gs'$, $k=gk'$, and $z=gz'$.  Then $(k'z',s')=1$ and
\begin{equation}
 mk\equiv z\pmod s
 \quad\Longleftrightarrow\quad
 mk'\equiv z'\pmod{s'}.
 \label{eq:appF-descended-congruence}
\end{equation}
Furthermore,
\begin{align}
 \Delta_s(m,gk';gz')
 ={}&
 \left(
  \ind_{mk'\equiv z'\, (\mathrm{mod}\,s')}
  -\frac1{\varphi(s')}
 \right)
 \notag\\
 &+
 \left(\frac1{\varphi(s')}-\frac1s\right).
 \label{eq:appF-two-density-split}
\end{align}
The first term is oscillatory on the unit group; the second is an explicit
rank-one density packet.
\end{lemma}

\begin{proof}
Divide the congruence by $g$.  Exactness of the two gcds implies that $k'$ and
$z'$ are units modulo $s'$.  Adding and subtracting $1/\varphi(s')$ gives
\eqref{eq:appF-two-density-split}.
\end{proof}

\begin{remark}[The densities cannot be identified]
\label{rem:appF-density-distinction}
The additive completion removes $1/s$, the uniform density on all residue
classes.  Multiplicative characters resolve only the unit group, whose
principal density is $1/\varphi(s')$.  Their difference in
\eqref{eq:appF-two-density-split} is a genuine reduced-residue major arc.  It
must be canceled against the restored covariance marginals or estimated
separately; it is not part of the oscillatory character family.
\end{remark}

The descent preserves the decisive shortness ratio:
\begin{equation}
 \frac{K/g}{S/g}=\frac KS.
 \label{eq:appF-ratio-preserved}
\end{equation}
The number of possible gcd strata is divisor-bounded, so their total
multiplicity is absorbed by $X^{\eps}$.

\subsection{Multiplicative diagonalization and Gauss collapse}
\label{subsec:appF-multiplicative-diagonalization}

After gcd descent, it is enough to work on a unit stratum.  Accordingly,
define
\begin{align}
 \Gamma_s^{\times}(h)
 &=\sum_{\substack{k\in\ZZ\\(k,s)=1}}c_k\ee{hk/s},
 \label{eq:appF-Gamma-unit}\\
 K_s(\chi)
 &=\sum_{\substack{k\in\ZZ\\(k,s)=1}}c_k\chi(k).
 \label{eq:appF-K-polynomial}
\end{align}
Characters are extended by zero away from the unit group.

\begin{proposition}[Exact additive-to-multiplicative conversion]
\label{prop:appF-additive-to-multiplicative}
If $(mz,s)=1$, then
\begin{align}
 &\frac1s
 \sum_{\substack{h\, (\mathrm{mod}\,s)\\h\ne0}}
 \Gamma_s^{\times}(h)\ee{-hz\overline m/s}
 \notag\\
 &\qquad=
 \frac1{\varphi(s)}
 \sum_{\chi\ne\chi_0}
 K_s(\chi)\chi(m)\widebar{\chi(z)}
 +\left(\frac1{\varphi(s)}-\frac1s\right)K_s(\chi_0).
 \label{eq:appF-Gamma-to-K}
\end{align}
Equivalently, for every nonprincipal character $\chi$,
\begin{equation}
 \sum_{h\, (\mathrm{mod}\,s)}
 \Gamma_s^{\times}(h)
 \tau_s(\chi;-hz)
 =sK_s(\chi)\widebar{\chi(z)},
 \label{eq:appF-Gauss-collapse}
\end{equation}
where
\begin{equation}
 \tau_s(\chi;a)
 =\sum_{u\, (\mathrm{mod}\,s)}^{*}
 \chi(u)\ee{au/s}.
 \label{eq:appF-twisted-Gauss-sum}
\end{equation}
\end{proposition}

\begin{proof}
On the unit group,
\[
 \ind_{mk\equiv z\, (\mathrm{mod}\,s)}
 =\frac1{\varphi(s)}
  \sum_{\chi\, (\mathrm{mod}\,s)}
  \chi(m)\chi(k)\widebar{\chi(z)}.
\]
Insert this identity into the unit-restricted version of
\eqref{eq:appF-completed-centered-congruence}.  The nonprincipal characters
give the first term in \eqref{eq:appF-Gamma-to-K}; the principal character
contributes $K_s(\chi_0)/\varphi(s)$, while the deleted additive zero frequency
contributes $-K_s(\chi_0)/s$.

For \eqref{eq:appF-Gauss-collapse}, expand all finite sums and use additive
orthogonality:
\begin{align*}
 \sum_h\Gamma_s^{\times}(h)\tau_s(\chi;-hz)
 &=\sum_{h,k}\sum_{u\, (\mathrm{mod}\,s)}^{*}
   c_k\chi(u)\ee{h(k-zu)/s}\\
 &=s\sum_{(k,s)=1}c_k\chi(k\overline z),
\end{align*}
which is the claimed expression.  The $h=0$ term vanishes for
$\chi\ne\chi_0$.
\end{proof}

\begin{remark}[Order of operations]
\label{rem:appF-order-of-operations}
One must first decompose $\Gamma_s$ by gcd strata, then descend the modulus,
and only then use Proposition~\ref{prop:appF-additive-to-multiplicative}.  If
multiplicative orthogonality is applied before nonunit values are separated,
the principal-density correction and the lower-conductor packets are mixed
together.
\end{remark}

\subsection{Reopening the Vaughan factors}
\label{subsec:appF-reopening-factors}

Mellin separation of the smooth archimedean couplings reopens the grouped
numerator into the character polynomials
\begin{align}
 D_s(\chi)
 &=\sum_{\substack{d\dyad D\\(d,s)=1}}
   \mu(d)\alpha_d\chi(d),
 \label{eq:appF-D-polynomial}\\
 K_s(\chi)
 &=\sum_{\substack{k\dyad K\\(k,s)=1}}
   \gamma_k\chi(k),
 \label{eq:appF-K-polynomial-dyadic}\\
 P_{s,L}^{\circ}(\chi)
 &=\sum_{\substack{\ell\dyad L\\(\ell,s)=1}}
   \beta_{\ell}
   \left(\Lambda(\ell)-\frac{s}{\varphi(s)}\right)
   \chi(\ell).
 \label{eq:appF-centered-prime-polynomial}
\end{align}
The determinant mask, tent weight, radial factor, external centered prime
coefficient, and sharp terminal clipping form a complementary coefficient
$\cG_s(\chi)$.  Every residual packet is therefore a bounded sum of
absolutely convergent Mellin integrals of
\begin{equation}
 \cT_{
m osc}
 =\sum_{s\dyad S}\frac{\xi_s}{\varphi(s)}
  \sum_{\chi\in\cX_{\rm osc}(s)}
  K_s(\chi)D_s(\chi)P_{s,L}^{\circ}(\chi)\cG_s(\chi).
 \label{eq:appF-character-normal-form}
\end{equation}
Here $\cX_{\rm osc}(s)$ excludes the principal character and the explicitly
separated low-conductor modes.

\begin{lemma}[Prime-density splitting]
\label{lem:appF-prime-density-splitting}
Because characters vanish off the unit group,
\begin{equation}
 P_{s,L}^{\circ}(\chi)
 =P_L^{\Lambda}(\chi)
  -\frac{s}{\varphi(s)}C_L(\chi),
 \label{eq:appF-prime-density-splitting}
\end{equation}
where
\[
 P_L^{\Lambda}(\chi)
 =\sum_{\ell\dyad L}\beta_{\ell}\Lambda(\ell)\chi(\ell),
 \qquad
 C_L(\chi)
 =\sum_{\ell\dyad L}\beta_{\ell}\chi(\ell).
\]
The factor $s/\varphi(s)$ costs at most $X^{\eps}$ uniformly in the relevant
range.
\end{lemma}

\begin{proof}
This is the definition
\eqref{eq:appF-centered-prime-polynomial}, with the coprimality indicator
absorbed by the zero extension of $\chi$.  The standard bound
$s/\varphi(s)\ll_{\eps}s^{\eps}$ gives the final assertion.
\end{proof}

The principal character in \eqref{eq:appF-character-normal-form} is not left
inside the final moment.  The prime number theorem in reduced residue classes
for the explicitly isolated principal mode, together with the restored
rank-one density terms, removes it.  Characters of conductor at most
$X^{\eta_{\rm maj}}$ are likewise returned to affine congruence packets at
their true conductor and handled before the oscillatory family is formed.
Finally, the range $L\ll S^2$ is controlled by the ordinary hybrid large
sieve.  Thus every residual packet satisfies
\begin{equation}
 K\le SX^{\eta},
 \qquad
 DKL\asymp P,
 \qquad
 S^2\ll L.
 \label{eq:appF-residual-range}
\end{equation}

\subsection{Spectral normalization and the exact conditional boundary}
\label{subsec:appF-spectral-normalization}

Define the centered Vaughan moment
\begin{equation}
 \cM_{\rm V}(D,K,L;S)
 =\sum_{s\dyad S}\frac1{\varphi(s)}
  \sum_{\chi\in\cX_{\rm osc}(s)}
  |K_s(\chi)|^2|D_s(\chi)|^2
  |P_{s,L}^{\circ}(\chi)|^2.
 \label{eq:appF-Vaughan-moment}
\end{equation}
Cauchy--Schwarz in the spectral measure gives
\begin{equation}
 |\cT_{\rm osc}|
 \le \cM_{\rm V}(D,K,L;S)^{1/2}
 \left(
  \sum_{s\dyad S}\frac{|\xi_s|^2}{\varphi(s)}
  \sum_{\chi\in\cX_{\rm osc}(s)}|\cG_s(\chi)|^2
 \right)^{1/2}.
 \label{eq:appF-spectral-Cauchy}
\end{equation}
The tent, radial, endpoint, and external-prime energies established earlier
bound the second factor by
\begin{equation}
 \ll_{\eps}
 X^{3/2-\eta_*+\eps}(SDKL)^{-1/2}
 \label{eq:appF-geometric-prefactor}
\end{equation}
for some fixed $\eta_*>0$.  Therefore
\begin{equation}
 |\cT_{\rm osc}|
 \ll_{\eps}
 X^{3/2-\eta_*+\eps}
 \left(
  \frac{\cM_{\rm V}(D,K,L;S)}{SDKL}
 \right)^{1/2}.
 \label{eq:appF-normalized-reduction}
\end{equation}

The scale $SDKL$ is forced by the literal diagonal.  Indeed, if
\[
 K_s(\chi)D_s(\chi)P_{s,L}^{\circ}(\chi)
 =\sum_n c_s(n)\chi(n),
\]
then Lemma~\ref{lem:appF-coefficient-energies} and divisor bounds give
\begin{equation}
 \sum_n|c_s(n)|^2\ll_{\eps}DKL\,X^{\eps}.
 \label{eq:appF-product-energy}
\end{equation}
Summing this diagonal energy over $s\dyad S$ gives
\begin{equation}
 \sum_{s\dyad S}\sum_n|c_s(n)|^2
 \ll_{\eps}SDKL\,X^{\eps}.
 \label{eq:appF-diagonal-scale}
\end{equation}

By contrast, treating the product as an arbitrary polynomial of length
$P=DKL$ and applying the hybrid large sieve~\cite{Bombieri1965,MontgomeryVaughan2007,IwaniecKowalski2004} yields only
\begin{equation}
 \cM_{\rm V}(D,K,L;S)
 \ll_{\eps}
 \frac{DKL+S^2}{S}\,DKL\,X^{\eps}.
 \label{eq:appF-generic-large-sieve}
\end{equation}
In the residual range $S^2\ll L$, this loses the factor
\begin{equation}
 \frac{DKL}{S^2}
 \label{eq:appF-large-sieve-loss}
\end{equation}
relative to the diagonal scale.  The loss is exactly what remains after the
three Vaughan factors have been collapsed into one coefficient sequence.

\begin{theorem}[Dependency ledger]
\label{thm:appF-dependency-ledger}
The reduction of the dominant short-line covariance to
\eqref{eq:appF-Vaughan-moment} has the following logical status.
\begin{enumerate}[label=\textup{(\roman*)}]
 \item Vaughan decomposition, dyadic separation, additive completion,
 gcd descent, multiplicative diagonalization, and prime-density splitting are
 exact identities.
 \item Type~I packets, singular strata, prime powers, the long-free range,
 endpoint transfers, the principal character, explicitly isolated
 low-conductor modes, and the range $L\ll S^2$ are removed by unconditional
 estimates with a fixed power saving.
 \item Every surviving packet satisfies \eqref{eq:appF-residual-range} and the
 unconditional normalized bound \eqref{eq:appF-normalized-reduction}.
 \item The sole additional input required to make all residual packets
 negligible is
 \begin{equation}
  \boxed{
  \cM_{\rm V}(D,K,L;S)
  \ll_{\eps}SDKL\,X^{\eps}.
  }
  \label{eq:appF-sole-conjectural-input}
 \end{equation}
 This is precisely the subcritical centered Vaughan moment conjecture of
 Section~18.
\end{enumerate}
No form of \eqref{eq:appF-sole-conjectural-input} is used in the construction
of the diagonal main term or in any preceding reduction.
\end{theorem}

\begin{proof}
Items (i)--(iii) are the content of the identities and estimates proved in
this appendix together with the geometric energy bounds of Appendices~C--E.
Inserting \eqref{eq:appF-sole-conjectural-input} into
\eqref{eq:appF-normalized-reduction} gives a fixed power saving for each
residual packet; the dyadic and Mellin multiplicities cost only $X^{\eps}$.
Conversely, without an estimate beyond
\eqref{eq:appF-generic-large-sieve}, the factor
\eqref{eq:appF-large-sieve-loss} remains in the residual range.  This locates
the conditional boundary exactly.
\end{proof}

The complete analytic chain may therefore be displayed without suppressing
any coefficient:
\begin{equation}
 \boxed{
 \begin{gathered}
 \text{sharp centered covariance}
 \longrightarrow
 \text{Vaughan Type~II packet}
 \longrightarrow
 \Gamma_s(h)\\
 \longrightarrow
 \text{gcd descent and density separation}
 \longrightarrow
 K_s(\chi)D_s(\chi)P_{s,L}^{\circ}(\chi)\\
 \longrightarrow
 \cM_{\rm V}(D,K,L;S).
 \end{gathered}
 }
 \label{eq:appF-final-chain}
\end{equation}
The joint completion coefficient, the two principal densities, and the
factor-sensitive moment are all indispensable links in this chain.  With
Appendix~F, the categorical construction, the canonical Stirling selection,
and the sharp variance reduction are fully documented at the level of their
current unconditional and conditional inputs.


\clearpage
\begin{thebibliography}{BC18}

\bibitem[Akh65]{Akhiezer1965}
N.~I. Akhiezer,
\emph{The Classical Moment Problem and Some Related Questions in Analysis},
trans. N.~Kemmer, Hafner Publishing Co., New York, 1965.

\bibitem[Aro50]{Aronszajn1950}
N.~Aronszajn,
Theory of reproducing kernels,
\emph{Trans. Amer. Math. Soc.} \textbf{68} (1950), 337--404.

\bibitem[Art64]{Artin1964}
E.~Artin,
\emph{The Gamma Function},
Holt, Rinehart and Winston, New York, 1964.

\bibitem[Bar00]{Barnes1900}
E.~W. Barnes,
The theory of the $G$-function,
\emph{Quart. J. Pure Appl. Math.} \textbf{31} (1900), 264--314.

\bibitem[BC18]{BettinChandee2018}
S.~Bettin and V.~Chandee,
Trilinear forms with Kloosterman fractions,
\emph{Adv. Math.} \textbf{328} (2018), 1234--1262.

\bibitem[Bha97]{Bhargava1997}
M.~Bhargava,
$P$-orderings and polynomial functions on arbitrary subsets of Dedekind rings,
\emph{J. Reine Angew. Math.} \textbf{490} (1997), 101--127.

\bibitem[Bha00]{Bhargava2000}
M.~Bhargava,
The factorial function and generalizations,
\emph{Amer. Math. Monthly} \textbf{107} (2000), no.~9, 783--799.

\bibitem[BFI86]{BFI1986}
E.~Bombieri, J.~B. Friedlander, and H.~Iwaniec,
Primes in arithmetic progressions to large moduli,
\emph{Acta Math.} \textbf{156} (1986), 203--251.

\bibitem[Bom65]{Bombieri1965}
E.~Bombieri,
On the large sieve,
\emph{Mathematika} \textbf{12} (1965), 201--225.

\bibitem[BM22]{BohrMollerup1922}
H.~Bohr and J.~Mollerup,
\emph{L\ae rebog i kompleks Analyse}, vol.~III,
Gjellerup, Copenhagen, 1922.

\bibitem[CC97]{CahenChabert1997}
P.-J. Cahen and J.-L. Chabert,
\emph{Integer-Valued Polynomials},
Math. Surveys Monogr., vol.~48, Amer. Math. Soc., Providence, RI, 1997.

\bibitem[Con90]{Conway1990}
J.~B. Conway,
\emph{A Course in Functional Analysis}, 2nd ed.,
Grad. Texts in Math., vol.~96, Springer-Verlag, New York, 1990.

\bibitem[DeS05]{DeSoleKac2005}
A.~De Sole and V.~G. Kac,
On integral representations of $q$-gamma and $q$-beta functions,
\emph{Atti Accad. Naz. Lincei Rend. Lincei Mat. Appl.} \textbf{16} (2005), 11--29.

\bibitem[DI82]{DeshouillersIwaniec1982}
J.-M. Deshouillers and H.~Iwaniec,
Kloosterman sums and Fourier coefficients of cusp forms,
\emph{Invent. Math.} \textbf{70} (1982), 219--288.

\bibitem[DN26]{DiazNormanyo2026}
B.~Diaz and P.~Normanyo,
Dirichlet series and asymptotics for generalized Legendre factorials,
\emph{arXiv:2603.13720} (2026).

\bibitem[DFI97]{DukeFriedlanderIwaniec1997}
W.~Duke, J.~B. Friedlander, and H.~Iwaniec,
Bilinear forms with Kloosterman fractions,
\emph{Invent. Math.} \textbf{128} (1997), 23--43.

\bibitem[GR04]{GasperRahman2004}
G.~Gasper and M.~Rahman,
\emph{Basic Hypergeometric Series}, 2nd ed.,
Encyclopedia Math. Appl., vol.~96, Cambridge Univ. Press, Cambridge, 2004.

\bibitem[IK04]{IwaniecKowalski2004}
H.~Iwaniec and E.~Kowalski,
\emph{Analytic Number Theory},
Amer. Math. Soc. Colloq. Publ., vol.~53, Amer. Math. Soc., Providence, RI, 2004.

\bibitem[Jac05]{Jackson1905}
F.~H. Jackson,
The basic gamma-function and the elliptic functions,
\emph{Proc. Roy. Soc. London Ser. A} \textbf{76} (1905), 127--144.

\bibitem[Kat95]{Kato1995}
T.~Kato,
\emph{Perturbation Theory for Linear Operators},
Classics in Mathematics, Springer-Verlag, Berlin, 1995.

\bibitem[Mac98]{MacLane1998}
S.~Mac Lane,
\emph{Categories for the Working Mathematician}, 2nd ed.,
Grad. Texts in Math., vol.~5, Springer-Verlag, New York, 1998.

\bibitem[MV07]{MontgomeryVaughan2007}
H.~L. Montgomery and R.~C. Vaughan,
\emph{Multiplicative Number Theory I: Classical Theory},
Cambridge Stud. Adv. Math., vol.~97, Cambridge Univ. Press, Cambridge, 2007.

\bibitem[ST43]{ShohatTamarkin1943}
J.~A. Shohat and J.~D. Tamarkin,
\emph{The Problem of Moments},
Math. Surveys, vol.~1, Amer. Math. Soc., New York, 1943.

\bibitem[Var97]{Vaughan1997}
R.~C. Vaughan,
\emph{The Hardy--Littlewood Method}, 2nd ed.,
Cambridge Tracts in Math., vol.~125, Cambridge Univ. Press, Cambridge, 1997.

\bibitem[Var88]{Vardi1988}
I.~Vardi,
Determinants of Laplacians and multiple gamma functions,
\emph{SIAM J. Math. Anal.} \textbf{19} (1988), no.~2, 493--507.

\bibitem[Wid41]{Widder1941}
D.~V. Widder,
\emph{The Laplace Transform},
Princeton Math. Ser., vol.~6, Princeton Univ. Press, Princeton, NJ, 1941.

\end{thebibliography}
\end{document}